\documentclass{memo-l}
\usepackage[utf8]{inputenc}
\usepackage{hyperref}
\usepackage[linesnumbered, norelsize]{algorithm2e}
\usepackage{
graphicx,amsmath,amssymb,enumerate,color,amsthm,psfrag}
\renewcommand{\P}{\mathbb{P}}

\newcommand{\W}{\mathbb{W}}
\newcommand{\F}{\mathcal{F}}

\newcommand{\G}{\mathcal{G}}
\newcommand{\R}{\mathbb{R}}
\newcommand{\N}{\mathbb{N}}
\newcommand{\Z}{\mathbb{Z}}
\newcommand{\E}{\mathbb{E}}
\newcommand{\1}{\mathbf{1}}
\newtheorem{teor}{Theorem}[chapter]
\newtheorem{lemaa}[teor]{Lemma}
\newtheorem{defin}[teor]{Definition}
\newtheorem{posl}[teor]{Corollary}
\newtheorem{propo}[teor]{Proposition}

\numberwithin{figure}{chapter}
\numberwithin{section}{chapter}
\numberwithin{equation}{chapter}

\newcommand{\be}{\begin{equation}}
\newcommand{\ee}{\end{equation}}
\newcommand{\mic}[1]{\left\lfloor #1\right \rfloor}

\newcommand{\defi}[1]{\defin{\sl #1}\rm}

\newcommand{\teo}[1]{\teor{\it #1}\rm}
\newcommand{\lem}[1]{\lemaa{\it #1}\rm}
\newcommand{\pos}[1]{\posl{\it #1}\rm}
\newcommand{\prop}[1]{\propo{\it #1}\rm}
\newcommand{\dok}[1]{\renewcommand\qedsymbol{$\blacksquare$} \begin{proof}[\underline{Proof}.] #1 \end{proof} }
\newcommand{\cM}{\mathcal{M}}
\newcommand{\M}{\mathcal{M}}%
\newcommand{\cal}{\mathcal }
\newcommand{\cN}{\mathcal{N}}
\newcommand{\cV}{\mathcal{V}}
\newcommand{\cE}{\mathcal{E}}
\newcommand{\cB}{\mathcal{B}}
\newcommand{\cK}{\mathcal{K}}

\newcommand{\cH}{\mathcal{H}}
\newcommand{\cG}{\mathcal{G}}
\newcommand{\cP}{\mathcal{P}}
\newcommand{\cS}{\mathcal{S}}
\newcommand{\cT}{\mathcal{T}}

\title{Time-like graphical models}
\date{May 11, 2014}

\subjclass[2010]{Primary  60G20,  60G60,  60H15, 60J65, 60J80, 62H05; Secondary 05C99}

\keywords{Stochastic processes indexed by graphs, graphical models, time-like graphs, martingales indexed by directed sets, stochastic heat equation}

\dedicatory{}

\author{Tvrtko Tadi\'c}
\address{University of Washington, Department of Mathematics,
Box 354350, Seattle, WA 98195-4350, USA \and
University of Zagreb, Department of Mathematics, Bijeni\v{c}ka cesta 30, 10000 Zagreb, Croatia}
\curraddr{Microsoft Corporation (City Center Plaza Bellevue), One Microsoft Way, Redmond, WA 98052, USA}
\email{tvrtko@math.hr}

\makeindex
\begin{document}

\begin{abstract}
We study continuous processes indexed by a special family of graphs. 
Processes indexed by vertices of graphs are known as {\it probabilistic graphical models}.
In 2011, Burdzy and Pal proposed a continuous version of graphical models indexed by graphs with an embedded time structure -- so called time-like graphs.
We extend the notion of time-like graphs and find properties of processes indexed by them. 
In particular, we solve the conjecture of uniqueness of the distribution for the process indexed by graphs with infinite number of vertices. 
We provide a new result showing the stochastic heat 
equation as a limit of the sequence of natural Brownian motions on time-like graphs. In addition, our treatment of time-like graphical models reveals
connections to Markov random fields, martingales indexed by directed sets and branching Markov processes.
\end{abstract}
\maketitle

\tableofcontents

\mainmatter

\chapter*{Introduction}
In the general theory of Markov processes (such as given in the book \cite[Blumenthal-Getoor]{blumenthal})
we have a process $X$ indexed by some \textbf{parameter set $T$}:
\be (X(t):t\in T).\ee
The set $T$ can be any set with some \textbf{order} $\preceq$. The book \cite{multparpr} by Khoshnevisan studies different 
cases of multiparameter processes ($T\subset \R^n$).  $T$ could, for instance, be
vertices of a directed graph with the order induced by the direction of edges.\vspace{0.1cm}

\begin{figure}[ht]
\begin{center}
 \includegraphics[width=5.5cm]{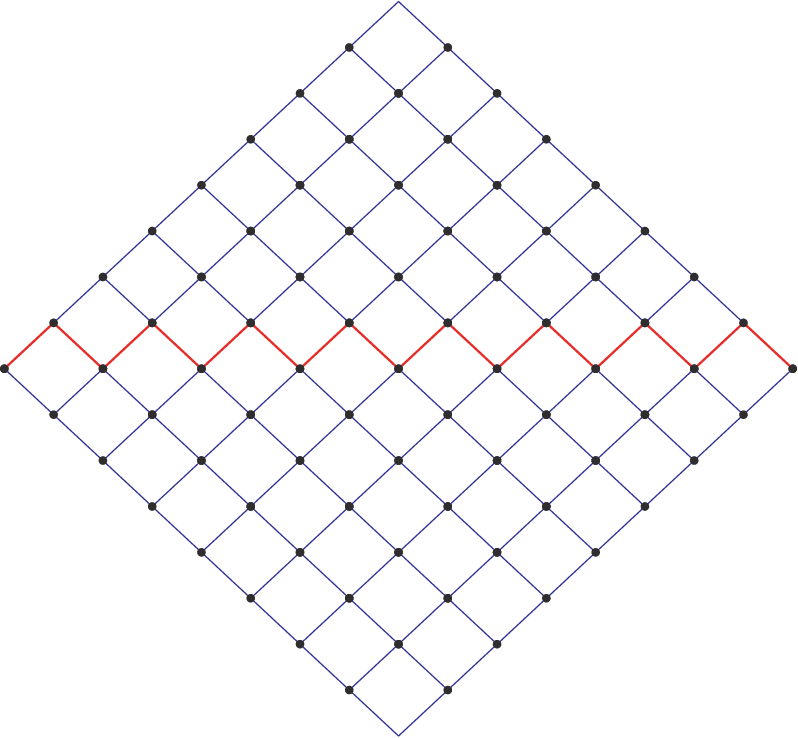}\quad \includegraphics[width=5cm]{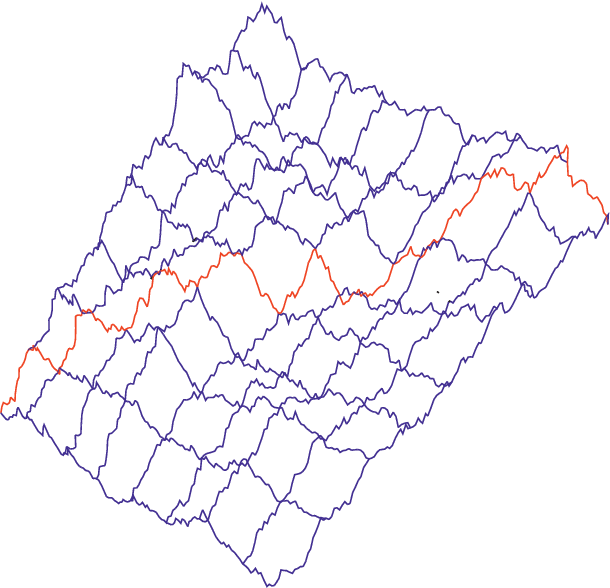}
\end{center}

\caption{Parameter set $T$ and the realization of the process indexed by $T$.}
\end{figure}

Processes indexed by vertices of graphs are well studied and are often
used in machine learning (\cite[Koller - Friedman]{prgphmdl}, \cite[Hastie et al.]{stat_learning}) and statistics (\cite[Lauritzen]{graph_models}, \cite[Studen\'{y}]{studeni}),  where they are called {\it probabilistic graphical models}.
In each of these models the \textbf{conditional independencies} can be read from the \textbf{structure of the graph}.  (A short introduction to 
undirected graphical models is given in \S \ref{sec:mrf}.) Graphical models have been intensively studied in the area of algebraic statistics (\cite[Drton et. al]{algebraicStatistics}), where techniques from algebraic geometry have been successfully used to study properties of conditional independence.
In probability, Markov processes indexed by trees have been studied (see \cite[Benjamini - Peres]{mchidxtr}), as well as Gibbs processes.\vspace{0.1cm}
\begin{figure}[ht]
\begin{center}

\includegraphics[width=6cm]{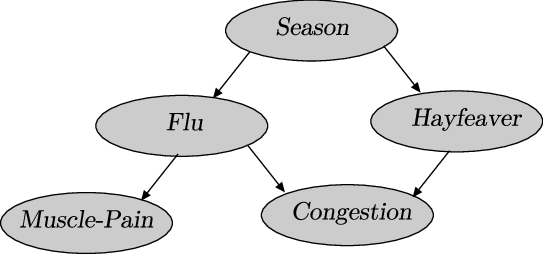}\quad \quad \includegraphics[width=6cm]{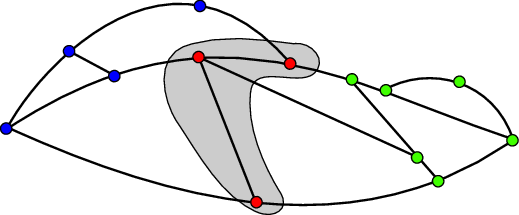}

\caption{In graphical models the structure of the graph induces conditional independencies.}\label{sl41}
\end{center}
\end{figure}

%

%
Similar continuous models such as the { branching  Brownian motion} (\cite{etheridge}), { Le Gall's Brownian snake} (\cite{etheridge}), { Brownian web} (\cite{bweb}) and {Brownian net} (\cite{bnet}) have been studied. In all
these models the underlying graph is a \textbf{random graph}.\vspace{0.1cm}

In 2011, Burdzy and Pal (\cite{tlg1}) introduced 
time-like graphs  (TLG's) and defined (Markov) processes on
graphs with no co-terminal cells (NCC-graphs). Compared to graphical models, these were {\bf continuous} processes (they have a random variable defined
at each point of the representation), and unlike the continuous 
models studied in probability,  the underlying graph was {\bf deterministic}. 
A number of properties (induced by the structure of the underlying graph) of these processes were proved. However, the model had 
 \textbf{strong restrictions} both on the degrees of vertices of the graph and the distribution of the process. \vspace{0.1cm}

In this paper we expand the definition of processes onto a wider 
family of graphs, answer open questions asked by Burdzy and Pal, and investigate
new properties and connections with some known processes.\vspace{0.1cm}

\begin{figure}[ht]
\begin{center}
\psfrag{0}{$\boldsymbol{0}$} 
\psfrag{1}{$\boldsymbol{1}$}
\psfrag{t}{\small $\boldsymbol{t}$}
\psfrag{a}{$\boldsymbol{1/3}$}
\psfrag{b}{$\boldsymbol{2/3}$}
\includegraphics[width=7cm]{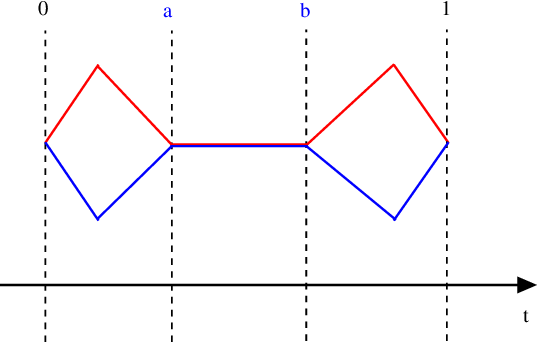}\quad\quad \includegraphics[width=4.5cm]{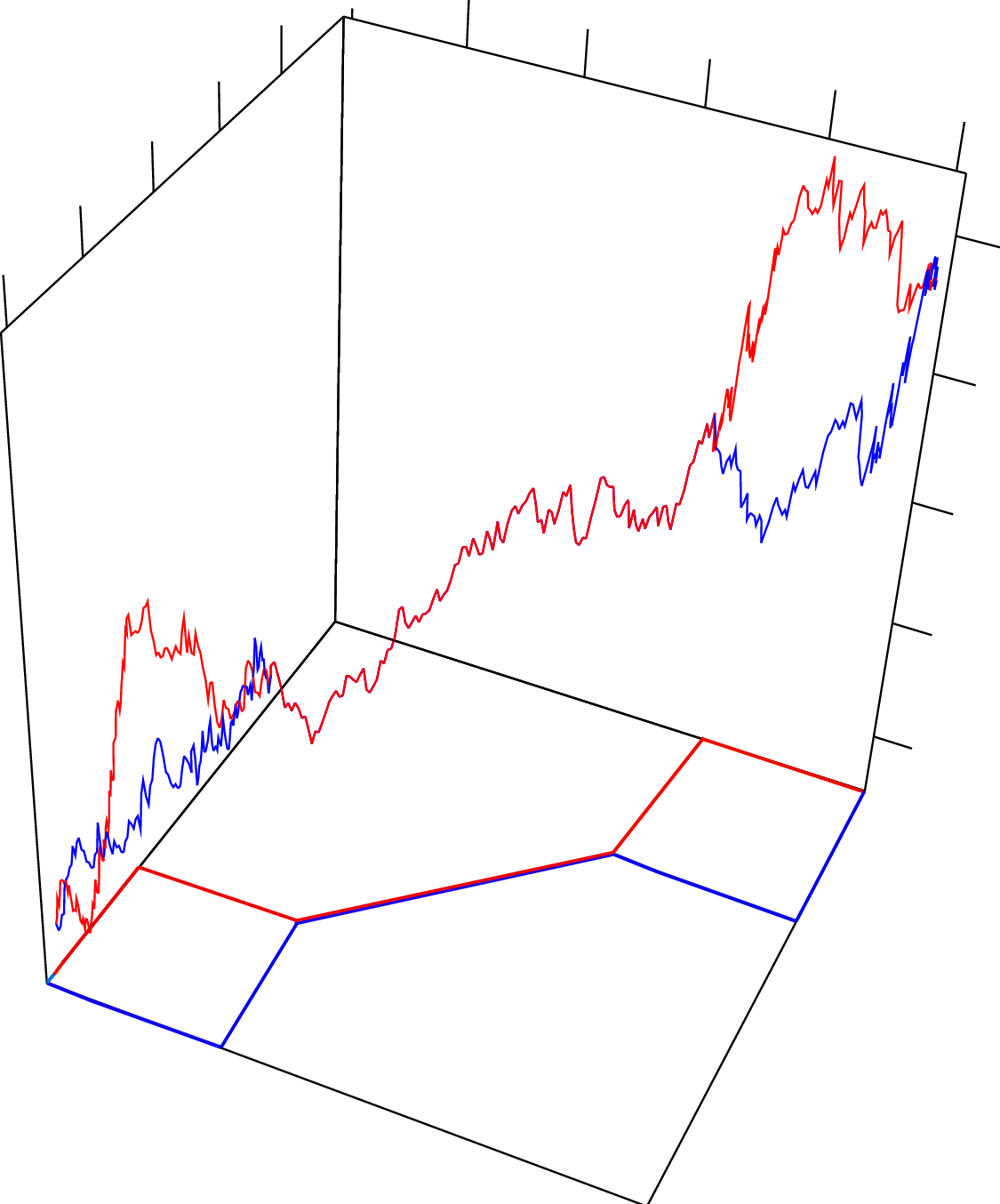} \\

\caption{ Time-like graph and a process indexed by it.} \label{pic:ex2pm4}
  \end{center}
\end{figure}

This paper has three main parts, and it ends with a list of open problems and
an appendix that contains definitions of some terms that might not be familiar to the reader.

\section*{Construction and properties}
In \S \ref{sec:01} we study the geometry of time-like graphs (TLG's). We are
focused on the TLG$^*$ family, since the processes
on this sub-family of time-like graphs can be well-defined. Many of the properties depend 
on the structure of the underlying time-like graph, so we investigate the \textbf{properties} 
and lastly give an \textbf{algorithm} for determining whether a graph belongs 
to the TLG$^*$ family.\vspace{0.15cm}

In \S \ref{sec:02} we give a very general criteria for constructing a process 
indexed by a TLG$^*$ $\cG$ (see \S \ref{condit}). Further, we show that the constructed 
process has the hereditary spine-Markovian property (see \S \ref{sec:spine2}) and we get that the distribution of the process
does not depend on its construction (see Theorem \ref{teo:uniq}). Burdzy and Pal (in \cite{tlg1}) \textbf{conjectured} 
that this holds for NCC graphs with infinitely many vertices. This is
proven here in a much more general setting (Theorem \ref{teo:burdzy-pal}).\vspace{0.15cm}

In \S\ref{sec:3a} we look into several properties of the constructed process induced by time and graph structure. Theorem \ref{thm:tm_mk} proves that a \textbf{generalized Markov property} holds, while Theorem \ref{thm:mrl_gp_mk} shows the connection between the 
constructed process and \textbf{Markov random fields}.\vspace{0.15cm}

Kurtz \cite{kurtz} studied \textbf{martingales} that are indexed by directed sets.  Theorem \ref{teo:tpltt} shows that every TLG$^*$ $\cG$
is a directed set, and under some conditions 
the process indexed by $\cG$ will be a martingale. In \S\ref{sec:3b} we develop stopping times and look at the properties
of filtrations to prove the \textbf{Optional Stopping Theorem} (Theorem \ref{teo:OST}) for martingales indexed by TLG$^*$'s.\vspace{0.15cm}

\section*{Natural Brownian motion and the stochastic heat equation}
In Part 2 we investigate another question from the original paper \cite{tlg1}. What happens when we have 
a process on a \emph{dense net} that covers (a subset of) the plane?  
In \S \ref{sec:04} 
we look at a rhombus grid that covers the whole plane and the two sided Brownian motion 
defined on this graph. We analyze 
what happens when the mesh size goes to zero, and study the connection with
the \textbf{stochastic heat equation} (Theorem \ref{teo:eul_grd}).\vspace{0.05cm}

In chapters  \S \ref{sec:03} and \S\ref{chp:heatrw} we develop tools to prove
the result about the stochastic heat equation. \S \ref{sec:03} reviews some results about maximums of Gaussian vectors and 
continuous Gaussian process. \S\ref{chp:heatrw} studies the approximation of the (stochastic) heat equation 
with one boundary and an initial value condition with the \textbf{Euler method} under 
very general conditions. The main tool for the analysis is the \textbf{simple random walk}.
\vspace{0.15cm}

  \begin{figure}[ht]
\begin{center}
\includegraphics[width=8cm]{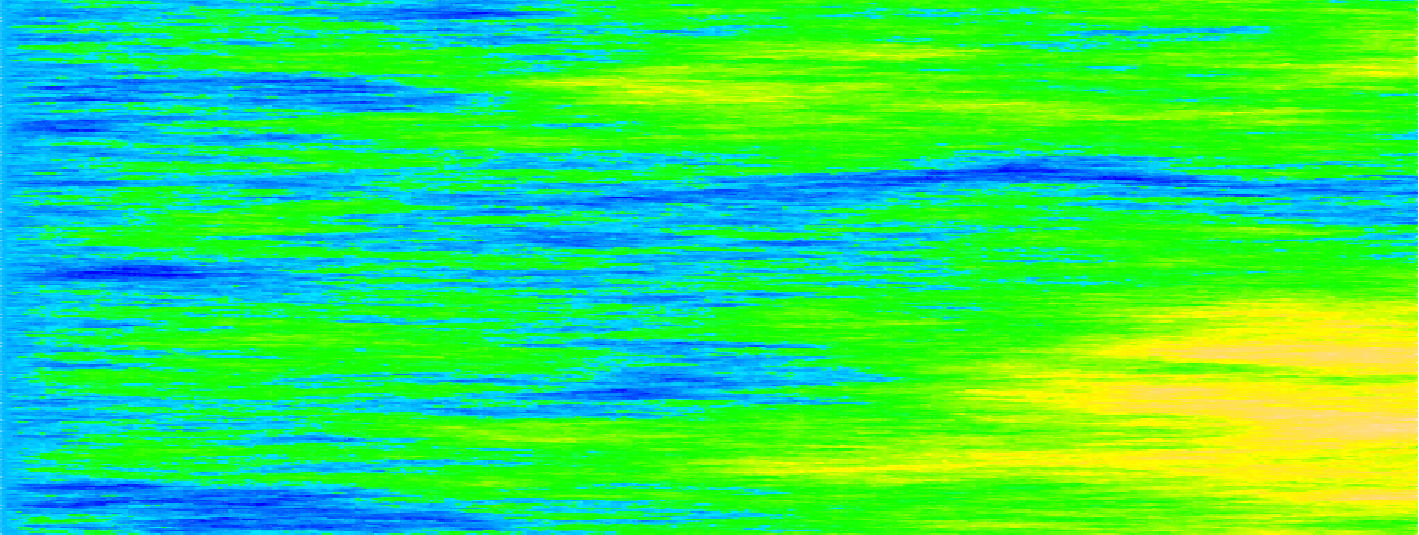}
\caption{Topographical image of the simulation of a process indexed by a dense rhombus grid}
\end{center}
\end{figure}

\section*{Processes on general and random time-like graphs}

The graphs used in Part 1 and 2 have one beginning and one end, so we can not define a process on a time-like tree. 
In \S\ref{sec:7} we  modify our approach to define a process with nice properties on a more general 
family of graphs - TLG$^{**}$'s. This family includes \textbf{trees}, and it turns out (see \S\ref{sec:ptlt}) that some properties which do not hold in general 
are true for time-like trees. We proceed to define \textbf{Galton-Watson time-like trees} (\S\ref{GWTLT}),
and investigate (\S\ref{sec:mbmpp}) what happens when we index the process by this type of random trees.\vspace{0.15cm}

\section*{Open questions and appendix}
This paper ends with several open questions: Under what conditions 
can we define a process on any TLG? If we know the process 
on some parts of the graph, what can we tell about the parts that are hidden from us?
Do we (under some conditions) have the strong Markov property? How would we model the 
evolution of the process on a graph over time?

 \begin{figure}[ht]
\begin{center}
\psfrag{a}{$\boldsymbol{-1}$}
\psfrag{b}{$\boldsymbol{2}$}
\psfrag{0}{$\boldsymbol{0}$}
\psfrag{1}{$\boldsymbol{1}$}

\includegraphics[width=6.8cm]{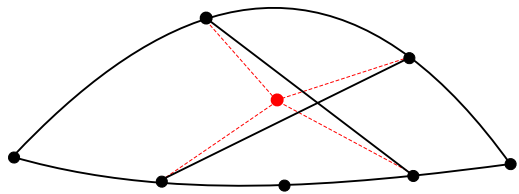}
\caption{Open question: We know about the black parts of the graph and the process on it, what can we say about 
the part of the structure that is \textcolor{red}{hidden}?}
  \end{center}
\end{figure}

The Appendix contains some definitions and known results that we will often use.

\part{Construction and properties}
\chapter{Geometry of time-like graphs}\index{Time-like graph (TLG)|(}\label{sec:01}

\begin{figure}[ht]
\begin{center}
\includegraphics[width=8cm]{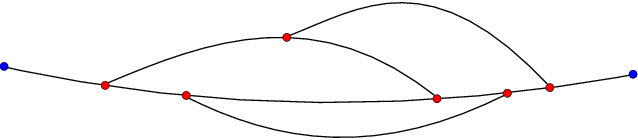}\\
\caption{ }\label{bptlg13:1}
\end{center}
\end{figure}

Most of the definitions presented in this chapter are modified from the original model presented in \cite{tlg1}.
The crucial difference is the Definition \ref{def:tlg} of time-like graphs. In the original model, Burdzy and Pal considered time-like graphs with the beginning and end vertex of degree 1,
and all other vertices of degree 3.
See Figure \ref{bptlg13:1}.

The rest is a deeper study of geometric properties of the special family TLG$^*$. These 
properties will later be vital for the construction of the processes and many of their 
properties.

\section{Basic definitions}
\defi{\label{def:tlg}A graph $\cG=(\cV,\cE)$ will be called a  \textbf{time-like graph (TLG)} if its sets 
of vertices $\cV$ and edges $\cE$ satisfy the following properties.
\begin{enumerate}[(i)]
 \item The set $\cV$ contains at least two elements, $\cV=\{t_0,t_1,\ldots,t_N\}$,
where $t_0=A$, $t_N=B$ and for $k=1,2,\ldots,N-2$, $A< t_k\leq t_{k+1}< B$.
 \item An edge between $t_j$ and $t_k$ will be denoted $E_{jk}$. We assume that there 
is no edge between $t_j$ and $t_k$ if $t_j=t_k$. $E_{jk}$ indicates that $t_j<t_k$. 
\item We assume that  all vertices
have a finite degree. 
\item We assume that for every vertex $t_k$ $k=1,\ldots, N-1$ there
exist edges $E_{jk}$ and $E_{kn}$ with $j<k<n$.

\end{enumerate}
We call TLG to be a \textbf{unit} TLG if $A=0$ and $B=1$.
}\index{Time-like graph (TLG)|textbf}

{\it Remarks.} (0) In our study of TLG's, we will assume that TLG is a unit TLG, unless
specified differently. (1) We do not exclude the case $\cV=\{t_0=A,t_N=t_1=B\}$.\par (2) The definition implies that TLG has no loops.  
\par (3) In (i) formally we should say that the elements have the form $(k,t_k)$,
so that $(k,t_k)$ and $(k+1,t_{k+1})$ are distinct even if $t_k=t_{k+1}$. This notation was simplified
to make writing easier.\par (4) An edge between $t_j$ and $t_k$ ($j<k$) will be denoted 
$E_{jk}$ (if it exists), and if we are using more of them we will use the notation $E_{jk}^1$, $E_{jk}^2$ \ldots
(or something similar).

 \begin{figure}[ht]
\begin{center}
\begin{minipage}{7cm}
\psfrag{0}{$\boldsymbol{t_0=0}$\quad} \psfrag{1}{$\boldsymbol{t_1}$}
\psfrag{2}{$\boldsymbol{t_2}$}\psfrag{3}{$\boldsymbol{t_3}$}
\psfrag{4}{$\boldsymbol{t_4}$} \psfrag{5}{$\boldsymbol{t_5}$}
\psfrag{6}{$\boldsymbol{t_6=1}$} \psfrag{a}{\textcolor{red}{$\boldsymbol{E_{36}^1}$}}
\psfrag{b}{\textcolor{red}{$\boldsymbol{E_{36}^2}$}} \psfrag{c}{\textcolor{blue}{$\boldsymbol{E_{02}}$}}
\includegraphics[width=7cm]{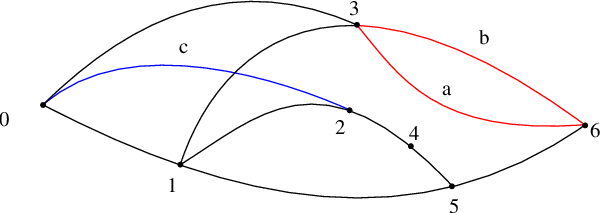}\\
\caption{TLG $\cG$}
\end{minipage}
\quad\quad
\begin{minipage}{7cm}
 \psfrag{0}{$\boldsymbol{0}$\quad} \psfrag{1}{$\boldsymbol{t_1}$}
\psfrag{2}{$\boldsymbol{t_2}$}\psfrag{3}{$\boldsymbol{t_3,t_4}$}
\psfrag{4}{$\boldsymbol{t_5}$} \psfrag{5}{$\boldsymbol{1}$}
\psfrag{R}{$\R^2$}\psfrag{t}{$\boldsymbol{t}$}
\includegraphics[width=7cm]{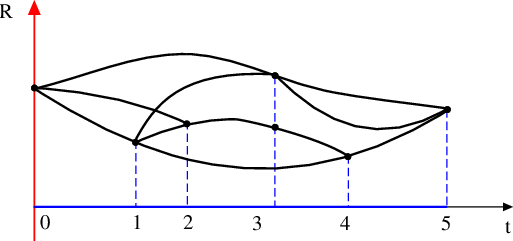}\\ \caption{Representation of a TLG $\cG$}
\end{minipage}

\end{center}
\end{figure}


The representation\index{Time-like graph (TLG)!representation} of a TLG in $\R^3$ is given by the following definition.

\defi{\label{def:rptlg}By abuse of notation let $E_{jk}:[t_j,t_k]\to \R^2$ denote a continuous function for all $E_{jk}\in \cE$.
Assume: 
\begin{enumerate}[(i)]
 \item That the images of the open sets $(t_j,t_k)$ under the maps
$t\mapsto (t,E_{jk}(t))$, where $E_{jk}\in \cE$ are disjoint.
\item That $E_{jk}(t_k)=E_{kn}(t_k)$ if $E_{jk}, E_{kn}\in \cE$; $E_{jk}(t_k)=E_{mk}(t_k)$
if $E_{jk},E_{mk}\in \cE$; and $E_{0k}(t_0)=E_{0j}(t_0)$ for $E_{0k},E_{0j}\in \cE$.
\end{enumerate}
We will call the set 
$$R(\cG)=\{(t,E_{jk}(t))\in [0,1]\times \R^2:E_{jk}\in \cE, t\in [t_j,t_k] \}$$
a \textbf{representation} of $\cG$. We will say that $\cG_1$ is a subgraph of $\cG_2$, and write
$\cG_1\subset \cG_2$  if there exist representations of the two such that $R(\cG_1)\subset R(\cG_2)$.
We will call $\cG$ \textbf{planar}\index{Time-like graph (TLG)!planar} if it has a representation $R(\cG)\subset \R^2$. 

Let $\bar{t}_j= (t_j,E_{jk}(t_j))$ for $j<N$ and $\bar{t}_N=(t_N,E_{N-1,N}(t_N))$. } 


\noindent {\it Remark.} There are many representations for a TLG, but there is a unique
TLG corresponding to a representation.


\begin{figure}[ht]
\begin{center}
\psfrag{a}{$\boldsymbol{E_{k_1k_2}}$}
\psfrag{b}{$\boldsymbol{E_{k_2k_3}}$}
\psfrag{c}{$\boldsymbol{E_{k_3k_4}}$}
\psfrag{d}{$\boldsymbol{\cdots}$}
\psfrag{e}{$\boldsymbol{E_{k_{n-1}k_n}}$}
\includegraphics[width=8cm]{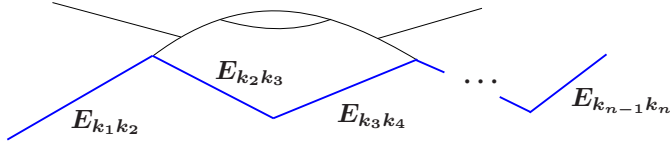}\\
\caption{A time-path}
\end{center}
\end{figure}

\defi{\label{def:tmpth}We will call a sequence of edges 
\be (E_{k_1k_2},E_{k_2k_3},\ldots, E_{k_{n-1}k_n}) \label{eq:tp1}\ee
a \textbf{time-path}\index{Time-like graph (TLG)!time-path}\index{Time-path|see{Time-like graph (TLG)}} if $E_{k_jk_{j+1}}\in\cE$ for every $j$. We will denote the set of all paths of the form $(\ref{eq:tp1})$
by $\sigma(k_1,k_2,\ldots, k_n)$. This time path is \textbf{full time-path}\index{Time-like graph (TLG)!time-path!full-time path}\index{Full time-path|see{Time-like graph (TLG)}} if $k_0=0$ and $k_n=N$.
We will denote the set of all full time-paths by $P_{0\to 1}(\cG)$.}

\noindent \emph{Remark.} Note that the notation $\sigma (k_1,k_2,\ldots,k_n)$ does not uniquely identify 
the path, since there can be more than one edge between the two vertices. 

 \begin{figure}[ht]
\begin{center}
\includegraphics[width=6cm]{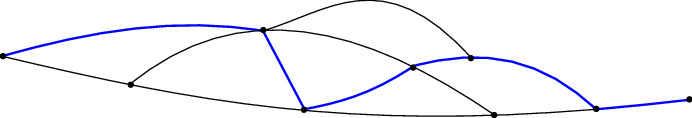}\quad \includegraphics[width=6.5cm]{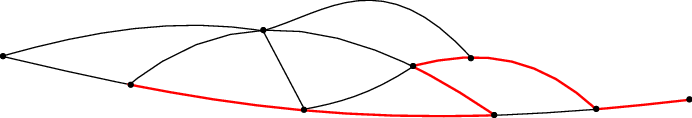}\\
\caption{Example of a \textcolor{blue}{full time-path} and an example of 
path that is \textcolor{red}{not a time-path}}

\end{center}
\end{figure}

\defi{\label{defct}\begin{enumerate}[(i)]
       \item Time paths $\sigma_j\in\sigma(j_1,j_2,\ldots,j_n)$ and $\sigma_k\in\sigma(k_1,k_2,\ldots,k_m)$
are \textbf{co-terminal} if $j_1=k_1$ and $j_n=k_m$.
      \item Co-terminal paths $\sigma_j\in\sigma(j_1,j_2,\ldots,j_n)$ and $\sigma_k\in\sigma(k_1,k_2,\ldots,k_m)$ will be
form a \textbf{cell} $(\sigma_j,\sigma_k)$ if $$\{j_2,j_3,\ldots,j_{n-1}\}\cap\{k_2,\ldots,k_{m-1}\}=\emptyset . $$
\item We will call a cell\index{Time-like graph (TLG)!cell} $(\sigma_j,\sigma_k)$ for $\sigma_j\in\sigma(j_1,j_2,\ldots,j_n)$ and $\sigma_k\in\sigma(k_1,k_2,\ldots,k_m)$ \textbf{simple} 
if if there does not exist a time path $\pi\in\sigma(i_1,i_2,\ldots,i_r)$ such that $i_1\in \{j_2,j_3,\ldots,j_{n-1}\}$
and $i_r\in \{k_2,\ldots,k_{m-1}\}$, or $i_1\in \{k_2,\ldots,k_{m-1}\}$ and $i_r\in \{j_2,j_3,\ldots,j_{n-1}\}$. 

      \end{enumerate}
}

\begin{figure}[ht]
\begin{center}

\includegraphics[width=9cm]{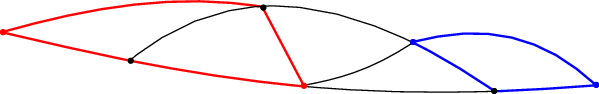}\\
\caption{\textcolor{red}{Non-simple cell} and a \textcolor{blue}{simple cell}}

\end{center}
\end{figure}


\section{TLG$^*$ family}\index{TLG$^*$ family|(}
We will now describe the family of TLG graphs that is generated from minimal graph 
by adding vertices and adding edges between vertices connected by a time-path.

\defi{\label{def:tlg*}The \textbf{TLG$^*$-family} is given in the following inductive way.
\begin{enumerate}[(i)]
 \item The minimal graph $\cG=(\cV,\cE)$, with $\cV=\{t_0=1,t_N=1\}$ and $\cE=\{E_{0N}\}$ is a 
TLG$^*$.

\begin{figure}[ht]
\begin{center}
\psfrag{0}{$\boldsymbol{t_0}$}
\psfrag{1}{$\boldsymbol{t_1}$}
\psfrag{E}{$\boldsymbol{E_{01}}$}
\includegraphics[height=0.5cm]{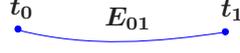}\\
\caption{The minimal graph}
\end{center}
\end{figure}

\item \label{def:tlg*:ii} Let $\cG_1=(\cV_1,\cE_1)$ be a TLG$^*$, where $\cV_1=\{t_0,t_2,\ldots,t_N\}$.
\begin{enumerate}[(1)]
 \item \label{def:tlg*:1}{\sc (adding a vertex)} If $\tau_k\in [0,1]$, and for some $E_{k_1k_2}\in \cE$ and
$t_{k_1}< \tau_k< t_{k_2}$ then set
$$\cV_2:=\cV_1\cup\{\tau_k\}\quad \textrm{and}\quad \cE_2:=\cE_1\cup \{E_{k_1k},E_{kk_2} \}\setminus \{E_{k_1k_2}\}.$$
$\cG_2:=(\cV_2,\cE_2)$ is also a TLG$^*$.

\begin{figure}[ht]
\begin{center}
\psfrag{a}{$\boldsymbol{t_{k_1}}$}
\psfrag{b}{$\boldsymbol{t_{k_2}}$}
\psfrag{c}{$\boldsymbol{\tau_{k}}$}
\psfrag{E}{$\boldsymbol{E_{k_1k_2}}$}
\psfrag{1}{$\boldsymbol{E_{k_1k}}$}
\psfrag{2}{$\boldsymbol{E_{kk_2}}$}
\includegraphics[height=2cm]{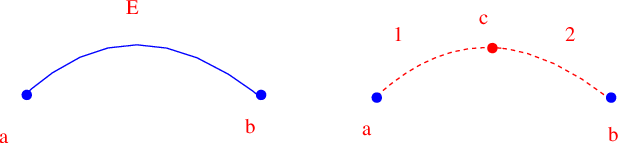}\\
\caption{Adding a vertex}
\end{center}
\end{figure}

\item {\sc (adding an edge)} Let $t_j,t_k\in \cV_1$ such that $t_j<t_k$, and assume that there exists a time-path $\sigma_{jk}\in\sigma(j,\ldots,k)$ 
between these vertices. Then set
$$\cV_2:=\cV_1\quad \textrm{and}\quad \cE_2:=\cE_1\cup \{E_{jk}^* \}.$$
$\cG_2:=(\cV_2,\cE_2)$ is also a TLG$^*$. ($E^{*}_{jk}$ is an new edge (not in $\cE_1$).)

\begin{figure}[ht]
\begin{center}
\includegraphics[height=1.3cm]{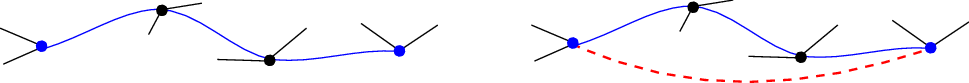}\\
\caption{Adding the edge \textcolor{red}{$E_{jk}^*$}}
\end{center}
\end{figure}

\end{enumerate}

\item \label{def:tlg*:iii} We will say that $(\cG_j)_{1\leq j\leq k}$ is a \textbf{tower of TLG$^*$'s} or \textbf{TLG$^*$-tower}\index{TLG$^*$ family!TLG$^*$-tower}\index{TLG$^*$-tower|see{TLG$^*$ family}} if 
for $j>1$, $\cG_j$ is constructed from $\cG_{j-1}$ as in (ii). 
\end{enumerate}
}
Remarks. 
(1) Clearly, all TLG$^*$'s are TLG's. (2) It is also clear that if $(\cG_j)_{1\leq j\leq k}$ 
is a tower of TLG$^*$'s and $G_k$ is planar that all the graphs in this tower of TLG$^*$'s are planar.\vspace{0.3cm}

We will turn our attention to the question which TLG's are TLG$^*$. The following is a generalization 
and a new proof of the result known to Burdzy and Pal (see Theorem 2.9 (ii) in \cite{tlg1}).

\teo{\label{thm:tlg*1}All planar TLG's are TLG$^*$'s.
\dok{Let $\cG$ be a planar TLG and $R(\cG)$ its representation in $\R^2$. We will prove
the claim in several steps.

(i) Denote time-paths from $t_0=0$ to $t_N=1$ in $\cG$ with $P_{0\to 1}(\cG)$.
For each $\sigma\in P_{0\to 1}(\cG)$ there exists a continuous function
$g_\sigma:[0,1]\to \R$ such that its graph $\Gamma_{g_{\sigma}}=\{(x,g_{\sigma}(x)):x\in [0,1]\}$ 
is the representation of $\sigma$ in $R(\cG)$. 
For two paths $\sigma'\neq \sigma''$ we have $g_{\sigma'}\neq g_{\sigma''}$, and there are three possibilities
\begin{itemize}
\item If $g_{\sigma'}\leq g_{\sigma''}$
or $g_{\sigma'}\geq g_{\sigma''}$. In the first case we say $\sigma' \leq \sigma''$
and in the second case we say $\sigma' \geq \sigma''$.
 \item If not, $\min\{g_{\sigma'},g_{\sigma''} \}$, $\max\{g_{\sigma'},g_{\sigma''} \} $ are also representations of  
paths from 0 to 1. (These paths use the same set of edges as paths $\sigma'$ and $\sigma''$.) 
\end{itemize}

We define $\sigma'\wedge \sigma''$ and $\sigma'\vee \sigma''$ to be the path represented by $\min\{g_{\sigma'}, g_{\sigma''}\}$ and $\max\{g_{\sigma'},g_{\sigma''} \} $ in $R(\cG)$.
This operation is closed, commutative and associative, and further $\sigma'\wedge \sigma''\leq \sigma' \leq \sigma'\vee \sigma''$ and $\sigma'\wedge \sigma''\leq \sigma''\leq \sigma'\vee \sigma''$.

(ii) We pick $\sigma_1$ to be $\wedge_{\sigma\in P_{0\to 1}(\cG)}\sigma$, and we set
$\cG_1=(\cV_1,\cE_1)$ such that all vertices and all edges of $\sigma_1$
are in $\cV_1$ and $\cE_1$. Clearly this is a planar TLG. Note that we choose
$\sigma_1$ such that there is no $\sigma'$ in $P_{0\to 1}(G)$ with $\sigma'\leq \sigma_1$.

Now we continue inductively. Let $\cG_{k-1}=(\cV_{k-1},\cE_{k-1})$ be a TLG obtained in the previous
step. If $\cE\setminus\cE_{k-1}=\emptyset$ clearly $\cG_{k-1}=\cG$. Otherwise,
choose $\sigma_{k}$ in $P_{0\to 1}(\cG)\setminus P_{0\to 1}(\cG_{k-1})$ such that
there is no $\sigma'$ in the same set with $\sigma'\leq \sigma_{k}$. (The set $P_{0\to 1}(\cG)\setminus P_{0\to 1}(\cG_{k-1})$ is nonempty
since every edge $E\in \cE\setminus\cE_{k-1}$ is part of a path from 0 to 1 in $\cG$. 
There is such minimal edge with respect to the given order, since
this is a finite set.) We now set $\cG_{k}=(\cV_{k},\cE_{k})$, where $\cV_{k}$ is the set of all  
vertices in $\cV_{k-1}$ and on $\sigma_k$ and $\cE_{k}$ is set of all  
edges in $\cE_{k-1}$ and that $\sigma_k$ is made of. Again, $\cG_{k}$ is a planar TLG. 

Since there is only a finite number of edges in $\cE$, at some step $K$ we will stop,
and we will have $\cG_K=\cG$.

(iii) Note that for each $k$ there is no edge $E_{jn}\in \cE\setminus\cE_{k}$ such that there exists $ \sigma\in P_{0\to 1}(\cG_{k})$
with $E_{jn}\leq g_{\sigma}|_{[t_j,t_n]}$. Otherwise, there would exist a $\sigma_l$
for some $l\leq k$ such that $E_{jn}\leq g_{\sigma_l}|_{[t_j,t_n]}$, and a path
$\sigma'\in P_{0\to 1}(\cG)$ that contains $E_{jn}$, but then $\sigma'\wedge \sigma_l\leq \sigma_l$,
and this contradicts the definition of $\sigma_l$.

(iv) From the definition in (i) it is clear that \be\sigma_{max}^k=\vee_{\sigma\in P_{0\to 1}(\cG_k)}\sigma \label{s_max}\ee
is also a path in $\cG_k$.

(v) Now we will show that all $\cG_k$ are TLG$^*$'s.
It is clear that $\cG_1$ can be obtained from the minimal 
graph $\cG_0$ by repeating step (\ref{def:tlg*:1}) in Definition \ref{def:tlg*}.

We assume that $\cG_{k-1}$ is a TLG$^*$. For an edge $E_{jn}$ in $\sigma_k$
that is not in $\cE_{k-1}$, we have by (iii) 
\be g_{\sigma_{max}^{k-1}}|_{[t_j,t_k]}\leq E_{jn}.\label{tlg*:j2}\ee

Further, $\sigma_{max}^{k-1}$ (see (\ref{s_max})) will have common vertices with $\sigma_k$ (at least in 0 and 1).
The set $T=\{t\in [0,1]: \sigma_{max}^{k-1}(t)\neq \sigma_k(t)\}$ has
at exactly one connected component. Otherwise, there would exist $t_{l_1}<t_{l_2}\leq  t_{l_3}<t_{l_4}$ in $\partial T$ and 
we would have two sub-paths $\sigma(l_1\ldots l_2)$ and $\sigma(l_3\ldots l_4)$ that 
start and end at vertices that are on $\sigma_{max}^{k-1}$, but since $(\ref{tlg*:j2})$
we have 
$$g_{\sigma_{max}^{k-1}}|_{[t_{l_1},t_{l_2}]\cup [t_{l_3},t_{l_4}]}\leq g_{\sigma_k}|_{[t_{l_1},t_{l_2}]\cup [t_{l_3},t_{l_4}]}$$
(their representations lie above
$R(\cG_{k-1})$).  
But, now $\sigma'$ is represented by $$g_{\sigma'}(t):=\left\{\begin{array}{cc}
                                                       g_{\sigma_k}(t) & t\in [t_{l_1},t_{l_2}]\\
                                                       g_{{\sigma_{max}^{k-1}}}(t) & t\in [t_{l_1},t_{l_2}]^c\\
                                                      \end{array}
 \right.$$
is also a path in $P_{0\to 1}(\cG)\setminus P_{0\to 1}(\cG_{k-1})$, such that $\sigma'\leq \sigma_k$.
This is a contradiction, with the definition of $\sigma_k$. Therefore, $T$
has only one connected component and  $\sigma_{max}^{k-1}$ and $\sigma_k$ have two common
vertices - $t_{l_1}$ and $t_{l_2}$. Since $t_{l_1}$ and $t_{l_2}$ are on the path
$\sigma_{max}^{k-1}$ by Definition \ref{def:tlg*} we can add an edge between them,
and after that add vertices  that are on the path that connects them. All the other edges of the path 
$\sigma_k$ (that are below the path $\sigma_{max}^{k-1}$ in the representation) are already
included in $\cG_{k-1}$ (by (iii)), so we get $\cG_{k}$.

}
}

\begin{figure}[ht]
\begin{center}
\psfrag{0}{$\boldsymbol{0}$} \psfrag{1}{}
\psfrag{2}{}\psfrag{3}{}
\psfrag{4}{} \psfrag{5}{}
\psfrag{6}{$\boldsymbol{1}$}
\includegraphics[width=6cm]{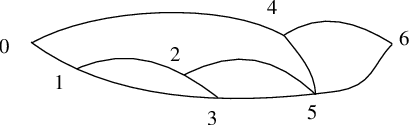}\quad \includegraphics[width=6cm]{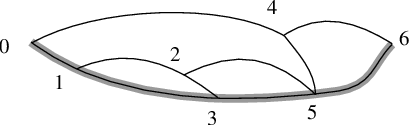}\\
\includegraphics[width=6cm]{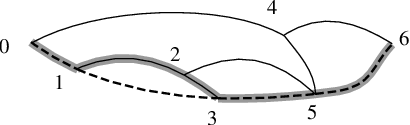}\quad \includegraphics[width=6cm]{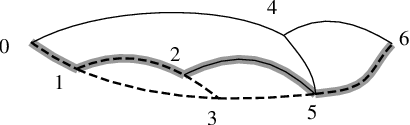}\\
\includegraphics[width=6cm]{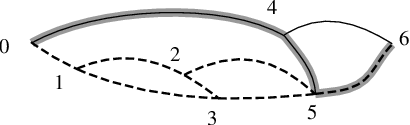}\quad \includegraphics[width=6cm]{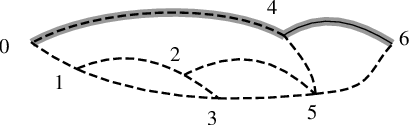}\\
\caption{Illustration of the proof of Theorem \ref{thm:tlg*1}. The the path colored in gray
represents $\sigma_{k}$, while dashed lines represent $\cG_{k-1}$. } \label{pic2}
  \end{center}
\end{figure}

\noindent\emph{Remark.} The proof gives us the following algorithm for constructing
a planar TLG $\cG$ as a TLG$^*$.\vspace{0.3cm}

\begin{algorithm}[H]
 $\sigma$ a minimal path with respect to $\leq$ in $P_{0\to 1}(\cG)$\;
$\cG^\#=(\cV^\#,\cE^\#)$ that consists of all vertices and all edges of $\sigma$ (in $\cG$)\;\label{alg:s2}

\While{$\cE\setminus \cE^\#\neq \emptyset$}{
      $\sigma$ a minimal path with respect to $\leq$ in $P_{0\to 1}(\cG)\setminus P_{0\to 1}(\cG^\#)$\;\label{alg:s6}
      add all edges and vertices that make $\sigma$ (in $\cG$) to $\cG^\#$\; \label{alg:lx}
} 
\caption{Constructing a planar TLG as a TLG$^{*}$.}\label{alg:pltlg}
\end{algorithm}\vspace{0.3cm}

We have shown that the step in line \ref{alg:lx} can be done by adding edges and vertices as described in 
Definition \ref{def:tlg*}. Since $\cG^\#$ is a TLG$^*$ in line \ref{alg:s2}, $\cG^\#$ remains a TLG$^*$
through the whole algorithm. The illustration of this algorithm is given in Figure \ref{pic2}.

\pos{For a planar TLG $\cG$ there exists a tower of planar TLG's (TLG$^*$'s) $(\cG_j)_{1\leq j\leq n}$.
such that $\cG_1=(\{t_0=0,t_N=1\},\{E_{0N}\})$ and $\cG_n=\cG$. Further, there exists a
sequence of representations $(R(G_j))_{1\leq j\leq n}$ such that $R(G_{j-1})\subset R(G_{j})$ for $j>1$. }

\begin{figure}[ht]
\begin{center}
\begin{minipage}{7cm}
 \psfrag{0}{$\boldsymbol{t_0}$} \psfrag{1}{$\boldsymbol{t_1}$}
\psfrag{2}{$\boldsymbol{t_2}$}\psfrag{3}{$\boldsymbol{t_4}$}
\psfrag{4}{$\boldsymbol{t_3}$} \psfrag{5}{$\boldsymbol{t_5}$}
\includegraphics[width=7cm]{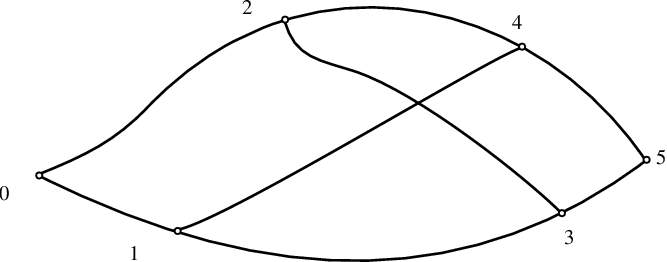}\\
\caption{ A TLG that is not a TLG$^*$} \label{pic1}
\end{minipage}\quad 
\begin{minipage}{7cm}
\includegraphics[width=7cm]{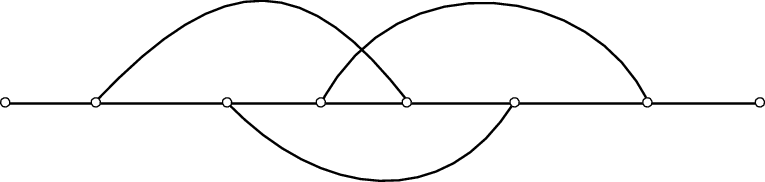}\\
\caption{Non-planar TLG$^*$.} \label{sl3}
 
\end{minipage}

  \end{center}
\end{figure}

\teo{\label{teo:tlg_exmpl}\begin{enumerate}[(i)]
      \item \label{teo:tlg_exmpl:1} There exists a TLG that is not a TLG$^*$.
      \item \label{teo:tlg_exmpl:2} There exists a non-planar TLG$^*$.
     \end{enumerate}
}
\dok{ We will show the claim using examples similar to those Burdzy and Pal gave in \cite{tlg1}.\vspace{0.2cm}

(i) Assume the TLG $\cG=(\cV,\cE)$, where $\cV=\{t_j=j/5:j=0,1,\ldots,5\}$ and
$$\cE=\{E_{01},E_{02}, E_{14}, E_{13},  E_{23},E_{24},E_{45},E_{35}\}$$ (on the Figure \ref{pic1}.) is a TLG$^*$. 
Then there exists a tower of TLG$^*$ $(\cG_j)_{1\leq j\leq n}$ such that $\cG_n=\cG$.
Let $E^*$ be the edge form the set $\cE^*=\{E_{14},E_{13},E_{24}, E_{23}\}$ with 
largest $j$ such that $E^*\in \cE_j\setminus \cE_{j-1}$. ($E^*$ is the last
edge from $\cE^*$ to be added to the graph.) 

In Definition \ref{def:tlg*}. we add 
edges in each step, so that their vertices lie on the same path from 0 to 1 and these
vertices will continue to be on the same path in future steps. Since, no three vertices
from the set $\{t_1,t_2,t_3,t_4\}$ are on the same path in $\cG$, in each step we can
add only one edge from the set $\cE^*$.

The graph $\cG_{j-1}$ contains the vertices $t_1$, $t_2$, $t_3$ and $t_4$,
since it contains three out of four edges from $\cE^*$ connecting them.

In order to obtain $\cG_j$ the endpoints of $E^*$ have to be connected by a time path.
It is clear that each element of the tower  $(\cG_j)_{1\leq j\leq n}$ the number of time
paths between the two vertices increases. This means that the number time paths
between the endpoints of $\cG$ will be at least two, but this is not true in
$\cG$. Hence, $\cG$ can not be a TLG$^*$.\vspace{0.2cm}

(ii) Let $\cG=(\cV,\cE)$, where  $\cV=\{t_j=j/7:j=0,1,\ldots,7\}$ and 
$$\cE=\{E_{01},E_{12},E_{23},E_{34},E_{45},E_{56},E_{67},E_{14},E_{25},E_{3,6}\}.$$
It is clear that this is a TLG$^*$ and it is not planar. See Figure \ref{sl3}
}

\index{TLG$^*$ family|)}
\section{Consistent representation of a TLG$^*$-tower, spines and (re)construction}
\begin{figure}[ht]
\begin{center}
\psfrag{a}{$\boldsymbol{a}$} \psfrag{b}{$\boldsymbol{b}$}
\includegraphics[width=7cm]{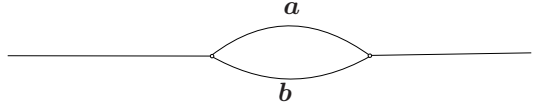}\\
\caption{The arc $a$ and $b$ might not always represent the same edges.} \label{pic11}

  \end{center}
\end{figure}

If $\cG$ is a TLG$^*$, then let $(\cG_j)_{j=0}^n$ be a TLG$^*$ tower.
In the corresponding sequence of representations $(R(\cG_j))_{j=0}^n$
we could have some inconsistencies which we would like to avoid.
For instance, let co-terminal edges $E^1=E^1_{m_1m_2}$ and $E^2=E^2_{m_1m_2}$ be present in
the whole tower and the graph in the Figure \ref{pic11} can represent part of 
each representation. The arcs $a$ and $b$ in representation $R(\cG_{j_1})$ might
represent $E^1$ and $E^2$, while in some other representation $R(\cG_{j_2})$ it might 
be the other way around. To avoid this we will only use 
\textit{consistent representations} of the TLG$^*$ tower $(\cG_j)_{j=0}^n$.

\defi{\label{def:consist_rep}We will call a sequence of representations $(R(\cG_j))_{j=0}^n$
a \textbf{consistent representation}\index{TLG$^*$ family!consistent representation} of the TLG$^*$-tower $(\cG_j)_{j=0}^n$ if:
\begin{enumerate}[(a)]
 \item If we add a new vertex $\tau_k$ to the TLG$^*$ $\cG_{j-1}$ 
to obtain $\cG_j$ by  removing an edge $E_{k_1k_2}$, and replacing it with
$E_{k_1k}$ and $E_{kk_2}$ (as in step (\ref{def:tlg*:1}) of Definition \ref{def:tlg*}.),  
then the representation of edges  $E_{k_1k}$ and $E_{kk_2}$ is the same as 
that of $E_{k_1k_2}$, i.e. $$E_{k_1k_2}([t_{k_1},t_{k_2}])=E_{k_1k}([t_{k_1},t_{k}])\cup E_{kk_2}([t_{k},t_{k_2}]).$$
\item All the edges that are in both $\cG_{j-1}$ and $\cG_{j}$, will have
the same representation in $R(\cG_{j-1})$ and $R(\cG_{j})$, i.e. for $E_{k_1k_2}\in \cE_{j-1}\cap \cE_{j}$
if $E_{k_1k_2}'$ is the representation in $R(\cG_{j-1})$ and $E_{k_1k_2}''$ is the representation in $R(\cG_{j})$
then $$E_{k_1k_2}'([t_{k_1},t_{k_2}])=E_{k_1k_2}''([t_{k_1},t_{k_2}]).$$

\end{enumerate}
}

The two following facts are true about consistent representations.

\prop{\label{prop:consis} \begin{enumerate}[(i)]
        \item If $(R(\cG_j))_{j=0}^n$ is a consistent representation of the TLG$^*$-tower $(\cG_j)_{j=0}^n$ then  $R(\cG_{j-1})\subset R(\cG_{j})$ for $j\geq 1$.
	\item If $(\cG_j)_{j=0}^n$ is a TLG$^*$-tower, for a fixed representation $R(\cG_n)$, there
exists a unique consistent representation  $(R(\cG_j))_{j=0}^n$ of this TLG$^*$ tower.
       \end{enumerate}
\dok{The claim (i) is clear from Definition \ref{def:consist_rep}. (ii) follows by induction on the number of edges.}  }

\defi{Let $\cG$ be a TLG$^*$ and fix its representation $R(G)$. 
By Definition \ref{def:tlg*}. of TLG$^*$'s there exists a TLG$^*$ tower
$(\cG_j)_{j=0}^n$, where $\cG_0$ is the minimal graph and
$\cG_n=\cG$. By Proposition \ref{prop:consis} there exists a consistent
representation $(R(\cG_j))_{j=0}^n$ where $R(\cG_n)=R(\cG)$.

It is easy to see that that $R(\cG_0)$ is the representation of
a full time-path $\sigma$ in $\cG$. We will call such a full
time-path a \textbf{spine}.\index{TLG$^*$ family!spine}}\index{Spine of a TLG$^*$|see{TLG$^*$ family}}\vspace{0.2cm} 

The question is each full time-path a spine? In other words, can we take any full time path, and by adding vertices
and edges as in the Definition \ref{def:tlg*} of TLG$^*$ get the TLG$^*$ $\cG$.

\teo{\label{thm:spine}Each full-time path\index{Time-like graph (TLG)!full-time path} in TLG$^*$ is a spine.
\dok{We will prove this claim by induction on the number of edges $m=|\cE|$ in $\cG$.

For $m=1$ the claim holds, since the spine is the whole $\cG$.\vspace{0.2cm}

Assume that the claim holds for $m\geq 1$. Let $\cG$ be a TLG$^*$
with $m+1$ edges. There exists a TLG$^*$ $\cG'$ such that by adding
a vertex or edge (as in step (\ref{def:tlg*:ii}) Definition \ref{def:tlg*}.)
we get $\cG$. (Note that in both cases $\cG'$ has $m$ edges.) \vspace{0.2cm}

If we added a new vertex to $\cG'$ there exists a representation of 
$R(\cG)$ that is the same as the one of $R(\cG')$. Now it is clear,
that if we pick any full time-path in $\cG$, there is a $\sigma'$ full
time-path in $\cG'$ with the same representation in $R(\cG')$.
We first construct $\cG'$, from $\sigma'$ and then we add $\cG$
to the tower describing that construction.

\begin{figure}[ht]
\begin{center}
\psfrag{1}{$\boldsymbol{\sigma^*}$} \psfrag{2}{$\boldsymbol{\sigma'}$}
\psfrag{E}{\small$ \boldsymbol{E^*_{h_1h_2}}$}
\includegraphics[width=12cm]{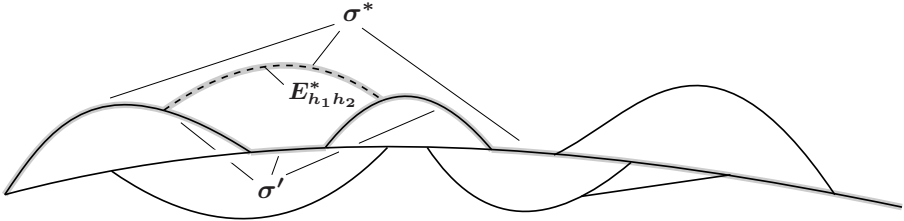}\\
\caption{Two spines $\sigma^*$ and $\sigma'$.} \label{pic12}

  \end{center}
\end{figure}

If we added a new edge $E^*_{h_1h_2}$. If we pick a full time-path $\sigma'$
that is in $\cG'$, then we first construct $\cG'$ from it and then add 
$\cG$ as the last member of the tower describing that construction.
If we pick a full time path $\sigma^*$ containing $E^*_{h_1h_2}$, let
$\sigma'$ be a full time path connecting $t_{h_1}$ and $t_{h_2}$,
such that $\sigma'$ and $\sigma^*$ are the same except between 
$t_{h_1}$ and $t_{h_2}$. We can construct $\cG'$ from $\sigma'$. To construct
$\cG$ from $\sigma^*$ we start with one edge representing $\sigma^*$,
and then add vertices $t_{h1}$ and $t_{h2}$ and an edge between them.
Now, we have a full time-path that has the same representation as 
$\sigma '$, and we keep adding edges and vertices in the same order as in the construction 
of $\cG'$ starting with $\sigma'$ (we skip the steps in which
$t_1$ and $t_2$ are added). At the end we have $\cG$. }}

We have an interesting consequence of the previous Theorem.

\pos{\label{pos:sped}If $E_{{h_1}{h_2}}$ is an edge between the two vertices connected by a 
time path (not containing that edge) in $\cG=(\cV,\cE)$, then $\cG'=(\cV,\cE\setminus\{E_{{h_1}{h_2}}\})$
is also a TLG$^*$.
\dok{We pick a full time-path containing that path. Now in the construction of 
$\cG$ from that time path we skip the step in which need to add the edge
$E_{{h_1}{h_2}}$ and we get $\cG'$.}}


\defi{\label{def:ptorder}A \textbf{point}\index{Time-like graph (TLG)!point} on $\cG=(\cV,\cE)$ is an element of the set 
$$\{(E_{jk},\tau)\ : \ E_{jk}\in \cE,  \tau\in [t_j,  t_k] \},$$
and the representation of the point $t=(E_{jk},\tau)$ is the 
point on $R(E_{jk})$ whose time coordinate is $\tau$. 
$t_1=(E^1,\tau_1)$ and $t_2=(E^2,\tau_2)$ are connected 
by a (time-)path if $E^1$ and $E^2$ are a part of some (time-)path. 
We will write $t_1\preceq t_2$ if  $\tau_1\leq \tau_2$ and $t_1$
and $t_2$ are connected by a time-path. }

\noindent \emph{Remark.} For a point $t$ on $\cG$ we will write $t\in G$. 
Note that vertices can be represented as several points, if they are endpoints 
to several edges, identify them as one point. The order '$\preceq $' introduced 
is the {\bf order induced by the structure of the graph}\index{Time-like graph (TLG)!order induced by a TLG $\preceq$}. We will write 
for the time of $t$, to simplify the notation, just $t$.\vspace{0.1cm}

We will give a criteria for connectedness of two points
by a time-path. This says that the two points are connected
by a time-path in $\cG$, if and only if 
their representations are connected from the moment that
these points exist in the TLG$^*$-tower (that leads to
the construction of $\cG$). A very similar result will hold 
for any path in $\cG$ with a given time frame.

\teo{\label{pre-con}Let $t_*$ and $t^*$ be two points on $\cG$ and let $(\cG_k)_{k=1}^n$
be a TLG$^*$-tower such that $\cG_n=\cG$ and $(R(\cG_k))_{k=1}^n$ its consistent representation. Assume $k_0$ is the smallest $k$
such that $\bar{t}_*$ and $\bar{t}^*$ are
on $R(\cG_k)$. Then $t_*$ and $t^*$ are connected by a time-path in
$\cG$ if and only if they are connected by a time-path in $\cG_{k_0}$.
\dok{If $t_*$ and $t^*$ are connected by a time-path in $\cG_{k_0}$, they
will remain connected by a time path in all $\cG_k$ for $k\geq k_0$.\vspace{0.2cm}

Let $k_*\geq k_0$ be the smallest $k$ such that $t_*$ and $t^*$
are connected in $\cG_k$. $k_*$ exists and is less or equal $n$. If $k_*>k_0$,
then $t_*$ and $t^*$ are points in $\cG_{k_*-1}$ but are not connected. This means that
an edge between two vertices $t_j$ and $t_h$ was added and $t_*$ and $t^*$ are
on some time-path. But since the points $t_j$ and $t_h$ need to be connected
in the previous step, this would not affect the connection between $t_*$ and $t^*$.
So $t_*$ and $t^*$ are connected in $\cG_{k_*-1}$. This contradicts the definition 
of $k_*$. Therefore, $k_0=k_*$.}}

From the last result we know that a simple cell will remain a simple cell
in the TLG$^*$-tower.

\pos{\label{pos:sm_cell}Let $(\cG_k)_{k=1}^n$ be a TLG$^*$-tower and $1\leq k<l\leq n$. 
If $(\sigma_1,\sigma_2)$ is a simple cell in $\cG_k$ then $(\sigma'_1,\sigma'_2)$
is a simple cell in $\cG_l$, where $(\sigma_1,\sigma_2)$ and $(\sigma'_1,\sigma'_2)$
have the same representation in the consistent representation of $(\cG_k)_{k=1}^n$.}

\defi{For any path $\rho$ in $\cG$ we say that the interval $I=[a,b]$
is its \textbf{time-frame} if $R(\rho)\subset I\times \R^2$.}

\teo{\label{pre-path-con}Let $t_*$ and $t^*$ be two points on $\cG$ and let $(\cG_k)_{k=1}^n$
be a TLG$^*$-tower such that $\cG_n=\cG$. Assume $k_0$ is the smallest $k$
such that $t_*$ and $t^*$ are points
on $\cG_k$. Then $t_*$ and $t^*$ are connected by a path $\rho$ within the time-frame $[a,b]$ in
$\cG$ if and only if they are connected by a path within the time-frame $[a,b]$
in $\cG_{k_0}$.
\dok{The proof is the same as in Theorem \ref{pre-con}. We look a the first 
member of the tower when $t_*$ and $t^*$ are connected by a path within the time frame $[a,b]$, if this
is not $k_0$, then the connection was established by adding an edge between 
some vertices $t_j$ and $t_k$, but these had to already be connected by a time-path. 
So the connection existed in the previous member of the tower. Which proves the claim.}
}

\section{Interval TLG$^*$'s }

In this section we will show the interval property of TLG$^*$'s.

\defi{Let $\cG$ be a TLG, and $\tau_1\leq \tau_2$ vertices on a TLG. We define 
$\cG[\tau_1,\tau_2]$ the \textbf{interval}\index{Time-like graph (TLG)!Interval TLG}\index{Interval TLG|see{Time-like graph (TLG)}} $[\tau_1, \tau_2]$ of $\cG$ to be
the graph $(\cV[\tau_1,\tau_2],\cE[\tau_1,\tau_2])$ such that
$\cV[\tau_1,\tau_2]$ are all the vertices $t_k$ such that there exist a time-paths
$\sigma_{\tau_1t_k}$ and $\sigma_{t_k\tau_2}$, and $\cE[\tau_1,\tau_2]$ are edges 
from $\cE$ that connect vertices from $\cV[\tau_1,\tau_2]$.}

\noindent\emph{Remark.} Note that if $\tau_1$ and $\tau_2$ are not connected by a time-path then $\cV[\tau_1,\tau_2]=\emptyset$.\vspace{0.3cm} 

The following result will show that interval TLG$^*$'s are TLG$^*$.

\begin{figure}[ht]
\begin{center}
\psfrag{a}{$\boldsymbol{\tau_1}$} \psfrag{b}{$\boldsymbol{\tau_2}$}
\psfrag{G}{$\boldsymbol{\cG[\tau_1,\tau_2]}$} 
\includegraphics[width=12cm]{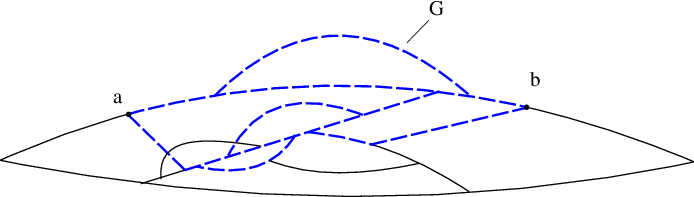}\\ 
\caption{The interval graph $\cG[\tau_1,\tau_2]$.} \label{pic35}

  \end{center}
\end{figure}

\teo{\label{teo:inttgl}Let $\cG$ be a TLG$^*$ and $\tau_1\leq \tau_2$ be two vertices connected 
by a time path. Then $\cG[\tau_1,\tau_2]$ is a TLG$^*$.
\dok{Pick a spine/full time-path $\sigma$ that contains $\tau_1$ and $\tau_2$. (It will exist since $\tau_1$
and $\tau_2$ are connected by a time-path.) Now, pick a TLG$^*$ tower $(\cG_j)_{j=1}^n$
that starts with $\sigma$ and ends with $\cG$.\vspace{0.2cm}

We will show by induction that $\cG_j[\tau_1,\tau_2]$ is a TLG$^*$ for all $j$. Without loss
of generality we can assume that $\cG_1$ contains all vertices on $\sigma$ in $\cG$.

It is clear that the claim holds for $j=1$. Assume it holds for $j\geq 1$, and let's prove it
for $j+1$. There are 4 cases to consider:
\begin{enumerate}[(1)]
\item If we added an vertex to $\cG_j$ to obtain $\cG_{j+1}$ the claim clearly holds.
 \item We added an edge that is not connecting vertices in $\cV_{j}[\tau_1,\tau_2]$. Then $\cG_{j+1}[\tau_1,\tau_2]$
is the same as $\cG_{j}[\tau_1,\tau_2]$. 
\item  We added an edge that is connecting vertices in $\cV_{j}[\tau_1,\tau_2]$, then these two vertices are connected
by a time-path in $\cG_j$, and hence they are connected by a time-path in $\cG_{j}[\tau_1,\tau_2]$. 
This is the same as if we added a new edge on $\cG_{j}[\tau_1,\tau_2]$ to obtain $\cG_{j+1}[\tau_1,\tau_2]$. 
\item We added an edge that is connecting a vertex in $\cV_{j}[\tau_1,\tau_2]$ and a vertex not in $\cV_{j}[\tau_1,\tau_2]$.
In this case $\cG_{j+1}[\tau_1,\tau_2]$
is the same as $\cG_{j}[\tau_1,\tau_2]$, because the vertex not in $\cV_{j}[\tau_1,\tau_2]$, 
by Theorem \ref{pre-con}, can't be in $\cV_{j+1}[\tau_1,\tau_2]$.  
\end{enumerate}
Since in all cases $\cG_{j+1}$ is either the same as $\cG_j[\tau_1,\tau_2]$, or obtained from $\cG_j[\tau_1,\tau_2]$
by adding and edge or a vertex, it is a TLG$^*$.

This proves that $\cG[\tau_1,\tau_2]$ is a TLG$^*$. }}

From this proof we can get the following conclusion.

\pos{When we erase the repeating elements the sequence $(\cG_j[\tau_1,\tau_2])_{j=1}^n$ is a TLG$^*$-tower
for $\cG[\tau_1,\tau_2]$.}

\pos{\label{cor:tlg*cnstint}For a TLG$^*$ $\cG$ and vertices $\tau_1$ and $\tau_2$ on a spine $\sigma$
we have that there exists a TLG$^*$-tower $(\cG_j)_{j=1}^n$ with consistent representation
$(R(\cG_j))_{j=1}^n$ such that for some $n_0\leq n$  
$$R(\cG_0)=R(\sigma), \quad R(\cG_{n_0})=R(\cG[\tau_1,\tau_2])\cup R(\sigma).$$
That is after the spine $\sigma$, we can construct $\cG[\tau_1,\tau_2]$, and then the rest of 
$\cG$.
\dok{We first construct the spine $\sigma$, and then construct TLG$^*$ $\cG[\tau_1,\tau_2]$.
Now, we apply steps from the proof of Theorem \ref{teo:inttgl}. that are using edges and vertices 
that haven't yet been constructed. In each of these steps when we add an edge time-path connectedness
is already guaranteed since the TLG$^*$ that we have is a sup-graph of the TLG$^*$ when the step was done
in the proof of  Theorem \ref{teo:inttgl}.} }

\section{Topology on TLG's}\index{Time-like graph (TLG)!topology|(}

For some things that follow we will need a notion of a limit 
of points on a TLG. In order to define a limit we need to define a 
topology. 

\defi{For a point $t$ on a TLG $\cG$, and $0<\delta <\min\{|t_k-t|:t_k\in \cV\setminus\{t\}\}$, 
we say that the \textbf{ball} $B_{\delta}(t)$ centered at $t$
with radius $\delta$ is the set of all points $s$ on a TLG, such that:
\begin{itemize}
 \item $t$ and $s$ are on a time-path;
 \item the absolute value of the time difference $|t-s|$ is less than $\delta$.
\end{itemize}}
\begin{figure}[ht]
\begin{center}
\begin{minipage}{7cm}
\psfrag{t}{\small$ \boldsymbol{t}$}
\psfrag{d}{\small$ \boldsymbol{\delta}$}
\includegraphics[width=7cm]{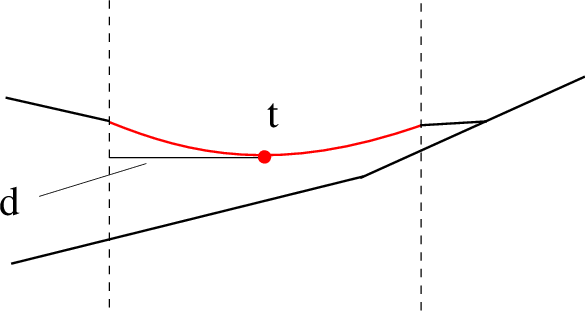}\\
\caption{Ball in a TLG} \label{pic38}
\end{minipage}\quad\quad 
\begin{minipage}{7cm}
\psfrag{t}{\small$ \boldsymbol{\bar{t}}$}
\psfrag{B}{\small$ \boldsymbol{\partial B(\bar{t},\delta_1)}$}
\psfrag{a}{\small$ \boldsymbol{\bar{t}_j}$}
\psfrag{b}{\small$ \boldsymbol{\bar{t}_k}$}
\psfrag{1}{\small$ \boldsymbol{\bar{i}_1}$}
\psfrag{2}{\small$ \boldsymbol{\bar{i}_2}$}
\psfrag{3}{\small$ \boldsymbol{\bar{i}_3}$}
\psfrag{4}{\small$ \boldsymbol{\bar{i}_4}$}
\psfrag{5}{\small$ \boldsymbol{\bar{i}_5}$}
\includegraphics[width=7cm]{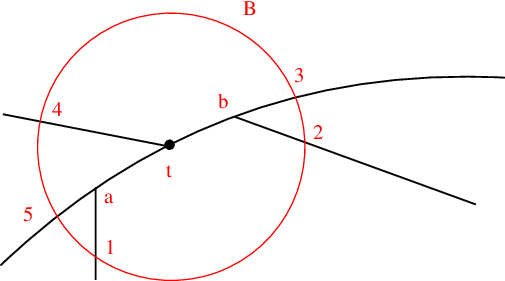}\\
\caption{The representation intersected by a sphere. In this case we have $I=\{i_1,\ldots, i_5\}$.} \label{pic40} 
\end{minipage}

\end{center}
\end{figure}


The following is a classical definition of open sets.
\defi{For a set $U$ of points on a TLG $\cG$ we say it is an \textbf{open set}, if 
for each 
$t\in U$ there exists a $\delta>0$ such that $B_{\delta}(t)\subset U$. 

We define ${\mathcal{T}}_{\G}$ to be the set of all open sets in TLG $\G$. 
}

\lem{Let $\cG$ be a TLG, and fix its representation $R(\cG)$. $U$ is an open
set in if and only if $R(U)$ is an open set in $R(\cG)$.
\dok{If $U\in \cT_{\cG}$, then pick arbitrary $\bar{t}\in R(U)$. There are only finitely many paths 
that don't pass through $t$, and the union of their representations is a compact set $K$ in $\R^3$.
Now, we pick $\delta_1=d(\bar{t},K)/2$ (where $d$ is the usual metric in $\R^3$). 
Also, we pick $\delta_2>0$ such that 
$B_{\delta_2}(t)\subset U$. For $\delta=\min\{\delta_1,\delta_2\}$, $\{s\in R(\cG):d(s,\bar{t})<\delta\}\subset R(U)$.
Hence $R(U)$ is opened.


If $R(U)$ is opened, then we pick $t\in U$. Pick $\delta_1>0$ such that $B(\bar{t},\delta_1)=\{s\in R(G): d(s,\bar{t})<\delta_1\}\subset R(U)$. 
There exists finitely many full time-paths 
$\pi_1,\ldots , \pi_k$ that contain $t$. Let $I$ be the points on $\cG$
whose representations are at the intersection 
of $R(\pi_1),\ \ldots\ , R(\pi_k)$  with $\partial B(\bar{t},\delta_1)$. (See Figure \ref{pic40}.) $I$ is finite, and now pick
$\delta = \min\{|t-z|:z\in I\cup (\cV\setminus \{t\})\}/2$. Since $t\notin I$, $\delta>0$.
Hence, $B_{\delta}(t)\subset U$.}}

\prop{$\cT_{\cG}$ is a topology on $\cG$.
\dok{Note that $t\mapsto \bar{t}$ is a bijection. Hence, if $(U_\alpha:\alpha\in A)$
is in $\cT_{\cG}$, then since
$$R(\bigcup_{\alpha\in A}U_{\alpha})=\bigcup_{\alpha\in A}R(U_{\alpha})$$
is an open set so is $\bigcup_{\alpha\in A}U_{\alpha}$. We can use the same 
approach for the finite intersection. }}

\pos{$t\mapsto \bar{t}$ is a homeomorphism (i.e. a continuous 
bijective function with a continuous inverse) from $\cG$ to $R(\cG)$.}

\pos{\label{pos:TGmtr}The topological space $(\cG,\cT_\cG)$ is metrizable\index{Time-like graph (TLG)!topology!metrizability}.
\dok{Fix the representation $R(\cG)$, and set $d_{\cG}(s,t):=d_{\R^3}(\bar{s},\bar{t})$.
$d_{\cG}$ is a metric and the topology induced by $d_{\cG}$ is $\cT_\cG$.}}

\pos{$(\cG,\cT_{\cG})$ is a Hausdorff space.
\dok{Follows from the fact that this space is metrizable.
%
%
%
%
%
%
}}

We define limit on TLG's in the following natural way. 

\defi{We say that the sequence of points $(t_n)$ converges\index{Time-like graph (TLG)!topology!Convergence of points} to the point $t$
in TLG if: 
\begin{itemize}
 \item there exists $n_0\in \N$ such that for each $n\geq n_0$ 
the points $t_n$ and $t$ are connected by a time-path;
 \item the absolute value of the time difference $|t_n-t|$
converges to $0$.
\end{itemize}}

\noindent\emph{Remark.} The time-path that connects $t_n$ and $t$ can depend on $n$.
and can be a different time-path for different $n$'s. (It will always contain $t$.)

We will show that this is also the limit in the topology that we defined.

\teo{\label{thm:convTLG}Let $\cG$ be a TLG, and $R(\cG)$ be its representation. A sequence of points
$(t_n)$ converges to $t$ in $\cG$ if and only if their representations $(\bar{t}_n)$ converge to $\bar{t}$.
\dok{If $t_n\to t$ in $\cG$. There are finitely many paths $\sigma_1,\ldots, \sigma_k$ 
going through $t$. In the representation each path $\sigma_j$ is represented by a graph 
of some continuous function $f_{\sigma_j}$. But now since 
$$(t_n,f_{\sigma_j}(t_n))\to (t,f_{\sigma_j}(t))=\bar{t},$$
and for each $\bar{t}_n$ there is $k_n$ such that $\bar{t}_n=(t_n,f_{\sigma_{k_n}}(t_n))$, the claim follows.\vspace{0.2cm}

Let $\bar{t}_n\to \bar{t}$ in $R(\cG)$. Now, there are only finitely many paths 
that don't pass through $t$, and the union of their representations is a compact set $K$.
Now we pick $\delta=d(\bar{t},K)/2$. Now, there exists $n_0$ such that for all 
$n\geq n_0$ $\bar{t}_n\in B_{\delta}(\bar{t})$, but this implies that all $t_n$
are connected by a time-path to $t$. It is clear that the absolute value of the time difference $|t-t_n|$
converges to $0$.    } }

\pos{$t_n\to t$ in $\cG$ if and only if $t_n\to t$ in $(\cG,\cT_{\cG})$. 
\dok{Fix a representation $R(\cG)$, and define a metric $d_{\cG}$ as in 
Corollary \ref{pos:TGmtr}. It is clear from Theorem \ref{thm:convTLG} that 
we have convergence if and only if $d_{\cG}(t_n,t)=d_{\R^3}(\bar{t}_n,\bar{t})\to 0$.}}\index{Time-like graph (TLG)!topology|)}

\section{TLG$^*$ as a topological lattice}\label{tlg*tltt}
In this section we will show that TLG$^*$'s are topological lattices.

\defi{\label{def:tplt}A Hausdorff space $X$ with some order '$\leq$' is called a \textbf{topological
lattice}\index{Topological lattice}\index{Topological lattice|seealso{TLG$^*$ family}} if for $x_1,x_2\in X$:
\begin{itemize}
 \item  there exists a unique element $x_1\wedge x_2$ such that
$$\{x\in X:x\leq x_1 \}\cap \{x\in X:x\leq x_2\} = \{x\in X:x\leq x_1\wedge x_2\};$$
\item  there exists a unique element $x_1\vee x_2$ such that
$$\{x\in X:x\geq x_1 \}\cap \{x\in X:x\geq x_2\} = \{x\in X:x\geq x_1\vee x_2\}.$$
\end{itemize}
and $x_1\wedge x_2$ and $x_1\vee x_2$ are continuous mappings of $X\times X$ (with product topology) onto $X$.}

\teo{\label{teo:tpltt}A TLG$^*$ $\cG$ is a topological lattice\index{TLG$^*$ family!topological lattice} with respect to the order $\preceq$ induced 
by the structure of $\cG$.
\dok{Let $(\cG_k)_{k=0}^n$ be a TLG$^*$-tower starting with 
the minimal graph $\cG_0$ and ending with $\cG_n=\cG$.

We will prove the claim by induction. Clearly, $\cG_0$ is a topological lattice.
Let's assume $\cG_k$ is a topological lattice. 

If we added a new vertex to $\cG_{k}$
in order to get $\cG_{k+1}$, then clearly $\cG_{k+1}$ is also a topological lattice.

If we added a new edge to $\cG_{k}$
in order to get $\cG_{k+1}$, then take two points $t,s\in \cG_{k+1}$. If $t,s\in \cG_k$,
then by assumption there exist $t\wedge s$ and $t\vee s$, the same is clear if 
$t,s$ are points of the new edge $E_{jk}^*$. The only case that remains to be checked is
when $t\in E_{jk}^*$ and $s\in \cG_k$. If $t\preceq s$, then $t\wedge s=t$ and $t\vee s=s$.
Similarly when  $s\preceq t$. Otherwise, we have $\{\tau \in \cG_{k+1}:\tau \preceq s\}$ 
is in $\cG_{k}$, so
\begin{align*}
 &\{\tau \in \cG_{k+1}:\tau \preceq s\}\cap \{\tau \in \cG_{k+1}:\tau \preceq t\}\\
=&\{\tau \in \cG_{k}:\tau \preceq s\}\cap \{\tau \in \cG_{k}:\tau \preceq t\}\\
=&\{\tau \in \cG_{k}:\tau \preceq s\}\cap \{\tau \in \cG_{k}:\tau \preceq t_j\}\\
=&\{\tau \in \cG_{k}:\tau \preceq s\wedge t_j\},
\end{align*}
therefore, we have  $s\wedge t=s\wedge t_j$. In the same way we can show 
that $s\vee t=s\vee t_k$. 
The uniqueness follows from the fact that if $u\preceq v$ and $v\preceq u$ we have $u=v$.

Let $(t^1_n)$ and $(t^2_n)$ be a sequence of points converging respectively to 
$t_1$ and $t_2$ on $\cG$. If $t_1=t_2$ both sequences converge to the same 
point, and so will $(t^1_n \vee t^2_n)$ and $(t^1_n \wedge t^2_n)$. If 
$t_1$ and $t_2$ are on the same time-path, assume $t_1\prec t_2$. Now,
by the definition of convergence, there will exist a $n_0$ such that 
for $n\geq n_0$ we have $t^1_n\prec t^2_n$, hence 
$$t^1_n\vee t^2_n=t^1_n\to t_1, \quad t^1_n\wedge t^2_n=t^2_n\to t_2 $$

If $t_1$ and $t_2$ are not connected by a time-path, let $\delta < \min\{|t_1-t|/2:t\in \cV\setminus\{t_1\}\}\wedge \min\{|t_2-t|/2:t\in \cV\setminus\{t_2\}\}$,
it is not hard to see that for $t'\in B_{\delta}(t_1)$ and $t''\in B_{\delta}(t_2)$,
we have $t'\vee t''=t_1\vee t_2$ and $t'\wedge t''=t_1\wedge t_2$. So for large $n$,
the sequences will have the values $t_1\vee t_2$ and $t_1\wedge t_2$. }}

It is not hard to see, that the TLG that is not a TLG$^*$ from Figure \ref{pic1} is not a topological
lattice -- there is no unique $t_1\vee t_2$ and $t_3\wedge t_4$. 

\lem{There exists a topological lattice TLG, that is not a TLG$^*$.\index{Topological lattice}
\dok{The TLG in the Figure \ref{pic33} is an example of a topological lattice TLG, that is not a 
TLG$^*$.

\begin{figure}[ht]
\begin{center}
\psfrag{0}{$\boldsymbol{t_0}$} 
\psfrag{1}{$\boldsymbol{t_1}$}
\psfrag{2}{$\boldsymbol{t_2}$}  
\psfrag{3}{$\boldsymbol{t_3}$}  
\psfrag{4}{$\boldsymbol{t_4}$}  
\psfrag{5}{$\boldsymbol{t_5}$}  
\psfrag{6}{$\boldsymbol{t_6}$}  
\psfrag{7}{$\boldsymbol{t_7}$}  
\psfrag{8}{$\boldsymbol{t_8}$}  
\psfrag{9}{$\boldsymbol{t_9}$}  
\includegraphics[width=12cm]{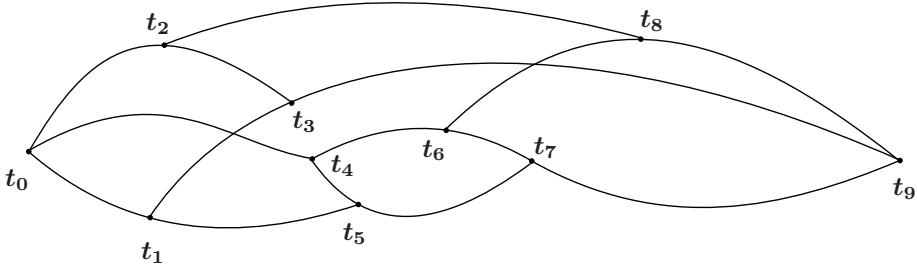}\\
\caption{Topological lattice TLG that is not a TLG$^*$.} \label{pic33}
\end{center}
\end{figure}

It is easy to see that $t_0\wedge t_j=t_0$ and $t_0\vee t_j=t_j$, and similarly
$t_9\wedge t_j=t_j$ and $t_9\vee t_j=t_9$. The following table will show what $t_k\wedge t_j$ and $t_k \vee t_j$ are.
In the table above the main diagonal (for $k< j$) $t_k\wedge t_j$ is calculated, and
below (for $k> j$) $t_k\vee t_j$. The diagonal is omitted, since $t_j\vee t_j=t_j\wedge t_j=t_j$.

$$\begin{array}{c||c|c|c|c|c|c|c|c}
    & t_1 &t_2 &t_3 &t_4 &t_5&t_6&t_7 &t_8\\ \hline \hline
  t_1 & \circ &t_0 &t_1 &t_0 &t_1&t_0&t_1 &t_0\\
  t_2  & t_3 &\circ &t_2 &t_0 &t_0&t_0&t_0 &t_2\\
  t_3& t_3 &t_3 &\circ &t_0 &t_1&t_0&t_1 &t_2\\
   t_4 & t_5 &t_8 &t_9 &\circ &t_4&t_4&t_4 &t_4\\
   t_5 & t_5 &t_9 &t_9 &t_5 &\circ&t_4&t_5 &t_4\\
  t_6 & t_7 &t_8 &t_9 &t_6 &t_7&\circ &t_6 &t_6\\
t_7 & t_7 &t_9 &t_9 &t_7 &t_7&t_7&\circ &t_6\\
t_8 & t_9 &t_8 &t_9 &t_8 &t_9&t_8&t_9 &\circ
  \end{array}
$$
This is not a TLG$^*$, since by applying the cell collapse transformation,
see Definition \ref{defi:sct} on the cell $(t_4-t_5-t_7,t_4-t_6-t_7)$ we
will no longer have a topological lattice, since $t_3\wedge t_8$ will no longer be unique.
Therefore the transformed graph is no longer a TLG$^*$ which contradicts Lemma \ref{lem:trtlt}.
if this is a TLG$^*$.
}}

A natural question that will be useful later is if we have a sequence of points
$(t_k)$ does there exist their infinitum and supreme, that is
$$\wedge_{k=1}^\infty t_k \quad \textrm{and}\quad \vee_{k=1}^\infty t_k.$$

\lem{\label{lem:kon:min}The order in which we take apply $\wedge$ and $\vee$ doesn't matter, that is
$$(t_1\wedge t_2)\wedge t_3 = t_1\wedge (t_2\wedge t_3)\quad \textrm{and} \quad (t_1\vee t_2)\vee t_3 = t_1\vee (t_2\vee t_3).$$
\dok{Let $t_*=(t_1\wedge t_2)\wedge t_3$ and $t^*=t_1\wedge (t_2\wedge t_3)$. 
It is clear that $t_* \preceq t_3$, and  $t_* \preceq t_1\wedge t_2$ implies 
$t_*\preceq t_2$ and $t_*\preceq t_1$. By definition it is clear that $t_*\preceq (t_2\wedge t_3)$, again using the same 
we have $t_*\preceq t_1\wedge (t_2\wedge t_3)=t^*$. In the same way, we can get $t^*\preceq t_*$,
and this implies $t_*=t^*$. Hence, the first equality follows. The second equality 
follows by similar arguments. These equalities imply the other statements. }}

\lem{\label{lem:bes:min}Let $(t_k)_{k=1}^{\infty}$ be a sequence of points in a TLG$^*$. We define the sequences
$(t_k^-)$ and $(t_k^+)$ by $t_1^-=t_1$, and $t_k^-=t_k\wedge t_{k-1}^-$, and $t_1^+=t_1$,
and $t_k^+=t_k\wedge t_{k-1}^+$. Sequence $(t_k^-)$ and $(t_k^+)$ will converge to limits 
$t_*$ and $t^*$. Further for any bijection $f:\N\to \N$ the sequences $(t_k^{f-})$
and $(t_k^{f+})$ obtained from $(t_{f(k)})$ in the same way will converge respectively to $t_*$
and $t^*$. 
\dok{By definition, for each $n$ the points $(t_k^-)_{k=1}^n$ there exists a full time-path 
$\sigma$, such that these points are all on $\sigma$. Further, the sequence of times $(t_k^-)$
converges to a time $t_*$. On the TLG$^*$ $\cG$ there are only finitely many points 
with that time, name them $t_{1*}$, \ldots, $t_{m*}$. Let $\varepsilon = \min\{|t_*-t_{j*}\vee t_{k*}| : k \neq j \}$
where the minimum is taken over the time distances. Now, if we pick $k_0$ such that
$|t_*-t_k^-|<\varepsilon$ (time distance) for $k\geq k_0$, then there will be only 
one $t_{j*}$ in the future of $t_k^-$'s for $k\geq k_0$. We set it to be $t_*$, and it is not hard to see 
that all the points are on the unique path between $t_*$ and $t_{k_0}^-$. Now it is 
clear, since the topology on that path is the same as the one on the open segment, that 
$t_k^-\to t_*$.

By what we have just proven $(t_k^{f-})$ converges to some point $t_*^f$. But then, we can show by definition, that
$t_*^{f}\preceq t_*$ and $t_*\preceq t_*^{f}$, which implies $t_*^{f}=t_*$.}}

\defi{For a finite sequence $(t_k)_{k=1}^n$ we define 
$$\wedge_{k=1}^n t_k:=t_{p(1)}\wedge (t_{p(2)} \wedge (\ldots (t_{p(n-1)} \wedge t_{p(n)}) )) \quad \textrm{and}\quad  \vee_{k=1}^n t_k=t_{p(1)}\vee (t_{p(2)} \vee (\ldots (t_{p(n-1)} \vee t_{p(n)}) )) .$$
where $p$ is any permutation of the set $\{1,2,\ldots, n\}$. For a sequence $(t_k)_{k=1}^\infty$ we define 
and any bijection $f:\N\to \N$ we define
$$\wedge_{k=1}^\infty t_k:=\lim_{n\to\infty}\wedge_{k=1}^n t_{f(k)}  \quad \textrm{and}\quad \vee_{k=1}^\infty t_k:=\lim_{n\to\infty}\vee_{k=1}^n t_{f(k)}.$$ }

\pos{The terms $\wedge_{k=1}^n t_k$, $\vee_{k=1}^n t_k$,  $\wedge_{k=1}^\infty t_k$ and $\vee_{k=1}^\infty t_k$ 
are well defined for any sequence $(t_k)$.
\dok{Follows from Lemma \ref{lem:kon:min}. and Lemma \ref{lem:bes:min}.}}

\section{Cell collapse transformation and the stingy algorithm}
Another property of TLG$^*$ will be
introduced in this section. This will be a transformation on TLG's that will map TLG$^*$'s
into TLG$^*$'s.

\defi{\label{defi:sct}We will call the map $\cG\mapsto \cG^{\circ}$ from TLG's into TLG's a \textbf{cell 
collapse transformation}\index{Time-like graph (TLG)!cell collapse transformation} \index{Cell collapse transformation|see{Time-like graph (TLG)}} if:\vspace{0.2cm}

Pick a cell $(\sigma_{uv},\sigma^1_{uv})$ (starting at $t_u$ and ending at $t_v$).
The transformation that we will describe,
basically, glues $\sigma^1_{uv}$ with its vertices to $\sigma_{uv}$, while keeping most
 of the connections between vertices in the graph.\vspace{0.2cm} 

We construct the graph $\cG^{\circ}=(\cV^{\circ},\cE^{\circ})$ in the following way:
\begin{itemize}
 \item In the first step we are maping the cell into a time-path. 

Let $t_u=t_{w_1}\leq \ldots\leq  t_{w_h}=t_v$ be the set of vertices on
the time-paths $\sigma_{uv}$ and $\sigma_{uv}^1$ ordered with respect to time.
We will map $t_{w_j}$ into $(t_{w_j})^{\circ}$ in $\cV^{\circ}$ so that the vertices
with the same time are mapped into same vertices, that is if $t_{w_{j_1}}=t_{w_{j_2}}$ 
then $(t_{w_{j_1}})^{\circ}\equiv (t_{w_{j_2}})^{\circ} $. 

We will use the notation $(t_{w_{j_1}})^{\circ}=t_{w_{j_1}^{\circ}}^{\circ}$. Note that
if $t_{w_{j_1}}=t_{w_{j_2}}$, then $w_{j_1}^{\circ}=w_{j_2}^{\circ}$.

We add an edge in $\cE^{\circ}$ between $t^{\circ}_{w_j^{\circ}}$ and $t^{\circ}_{w_{j+1}^{\circ}}$ if their times are different.
(Note that in this way all the vertices in $\{(t_{w_j})^\circ:j=1,2,\ldots n\}$ are on the same time-path.)
\item Every other vertex $t_j$ from $\cV$ not contained on the paths $\sigma_{uv}$ and $\sigma^1_{uv}$ 
is mapped into $(t_j)^{\circ}$ in $\cV^{\circ}$ so that the time is preserved, and these vertices are 
mapped into different vertices and disjoint from where the vertices on $\sigma_{uv}$
and $\sigma^1_{uv}$ were maped.
\item For each edge $E$ in $\cE$ not a part of $\sigma_{uv}$ or $\sigma^1_{uv}$ we add a 
$E^{\circ}$ in $\cE^{\circ}$ between the corresponding vertices. We color $E^{\circ}$ in red if $E$ is adjacent to
a vertex from $\sigma_{uv}$, or in blue if it is adjacent to the vertex from $\sigma^1_{uv}$. 
\end{itemize}
}

\begin{figure}[ht]
\begin{center}
\psfrag{1}{$\boldsymbol{\sigma^1_{uv}}$} 
\psfrag{2}{$\boldsymbol{\sigma_{uv}}$}
\includegraphics[width=12cm]{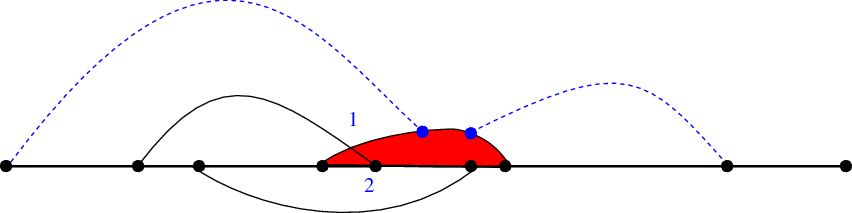}\\ 
\includegraphics[width=12cm]{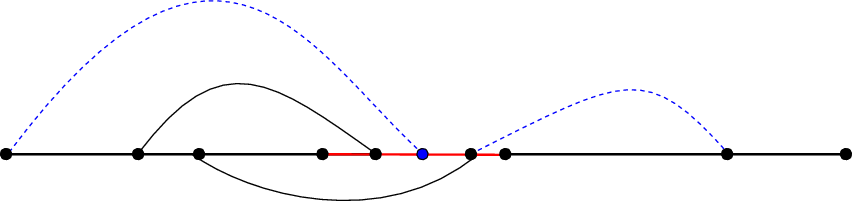}\\
\caption{Transformation from $\cG$ to $\cG^\circ$.} \label{pic31}

  \end{center}
\end{figure}

\defi{\label{defi:msct} A cell $(\sigma_1,\sigma_2)$ in TLG $\cG$ starting at $t_{k_1}$ and ending at 
$t_{k_2}$ will be called \textbf{truly simple}\index{Time-like graph (TLG)!truly simple cell} \index{Truly simple cell|see{Time-like graph (TLG)}}, if there is no path in $\cG[t_{k_1},t_{k_2}]$ connecting
the interior of $\sigma_1$ and $\sigma_2$. }

\noindent \emph{Remark.} The path in question does not have to be a time-path. If there exists a
a time path between $\sigma_1$ and $\sigma_2$ then it will be in $\cG[t_{k_1},t_{k_2}]$, so a truly simple cell 
is a simple cell.

Before we prove the main result of this section we will prove the following lemma.

\lem{\label{lem:dcmptlg*}\begin{enumerate}[(a)]
      \item Let $(\sigma_1,\sigma_2)$ starting at $t_{k_1}$ and ending at 
$t_{k_2}$ in in TLG$^*$ $\cG$ be a truly simple cell. Then 
$$R(\cG[t_{k_1},t_{k_2}])\setminus \{\bar{t}_{k_1},\bar{t}_{k_2}\}$$
has at least two connected components.
\item Let $\cG$ be a TLG$^*$, and let $R(\cG)\setminus \{\bar{t}_0,\bar{t}_N\}$ have 
two connected components. Closure  of each of these components, is a representation of a TLG$^*$.
     \end{enumerate}
\dok{(a) Since $(\sigma_1,\sigma_2)$ is a truly simple cell, there is no path between
the interior vertices of $\sigma_1$ and $\sigma_2$. So $R(\sigma_1)$ and $R(\sigma_2)$ are 
connected only through $\bar{t}_{k_1}$ and $\bar{t}_{k_2}$. Therefore, $R(\sigma_1)\setminus \{\bar{t}_{k_1},\bar{t}_{k_2}\} $
and $R(\sigma_2)\setminus \{\bar{t}_{k_1},\bar{t}_{k_2}\} $ are in two different connected components
of $R(\cG[t_{k_1},t_{k_2}])\setminus \{\bar{t}_{k_1},\bar{t}_{k_2}\}$. \vspace{0.2cm}

\begin{figure}[ht]
\begin{center}
\psfrag{a}{$\boldsymbol{0}$} \psfrag{b}{$\boldsymbol{1}$}
\psfrag{1}{$\textcolor{blue}{\boldsymbol{\cH}}$} \psfrag{2}{$\textcolor{red}{\boldsymbol{\cK}}$} 
\includegraphics[width=9cm]{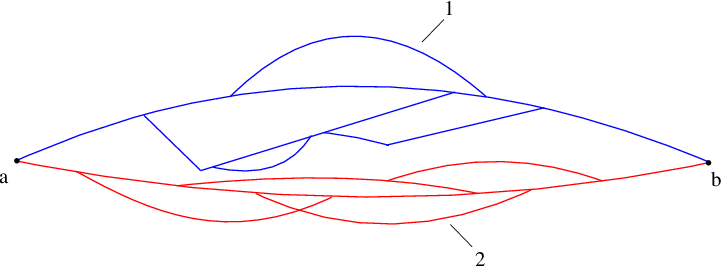}\\ 
\caption{$\cH$ and $\cK$ are TLG$^*$'s.} \label{pic46}

  \end{center}
\end{figure}

(b) Pick a component, and let $\cH$ be the sub-graph of $\cG$ that represents this
component and the union of $\{\bar{t}_1,\bar{t}_N\}$. Pick a TLG$^*$-tower $(\cG_j)_{j=0}^n$ that starts with a minimal edge
and ends with $\cG$. Let $(\cG_{j_k})_{k=1}^{n_1}$ be the subsequence of all members of 
$(\cG_j)_{j=0}^n$ such that an edge or a vertex whose representation intersects $R(\cH)\setminus \{\bar{t}_1,\bar{t}_N\}$ has been added to $\cG_{j_k-1}$ to
obtain $\cG_{j_k}$.\vspace{0.1cm}

By the definition of the sequence $(\cG_{j_k})$, an edge has been added to $\cG_{j_1-1}$
in order to obtain $\cG$. Since the representation of that edge intersects $R(\cH)\setminus \{\bar{t}_1,\bar{t}_N\}$
which is a disconnected component of $R(\cG)\setminus \{\bar{t}_1,\bar{t}_N\}$. So therefore that
edge needs to be between $t_1$ and $t_N$. Set $$\cH_1:=(\cV_{j_1}\cap \cV_{\cH}, \cE_{j_1}|_{\cV_{j_1}\cap \cV_{\cH}}),$$
where $\tilde{\cE}|_{\tilde{\cV}}$ represents the subset of edges in $\tilde{\cE}$ that are connecting
vertices in $\tilde{\cV}$. It is clear that $\cH_1$ is a minimal graph.\vspace{0.1cm} 

Further, define $\cH_k=(\cV_{j_k}\cap \cV_{\cH}, \cE_{j_k}|_{\cV_{j_k}\cap \cV_{\cH}})$ for $k=2,\ldots, n_1$. 
We will show that $(\cH_k)_{k=1}^{n_1}$ is a TLG$^*$-tower. $\cH_{n_1}$ by construction 
equals $\cH$. $\cH_1$ is a TLG$^*$. Let's assume $\cH_{k}$ is a TLG$^*$ (for $k
\geq 1$) and show that $\cH_{k+1}$ is a TLG$^*$. If a new vertex has been added to $\cG_{j_{k+1}-1}$
 to obtain $\cG_{j_{k+1}}$, this is, by construction, the same as if we added a new vertex
to $\cH_k$ in order to obtain $\cH_{k+1}$. If we added a new edge, the representation of that edge 
intersects $R(\cH)\setminus \{\bar{t}_1,\bar{t}_N\}$, and therefore is in that component.
Since the new edge is connecting two vertices connected by a time-path in $R(\cH)\cap R(\cG_{j_{k+1}-1})$
these vertices are in $\cH$, and they are connected in $\cH_{k}$. Hence, we added an edge 
to $\cH_k$ between two vertices connected by a time-path. In both cases $\cH_{k+1}$ is a TLG$^*$
obtained from $\cH_{k}$. Hence, $\cH$ is a TLG$^*$.  }}

\teo{\label{lem:trtlt}If $\cG$ is a TLG$^*$ and $\circ$ is collapsing a truly simple cell, then $\cG^\circ$ is also a TLG$^*$.
Further, if $\circ$ is collapsing a simple cell, then $\cG^\circ$ doesn't have to be a TLG$^*$.
\dok{Pick a spine $\sigma$ that contains $\sigma_{uv}$ side of the chosen cell
$(\sigma_{uv},\sigma^1_{uv})$. We know from Theorem \ref{teo:inttgl} that 
$\cG[t_u,t_v]$ is a TLG$^*$. By Lemma \ref{lem:dcmptlg*} $\cG[t_u,t_v]$ is a union 
of two or more TLG$^*$'s that only have vertices $t_u$ and $t_v$ in common, and
are otherwise disconnected. $\sigma_{uv}$ and $\sigma_{uv}^1$ are in two different
TLG$^*$ components and they are also spines in these components. Therefore, we can
construct $\cG[t_u,t_v]$, from the TLG$^*$ whose representation is $R(\sigma_{uv})\cup R(\sigma^1_{uv})$.
We start with the two edges that have the representation $R(\sigma_{uv})$ and $R(\sigma^1_{uv})$, and then we first construct the component
that contains $\sigma_{uv}$, then the one that contains $\sigma^1_{uv}$, and possible other components.
At the end we get $\cG[t_u,t_v]$. But then, for any full-time path $\sigma$ 
that contains $t_u$ and $t_v$ we  can construct the TLG$^*$ whose representation is
$R(\sigma)\cup R(\cG[t_u,t_v])$ starting with the TLG$^*$
$$\cG_1=(\{t_0,t_u,t_v,t_N\},\{E_{0u},E^1_{uv},E^2_{uv}, E_{vN}\}),$$
and later, by Corollary \ref{cor:tlg*cnstint}, we can construct 
$\cG$. Hence, there exists a TLG$^*$-tower $(\cG_{j})_{j=1}^n$ such that 
ends with $\cG$, and its consistent representation has the representation of 
the cell $(\sigma_{uv},\sigma^1_{uv})$ at each level. Now, we define 
$\circ$-transformation to collapse the cell whose representation is $R(\sigma_{uv},\sigma^1_{uv})$.
We will show that $(\cG_j^{\circ})_{j=1}^n$ is a TLG$^*$-tower.\vspace{0.2cm} 

It is clear that $\cG_1^{\circ}$ is a TLG$^*$ and that images of all
points connected by a time path in $\cG_1$ are connected in $\cG_1^{\circ}$. Let's assume $\cG_k^{\circ}$ 
is a TLG$^*$ and that images of all
points connected by a time path in $\cG_k$ are connected in $\cG_k^{\circ}$. 

If we added a vertex to $\cG_k$ in order to obtain $\cG_{k+1}$,
then $\cG_{k+1}^{\circ}$ is either the same as $\cG_k$ or it has an added vertex. It is clear
in this case that images of all the points that are connected in $\cG_{k+1}$ by a time-path
are connected by a time-path in $\cG_{k+1}^{\circ}$.

If we added an edge to $\cG_k$ in order to obtain $\cG_{k+1}$,
then $\cG_{k+1}^{\circ}$ is the same as if added an edge to $\cG_{k}^{\circ}$. Since this edge
is connecting image of two points  in $\cG_k$ that are connected by a time-path,
they are also connected by a time-path in $\cG_k^{\circ}$. Hence, $\cG_{k+1}^{\circ}$
is also a TLG$^*$. Images of all the time-path connected points in $\cG_{k+1}$ that are not on the 
edge added, are connected by a time path in $\cG_{k+1}^{\circ}$. (This is inherited from $\cG_k$.)
The points on the edge are connected through the endpoints, and since the image of the edge
is connected through the image of the endpoints, the connectedness follows.\vspace{0.1cm}

Hence $(\cG_k^{\circ})$ is a TLG$^*$-tower ending with $\cG^\circ$.

\vspace{0.2cm}
  
The example when we collapse a simple cell in a TLG$^*$ an don't obtain a TLG$^*$ is given on Figure \ref{pic34}.
The second figure is not a topological lattice, so it is not a TLG$^*$.
\begin{figure}[ht]
\begin{center}
\includegraphics[width=6cm]{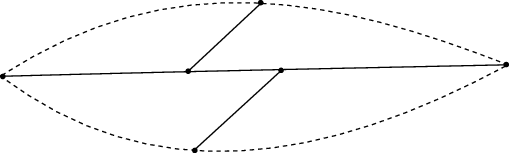}\ \includegraphics[width=6cm]{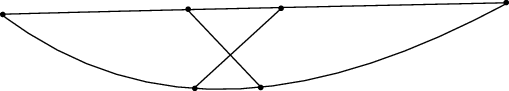}\\
\caption{Planar TLG$^{(*)}$, when we collapse the dashed (simple) cell we no longer have a TLG$^*$.} \label{pic34}
\end{center}
\end{figure}

}}

We obtain following from the previous proof. 

\pos{\label{circ-tlgt}Let $(\cG_j)_{j=1}^n$ be a TLG$^*$-tower. If there exists a truly simple cell,
in $\cG_n$ such that the representation of this cell is truly simple in each member of the tower,
then for $\circ$ the cell collapsing transformation of this cell $(\cG_j^{\circ})_{j=1}^n$ is also a TLG$^*$-tower.
\dok{We first construct a pre-tower $(\cG_{j})_{j=-m}^1$, where $\cG_{-m}$ 
is the minimal TLG$^*$ containing the cell. Now, $(\cG_{j})_{j=-m}^n$ is a TLG$^*$-tower, and in the same way as in the previous proof
we can show that $(\cG_{j}^{\circ})_{j=-m}^n$ is a tower, and the claim follows. }}

%

In what follows we will define an algorithm which will
give us the criteria for understanding is something a
TLG$^*$ or not.
\defi{For a TLG $\cG$ and a full-time path $\sigma$ in $\cG$, the following algorithm \index{TLG$^*$ family!stingy algorithm}\index{Stingy algorithm|see{TLG$^*$ family}}
will be called the \textbf{stingy algorithm for the TLG $\cG$ with respect to the full-time path $\sigma$}:
\begin{center}
\begin{algorithm}[H]
 $\sigma$ a full time-path (i.e. in $P_{0\to 1}(\cG)$)\;
$\cG^\#=(\cV^\#,\cE^\#)$ a TLG that consists of all vertices and all edges of $\sigma$ (in $\cG$)\;\label{alg2:s2}

\While{$\cE\setminus \cE^\#\neq \emptyset$\label{alg2:lwhl}}{
      $\sigma_{kl}$ a time-path in $\cG$ and not in $\cG^{\#}$ between $t_k$ and $t_l$ in $\cV^{\#}$ such that $|t_l-t_k|$ is minimal\;\label{alg2:s6}
      add all edges and vertices that make $\sigma_{kl}$ (in $\cG$) to $\cG^\#$\; \label{alg2:lx}
} 
\caption{Stingy algorithm for constructing $\cG$ with respect to $\sigma$}\label{alg:tpltt}
\end{algorithm}
\end{center}}
\lem{The stingy algorithm for any $\cG$ and any full-time path $\sigma$ in $\cG$
terminates in finitely many steps.
\dok{It is clear that as long the condition $\cE\setminus \cE^\#\neq \emptyset$ is satisfied, we
can find $\sigma_{kl}$ as in the line \ref{alg2:s6} (it may not be unique, but it will exist).
Since in each {\bf while} loop execution we add at least one edge, eventually we will have $\cE=\cE^{\#}$. Clearly,
at that point we have $\cG=\cG^{\#}$.}}

The key to answering is $\cG$ a TLG$^*$ is in line \ref{alg2:s6}. We claim that
if $\cG$ is a TLG$^*$, then for the chosen $\sigma_{kl}$ the vertices
$t_k$ and $t_l$ are connected by a time path in $\cG^\#$ 
(constructed before we picked $\sigma_{kl}$).

\teo{If $\cG$ is a TLG$^*$ and $\sigma$ a spine in $\cG$, then in the stingy algorithm for $\cG$ with
respect to $\sigma$, each time line \ref{alg2:s6} is executed 
we pick a time-path between two points connected by a time-path in $\cG^{\#}$.}
\dok{
%
%
%
%
Let $n$ be the sum of degrees of vertices in $\cG$ whose degree is at least 3, that is
$$n(\cG)=\sum_{v\in \cV, d(v)\geq 3}d(v).$$ We will prove the following claim
by induction on $n$:\vspace{0.2cm}

{\sl For a TLG$^*$ $\cG$ where $n(\cG)=n$, when we run the algorithm
on $\cG$ for any spine $\sigma$ in line \ref{alg2:s6}  the chosen $\sigma_{kl}$ is such that 
$t_k$ and $t_l$ are connected by a time path in $\cG^\#$ from the previous iteration.}\vspace{0.2cm}

%

For $n=0$ this claim is clearly true (then we have a TLG$^*$ with one spine). 
Assume that this claim holds for all $n\leq m$ where $m\geq 0$.\vspace{0.2cm}

Let's show that this claim holds for $n=m+1$. If there is no such TLG$^*$ $\cG$,
then we say that the claim holds trivially. Otherwise, let $\cG$ be such a TLG$^*$, and 
$\sigma$ its arbitrary spine from $P_{0\to 1}(\cG)$.\vspace{0.2cm} 

We pick $t_u$ and $t_v$ on $\sigma$ that are connected by a time-path $\sigma^1_{uv}$ in $\cG$
outside of $\sigma$ such that $|t_u-t_v|$ is minimal. Let $\sigma_{uv}$ be the time-path between
$t_u$ and $t_v$ on $\sigma$. Note that, by the construction, the cell $(\sigma_{uv},\sigma^1_{uv})$ 
is truly simple. (Otherwise, if the sides $\sigma_{uv}$ and $\sigma^1_{uv}$ are connected
by a path in $\cG[t_u,t_v]$ that would contradict the minimality of $t_v-t_u$.)\vspace{0.2cm}

The graph  constructed by the simple cell collapsing transformation with respect to $(\sigma_{uv},\sigma^1_{uv})$ 
- $\cG^{\circ}$ is by Lemma \ref{lem:trtlt} a TLG$^*$ and we have $n(\cG^\circ)<n(\cG)$ (it is clear that $d((t_u)^\circ)<d(t_u)$ and $d((t_v)^\circ)<d(t_v)$). 
So by induction assumption we can 
apply the algorithm to $\cG^\circ$ and in this way show that it is a TLG$^*$.\vspace{0.2cm}

%
%
%

We will parallely run the algorithm on $\cG^\circ$ and $\cG$ with the given spine $(\sigma)^\circ$
and the corresponding spine $\sigma$. 

We will assume that in the first
iteration of the {\bf while} loop in line \ref{alg2:lwhl} (of the algorithm on $\cG$) time-path $\sigma^1_{uv}$ 
was chosen.\vspace{0.1cm} 

Let $p$ denote the number of iterations of the {\bf while} loop in line \ref{alg2:lwhl}, and 
$\cG^{\#}_p$ the graph constructed until that point when we run the algorithm on $\cG$.\vspace{0.2cm}


Now, we will show that if $\sigma_{kl}^\circ$ was chosen in the $p$-th iteration of the while loop
on $\cG^\circ$, then we can choose $\sigma_{kl}$ in $p+1$-st iteration of the {\bf while} loop on 
$\cG$.\vspace{0.2cm} 

For $p=1$ this holds, $\sigma_{kl}^\circ$ is connecting $(t_k)^\circ$ and $(t_l)^\circ$, and by the construction
of $\cG^\circ$, $t_k$ and $t_l$ are connected by a time path in $\cG$. (Otherwise, we $t_k$ and $t_l$ would be
points on different sides of the cell, connected by a the time path $\sigma_{kl}$, and the cell $(\sigma_{uv},\sigma^1_{uv})$ wouldn't be minimal.) 

Assume this holds for $p=r\geq 1$. 

For $p=r+1$ let $\sigma_{kl}^\circ$, be chosen. By assumption we know that $(t_k)^\circ$ and $(t_l)^\circ$
are connected by a time path in $(\cG')^{\#}_{p-1}$, we know, that $\sigma_{kl}$ is a path connecting $t_k$
and $t_l$ in $\cG$, there can't be a path whose time difference is smaller, because such would exist
in $\cG^\circ$ also. The only thing that we need to show is that $t_k$ and $t_l$ are connected by a time-path
in $\cG^{\#}_p$. 

\begin{figure}[ht]
\begin{center}
\psfrag{1}{$\boldsymbol{t_k}$} \psfrag{2}{$\boldsymbol{t_{k'}}$}
\psfrag{5}{$\boldsymbol{t_u}$} \psfrag{3}{$\boldsymbol{t_{v}}$}
\psfrag{6}{$\boldsymbol{t_l}$} \psfrag{4}{$\boldsymbol{t_{l'}}$}
\psfrag{7}{$\boldsymbol{t_{k'} \wedge t_{l}}$}
\includegraphics[width=8cm]{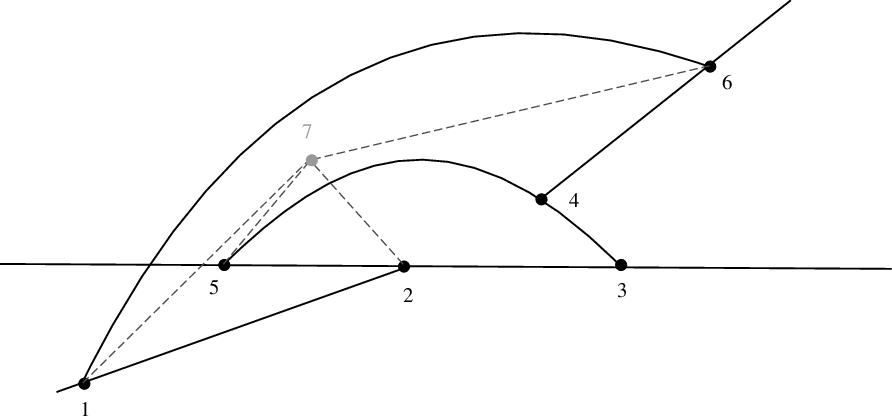}\\ 
\caption{} \label{pic32}

  \end{center}
\end{figure}

Assume the opposite. This would mean, by construction that there exists 
$t_{k'}\in \sigma_{uv}$ and $t_{l'}\in \sigma_{uv}^1$ on different sides of the cell $(\sigma_{uv}^1,\sigma_{uv})$,
such that $t_k\leq t_{k'}\leq t_{l'}\leq t_l$ (time order), and there exist paths $\sigma_{kk'}$ and $\sigma_{l'l}$.
Now this would mean, since $\cG$ is a TLG$^*$, and therefore by Theorem \ref{teo:tpltt} a topological lattice , that
$$\{t\in \cG: t\preceq t_{k}\}\cap \{t\in \cG: t\preceq t_{u}\}\subset \{t\in \cG: t\preceq t_{k'}\}\cap \{t\in \cG: t\preceq t_{l}\}=\{t\in \cG: t\leq t_{k'} \wedge t_{l}\}, $$
Hence, in $\cG$ there exists a vertex $t_{k'} \wedge t_{l}$ connected by time paths to $t_k$, $t_u$, $t_{k'}$ and $t_{l'}$. 
Now, $t_{k'} \wedge t_{l}$ has to be on $\sigma$, or otherwise $t_u$, $t_{k'} \wedge t_{l}$ and $t_{k'}$
form a cell, that will be a smaller cell whose on side is on the spine $\sigma$ in $\cG$. But this contradicts 
the choice of $t_u$ and $t_v$. Now, if $t_{k'} \wedge t_{l}$ is on $\sigma$, this contradicts the choice of $t_k$
and $t_l$, since $t_k<t_{k'} \wedge t_{l}<t_l$, because $t_{k'} \wedge t_{l}$ and $t_l$
are in $\cG^{\#}_p$, they are connected in $\cG$ and their time difference is less than 
$t_l-t_k$. 

Hence $t_k$ and $t_l$ have to be connected in $\cG^{\#}_p$.\vspace{0.2cm}

This shows that the algorithm will be making a connection between two connected vertices in each step.\vspace{0.2cm}

Finally, this proves the claim.}

\pos{The following algorithm determines is a TLG $\cG$ a TLG$^*$ or not:
\begin{center}
\begin{algorithm}[H]
 $\sigma$ a full time-path (i.e. in $P_{0\to 1}(\cG)$)\;
$\cG^\#=(\cV^\#,\cE^\#)$ a TLG that consists of all vertices and all edges of $\sigma$ (in $\cG$)\;\label{alg3:s2}

\While{$\cE\setminus \cE^\#\neq \emptyset$\label{alg3:lwhl}}{
      $\sigma_{kl}$ a time-path in $\cG$ and not in $\cG^{\#}$ between $t_k$ and $t_l$ in $\cV^{\#}$ such that $|t_l-t_k|$ is minimal\;\label{alg3:s6}
       \eIf{$t_k$ and $t_l$ are connected by a time-path in $\cG^{\#}$}{add all edges and vertices that make $\sigma_{kl}$ (in $\cG$) to $\cG^\#$\; \label{alg3:lx}}{
	{\bf return} This is not a TLG$^*$\;}
      	
} 
{\bf return} This is a TLG$^*$\;
\caption{Determine is TLG $\cG$ a TLG$^*$ or not.}\label{alg3:tpltt}
\end{algorithm}
\end{center}}

%
%

\section{TLG's with infinitely many vertices}\label{inf_vert}
\index{Time-like graph (TLG)!with infinitely many vertices}
We will allow $t_0$ and $t_N$ to take values in $\R\cup \{-\infty,\infty\}$.

\defi{\label{def:intlg*}\begin{enumerate}[(i)]
       \item Suppose that the vertex set of a graph $\cG=(\cV,\cE)$ is infinite. We will call
$\cG$ a time-like graph (TLG) if it satisfies the following conditions.
\begin{enumerate}[(a)]
 \item There is a sequence of TLG's $\cG_n=(\cV_n,\cE_n)$ with finite vertex
set $\cV_n$, $n\geq 1$, and for some representations of $\cG_n$'s and $G$ we have
$$\bigcup_{n=1}^{\infty}R(\cG_n)=R(\cG).$$ 
\item \label{def:intlg*:b} The graph $\cG$ is locally finite, i.e. it has a representation $R(\cG)$
such that for any compact $K\subset \R^3$ a finite number of edges intersects $K$.
\end{enumerate}
\item A TLG $\cG$ with infinite vertex set will be called an TLG$^*$ if it satisfies
the following conditions.
\begin{enumerate}[(a)]
  \item We can choose a sequence of TLG$^*$'s $\cG_n$ in (i).
(In the sense of the Definition \ref{def:tlg*}.(\ref{def:tlg*:iii}), i.e. $(\cG_j)_{1\leq j\leq n}$
is a tower of TLG$^*$'s for all $n$.)
  \item Let $\cV_n=\{t_{0,n},t_{1,n},\ldots,t_{N_n,n}\}$. The initial 
vertices $t_{0,n}\in \cV_n$ and $t_{N_n,n}\in \cV_n$ are the same 
for all $\cG_n$, i.e. for all $m,n\geq 1$
$$t_{0,n}=t_{0,m}\quad \textrm{and}\quad t_{N_n,n}=t_{N_m,m}.$$
\item The initial and terminal edges form a decreasing sequence in the representations
of $\cG_n$'s, i.e. if $n>m$
$$E_{t_{0,n},t_{1,n}}((t_{0,n},t_{1,n}))\subset E_{t_{0,m},t_{1,m}}((t_{0,m},t_{1,m}))$$ and
$$ E_{t_{N_n-1,n},t_{N_n,n}}((t_{N_n-1,n},t_{N_n,n}))\subset E_{t_{N_m-1,m},t_{N_m,m}}((t_{N_m-1,m},t_{N_m,m})). $$ 

\end{enumerate}
      \end{enumerate}
 }

The following lemma will be useful for the construction of processes.

\lem{\label{lem:subgphfcv}Let $(\cG_n)$ and $(\cG_n')$ be two TLG$^*$-towers that lead to the construction 
of $\cG$. Let $\cH$ be a sub-graph (not necessarily a TLG$^*$) of some $\cG_{n_0}$
whose vertices have a finite time. Then there exists
$\cG_{n_1}'$ such that $R(\cH)\subset R(\cG_{n_1}')$ and all the vertices of $\cH$ 
are contained in $\cG_{n_1}'$.
\dok{Since $\cG$ is locally finite, there are finitely many vertices with representation 
on $R(\cH)$, also these vertices are of finite degree. For each such vertex $v$, by same argument,
there has to be $\cG_{n_v}'$ such $v$ in $\cG_{n_v}'$ has that degree. Now if $n_1$
is the maximum of $n_v$ over each such vertex $v$ the claim follows. }}
\index{Time-like graph (TLG)|)}

\chapter{Processes indexed by time-like graphs}\label{sec:02}

Let
$\cG=(\cV,\cE)$
be a TLG$^*$. In this chapter we construct a stochastic process on $\cG$
in such a
way that we have a random variable defined
at every point of the representation. (See Figure \ref{pic:sl9a:2}. for illustration.)
\begin{figure}[ht]
\begin{center}
\psfrag{0}{$\boldsymbol{0}$} 
\psfrag{1}{$\boldsymbol{1}$}
\psfrag{t}{\small $\boldsymbol{t}$}
\psfrag{a}{$\boldsymbol{1/3}$}
\psfrag{b}{$\boldsymbol{2/3}$}
\includegraphics[width=7cm]{ex2pm_4.eps}\quad\quad \includegraphics[width=4.5cm]{prez2pr3.eps} \\

\caption{ Time-like graph $\cG$ and a process indexed by it.} \label{pic:sl9a:2}
  \end{center}
\end{figure}

\vspace{0.2cm}
%

\defi{\label{def:pcpath}We define  $X=(X(t):t\in \cG)$\index{Process indexed by a TLG|textbf} as a collection of random variables with
$$X=(X_{E}(t):E=E_{jk}\in \cE,t\in [t_j,t_k]).$$
We will assume the following things.
\begin{itemize}
 \item If $E_{jk},E_{kn}\in \cE$ then $X_{E_{jk}}(t_k)=X_{E_{kn}}(t_k)$.
 \item If $E_{jk},E_{nk}\in \cE$ then $X_{E_{jk}}(t_k)=X_{E_{nk}}(t_k)$.
 \item Finally, if $E_{0j},E_{0k}\in\cE$ then $X_{0j}(t_0)=X_{0k}(t_0)$.
\end{itemize}
For a path $\sigma_1\in\sigma(k_1,k_2,\ldots,k_n)$ we use the notation
$$X_{\sigma_1}(t)=X_{E_{k_{j-1}k_{j}}}(t),$$
for all $j=2,3,\ldots,n$ and $t\in [t_{k_{j-1}},t_{k_j}]$.
}\vspace{0.2cm}

\noindent Remark. (1) If there are two edges $E_{jk}^q$ and $E_{jk}^p$ with the same endpoints 
we will denote processes on them by $X_{jk}^q$ and $X_{jk}^p$.\par (2) We will write $X(t)$ instead of $X_{jk}(t)$ or $X_\sigma$ when this will not cause any confusion.
\par
(3) In an infinite graph case we will do the same thing, but we will not define the process
at $t_0$ and $t_N$, if they are not in $\R$.\vspace{0.3cm}

If $\cP$ is the distribution of a Markov process $(Y(t):t\in [t_0,t_N])$, note that for every TLG there exists a $\cP$-process on $\cG$. 
Trivial example
of a $\cP$-process on a TLG can be constructed by taking a Markov process $(Y(t):t\in [t_0,t_N])$ with
distribution $\cP$ and then letting $X_{\sigma}(t)=Y(t)$ for all full time-paths $\sigma\in P_{0\to 1}(\cG)$.\vspace{0.2cm}

We will require some properties to hold for the process to be non-trivial.

\section{Spine-Markovian property}\index{Spine-Markovian property|(}\label{sec:spine}

\defi{\label{def:spi_mrk} Let $\sigma$ be any
full-time path  (from 0 to 1) in the TLG $\cG=(\cV,\cE)$.
Let $\cG_-$ be a subgraph (not necessarily a TLG) of $\cG$ whose representation is a closure of a connected 
component of $R(\cG)\setminus R(\sigma)$. Let $W$ be the set
of vertices - \textbf{roots} connecting $\cG_-$ to $\sigma$ and let
$\cG_+$ denote the graph represented by $R(\cG)\setminus R(\cG_-)$.

We say that the process $X$ on a TLG $\cG$ is  \textbf{spine-Markovian} if for each such $\sigma$ and $\cG_-$ the processes $(X(t):t\in \cG_-)$ and $(X(t):t\in \cG_+)$ given $(X(t):t\in W)$ are independent.}

\begin{figure}[ht]
\begin{center}
\psfrag{s}{\textcolor{green}{$\boldsymbol{\sigma}$}}
\psfrag{G}{\textcolor{blue}{$\boldsymbol{\cG_-}$}}
\psfrag{H}{$\boldsymbol{\cG_+}$}
\includegraphics[width=10cm]{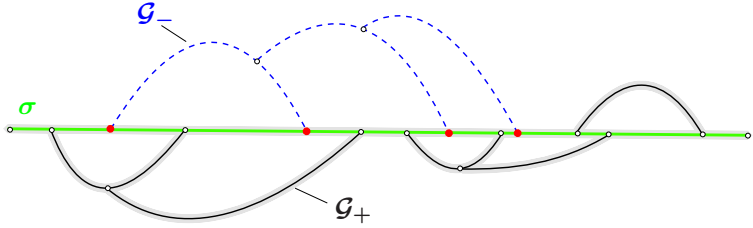}\\
\caption{Spine-Markovian property: The set of roots $\textcolor{red}{W}$ is illustrated by bullet points (\textcolor{red}{$\bullet$}).} \label{pic:sl9a}
  \end{center}
\end{figure}

\prop{\label{prop:spi_mrk}Let $\sigma$, $\cG_-$, $\cG_+$ and $W$ in a TLG $\cG$ be as in the Definition \ref{def:spi_mrk}. Then
for any $\sigma$-algebra $\F$ such that $\sigma(X(t):t\in W)\subset \F\subset \sigma(X_{\sigma})$,
If the process $X$ on $\cG$ is  spine-Markovian then the processes $(X(t):t\in \cG_-)$ and $(X(t):t\in \cG_+)$ given $\F$ are independent.
\dok{Let $Y_-$ and $Y_+$ denote bounded random variables respectively measurable 
in $\sigma(X(t):t\in \cG_-) $ and $\sigma(X(t):t\in \cG_+)$. For $A\in \F$,
 $Y_+\1_A$ is a bounded $\sigma(X(t):t\in \cG_+)$-measurable random variable, and
we have 
\begin{align*}
\E(Y_-Y_+\1_A) & =  \E(\E(Y_-Y_+\1_A|(X(t):t\in W)))\\
&= \E(\E(Y_-|(X(t):t\in W))\E(Y_+\1_A|(X(t):t\in W)))\\
&= \E(\E(Y_-|(X(t):t\in W))\E(\E(Y_+\1_A|\F)|(X(t):t\in W)))\\
&=\E(\E(Y_-|(X(t):t\in W))\E(\E(Y_+|\F)\1_A|(X(t):t\in W)))\\
&=\E(\E(Y_-\E(Y_+|\F)\1_A|(X(t):t\in W)))=\E(Y_-\E(Y_+|\F)\1_A)\\
&=\E(\E(Y_-\E(Y_+|\F)\1_A|\F))=\E(\E(Y_-|\F)\E(Y_+|\F)\1_A).
\end{align*}
}}

\noindent\emph{Remark.} Note that $\cG_+$ is a TLG while $\cG_-$ does not have to be (it is still a connected graph). 
Also, $\cG_+$ contains $\sigma$, so we can find $\cG_-^2$ a connected component
of $R(\cG_+)\setminus R(\sigma)$, and so on\ldots So, the TLG $\cG$ can be decomposed
into $\cG_-^1$, \ldots, $\cG_-^n$ that are connected components of $R(\cG)\setminus R(\sigma)$
and the spine $\sigma$.

\defi{We will call $(\sigma; \cG_-^1, \ldots, \cG_-^n)$ the \textbf{decomposition} of the TLG $\cG$
with respect to $\sigma$. The elements of this decomposition 
(including $\sigma$) we will call \textbf{components}.}

\noindent \emph{Remark.} Notice that the decomposition, given $\sigma$, is unique up to an order of $\cG_-^1, \ldots, \cG_-^n$.

\begin{figure}[ht]
\begin{center}
\psfrag{s}{$\boldsymbol{\sigma}$}
\psfrag{G}{$\boldsymbol{\cG_-^1}$}
\psfrag{H}{$\boldsymbol{\cG_-^2}$}
\psfrag{K}{$\boldsymbol{\cG_-^3}$}
\psfrag{L}{$\boldsymbol{\cG_-^4}$}
\includegraphics[width=12cm]{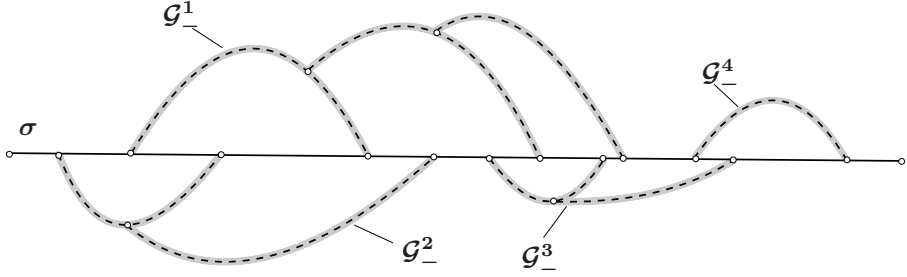}\\
\caption{The decomposition of $\cG$ with respect to $\sigma$.  } \label{pic:sl9b}
  \end{center}
\end{figure}

\prop{\label{thm:decm_tlg}Let $\cG'$ be a graph that is  the union of some graphs $\cG_-^{j_1}$, \ldots, $\cG_-^{j_k}$ in the decomposition,
and $\cG''$ the union of all the other components in the decomposition. Then the process
$(X(t):t\in\cG')$ and $(X(t):t\in \cG'')$ are independent given $(X(t):t\in W_{j_1}\cup \ldots \cup W_{j_k})$,
where $W_{j}$ is the set of roots of $\cG_{j}$. 
\dok{For $l=1,\ldots,k$ let $Y_l$  be a bounded $\sigma (X(t):t\in \cG_{j_l})$-measurable random 
variable, $Z$ a bounded $\sigma(X(t):t\in \cG'')$-measurable random variable and $A\in\sigma (X(t):t\in W_{j_1}\cup \ldots \cup W_{j_k})$.
Using the spine-Markovian property for each $\cG_{j_l}$ $l=1,2,\ldots,k$ at a time with respect to $\sigma$
we get
\begin{align*}
 \E(Y_1Y_2\ldots Y_kZ\1_A) &= \E(\E(Y_1|X_{W_{j_1}})Y_2\ldots Y_kZ\1_A)\\
&= \E(\E(Y_1|X_{W_{j_1}})\E(Y_2|X_{W_{j_2}})\ldots Y_kZ\1_A)\\
& \ \ldots\ \\
&= \E(\E(Y_1|X_{W_{j_1}})\E(Y_2|X_{W_{j_2}})\ldots \E(Y_k|X_{W_{j_k}})Z\1_A).
\end{align*}
Now, taking the conditional expectation with respect to $\sigma(X_{W_{j_1}}\ldots X_{W_{j_k}})$
\begin{align*}
 &= \E(\E(\E(Y_1|X_{W_{j_1}})\ldots \E(Y_k|X_{W_{j_k}})Z\1_A|X_{W_{j_1}}\ldots X_{W_{j_k}}))\\
&= \E(\E(Y_1|X_{W_{j_1}})\ldots \E(Y_k|X_{W_{j_k}})\E(Z|X_{W_{j_1}}\ldots X_{W_{j_k}})\1_A)
\end{align*}
Now, again using the spine-Markovian property on each graph in the union we get
\begin{align*}
 &= \E(Y_1\ldots \E(Y_k|X_{W_{j_k}})\E(Z|X_{W_{j_1}}\ldots X_{W_{j_k}})\1_A)\\
&\ldots\\
&= \E(Y_1\ldots Y_k\E(Z|X_{W_{j_1}}\ldots X_{W_{j_k}})\1_A).\\
\end{align*}
Which, finally, gives us
\begin{align*}
 &= \E(\E(Y_1\ldots Y_k\E(Z|X_{W_{j_1}}\ldots X_{W_{j_k}})\1_A|X_{W_{j_1}}\ldots X_{W_{j_k}}))\\
 &= \E(\E(Y_1\ldots Y_k|X_{W_{j_1}}\ldots X_{W_{j_k}})\E(Z|X_{W_{j_1}}\ldots X_{W_{j_k}})\1_A).\\
\end{align*}
Now from the Monotone Class Theorem the claim follows.  }}

We will need a stronger property for some proofs.

\defi{\label{def:S*}For a TLG$^*$ $\cG$ we define $S^*(\cG)$ to be the set of all TLG$^*$'s 
$\cH$ such that there exists a TLG$^*$-tower $(\cK_k)_{k=0}^n$ that starts with 
$\cK_0=\cH$ and ends with $\cK_n=\cG$.}

\defi{\label{def:her_mp}The process $(X(t):t\in \cG)$ has a \textbf{hereditary spine-Markovian property}\index{Spine-Markovian property!hereditary} 
if $(X(t):t\in \cH)$ is a spine-Markovian process
for each $\cH\in S^*(\cG)$. }
\index{Spine-Markovian property|)}

\section{Consistent distributions on paths}\label{cdst_pth}
\defi{Let $\cG$ be a TLG, for a family of distributions of stochastic processes 
on $[0,1]$
$$\{\mu_{\sigma}:\sigma \in H\},$$
where $H\subset P_{0\to 1}(\cG)$ (a subset of the set of full time-paths),  we say
that it is \textbf{consistent}\index{Consistent family of distributions along time-paths|see{Process indexed by a TLG}}
\index{Process indexed by a TLG!consistent family of distributions along time-paths}
 if for $\sigma_1,\sigma_2 \in H$
$$\mu_{\sigma_1}\circ\pi_T^{-1}=\mu_{\sigma_2}\circ\pi_T^{-1},$$
where $T=\{t:t\in E,E\in\sigma_1\ \&\ E\in\sigma_2\}$.}

\prop{If $\mu$ is the distribution of the process $X$ on a TLG $\cG$, then 
\be\{\mu_\sigma= \P\circ X_\sigma^{-1}:\sigma\in P_{0\to 1}(\cG)\} \label{mfam}\ee
is a consistent family.}

\begin{figure}[ht]
\begin{center}
 \includegraphics[width=5cm]{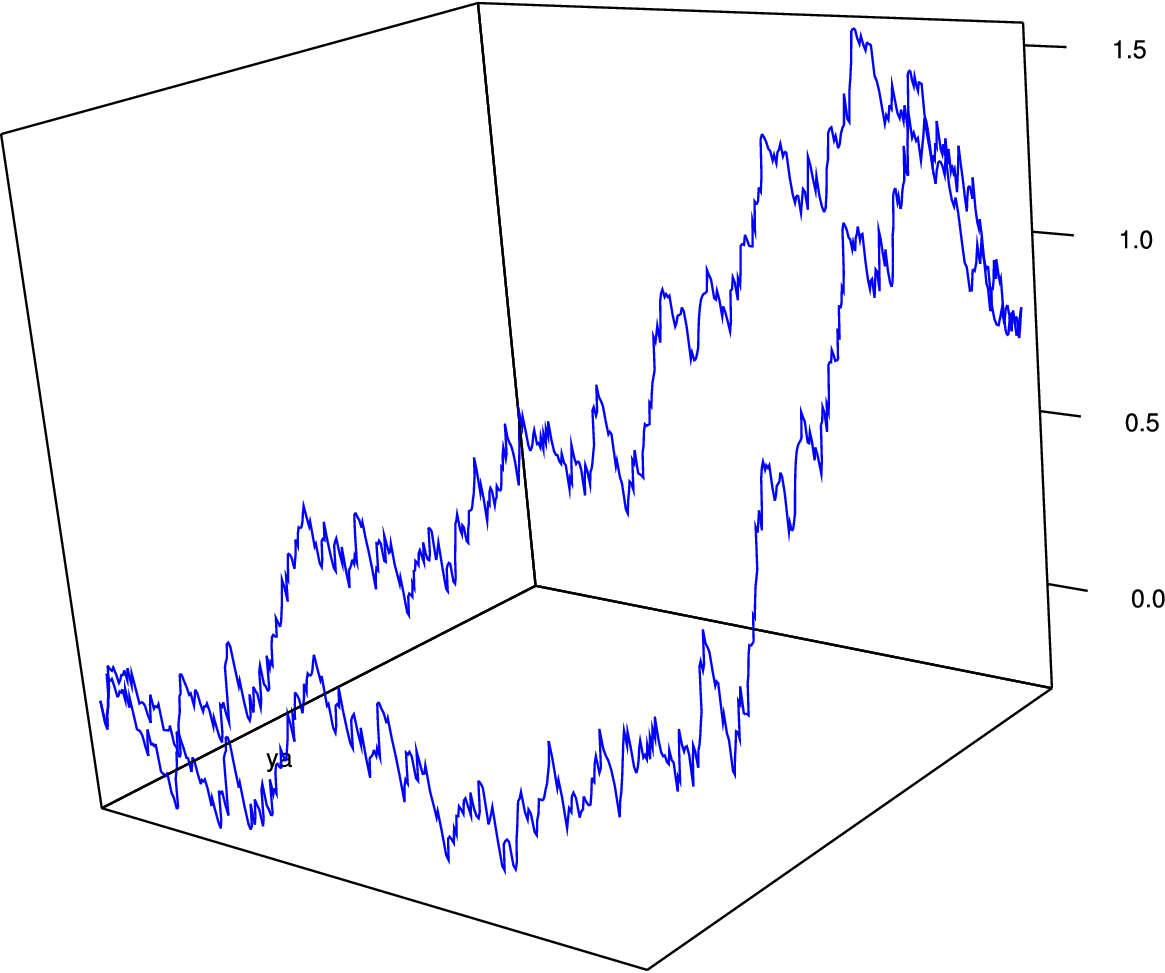}\quad \includegraphics[width=5cm]{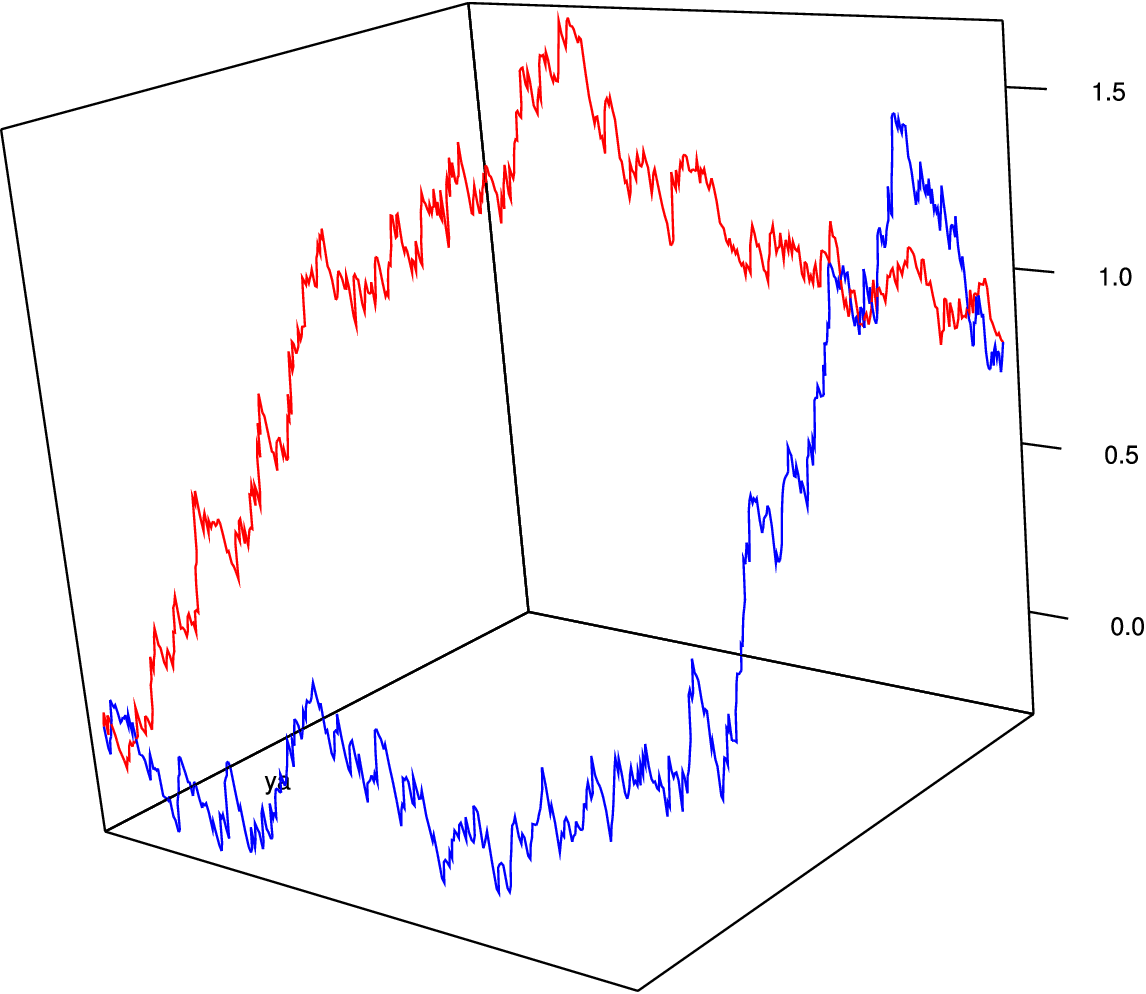}\\
\caption{Processes $X^1$ and $X^2$} \label{pic:prez1ab}
  \end{center}
\end{figure}

\noindent \emph{Remark.} It is not hard to see that the family of distributions given by $(\ref{mfam})$
does not uniquely determine $\mu$ - the distribution on $\cG$. For example if we take 
a Markov process $\cP$ on $[0,1]$, and we take the TLG graph $\cG=(\cV,\cE)$ where 
$\cV=\{0,1\}$ and $\cE=\{E_{01}^1, E_{01}^2\}$. Let $Y^1$ be a Markov process
on $[0,1]$ with distribution $\cP$, and $Y^2$ a $\cP$-Markov bridge
starting at $Y^1(0)$ and ending at $Y^1(1)$ conditionally independent 
given $Y^1(0)$ and $Y^1(1)$. (This can be done as in Theorem \ref{2mrk_pro}.)
Now, the process $X^1$ such that $X^1_{E^1_{01}}=Y^1$ and $X^1_{E^2_{01}}=Y^1$,
has the same distributions along the full-time paths as $X^2$ given by $X^2_{E^1_{01}}=Y^1$ 
and $X^2_{E^2_{01}}=Y^2$. But, these two processes are clearly different in distribution. (See Figure \ref{pic:prez1ab}.)

\pos{Let $\cP$ be a distribution of some process on $[0,1]$. If $\mu_\sigma=\cP$ for each full-time path $\sigma$ in TLG $\cG$, then 
$$\{\mu_\sigma : \sigma\in P_{0\to 1}(\cG)\}$$ is a consistent family.}\vspace{0.3cm}

\section{Construction from a consistent family}\label{constr_cnsfm}\index{Process indexed by a TLG!construction|(}\label{condit}
The interesting question is if we have a consistent family 
$$\cM:=\{\mu_{\sigma}:\sigma\in P_{0\to 1}(\cG)\},$$
under what conditions can we construct a process $X$ on $\cG$ such that $X_\sigma$ has the distribution
$\mu_{\sigma}$. We will call $X$ to be an \textbf{$\cM$-process}\index{Process indexed by a TLG!$\cM$-process}\index{$\cM$-process|see{Process indexed by a TLG}}\label{def:Mprc}\vspace{0.3cm}

We will show that such a process exists under the following assumptions\index{Process indexed by a TLG!construction!(3T) conditions}:
\begin{description}
  \item[T1] $\cG=(\cV,\cE)$ is a TLG$^*$.
 \item[T2] $\cM$ is a consistent family of measures that induce continuous or RCLL processes.
 \item[T3] For each simple cell $(\sigma_1,\sigma_2)$ in $\cG$ starting at $t_*$ and ending
at $t^*$, if $\sigma$ is a full-time path containing $\sigma_1$ (or $\sigma_2$) then the $\mu_\sigma$-distributed 
process
\be(Y(t):t\in [0,1])\label{wptm}\ee
has the property that $(Y(t):t\in [0,t_*]\cup [t^*,1])$ and $(Y(t):t\in [t_*,t^*])$ given $Y(t_*)$ and $Y(t^*)$
are independent.
\end{description}
Conditions (T1)-(T3) we will call \textbf{(3T) conditions}.

\noindent\emph{Remark.} Condition (T2) is needed so that we could define a conditional 
distribution when needed. So other $\cM$ can be a family of other types of processes for which
this would be possible (for example all the arguments would work for discrete processes). \vspace{0.3cm}

The condition (T3) can be rewritten in a different way.
\lem{\label{lem:wptm}The process given by $(\ref{wptm})$ has the property
that the distribution $(Y(t):t\in [t_*,t^*])$ given $(Y(t):t\in [0,t_*]\cup [t^*,1])$
depends only on $(Y(t_*),Y(t^*))$, in other words if $Z$ is a bounded 
$\sigma(Y(t):t\in [t_*,t^*])$-measurable random variable then
$$\E(Z|Y(t):t\in [0,t_*]\cup [t^*,1])=\E(Z|Y(t_*),Y(t^*)).$$
\dok{Let $A\in \sigma(Y(t):t\in [0,t_*]\cup [t^*,1])$ and define $U:=\E(Z|Y(t):t\in [0,t_*]\cup [t^*,1])$. 
Using the definition of the conditional expectation, and the property of $Y$ 
\begin{align*}
 &\E(U\1_A)\\
&=\E(Z\1_A)=\E(\E[Z\1_A|Y(t_*),Y(t^*)])\\
&=\E(\E[Z|Y(t_*),Y(t^*)]\E[\1_A|Y(t_*),Y(t^*)])\\
&=\E(\E[\E[Z|Y(t_*),Y(t^*)]\1_A|Y(t_*),Y(t^*)])\\
&=\E(\E[Z|Y(t_*),Y(t^*)]\1_A).
\end{align*}
The claim follows from the a.s. uniqueness of the conditional expectation.
}}

\subsection{Construction}\label{constr}
We will define a $\cM$-process on a TLG$^*$ $\cG$
with finite sets $\cV$ and $\cE$, where 
$t_0=0$ and $t_N=1$.

\defi{\label{def:constr}Let $(\cG_l)_{0\leq l\leq n}$ be a tower of TLG$^*$ where $\cG_0$
is a minimal graph $\cV_0=\{t_0=0,t_N=1\}$, $\cE_0=\{E_{0N}\}$ and $\cG_n=\cG$. Further let 
$\cM$ be a family of distributions satisfying (3T) conditions.
\begin{itemize}
 \item On $\cG_0$ we define a process $X^0$ with $\mu_{E_{0n}}$ distribution.
 \item If we have already defined $X^l$ on $\cG_l$ (for some $l<n$), then
we define $X^{l+1}$ on $\cG_{l+1}$ in the following way depending how we
constructed $\cG_{l+1}$ from $\cG_l$ (recall part (\ref{def:tlg*:ii}) of Definition \ref{def:tlg*}.).
\begin{enumerate}[(1)]
 \item In the construction a new vertex $\tau_l\in [0,1]\setminus \cV_{l}$ was added to graph $\cG_l$, by subdividing
some $E_{jk}$ such that $t_{j}<\tau_l<t_{k}$, into $E_{jl}$ and $E_{lk}$ to get $G_{l+1}$.  In this case, the two graphs $\cG_l$ and $\cG_{l+1}$ have a common representation,
$R(\cG_l)=R(\cG_{l+1})$, and we can define $X^{l+1}$ on $\cG_{l+1}$ to have the same values on this
representation as $X^{l}$.
\begin{figure}[ht]
\begin{center}
\psfrag{1}{$\boldsymbol{1}$}
\psfrag{0}{$\boldsymbol{0}$}
\psfrag{t}{$\boldsymbol{\tau_l}$}
\psfrag{2}{$\boldsymbol{E_{jl}}$}
\psfrag{3}{$\boldsymbol{E_{lk}}$}
\psfrag{4}{$\boldsymbol{E_{jk}}$}
\includegraphics[width=5cm]{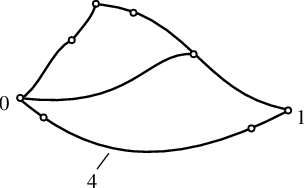} $\to$ \includegraphics[width=5cm]{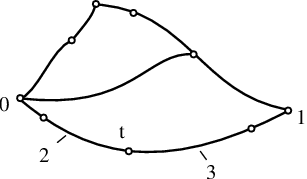}\\
\caption{Case (1) in the construction.} \label{pic:sl5a2}

\psfrag{1}{$\boldsymbol{1}$}
\psfrag{0}{$\boldsymbol{0}$}
\psfrag{4}{$\boldsymbol{E_{jk}}$}
\psfrag{5}{$\boldsymbol{E_{jk}^*}$}
\includegraphics[width=5cm]{sl5a1.eps} $\to$ \includegraphics[width=5cm]{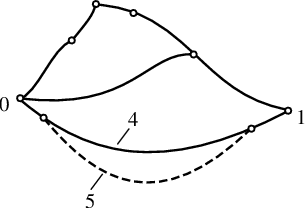}\\
\caption{Case (2) in the construction.} \label{pic:sl5b}

  \end{center}
\end{figure}

\item In the construction a new edge between two vertices $t_j<t_k$ in $\cV_{l}$ that are connected by a time path
in $\cG_l$ , was added to get $G_{l+1}$. So, $\cG_{l+1}$ has a new edge $E^*=E^*_{jk}$. 
Let $Z_j=X^l(t_{j})$ and $Z_k=X^l(t_{k})$. 

Now we pick a full-time path $\sigma$ that contains $E^*$. Now we define 
$\mu_*(\cdot |x,y)$ to be the conditional probability of the process with
the distribution $\mu_{\sigma}\circ \pi_{[t_j,t_k]}^{-1}$ conditioned to
have values $x$ at $t_j$ and $y$ at $t_k$. So we construct the process
$X^{l+1}$ in such a way that $X^{l+1}$ on $R(\cG_l)$ is equal to $X^l$
and $X^{l+1}_{E^*}$ is the process with distribution $\mu_*(\cdot |Z_j,Z_k)$
and otherwise independent of $X^l$ given $Z_j$ and $Z_k$.

\end{enumerate}

\end{itemize}
Since $n$ is finite this procedure will end and we will have a process $X=X^n$ defined on $\cG$.}\vspace{0.2cm}

\noindent \emph{Remark.} Note that this construction, i.e. the definition of $X$ on $\cG$ \textsl{depends} on the choice of the TLG$^*$ tower $(\cG_l)_{0\leq l\leq n}$.
\index{Process indexed by a TLG!construction|)}

\subsection{Constructed process is an $\cM$-process}\index{Process indexed by a TLG!$\cM$-process|(}

\defi{If $(\cG_k)_{k=0}^n$ is a TLG$^*$-tower where $\cG_n=\cG$.
If $\cM$ is a family of distributions on full time-paths of
$\cG$. This naturally induces a family $\cM(\cG_k)$ of distributions on full time-paths of
$\cG_k$.}

\noindent \emph{Remark.} This is well-defined since a representation of every full time-path in $\cG_k$,
is a representation of a full time-path in $\cG$ (in the consistent representation of the TLG$^*$-tower $(\cG_k)_{k=0}^n$).

The only question remains will the family induced  by $\cM$ have the same 
properties as $\cM$. This is shown to be true.

\lem{\label{lem:inr:A13}If $\cM$ is a family of distributions on full time-paths of a TLG$^*$ $\cG$
satisfying properties (T1)-(T3), then for any $\cH\in S^*(\cG)$ the 
family $\cM(\cH)$ also satisfies properties (T1)-(T3).
\dok{(T1) is clearly satisfied since $\cH$ is a TLG$^*$. (T2) is satisfied
since in the consistent representation all the full time paths in $\cH$
are full time paths in $\cG$. By Corollary \ref{pos:sm_cell}, in a consistent representation
a representation of a simple cell in $\cH$ is a representation of a simple cell
in $\cG$. Therefore (T3) holds.}}

\lem{\label{lem:P-proc}The process $X$ on $\cG$ defined in \ref{constr} is an $\cM$-process.}
\begin{proof}
 It is clear that $X^0$ is a $\cM(\cG_0)$-process on the minimal graph
$\cG_0$.\vspace{0.2cm}

For, $l<n$ we assume $X^l$ is a $\cM(\cG_l)$-process on $\cG_l$. If we got 
$X^{l+1}$ using step (1) in the construction, then we inherited this property
from $X^{l}$, since $\cM(\cG_l)=\cM(\cG_{l+1})$. If we got $X^{l+1}$ using step (2), recall that $\cG_{l}$ contains a time-path $\sigma_{jk}$ connecting $t_j$ and $t_k$, so there
is a full path $\sigma'$ in $\cG_{l+1}$ that starts with a time-path $\sigma_{0j}$ from $t_0$ to $t_j$, contains $\sigma_{jk}$, and ends with a time-path $\sigma_{kN}$. 
\begin{figure}[ht]
\begin{center}
\psfrag{1}{$\boldsymbol{1}$}
\psfrag{0}{$\boldsymbol{0}$}
\psfrag{j}{$\boldsymbol{t_j}$}
\psfrag{k}{$\boldsymbol{t_k}$}
\psfrag{3}{$\boldsymbol{E_{jk}^*}$}
\psfrag{2}{$\boldsymbol{\sigma'}$}
\psfrag{4}{$\boldsymbol{\sigma^*}$}
\includegraphics[width=7cm]{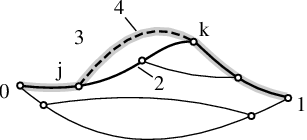}\\
\caption{} \label{pic:sl6}
  \end{center}
\end{figure}

Now for every full time-path $\sigma$ that contains the new edge $E^*=E_{jk}^*$ that was added 
in the construction, $\mu_{E^*}=\mu_{\sigma}\circ \pi_{[t_j,t_k]}^{-1}$ is well
defined since $\cM$ is a consistent family, and $\mu_{E^*}$ doesn't depend on the choice of $\sigma$.\vspace{0.2cm} 

Since, $(E^*,\sigma_{jk})$ is a simple cell, by property (T3) of $\cM$ we have
that $\mu_*(\cdot|x,y)$ is the conditional distribution of a $\mu_{E^*}$-distributed process on $[t_j,t_k]$ conditioned to have value $x$ at $t_j$ and $y$ at $t_k$.\vspace{0.2cm}

The process $X^{l+1}_{E^*}$, by construction, is independent given $(t_j,Z_j)$ and $(t_k,Z_k)$. 
By property (T3) of $\cM$ the distribution of $X^{l+1}_{\sigma'}$  where 
$\sigma'$ is the union of $\sigma_{0j}$, $E^*$, $\sigma_{k1}$ is
given by $\mu_{\sigma'}$: Let $A_0$ be an event in the path $\sigma$-algebra
on $[0,t_j]$, 
$A_1$ an event in the path $\sigma$-algebra
on $[t_k,1]$, and $B$ an event in the path $\sigma$-algebra
on $[t_j,t_k]$ we have:
\begin{align*}
 \P\circ (X^{l+1}_{\sigma'})^{-1}(A_0\times B\times A_1)&=\E(\1_{A_0}(X^{l+1}_{\sigma_{0j}})\1_{B}(X^{l+1}_{E^*})\1_{A_1}(X^{l+1}_{\sigma_{k1}}))\\ 
&= \E(\E(\1_{A_0}(X^{l+1}_{\sigma_{0j}})\1_{B}(X^{l+1}_{E^*})\1_{A_1}(X^{l+1}_{\sigma_{k1}})|Z_j,Z_k))\\
&=\E(\1_{A_0}(X^{l+1}_{\sigma_{0j}})\E(\1_{B}(X^{l+1}_{E^*})|Z_j,Z_k)\1_{A_1}(X^{l+1}_{\sigma_{k1}}))\\
&\stackrel{(T3)}{=} \int_{A_0\times A_1}\mu_{*}(B|\pi_{t_j}(x),\pi_{t_k}(x))\mu_{\sigma}\circ \pi_{[0,t_j]\cup [t_k,1]}^{-1}(dx)\\
&\stackrel{(T2)}{=}\int_{A_0\times A_1}\mu_{*}(B|\pi_{t_j}(x),\pi_{t_k}(x))\mu_{\sigma'}\circ \pi_{[0,t_j]\cup [t_k,1]}^{-1}(dx)\\
&\stackrel{(T3)}{=} \mu_{\sigma'}(A_0\cap B\cap A_1).
\end{align*}
By Monotone Class Theorem $\P\circ (X^{l+1}_{\sigma'})^{-1}=\mu_{\sigma'}$.
\end{proof}

\noindent \emph{Remark.} Note that just for the  existence of an $\cM$-process
on the TLG$^*$ we could weaken condition (T3). If we fix a construction to a
TLG$^*$-tower, then only some full time-paths need to have the described 
property, but then we would lose some properties of the constructed process.

\index{Process indexed by a TLG!$\cM$-process|)}

\subsection{The constructed process is a spine-Markovian process}\index{Process indexed by a TLG!spine-Markovian property|(}\label{sec:spine2}

\lem{\label{lem:sp-mk}The process $X$ on $\cG$ defined in \ref{constr} is a spine-Markovian process.}
\begin{proof}
$X^0$ is trivially an spine-Markovian process. Let's assume that $X^{l}$
is spine-Markovian. We have two cases to study to show that $X^{l+1}$ is spine-Markovian.\vspace{0.3cm}

$(\bullet 1)$ If we added a new vertex to the graph $\cG_l$ to obtain $\cG_{l+1}$.
Then the spine-Markovian property is directly inherited from the process $X^l$, since 
$W$ can't contain the new vertex.\vspace{0.3cm}

$(\bullet 2)$ We added a new edge $E^*$  to the graph $\cG_l$ between two existing
time-path connected vertices to obtain $\cG_{l+1}$. Pick a full time path $\sigma$, and the subgraphs $\cG_-$ and $\cG_+$ in
the graph $\cG_{l+1}$ (in the sense of the Definition \ref{def:spi_mrk}). First, note 
that from the construction the process $X^{l+1}_{E^*}$ is
independent from $X^{l+1}_{\cG_+}$ given the values of the process 
at the endpoints of $E^*$. We will call this property {\it edge-Markovian for the edge $E^*$ (in $\cG_{l+1}$)}.
(This property does not need to hold for other edges.) This will be used often during
the this proof.  We have the following cases.\vspace{0.15cm}

$(\circ 1)$ If the new edge is the only edge in $\cG_-$, i.e. $E^*$
is connecting two vertices on $\sigma$. The claim follows from the edge-Markovian property
for $E^*$.\vspace{0.15cm}

$(\circ 2)$ The new $E^*=E^*_{t_1^*t_2^*}$ edge is in $\cG_-=(\cV_-,\cE_-)$, but one of the vertices that $E^*$ is connecting 
is on $\sigma$. (See Figure \ref{pic9c}.) 
Let $\cG_{-}^*$ be the graph in $\cG_{l}$ that has the edges $\cE_-\setminus \{E^*\}$.  From Proposition \ref{prop:spi_mrk}. and the spine-Markovian property of $X^l$
we know $(X^{l+1}(t):t\in \cG_{-}^*)$ and  $(X^{l+1}(t):t\in \cG_+)$ given $(X(t):t\in W)$
are independent. (Note that one vertex in $W$ may not be in $\cG_{-}^*$.)
Now, let $Y_{-}^*$ be a bounded $\sigma(X^{l+1}(t):t\in \cG_{-}^*)$-measurable,
$Y_*$ a bounded $\sigma(X^{l+1}(t):t\in E^*)$-measurable, and $Y_+$ 
a bounded $\sigma(X^{l+1}(t):t\in \cG_+)$ measurable random variable. For
$A\in \sigma (X^{l+1}(t):t\in W)$ we have using edge-Markov property for $E^*$:

$$ \E(Y^*_-Y_*Y_+\1_A)=\E(Y^*_-\E(Y_*|X_{t_1^*},X_{t_2^*})Y_+\1_A).$$
Now, since $Y^*_-\E(Y_*|X_{t_1^*},X_{t_2^*})\in \sigma(X^{l+1}(t):t\in \cG_{-}^*)\vee \sigma (X(t):t\in W)$,
and this is independent of  $(X^{l+1}(t):t\in \cG_+)$ given $(X(t):t\in W)$. So,
\begin{align*}
\E(Y^*_-\E(Y_*|X_{t_1^*},X_{t_2^*})Y_+\1_A)&=\E(\E(Y^*_-\E(Y_*|X_{t_1^*},X_{t_2^*})Y_+\1_A|(X(t):t\in W)))\\
&= \E(\E(Y^*_-\E(Y_*|X_{t_1^*},X_{t_2^*})Y_+|(X(t):t\in W))\1_A)\\
&=\E(\E(Y^*_-\E(Y_*|X_{t_1^*},X_{t_2^*})|(X(t):t\in W))\E(Y_+|(X(t):t\in W))\1_A)\\
&=\E(Y^*_-\E(Y_*|X_{t_1^*},X_{t_2^*})\E(Y_+|(X(t):t\in W))\1_A)
\end{align*}

\begin{figure}[ht]
\begin{center}
\psfrag{s}{\textcolor{green}{$\boldsymbol{\sigma}$}}
\psfrag{G}{\textcolor{blue}{$\boldsymbol{\cG_-}$}}
\psfrag{H}{$\boldsymbol{\cG_+}$}
\psfrag{E}{$\textcolor{red}{E^*}$}
\includegraphics[width=9cm]{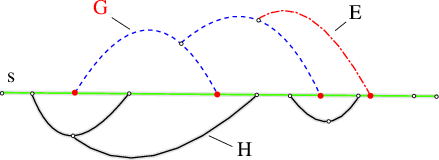}\\
\caption{}\label{pic9c}
  \end{center}
\end{figure}

Again using edge-Markovian property for $E^*$ we get
\begin{align*}
&\E(Y^*_-\E(Y_*|X_{t_1^*},X_{t_2^*})\E(Y_+|(X(t):t\in W))\1_A)\\
&=\E(Y^*_-Y_*\E(Y_+|(X(t):t\in W))\1_A)\\
&=\E(\E(Y^*_-Y_*\E(Y_+|(X(t):t\in W))\1_A|(X(t):t\in W)))\\
&=\E(\E(Y^*_-Y_*|(X(t):t\in W))\E(Y_+|(X(t):t\in W))\1_A).
\end{align*}
This proves the claim for $(\circ 2)$.\vspace{0.15cm}

$(\circ 3)$ The new $E^*$ edge is in $\cG_-=(\cV_-,\cE_-)$, both of the vertices that $E^*$ is connecting 
are not on $\sigma$ but are on $\cG_-$. In this case we fist use the edge-Markov property for $E^*$
and then in the similar way as in $(\circ 2)$ we use the spine-Markovian property
or Theorem \ref{thm:decm_tlg} if the graph $(\cV_-,\cE_-\setminus\{E^*\})$ is made of
two components.\vspace{0.15cm}

$(\circ 4)$ The new $E^*$ edge is in $\cG_+=(\cV_+,\cE_+)$ and not a part of  $\sigma$. Using the spine-Markov
property of $X^l$ we know that $(X^{l+1}(t):t\in\cE_+\setminus\{E^*\})$ and $(X^{l+1}(t):t\in \cE_-)$
are independent given $(X^{l+1}(t):t\in W)$. Using the edge-Markovian property for 
$X^{l+1}$ we get that $(X^{l+1}(t):t\in\cE_+)$ and $(X^{l+1}(t):t\in \cE_-)$
are independent given $(X^{l+1}(t):t\in W)$. (This is proven similar as in $(\circ 2)$.)\vspace{0.15cm}

$(\circ 5)$ If $E^*=E_{t_1^*t_2^*}$ is a part of the spine $\sigma$. By the construction
of $E^*$ we know that there exists a time-path going through vertices 
$t_1^*$ and $t_2^*$, and therefore there is a full time-path $\sigma'$ which
contains whole of $\sigma$ except $E^*$. Let $\sigma_{12}'$ be the part of $\sigma'$
connecting $t_1^*$ and $t_2^*$. We will use the spine-Markov property for $\sigma'$ on $\cG_l$
to prove the one for $\sigma$ on $\cG_{l+1}$. Take $\cG_-$ and $W$ in $\cG_{l+1}$ relative to
$\sigma$. Clearly, none of the vetrices in $W$ are on $E^*$. 
If none of them are on
$\sigma_{12}'$ (except maybe $t_1^*$ and $t_2^*$), we can apply the spine-Markovian property
relative to $\sigma'$ in the case ($\circ 3$), and we are done.

\begin{figure}[ht]
\begin{center}
\psfrag{z}{$\boldsymbol{\sigma}$}
\psfrag{s}{$\boldsymbol{\sigma'}$}

\psfrag{E}{$\boldsymbol{E^*}$}
\includegraphics[width=12cm]{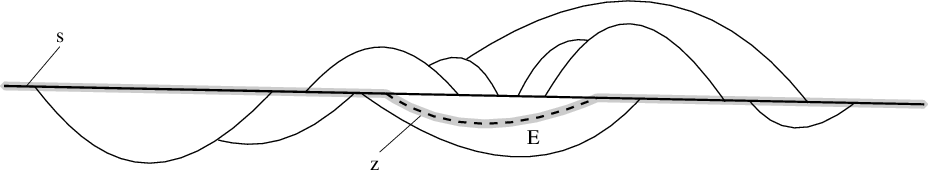}\\
\caption{The two spines $\sigma$ and $\sigma'$.} \label{pic:sl10a}
\psfrag{z}{$\boldsymbol{\sigma}$}
\psfrag{G}{$\boldsymbol{\cG_-}$}
\psfrag{H}{$\boldsymbol{\cG_+}$}
\psfrag{E}{$\boldsymbol{E^*}$}
\includegraphics[width=12cm]{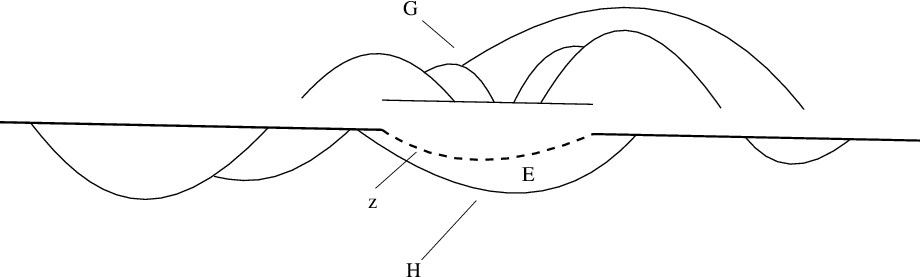}\\
\caption{$\cG_-$ and $\cG_+$ with respect to $\sigma$.} \label{pic:sl10b}
\psfrag{1}{$\boldsymbol{\cG_-^1}$}
\psfrag{2}{$\boldsymbol{\cG_-^2}$}
\psfrag{3}{$\boldsymbol{\cG_-^3}$}
\psfrag{4}{$\boldsymbol{\cG_-^4}$}
\psfrag{5}{$\boldsymbol{\cG_-^5}$}
\psfrag{A}{$\boldsymbol{A}$}
\psfrag{B}{$\boldsymbol{B}$}
\psfrag{C}{$\boldsymbol{C}$}
\psfrag{D}{$\boldsymbol{D}$}
\psfrag{E}{$\boldsymbol{E}$}
\includegraphics[width=12cm]{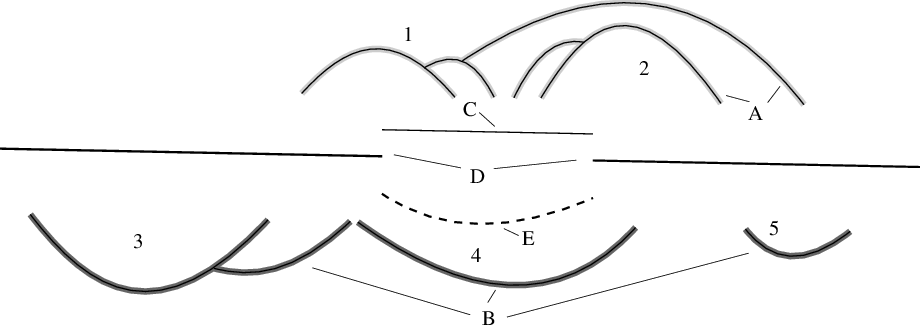}\\
\caption{$A$, $B$, $C$, $D$ and $E$ parts of $\cG$.} \label{pic:sl10c}
  \end{center}
\end{figure}

If some of the vertices in $W\setminus\{t_1^*,t_2^*\}$ are on $\sigma_{12}'$, then the whole $\sigma_{12}'$
is in $\cG_-$. Let's decompose $\cG_l$ with the respect to $\sigma'$. Now, the graph $\cG_-$
is a union of some components $\cG_{l-}^{1}$, \ldots, $\cG_{l-}^{k}$ and $\sigma_{12}'$.
$\cG_+$ is a union of  some other components $\cG_{l-}^{k+1}$, \ldots, $\cG_{l-}^{h}$ and $\sigma$.
Now we look a the following parts of $\cG$ (for a illustration see Figure \ref{pic:sl10c}.)
\begin{itemize}
 \item $A=\cG_{l-}^{1}\cup \ldots\cup\cG_{l-}^{k}$, $W_A=W_1\cup\ldots\cup W_k$.

\item $B=\cG_{l-}^{k+1}\cup \ldots\cup\cG_{l-}^{h}$, $W_B=W_{k+1}\cup\ldots \cup W_h$.
Note that $W_B\setminus\{t_1^*,t_2^*\}$  contains no vertices on the path $\sigma_{12}$
\item $C$ is the graph containing $\sigma_{12}$.
\item $D$ is the graph containing $\sigma'$ without $\sigma_{12}$.
This is the same as a graph containing $\sigma$ without $E^*$.
\item $E$ is the graph containing $E^*$.
\end{itemize}
Let's review which parts we have in the important graphs 
$$\begin{array}{c|c|c|c|c}
   \cG_-&\cG_+ &\sigma & \sigma'&E^*\\ \hline
A,C& B,D,E& D,E & C,D& E
  \end{array}.
$$
Note that $W$ the roots of $\cG_-$ are $(W_A\setminus (t_1^*,t_2^*))\cup \{t_1^*,t_2^*\}$.
Let $Y_H$ be a bounded $\sigma(X(t):t\in H)$-measurable random variable, for $H=A,B,C,D,E$, and let 
$T\in \sigma(X_W)$.
Now, we go step by step, using the right Markovian properties. First we use the edge-Markovian
property for edge $E^*=E_{t_1^*t_2^*}$, hence $Y_E$ is independent of the rest of the
$Y$-variables given $X(t_1^*)$ and $X(t_2^*)$
$$\E(Y_AY_BY_CY_DY_E\1_T)=\E(Y_AY_BY_CY_D\E(Y_E|X(t_1^*),X(t_2^*))\1_T).$$
Now, we apply the spine-Markovian property for on $A$ and $B$ relative to $\sigma'$:
\begin{align}
 =&\E(\E(Y_A|X_{W_A})Y_BY_CY_D\E(Y_E|X(t_1^*),X(t_2^*))\1_T) \nonumber\\
=&\E(\E(Y_A|X_{W_A})\E(Y_B|X_{W_B})Y_CY_D\E(Y_E|X(t_1^*),X(t_2^*))\1_T). \label{sp_dk:1}
\end{align}
Note, that $X_{W_A}$, $X_{W_B}$, $Y_C$, $Y_D$, and $X(t_1^*)$, $X(t_2^*)$, 
are all $\sigma(X_{\sigma'})$-measurable, and we can use the (T3) property
of the process $X_{\sigma'}$ ($X$ on the path $\sigma'$). Let $\F_{\sigma'}(t_1^*)=\sigma\{X_{\sigma'}(t):t\leq t_1^*\}$
and $\cG_{\sigma'}(t_2^*)=\sigma\{X_{\sigma'}(t):t\geq t_2^*\}$. Now, we take the conditional expectation
in $(\ref{sp_dk:1})$ with respect to $\F_{\sigma'}(t_1^*)\vee \cG_{\sigma'}(t_2^*)$. Note that
$Y_D $, $\1_T$ and $X_{W_B}$ are $\F_{\sigma'}(t_1^*)\vee \cG_{\sigma'}(t_2^*)$-measurable. Hence,
 
\begin{align}
 =&\E(\E[\E(Y_A|X_{W_A})\E(Y_B|X_{W_B})Y_CY_D\E(Y_E|X(t_1^*),X(t_2^*))\1_T|\F_{\sigma'}(t_1^*)\vee \cG_{\sigma'}(t_2^*)]), \nonumber\\ 
=&\E(\E[\E(Y_A|X_{W_A})Y_C|\F_{\sigma'}(t_1^*)\vee \cG_{\sigma'}(t_2^*)]\E(Y_B|X_{W_B})Y_D\E(Y_E|X(t_1^*),X(t_2^*))\1_T). \label{sp_dk:2}
\end{align}
Using, the spine-Markovian property of $B$ with respect to $\sigma'$ and the edge-Markovian
property of $E^*$, respectively we get
\begin{align}
 =&\E(\E[\E(Y_A|X_{W_A})Y_C|\F_{\sigma'}(t_1^*)\vee \cG_{\sigma'}(t_2^*)]Y_BY_D\E(Y_E|X(t_1^*),X(t_2^*))\1_T) \nonumber \\
=&\E(\E(\E(Y_A|X_{W_A})Y_C|\F_{\sigma'}(t_1^*)\vee \cG_{\sigma'}(t_2^*))Y_BY_DY_E\1_T) \label{sp_dk:2a}
\end{align}
It remains to show that $\E[\E(Y_A|X_{W_A})Y_C|\F_{\sigma'}(t_1^*)\vee \cG_{\sigma'}(t_2^*)]$ is $\sigma(X_W)$ measurable.
Let $W_A=W_A'\cup W_A^*$, where $W_A'\subset [0,t_1^*]\cup [t_2^*,1]$, and $W_A^*=W_A\setminus W_A'$.
We can assume $X_{W_A}=(X_{W_A'}, X_{W_A^*})$
If $f_{Y_A}(x_{W_A'},x_{W_A^*})=\E(Y_A|X_{W_A}=(x_{W_A'},x_{W_A^*}))$, then 
\begin{align}
 &\E[\E(Y_A|X_{W_A})Y_C|\F_{\sigma'}(t_1^*)\vee \cG_{\sigma'}(t_2^*)] \nonumber\\
=&\E[f(X_{W_A})Y_C|\F_{\sigma'}(t_1^*)\vee \cG_{\sigma'}(t_2^*)] \nonumber\\
=&\int f(X_{W_A'},x_{W_A^*})y_C\P(X_{W_A^*}\in dx_{W_A^*},Y_C\in dy_C |\F_{\sigma'}(t_1^*)\vee \cG_{\sigma'}(t_2^*))\label{sp_dk:3}
\end{align}
Now, since $W_A^*$ and $Y_C$ are $\sigma(X_{\sigma}(t):t\in[t_1^*,t_2^*])$-measurable, using the
(T3) (note that $C$ and $E$ form a simple cell) and Lemma \ref{lem:wptm}, we have
$$\P(X_{W_A^*}\in dx_{W_A^*},Y_C\in dy_C |\F_{\sigma'}(t_1^*)\vee \cG_{\sigma'}(t_2^*))=\P(X_{W_A^*}\in dx_{W_A^*},Y_C\in dy_C |X(t_1^*),X(t_2^*)).$$
This and $(\ref{sp_dk:3})$ implies that $\E[\E(Y_A|X_{W_A})Y_C|\F_{\sigma'}(t_1^*)\vee \cG_{\sigma'}(t_2^*)]$ is $\sigma(X_W)$-measurable,
since $W=W_A'\cup\{t_1^*,t_2^*\}$.
Now taking the conditional expectation in $(\ref{sp_dk:2a})$ with respect to $X_W$ we get
\begin{align*}
 &\E(\E(\E(\E(Y_A|X_{W_A})Y_C|\F_{\sigma'}(t_1^*)\vee \cG_{\sigma'}(t_2^*))Y_BY_DY_E\1_T|X_W))\\
=&\E(\E(\E(Y_A|X_{W_A})Y_C|\F_{\sigma'}(t_1^*)\vee \cG_{\sigma'}(t_2^*))\E(Y_BY_DY_E|X_W)\1_T)
\end{align*}
Using the fact that $X_W$ is $\F_{\sigma'}(t_1^*)\vee \cG_{\sigma'}(t_2^*)$-measurable, we have
\begin{align*}
 &\E(\E(\E(Y_A|X_{W_A})Y_C\E(Y_BY_DY_E|X_W)\1_T|\F_{\sigma'}(t_1^*)\vee \cG_{\sigma'}(t_2^*)))\\
=&\E(\E(Y_A|X_{W_A})Y_C\E(Y_BY_DY_E|X_W)\1_T)
\end{align*}
Applying the spine-Markovian property to $A$ with respect to $\sigma'$ we get
\begin{align*}
 &\E(\E(Y_A|X_{W_A})Y_C\E(Y_BY_DY_E|X_W)\1_T)\\
=&\E(Y_AY_C\E(Y_BY_DY_E|X_W)\1_T).
\end{align*}
Finally, taking the conditional expectation with respect to $X_W$
we get
\begin{align*}
 &\E(\E(Y_AY_C\E(Y_BY_DY_E|X_W)\1_T|X_W))\\
 =&\E(\E(Y_AY_C|X_W)\E(Y_BY_DY_E|X_W)\1_T).
\end{align*}
From the Monotone Class Theorem the claim follows.
\end{proof}

\subsection{The constructed process is a hereditary spine-Markovian process}\label{sec:hsmrk}\index{Process indexed by a TLG!spine-Markovian property!hereditary|(}
Recall how we defined $S^*(\cG)$ and the hereditary spine-Markovian property. 
(See Definition \ref{def:S*} and Definition \ref{def:her_mp}. on page \pageref{def:her_mp}.)

\prop{\label{prop:hrdspm}The process $X$ on $\cG$ defined as in \ref{constr}. is hereditary spine-Markovian.
\begin{proof}
Fix an arbitrary TLG$^*$ $\cG$ and an TLG$^*$ tower $(\cG_k)_{k=0}^n$ such that $\cG_0$ is the minimal
graph and $\cG_n=\cG$.

Clearly, $X^0$ is spine-Markovian, and the claim holds since $S^*(\cG_0)=\{\cG_0\}$.
Now, we will show that if the process $X^{k-1}$ on $\cG_{k-1}\neq \cG$ is hereditary spine-Markovian,
so is $X^{k}$ on $\cG_{k}$.

$\bullet$ If we got $\G_k$ by adding a new vertex to $\cG_{k-1}$ then we are done, since the
distribution of the process didn't change on the joint representation of these two TLG$^*$.\vspace{0.3cm}

$\bullet$ Let's view the case when we added a new edge $E^*$ (between the existing vertices)
to $\cG_{k-1}$ to obtain $\cG_k$. Take any $\cH\in S^*(\cG_{k})$. If $\cH\in S^*(\cG_{k-1})$,
then we are done. Otherwise, $\cH=(\cV_H,\cE_H)$ contains the new edge $E^*$, i.e. $E^*\in \cE_H$.
Let $E^*=E_{t_1^*t_2^*}$. \vspace{0.2cm}

$(\circ 1)$ If there exists a path $\sigma_{12}$ connecting $t_1^*$
and $t_2^*$ (not containing $E^*$), then $\cH'=(\cV_H,\cE_H\setminus\{E^*\})$ is a TLG$^*$ (Corollary \ref{pos:sped}.) and in $S^*(\cG_{k-1})$.
This implies that $(X^k(t):t\in \cH')$ is spine-Markovian, and in the same way as in the Lemma \ref{lem:sp-mk},
we can show that  $(X^k(t):t\in \cH)$ is spine-Markovian.\vspace{0.2cm}

$(\circ 2)$ If a path $\sigma_{12}$ connecting $t_1^*$ and $t_2^*$ does not exist, then
take any tower $(\cK_l)_{l=0}^n$ such that $\cK_0=\cH$ and $\cK_n=\cG$.

Let $k$ be a minimum $l$ such $t_1^*$ and $t_2^*$ are connected in $K_l$ by some path
not containing $E^*$. Such a $k$ exists, because for the construction of $E^*$
$t_1^*$ and $t_2^*$ need to be connected by a time-path in $\cG_{k-1}$, so this is
also true in $\cG_k=\cK_n$. But then, we just added a new edge $E_{t_1^*t_2^*}'$ 
to $\cK_{m-1}$. Now we can  first add an edge 
$E_{t_1^*t_2^*}'$ to $\cK_0$, and after that add vertices and edges in the order we added them to obtain $\cK_{m-1}$
from $\cK_0$. 

In this way, we would still get $\cK_m$ at the end. This shows that a TLG$^*$ 
$(\cV_H,\cE\cup\{E'_{t_1^*t_2^*}\})$ 
(the TLG$^*$ that we get when we add a new edge connecting $t_1^*$ and $t_2^*$ to $\cK_0$) is in $S^*(\cG_k)$.

Now, we are previous case $(\circ 1)$: $X^k$ on $(\cV_H,\cE\cup\{E'_{t_1^*t_2^*}\})$ is spine-Markovian.

To prove that $X^k$ on $\cH$ is spine-Markovian we need to consider two cases:
If a spine $\sigma$ in $\cH$ contains $E^*$, then $E'_{t_1^*t_2^*}$ is just
one of the components (disjoint from others) in $(\cV_H,\cE\cup\{E'_{t_1^*t_2^*}\})$
with respect to $\sigma$. For any other spine $\sigma$ not containing $E^*$, since
$E'_{t_1^*t_2^*}$ will be an extra part of some component in $(\cV_H,\cE\cup\{E'_{t_1^*t_2^*}\})$
with respect to $\sigma$. This shows that $X^k$ on $\cH$ is spine-Markovian. 
\end{proof}
}
\index{Process indexed by a TLG!spine-Markovian property!hereditary|)}\index{Process indexed by a TLG!spine-Markovian property|)}
\index{Spine-Markovian property|seealso{Process indexed by a TLG}}
\subsection{Uniqueness in law of hereditary spine-Markovian $\cM$-processes}\index{Process indexed by a TLG!uniqueness of distribution|(}
Lemmas \ref{lem:P-proc} and \ref{lem:sp-mk} give the following proposition.

\prop{The process $X$ on $\cG$ defined in \ref{constr} is a hereditary spine-Markovian $\cM$-process.
}\vspace{0.3cm}

We will finish this discussion by showing uniqueness in law of hereditary spine-Markovian $\cM$-processes.\vspace{0.3cm}

As we noticed in the Remark after the Definition \ref{def:constr}, the definition of the 
process $X$ on $\cG$ depends on the choice of the TLG$^*$ tower, on which we
inductively define the process.  It turns out, that the distribution of the
process $X$ is unique, and therefore it doesn't depend on the choice of the TLG$^*$ 
tower.

First, let's prove the following lemma. 
\lem{\label{lem:uniq}Let $X$ be a hereditary spine-Markovian $\cM$-process on a TLG$^*$ $\cG$.
If $\cG$ can be obtained from a TLG$^*$ $\cG'$, by adding a new edge or vertex as in Definition \ref{def:tlg*}.(ii),
then $X'$ a restriction of $X$ to $\cG'$ is also a hereditary spine-Markovian $\cM$-process.
\dok{Any full-time path in $\cG'$ is also a full-time path in $\cG$. 
Since, $S^*(\cG')\subseteq S^*(\cG)$, it is clear that $(X(t):t\in \cG')$
is hereditary spine-Markovian.}}

\teo{\label{teo:uniq}A hereditary spine-Markovian $\cM$-process (satisfying (3T) properties) on a TLG$^*$ $\cG$ has a unique 
distribution.
\dok{We will prove this using the induction on the number of edges $n$ of the TLG$^*$.\vspace{0.2cm} 

For $n=1$, we have a minimal graph and its distribution is clearly uniquely given.\vspace{0.2cm} 

For $n>1$, suppose $\cG$ can be obtained from $\cG'$  by adding a new edge or vertex as in Definition \ref{def:tlg*}.(ii).
If we just added a vertex to $\cG'$ in order to obtain $\cG$, then we are
done since these two graphs have the same representation $R(\cG)=R(\cG')$.
Since $\cG'$ has $n-1$ edge, the distribution on it is unique, and so is
on $\cG$.\vspace{0.2cm} 

If we added a new edge between the existing $E^*$ two vertices $t_1$ and $t_2$ on $\cG'$.
We are done since, there has to exist a full time-path $\sigma$ in $\cG'$  
containing $t_1$ and $t_2$. But now, $\cG_-=E^*$ is a component in the 
decomposition of $\cG$ with 
the respect to $\sigma$ and $\cG_+=\cG'$ is the rest of $\cG$. 
Now, the processes
$(X(t):t\in \cG')$ and  $(X(t):t\in E^*)$ are independent given $X(t_1)$ and $X(t_2)$.
By Lemma \ref{lem:uniq} $(X(t):t\in \cG')$ is a  hereditary spine-Markovian $\cM$-process, so its 
distribution is unique. 
The distribution of $(X(t):t\in E^*)$ given $X(t_1)$ and $X(t_1)$
is also uniquely given because of the consistency (i.e. (T2) property) of $\cM$.\vspace{0.2cm} 

Hence, the distribution of $X$ on $\cG$ is unique.}}

\defi{\label{def:nmprc}We define the process constructed in \S\ref{constr} to be the natural $\cM$-process on the TLG$^*$ $\cG$.}\index{Process index by a TLG!$\cM$-process!natural|textbf}
\index{Process indexed by a TLG!uniqueness of distribution|)}

\section{Processes on TLG's with infinite number of vertices}
\index{Process indexed by a TLG!with infinite number of vertices|(}
In Section \ref{inf_vert} (see Definition \ref{def:intlg*}) we introduced
TLG's and TLG$^*$ with infinitely many vertices. As in the case where 
we had only a finite number of vertices, here also we will construct 
a process on TLG$^*$ graphs.

\subsection{Construction}\index{Process indexed by a TLG!with infinite number of vertices!construction}

Let $\cG=(\cV,\cE)$ a TLG$^*$ such that $\cV$ is infinite. According to the definition,
there exists a tower of TLG$^*$'s $\cG_n=(\cV_n,\cE_n)$, $n\geq 1$, such that 
$\cV_n$ is finite, where $\cV=\bigcup_{n\geq 1}\cV_n$.\vspace{0.3cm}

Let 
\begin{equation}
\cM=\{\mu_{\sigma}: \sigma\in P_{0\to 1}(\cG)\} \label{kf:1} 
\end{equation}

be a family of distributions of Markov processes along
full-time paths in $\cG$ satisfying conditions (T1)-(T3) given in Section \ref{condit}. 
(Although $0$ and $1$ don't have to be the start and the end
of time in $\cG$, we will still use the notation $P_{0\to 1}(\cG)$ for full-time paths
in $\cG$.)\vspace{0.3cm}

Since $$\cM(\cG_n)=\{\mu_{\sigma}: \sigma\in P_{0\to 1}(\cG_n)\}$$ 
is well-defined, and we can show similarly as in Lemma \ref{lem:inr:A13} 
that $\cM(\cG_n)$ satisfies (T1)-(T3), we can define a hereditary spine-Markovian process $X^n$ on $\cG_n$, such that
for each $\sigma\in P_{0\to 1}(\cG_n)$ the process $X^n_{\sigma}$ has the distribution $\mu_{\sigma}$. 
Further, 
the restriction of this process to $\cG_k$ ($k\leq n$)  has the same distribution as
the $\cM(\cG_k)$-process $X^k$ defined on $\cG_k$ in the similar manner.\vspace{0.3cm}

Now, Kolomogorov's consistency theorem shows, that there
exists a process $X$ on $\cG$ such that the restriction of $X$ to any $\cG_k$
has same distribution as $X^k$. Note, that since each $\sigma\in P_{0\to 1}(\cG)$ is in some of the
$\cG_k$'s we have $X_\sigma$ has the distribution $\mu_{\sigma}$.

\subsection{Uniqueness of the distribution}\label{subsec:undistr}\index{Process indexed by a TLG!with infinite number of vertices!uniqueness of distribution}

\lem{\label{lem:2seqcrt}
Let $\cG_0$, $\cH$ and $\cG_1$ be TLG$^*$'s with the following
properties:
\begin{enumerate}[(1)]
 \item $\cG_0\in S^*(\cG_{1})$;
 \item $\cV_{\cG_{0}}\subset \cV_{\cH}\subset \cV_{\cG_{1}}$;
 \item $R(\cG_0)\subset R(\cH)\subset R(\cG_{1})$.
\end{enumerate}
Then $\cG_{0}\in S^*(\cH)$.

\begin{figure}[ht]
\begin{center}
\psfrag{0}{\textcolor{blue}{$\boldsymbol{\cG_0}$}}
\psfrag{1}{$\boldsymbol{\cG_1}$}
\psfrag{h}{$\boldsymbol{\cH}$}
\includegraphics[width=7cm]{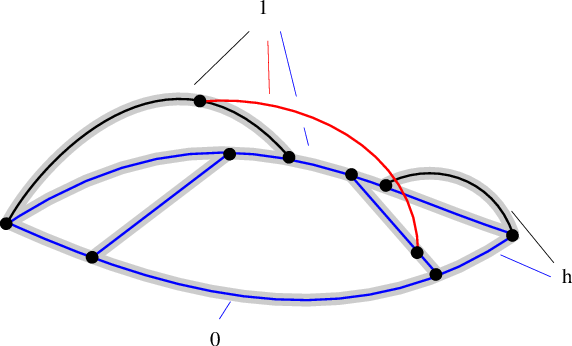}\\ 
\caption{} \label{sl47}

  \end{center}
\end{figure}

\dok{

We will show the claim by induction on 
$$n(\cG_0)=\sum_{v\in \cV_{\cG_0}, d(v)\geq 3}d(v). $$
If $n=0$ the claim is $\cG_0$ represents a spine in $\cH_0$, 
and the claim follows by Theorem \ref{thm:spine}.\vspace{0.2cm}

Assume the claim holds for $n\leq k$. We will prove the claim 
for $n=k+1$. Pick a truly simple cell $(\pi_{uv},\pi^1_{uv})$ (recall Definition \ref{defi:msct}) in $\cG_0$ (for example 
pick a spine $\pi$ and then a time path $\pi^1_{uv}$ not contained in $\pi$ connecting $t_u$ and $t_v$
such that $|t_u-t_v|$ is minimal), the representation of this cell will 
remain a truly simple cell in $\cG_1$ (by Theorem \ref{pre-path-con}.) and 
therefore also in $\cH$.\vspace{0.2cm}

Hence,  we define a cell collapsing transformation $\circ$ that is collapsing this cell.\vspace{0.2cm}

For the TLG$^*$'s $\cG_0^{\circ}$, $\cH^{\circ}$ and $\cG_1^{\circ}$ property (1) holds by Corollary \ref{circ-tlgt},
while (2) and (3) are clear.
Now, since $n(\cG_0^{\circ})<n(\cG_0)$, by induction assumption $\cG_0^{\circ}\in S^*(\cH^\circ)$. 

We follow the construction from $\cG_0^\circ$ to $\cH^{\circ}$, to obtain a
TLG$^*$-tower going from $\cG_0$ to $\cH$. Let $(\cK_j')_{j=0}^n$ be TLG$^*$-tower starting with
$\cK_0'=\cG_0^{\circ}$ and $\cK_n'=\cH^{\circ}$. Now we construct a TLG$^*$-tower
$(\cK_l)$ staring with $\cK_0=\cG_0$. The idea of the construction is the following:
if  on $\cK_j'$ to obtain $\cK_{j+1}'$ we added
\begin{itemize}
 \item a vertex, then add an appropriate vertex to $\cK_j$ to obtain $\cK_{j+1}$;
 \item an edge, then connect two appropriate vertices in $\cK_j$ by an edge to obtain $\cK_{j+1}$.
\end{itemize}

The main question is: When we add an edge, are we connecting two vertices that are 
connected by a time-path? That means that in $K_j$ two vertices $t_k$ and $t_l$ are 
not connected by a time path, but $(t_k)^{\circ}$ and  $(t_l)^{\circ}$ are connected by
a time-path in $\cK_j'$.
So we have a situation like on the Figure \ref{pic32:2}. (Other situations are similar.)

\begin{figure}[ht]
\begin{center}
\psfrag{1}{$\boldsymbol{t_k}$} \psfrag{2}{$\boldsymbol{t_{k'}}$}
\psfrag{5}{$\boldsymbol{t_u}$} \psfrag{3}{$\boldsymbol{t_{v}}$}
\psfrag{6}{$\boldsymbol{t_l}$} \psfrag{4}{$\boldsymbol{t_{l'}}$}
\psfrag{7}{$\boldsymbol{t_{k'} \wedge t_{l}}$}
\psfrag{8}{$\boldsymbol{(t_{k'} \wedge t_{l})\vee t_{l'}}$}
\includegraphics[width=11cm]{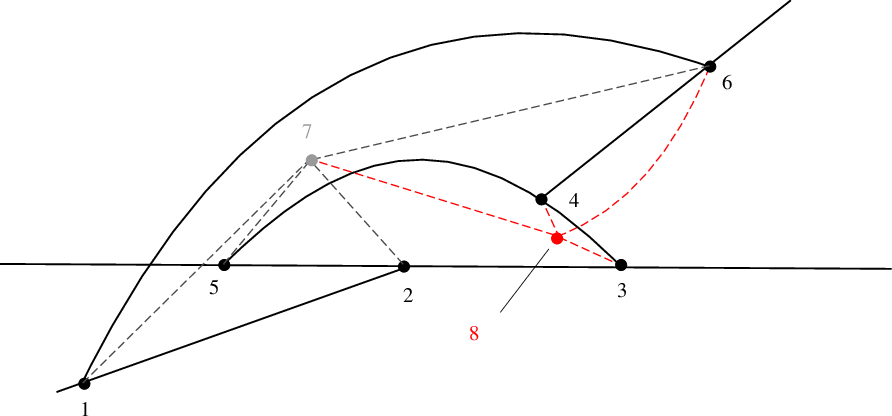}\\ 
\caption{} \label{pic32:2}

  \end{center}
\end{figure}

By Theorem \ref{teo:tpltt} $\cH$ is a topological lattice, there exists a unique vertex
$t_{u'}:=t_{k'} \wedge t_{l}$. By definition since $t_u$ and $t_k$ are in the past of 
$t_l$ and $t_k'$, $t_u$ and $t_k$ are connected by a time-path to $t_{u'}$.

Using this property, again, we know that 
in $\cH$ there exists $t_{v'}=t_{u'}\vee t_{l'}$. And know by the same argumentation
$t_{v'}$ is connected by a time path to $t_l$ and $t_v$.

Note that $t_{u'}$ and $t_{v'}$ are in the time frame $[t_u,t_v]$. 
In order for the cell $(\pi_{uv},\pi^1_{uv})$ to remain truly simple,
$t_{u'}=t_u$ or $t_{v'}=t_v$ (otherwise the path $t_{k'}-t_{u'}-t_{v'}-t_{l'}$
will go from one side of the cell to the other within time frame $[t_u,t_v]$).\vspace{0.2cm}

But, since $t_v$ and $t_l$ or $t_u$ and $t_k$ are not connected  by a time path
in $\cK_j$ (since $t_k$ and $t_l$ are not), it follows that their images under the transformation 
are not connected in $\cK_j'$. Hence this is a contradiction.\vspace{0.2cm}

Therefore, in our procedure we construct a TLG$^*$-tower.}
}
 
\begin{figure}[ht]
\begin{center}
\includegraphics[width=7cm]{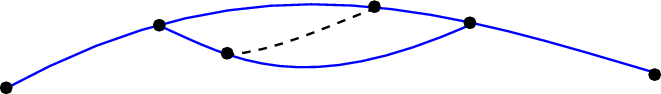}\\ 
\caption{} \label{sl43}

  \end{center}
\end{figure}

\noindent {\it Remark}. The conditions (2) and (3) are not sufficient to imply the conclusion 
of the Lemma. The example is given on Figure \ref{sl43}. The whole line graph with vertices, and the whole graph
are TLG$^*$'s (since they are planar), but we can't construct the second from the first, since a 
simple cell is not a simple cell in the second.

\lem{Let $\cG$ be a TLG$^*$ with infinitely many vertices and $(\cG_j^1)$ and $(\cG_j^2)$ two TLG$^*$-towers 
that construct $\cG$. For any points $\tau_1\prec \tau_2$ on $\cG_1$ with finite times and $j_1\geq 1$, 
the distribution of the natural $\cM$-processes $X^1$ and $X^2$ restricted on $R(\cG_{j_1}^1[\tau_1,\tau_2])$
is unique.
\dok{First, we know that $\cG_{j}^{h}[\tau_1,\tau_2]$ ($h=1,2$) is a TLG$^*$ (see Theorem \ref{teo:inttgl}), also note that
all of its vertices have finite time. \vspace{0.2cm}

By Lemma \ref{lem:subgphfcv}, there is a $\cG_{k_1}^2$ such that $R(\cG_{j_1}[\tau_1,\tau_2])\subset R(\cG_{k_1}^2)$, and $\cV_{j_1}^1\subset \cV_{k_1}^2$.
Further, note that  $R(\cG_{j_1}^1[\tau_1,\tau_2])\subset R(\cG_{k_1}^2[\tau_1,\tau_2])$. Using the same idea, we can find 
$j_2$ such that $R(\cG_{k_1}^2[\tau_1,\tau_2]) \subset R(\cG_{j_2}^1[\tau_1,\tau_2])$ and $\cV_{k_1}^2\subset \cV_{j_2}^1$.
In this way  $\cG_{j_1}^1[\tau_1,\tau_2]$, $\cG_{k_1}^2[\tau_1,\tau_2]$ and $\cG_{j_2}^1[\tau_1,\tau_2]$ satisfy the properties
of Lemma \ref{lem:2seqcrt}. Therefore, we can construct 
$\cG_{k_2}^2[\tau_1,\tau_2]$ from $\cG_{j_1}^1[\tau_1,\tau_2]$.\vspace{0.2cm}

By Corollary \ref{cor:tlg*cnstint},
we can construct a spine $\pi$ going through $\tau_1$ and $\tau_2$, then 
$\cG_j^{h}[\tau_1,\tau_2]$ ($h=1,2$) on that spine, and after that the rest of $\cG_j^{h}$. Since, 
$\cM(\cG_{j}^{h}[\tau_1,\tau_2])$ --
the restriction of the family $\cM$ on $\cG_{j}^{h}[\tau_1,\tau_2]$, is a (3T) family, 
$X^{h}$ restricted on $G_{j}^{h}[\tau_1,\tau_2]$ is a natural $\cM(\cG_{j}^{h}[\tau_1,\tau_2])$-process.\vspace{0.2cm}

Hence, $X^2$ on $\cG_{k_1}^2([\tau_1,\tau_2])$ is distributed as a natural $\cM(\cG_{k_1}^2[\tau_1,\tau_2])$-process.
Since $\cG_{k_1}^2[\tau_1,\tau_2])$ can be constructed from $\cG_{j_1}^1[\tau_1,\tau_2]$, 
$X^2$ restricted on $R(\cG_{j_1}^1[\tau_1,\tau_2])$ is a natural $\cM(\cG_{j_1}^1[\tau_1,\tau_2])$-process. 
Therefore, $X^2$ has the same distribution as $X^1$ on $R(\cG_{j_1}^1[\tau_1,\tau_2])$.}}

 Burdzy and Pal were able to prove the uniqueness only in the case of planar NCC TLG's with 
infinite vertex set. The following proves their conjecture (see the sentence before Theorem 3.9. in \cite{tlg1})
that this is true in general case (including the non-planar case).

\teo{\label{teo:burdzy-pal}Let $\cG=(\cV,\cE)$ be a TLG$^*$'s with infinitely many vertices in $\cV$,
and let $X^1$ and $X^2$ be two $\M$-processes constructed using the TLG$^*$-towers 
$(\cG_n^1)$ and $(\cG_n^2)$, then  $X^1$ and $X^2$ have the same distribution.
\dok{Pick points $\tau_1^{(n)}\prec \tau_2^{(n)}$ on $\cG_1^1$ with finite time such that 
$\tau_1^{(n)}\downarrow -\infty$ and $\tau_2^{(n)}\uparrow +\infty$ (in time). 
Now, the distributions of $X^1$ and $X^2$ on $R(\cG_n[\tau_1^{(n)},\tau_2^{(n)}])$
are the same, and since $$\bigcup_{n=1}^{\infty}R(\cG_n[\tau_1^{(n)},\tau_2^{(n)}])=R(\cG),$$
by Kolmogorov's consistency theorem we have that $X^1$ and $X^2$ have the same distribution.}}
\noindent {\it Remark.} To use the Kolmogorov's consistency theorem we need to look at finite dimensional 
vectors $(X^1(t_1),\ldots, X^1(t_m))$ and $(X^2(t_1),\ldots, X^2(t_m))$ for a finite 
number of points $t_1,\ldots, t_m \in \cG$ with finite time. Since each point is in some subgraph of $\cG$, 
there exists a $n$ such that $$\{t_1,\ldots, t_m\}\subset \cG_n[\tau_1^{(n)},\tau_2^{(n)}]$$
and hence the random vectors have the same distribution.
\index{Process indexed by a TLG!with infinite number of vertices|)}

\chapter{Markov properties of processes indexed by TLG's}\label{sec:3a}
From \S\ref{sec:spine2}. and \S\ref{sec:hsmrk}. we know that the constructed process 
has a (hereditary) spine-Markovian property\index{Spine-Markovian property}. This property is induced by the graph 
structure and as we will see there is one more property this process has when 
$\cM$ is a (3T)-family. If $\cM$ has some additional properties we will have some 
additional properties of the process on the TLG$^*$ $\cG$. 

\index{Cell-Markovian property|(}\section{Cell-Markov properties}
Recall, truly simple cell has been defined in Definition \ref{defi:msct}.

\defi{\label{def:cllmk}We will say that a process $X$ on a TLG $\cG$ is \textbf{cell-Markovian} if for any truly simple cell $(\sigma_1,\sigma_2)$
starting at $t_*$ and ending at $t^*$ the processes $X_{\sigma_1}$ and $X_{\sigma_2}$
are conditionally independent, given the values $X(t_*)$ and $X(t^*)$.}

\defi{\label{def:scllmk}We will say that a process $X$ on a TLG $\cG$ is \textbf{strong cell-Markovian}\index{Cell-Markovian property!strong} if for for any truly simple cell $(\sigma_1,\sigma_2)$
starting at $t_*$ and ending at $t^*$ the processes is cell-Markovian and 
$(X(t):t\in \cG[t_*,t^*])$ and $(X(t):t\in \cG[0,t_*]\cup \cG[t^*,1])$ are independent, given the values $X(t_*)$ and $X(t^*)$.}




\begin{figure}[ht]
\begin{center}
\psfrag{a}{$\boldsymbol{t_*}$}
\psfrag{b}{$\boldsymbol{t^*}$}
\psfrag{1}{$\boldsymbol{1}$}
\psfrag{0}{$\boldsymbol{0}$}
\includegraphics[width=8cm]{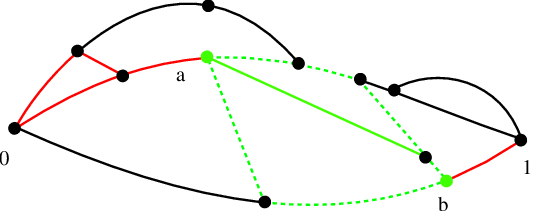}\\
\caption{Strong cell-Markovian property: $X_{\textcolor{green}{\cG[t_*,t^*]}}\perp X_{\textcolor{red}{\cG[0,t_*]\cup \cG[t^*,1]}} | (X(t_*),X(t^*))$ } \label{pic:sl51}
  \end{center}
\end{figure}

Before we prove the that the strong cell-Markovian property holds,
we will prove the following lemma. 

\lem{\label{lem:ccind}Let $T=A\cup B\cup C\cup D$, and $X=(X(t):t\in T)$ a stochastic process.
If \begin{enumerate}[(1)]
    \item $X_A=(X(t):t\in A)$ and $X_C=(X(t):t\in C)$ are independent given $X_B=(X(t):t\in B)$
    \item for some subset $C'\subset C$ $X_{A\cup B\cup C}$ and $X_D$ are independent
given $X_{C'}$ 
   \end{enumerate}
then $X_A$ and $X_{C\cup D}$ are independent given $X_B$.
\dok{Let $Y_S$ be a bounded $\sigma(X_S)$-measurable  function, for $S=A,B,C,D$, and 
$U$ be an element in $\sigma(X_B)$. Then, using (2) we have
$$ \E(Y_AY_BY_CY_D\1_U)= \E(Y_AY_C\E(Y_D|X_{C'})\1_U).$$
Using (1) we get 
$$\E(Y_AY_C\E(Y_D|X_{C'})\1_U)=\E(\E(Y_A|X_B)Y_C\E(Y_D|X_{C'})\1_U),$$
and again using (2) and the fact that $\E(Y_A|X_B)Y_C\1_U$ is a bounded $\sigma(X_{A\cup B\cup C})$-measurable
random variable we get 
 $$\E(\E(Y_A|X_B)Y_C\E(Y_D|X_{C'})\1_U)=\E(\E(Y_A|X_B)Y_CY_D\1_U).$$
Finally, conditioning everything (under the expectation) on $X_B$ we get 
 $$\E(\E(Y_A|X_B)Y_CY_D\1_U)= \E(\E[\E(Y_A|X_B)Y_CY_D\1_U|X_B])=\E(\E(Y_A|X_B)\E(Y_CY_D|X_B)\1_U).$$
Now, using the Monotone Class Theorem the claim follows.
}} 

\teo{\label{thm:strcell}The process $X$ on $\cG$ defined in \S\ref{constr} is strong cell-Markovian process.}
\dok{It is known from Corollary \ref{cor:tlg*cnstint} that there exists a TLG$^*$-tower $(\cG_k)_{k=0}^n$ 
that starts with a spine $\pi$ containing $t_*,t^*$, there exists $n_0$ such that 
$\cG_{n_0}$ such that $R(\cG[t_*,t^*])\cup R(\pi)=R(\cG_{n_0})$, and then we can construct 
the rest of $\cG$ (i.e. $\cG_n=\cG$). Since, $\cG_{n_0}$ is a TLG$^*$, we know by Theorem 
that $X^{n_0}$ the natural $\cM(\cG_{n_0})$-process on $\cG_{n_0}$ is the same as the restriction
of the process $X$ on $\cG_{n_0}$.\vspace{0.2cm} 

Assume that $\pi$ is the spine that contains $\sigma_1$. 
Since $\sigma_2$ will in a decomposition component $\cG_-$ with roots $t_*$
and $t^*$,  by the spine-Markovian property, $X_{\sigma_1}^{n_0}$ is independent of $X_{\pi}^{n_0}$ 
given $X^{n_0}(t_*)$ and $X^{n_0}(t^*)$. This proves the cell-Markovian property.

We use induction to show that $(X(t):t\in \cG_k[t_*,t^*])$ is independent of 
$(X(t):t\in \cG_{k}[0,t_*]\cup \cG_{k}[t^*,1])$. For $n=0$ the claim 
follows from (T3) property. For $k=1,\ldots, n_0$ the process on every edge that we add 
will depend only on the value of the process $(X(t):t\in \cG_{k-1}[t_*,t^*])$
at its endpoints, so the claim will follow by Lemma \ref{lem:ccind}. For $k>n_0$ we have the following
cases:
\begin{itemize}
 \item We added an vertex - nothing changes since the representation is the same.
 \item We added an edge not in $\cG[0,t_*]\cup \cG[t^*,1])$ - this has no impact.
 \item We added an edge $E$ that connects two vertices in $\cG_{k-1}[0,t_*]\cup \cG_{k-1}[t^*,1])$. Then
the process depends only on the values of $X$ at the endpoints, the claim is true by Lemma \ref{lem:ccind}.
\end{itemize}
Since the distribution of the process, by Theorem \ref{teo:uniq}, doesn't depend on the construction
the claim follows.}

\pos{\label{pos:t-mk}For the process $X$ on $\cG$ defined in \S\ref{constr}, if $(\sigma_1,\sigma_2)$ is a truly simple cell
starting at $0$ and ending at $t^*$, then the processes $(X(t):t\in \cG[0,t^*])$ and $(X(t):t\in \cG[t^*,1])$ 
are independent given the values of $X(0)$ and $X(t^*)$.}

\index{Cell-Markovian property|)}

\section{Graph-Markovian and time-Markovian property}\label{grph_tm_mrk}


First, we introduce the graph-Markovian property\index{Graph-Markovian property|textbf}, a version of the global Markov property\index{Global Markov property} in graphical models (see Definition 
\ref{def:gr_mpr} (c)). 

\defi{\label{def:grpMrk}Suppose that $W\subset R(\cG)$ is a finite non-empty set such that $R(\cG)\setminus W$
is disconnected. Some edges of $\cG$ are cut by $W$ into two or more components. Let
us call this new collection of edges $\cE_0$. Suppose that $\cE_1$  and $\cE_2$ are disjoint sets
of edges with the union equal $\cE_0$. We will call a process $X$ on a TLG  graph $\cG$ a
\textbf{graph-Markovian process} if for all $W$, $\cE_1$, $\cE_2$, the conditional distribution
of $(X_t:t\in E,E\in \cE_1)$ given $(X_t:t\in E,E\in \cE_2)$ depends only on $(X_t:t\in W)$. }

\begin{figure}[ht]
\begin{center}
\begin{minipage}{7cm}
\psfrag{1}{$\boldsymbol{\cE_1}$}
\psfrag{2}{$\boldsymbol{\cE_2}$}
\includegraphics[width=7cm]{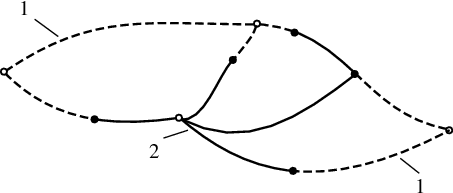}\\
\caption{Graph-Markovian property: Black points ($\bullet$) represent $W$, dashed 
lines $\cE_1$, and full lines $\cE_2$. } \label{pic:gr_mk}

\end{minipage}\ \ \begin{minipage}{7cm}
\psfrag{1}{$\boldsymbol{1}$}
\psfrag{0}{$\boldsymbol{0}$}
\psfrag{t}{$\boldsymbol{t}$}
\psfrag{P}{$\boldsymbol{P(t)}$}
\psfrag{F}{$\boldsymbol{F(t)}$}
\includegraphics[width=7cm]{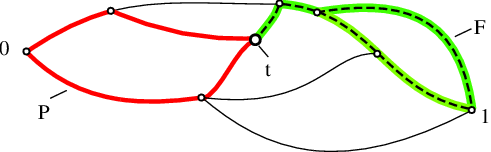}\\
\caption{Time-Markovian property: The past $P(t)$ - of $t$, and the future - $F(t)$ of $t$.} \label{pic:tm_mk}
\end{minipage}

  \end{center}
\end{figure}


The second property is the time-Markovian property.\index{Time-Markovian property|textbf}

\defi{\label{def:t-mrkov}(a) Let $t$ be a point in $\cG$. 
\begin{enumerate}[(i)]
       \item ({\sc the future}) $F(t)=\{s\in \cG: s\succeq  t\}$  is the set of all 
points with times $s\geq t$, such that there is a full path passing through $t$ and $s$.
\item ({\sc the past}) $P(t)=\{s\in \cG: s\preceq  t\}$  is the set of all 
points with times $s\leq t$, such that there is a full path passing through $t$ and $s$.
      \end{enumerate}


(b) We will say a process $X$ on a TLG  graph $\cG$ is a
\textbf{time-Markovian process} if for every $t$, if the conditional distributions
of $(X(s):s\in P(t))$ and $(X(s):s\in F(t))$ given $X(t)$  are independent.
 }

\noindent\emph{Remark.} Note that if $(X(t):t\in \cG)$ is time-Markovian, then for every 
full time-path $\pi$ the process $X_{\pi}=(X(t):t\in \pi)$ is a Markov process.


\section{Processes on TLG's for Markov family $\cM$}


Some additional properties will hold if the distributions in the family
$\cM$ are all distributions of Markov processes.\vspace{0.2cm}

Note that in this case the property (T3) is automatically satisfied,
so the only thing that we need for the construction is the fact 
that $\cM$ is a consistent family of distributions of Markov processes
that are continuous or RCLL (or any other that we can define conditional
distributions on) on a TLG$^*$ $\cG$.\vspace{0.3cm}

In the next few subsections we will show that in this case we have additional properties
- edge-Markovian and time-Markovian properties.

\subsection{The constructed process is a time-Markovian process}\index{Time-Markovian property|(}
\teo{\label{thm:tm_mk}The process $X$ defined on $\cG$ defined in \S\ref{constr} for a Markov family $\cM$ is a time-Markovian process.}

\begin{figure}[ht]
\begin{center}
\psfrag{1}{$\boldsymbol{1}$}
\psfrag{0}{$\boldsymbol{0}$}
\psfrag{t}{$\boldsymbol{t}$}
\psfrag{e}{$\boldsymbol{t_{-\varepsilon}}$}
\psfrag{G}{$\boldsymbol{\cG_{\varepsilon}[t_{-\varepsilon},t]}$}
\psfrag{P}{\textcolor{red}{$\boldsymbol{P(t)}$}}
\psfrag{F}{$\boldsymbol{F(t)}$}
\includegraphics[width=10cm]{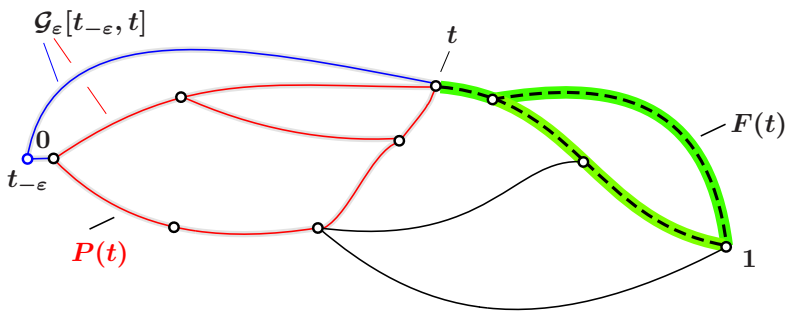}\\
\caption{Idea of the proof of time Markovian property} \label{pic:tm_mk1}
  \end{center}
\end{figure}

\dok{Let $t$ be a point on $\cG$. We can assume it is a vertex in $\cV$. We will expand the vertex set $\cV$, by adding 
the vertex $t_{-\varepsilon}=-\varepsilon$, i.e. $\cV_{\varepsilon}=\cV\cup\{t_{-\varepsilon}\}$.
Further, we will expand the edge set by adding edges connecting $t_{-\varepsilon}$
with $0$ and $t$, i.e. $\cE_{\varepsilon}=\cE\cup\{E_{-\varepsilon,0},E_{-\varepsilon,t}\}$. (See Figure \ref{pic:tm_mk1}.)
It is not hard to see that $\cG_{\varepsilon}=(\cV_{\varepsilon},\cE_{\varepsilon})$ is a TLG$^*$. We define 
$X(-\varepsilon)=0$ and
$X_{E_{-\varepsilon,0}}$ and $X_{E_{-\varepsilon,t}}$ to be
interpolations between the values of the processes at the end points. 
$X_{\cG_{\varepsilon}}$ is a continuous or RCLL process with Markov processes along full time-paths,
and since $X_{\cG}$ is a hereditary time-Markovian, so is $X_{\cG_{\varepsilon}}$.

For any path $\pi$ between $0$ and $t$, $(E_{-\varepsilon,0} \pi ,E_{-\varepsilon,t})$
is a truly simple cell. Now, using the strong cell-Markovian property, we have 
that $(X(t):t\in \cG_{\varepsilon}[-\varepsilon, t]))$ and $(X(t):t\in \cG_{\varepsilon}[t,1]))$ 
are independent given $X(t)$ and $X(-\varepsilon)$. Since $X(-\varepsilon)$ is deterministic, $\cG_{\varepsilon}[t,1]=\cG[t,1]=F(t)$
and $P(t)=\cG[0, t]\subset \cG_{\varepsilon}[-\varepsilon, t])$, the claim follows.}
\index{Time-Markovian property|)}
\subsection{Moralized graph-Markovian property}\index{Graph-Markovian property!moralized|(}

In graphical models when we turn Bayes nets into Markov random fields\index{Markov random field (MRF)},
we \textit{moralize}\index{Moralization} the graph (see \S 4.5 \cite[Koller, Friedman]{prgphmdl}). It turns out that the Markov processes on TLG$^*$'s,
in general, don't satisfy the graph-Markovian property described in Section 
\ref{grph_tm_mrk} (see discussion given in Subsection \ref{not_grph}).

But under the modification of the graph, that we will call \textit{moralization}, 
we will have a similar property.

\defi{Let $\cG=(\cV,\cE)$ be a TLG. The graph $\cG^{\heartsuit}=(\cV^{\heartsuit},\cE^{\heartsuit})$ 
given by $\cV^{\heartsuit}=\cV$ and 
$$\cE^{\heartsuit}=\cE\cup \{E^{\heartsuit}_{ij}\ :\ i\ \textrm{and} \ j\ \textrm{are begining and end of a truly simple cell in}\ \cG\}$$ 
will be called a \textbf{moralized graph}. }

\noindent \emph{Remark.} Note that for a TLG$^*$ $\cG$, $\cG^{\heartsuit}$ is also a 
TLG$^*$ - we are adding edge between points that are connected by a time-path.

\begin{figure}[ht]
\begin{center}
\includegraphics[width=6cm]{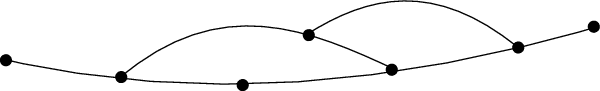}\quad \quad\includegraphics[width=6cm]{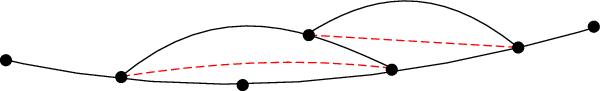}\\
\caption{Moralization of a TLG $\cG$ into $\cG^{\heartsuit}$. } \label{pic:morl}
  \end{center}
\end{figure}

\defi{\label{def:mrlgrpMrk}\index{Graph-Markovian property!moralized|textbf}\index{Moralized graph-Markovian property|see{Graph-Markovian property}}Let $\cG$ be a TLG, and $\cG^{\heartsuit}$ its moralization. 
Suppose that $W\subset R(\cG)\subset R(\cG^{\heartsuit})$ is a finite non-empty set such that $R(\cG^{\heartsuit})\setminus W$
is disconnected. Some edges of $\cG$ are cut by $W$ into two or more components. Let
us call this new collection of edges $\cE_0$. Suppose that $\cE_1$ and $\cE_2$ are disjoint sets
of edges with the union equal $\cE_0$. We will call a process $X$ on a TLG  graph $\cG$ a
\textbf{moralized graph-Markovian process} if for all $W$, $\cE_1$, $\cE_2$, the conditional distribution
of $(X_t:t\in E,E\in \cE_1)$ given $(X_t:t\in E,E\in \cE_2)$ depends only on $(X_t:t\in W)$. } 

\begin{figure}[ht]
\begin{center}
\psfrag{1}{\textcolor{blue}{$\boldsymbol{\cE_1}$}}
\psfrag{2}{\textcolor{magenta}{$\boldsymbol{\cE_2}$}}
\psfrag{W}{\textcolor{green}{$\boldsymbol{W}$}}
\includegraphics[width=11cm]{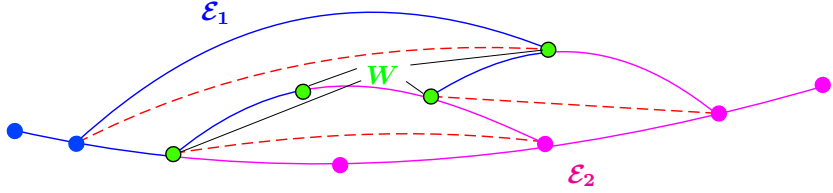}\\
\caption{$(X(t):t\in \cE_1)$ is independent of $(X(t):t\in \cE_2)$ given $X_W$.} \label{pic:49}
  \end{center}
\end{figure}

Before, we prove the moralized graph-Markovian property, we will prove the following lemma.

\lem{\label{lem:sepind}Let $T=A\cup B$, and a stochastic process $X=(X(t):t\in T)$ such that 
\begin{enumerate}[(1)]
 \item there exist $A_1$ and  $A_2$ subsets of $A$ such that
$X_{A_1}=(X(t):t\in A_1)$ is independent of $X_{A_1^c}=(X(t):t\in A\setminus A_1)$
given $X_{A_2}=(X(t):t\in A_2)$;
 \item there exists $A_b$ subset of $A\setminus A_1$ such that $X_{A}$ is independent 
of $X_{B}$ given $X_{A_b}$;
\end{enumerate}
then $X_{B\cup A_1^c}$ is independent of $X_{A_1}$ given $X_{A_2}$.
}
\dok{Let $Y_S$ be a bounded $\sigma(X_{S})$-measurable random variable, and $U\in \sigma(X_{A_2})$.
Now, using (2) we have 
$$\E(Y_BY_{A_1}Y_{A_1^c}\1_U)=\E(\E(Y_B|X_{A_b})Y_{A_1}Y_{A_1^c}\1_U),$$
and using (1) we get 
$$\E(\E(Y_B|X_{A_b})Y_{A_1}Y_{A_1^c}\1_U)=\E(\E(Y_B|X_{A_b})\E(Y_{A_1}|X_{A_2})Y_{A_1^c}\1_U).$$
Using, (2) once more we have 
$$\E(\E(Y_B|X_{A_b})\E(Y_{A_1}|X_{A_2})Y_{A_1^c}\1_U)=\E(Y_B\E(Y_{A_1}|X_{A_2})Y_{A_1^c}\1_U),$$
and now conditioning everything under the expectation on $X_{A_2}$ we get
$$\E(Y_B\E(Y_{A_1}|X_{A_2})Y_{A_1^c}\1_U)=\E(\E(Y_BY_{A_1^c}|X_{A_2})\E(Y_{A_1}|X_{A_2})\1_U).$$
From the Monotone Class Theorem the claim follows.}

\teo{\label{thm:mrl_gp_mk}For a Markov family $\cM$, the natural $\cM$-process on a TLG$^*$ $\cG$
is a moralized graph-Markovian process.
\dok{We use induction on the number of edges $|E|$ for a TLG$^*$ $\cG=(\cV,\cE)$.
For $|E|=1$, the claim is clearly true. Assume that the claim is true for 
$|E|=k\geq 1$. Let's show the claim for $|E|=k+1$. Pick $\cG$ and $W$ a set of points $\cG$,
such that $R(\cG^{\heartsuit})\setminus R(W)$ is disconnected. We need to consider the following cases:
 
If we got $\cG$ by adding a new vertex to some TLG$^*$ $\cH$. In that case, 
since the representation of $\cH$ and $\cG$ is the same, the claim follows.

If we got $\cG$ by adding a new edge $E_*$ between the vertices $t_*$ and $t^*$ in some TLG$^*$ $\cH$, 
we first have to note that $t_*$ and $t^*$ are the begining and the end of a (truly) simple cell whose one side is 
$E_*$. Hence, $t_*$ and $t^*$ are both in one of the following $\cE_1\cup W$ or $\cE_2\cup W$.

We have the following cases to consider:
\begin{itemize}
 \item $R(E_*)\cap R(W)=\emptyset$ then $E_*$ will entierly be in one of 
$\cE_1$ or $\cE_2$. We will assume $E_*\in \cE_1$, and let $\cE_1'=\cE_1\setminus\{E_*\}$.
In 
we use the spine-Markovian property with roots $t_*$ and $t^*$, hence $X_{E_*}$ 
is independent of $X_{\cE_1'\cup \cE_2}$ given $X(t_*)$ and $X(t^*)$. Now since $X_{\cE_2}$
is independent of $X_{\cE_1'}$ given $X_W$, by Lemma \ref{lem:ccind}. it follows that $X_{\cE_1}$
is independent of $X_{\cE_2}$ given $X_W$.


\begin{figure}[ht]
\begin{center}
\psfrag{1}{\textcolor{blue}{$\boldsymbol{\cE_1}$}}
\psfrag{2}{\textcolor{magenta}{$\boldsymbol{\cE_2}$}}
\psfrag{W}{\textcolor{green}{\small $\boldsymbol{W}$}}
\psfrag{E}{$\boldsymbol{E_*}$}
\psfrag{u}{$\boldsymbol{t_*}$}
\psfrag{v}{$\boldsymbol{t^*}$}
\includegraphics[width=12cm]{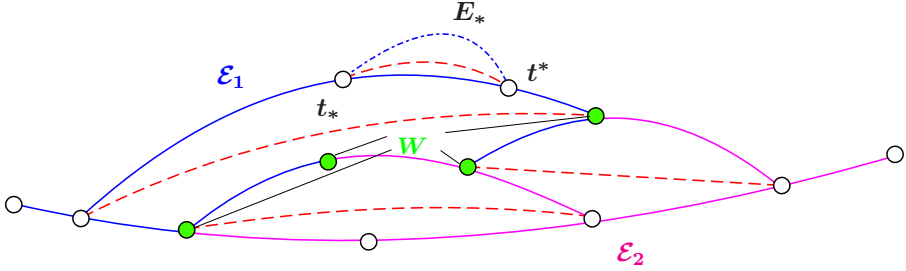}\\
\caption{The new edge $E^*$ doesn't contain points from $W$. } \label{pic:49b}
  \end{center}
\end{figure}

\item $R(E_*)\cap R(W)\neq \emptyset$, then we assume $t_*, t^*\in \cE_1\cup W$.
Denote, $W_{\cH}$ the points represented by $R(W)\cap R(\cH)$. Note that these 
points separate $\cH$.

\begin{figure}[ht]
\begin{center}
\psfrag{1}{\textcolor{blue}{$\boldsymbol{\cE_1}$}}
\psfrag{2}{\textcolor{magenta}{$\boldsymbol{\cE_2}$}}
\psfrag{W}{\textcolor{green}{\small $\boldsymbol{W_{\cH}}$}}
\psfrag{T}{\textcolor{green}{\small $\boldsymbol{W_*}$}}
\psfrag{E}{$\boldsymbol{E_*}$}
\psfrag{a}{\small $\boldsymbol{E_{*,1}^1}$}
\psfrag{b}{\small $\boldsymbol{E_{*,1}^2}$}
\psfrag{c}{\small $\boldsymbol{E_{*,2}^1}$}
\psfrag{d}{\small $\boldsymbol{E_{*,2}^2}$}
\psfrag{u}{$\boldsymbol{t_*}$}
\psfrag{v}{$\boldsymbol{t^*}$}
\includegraphics[width=12cm]{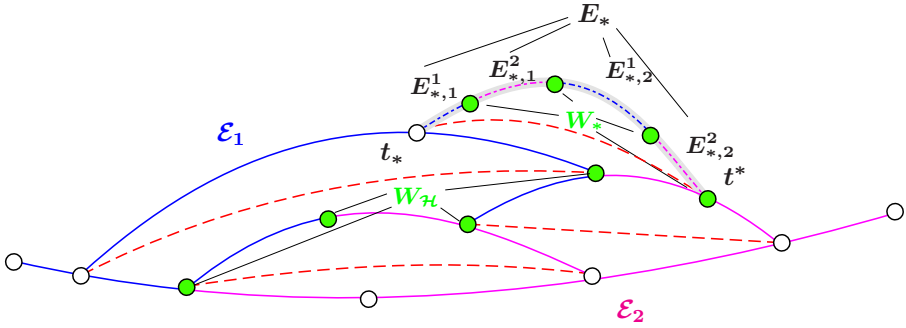}\\
\caption{The new edge $E^*$ contains points from $W$.} \label{pic:49a}
  \end{center}
\end{figure}

Let $E_{*,k}^j$, $k=1,\ldots, n_j$ denotes the edges in $\cE_j$ that cover the edge 
$E_*$. It is not hard to see, since the process along $X_{E^*}$ is Markov 
that 
\be X_{E_{*,k}^j}\perp X_{R(\cG)\setminus R(E_{*,k}^j)}| X_{\partial E_{*,k}^j}.\label{mrl:mp1}\ee 

The endpoints of at least one of the sequences $(E_{*,k}^1)$ or $(E_{*,k}^2)$
will be only in $W$. Otherwise, $t_*$ and $t^*$ won't be in $W$, and
they won't be both in $\cE_1$. Under the assumption that $t_*$ and $t^*$ are
in $\cE_1\cup W$, it follows that $(E_{*,k}^2)$ has all its 
endpoints in $W$, and call that set $W_*$.

Let $Y_1$ be a bounded $\sigma (X(t):\bar{t} \in R(\cE_1)\setminus R(E_*))$, $Y_2$ a bounded $\sigma(X(t):\bar{t} \in R(\cE_2)\setminus R(E_*))$-measurable,
and $Y_{*,k}^j$ a bounded $\sigma(X(t):t\in E_{*,k}^j)$-measurable random variable, for $j=1,2$, $k=1,\ldots,n_j$,
$$Y_1'=\prod_{k=1}^{n_1}Y_{*,k}^1,\quad Y_2'=\prod_{k=1}^{n_1}Y_{*,k}^2. $$

First, we will show that $Y_1'$ is independent of $Y_2'$ given $X_{W_*}$. 
Let $A\in \sigma(X_{W_*})$. Using $(\ref{mrl:mp1})$ we get 
\begin{align*}
 \E(Y_1'Y_2'\1_A)&=\E(Y_1'\E(Y_{*,2}^1|X_{\partial E_{*,k}^1})Y_{*,2}^2\ldots Y_{*,2}^{n_2}\1_A )\\
 &=\E(Y_1'\E(Y_{*,1}^2|X_{\partial E_{*,1}^2})\E(Y_{*,2}^2|X_{\partial E_{*,2}^2})\ldots Y_{*,2}^{n_2}\1_A ) \\
 &\vdots \\
 &=\E(Y_1'\E(Y_{*,1}^2|X_{\partial E_{*,1}^2})\E(Y_{*,2}^2|X_{\partial E_{*,2}^2})\ldots \E(Y_{*,2}^{n_2}|X_{\partial E_{*,n_2}^2})\1_A)
\end{align*}
Now, we condition everything under the expectation with respect to $X_{W_*}$:
\begin{align*}
 &=\E(\E[Y_1'\E(Y_{*,1}^2|X_{\partial E_{*,1}^2})\E(Y_{*,2}^2|X_{\partial E_{*,2}^2})\ldots \E(Y_{*,2}^{n_2}|X_{\partial E_{*,n_2}^2})\1_A|X_{W_*}])\\
 &=\E(\E[Y_1'|X_{W_*}]\E(Y_{*,1}^2|X_{\partial E_{*,1}^2})\E(Y_{*,2}^2|X_{\partial E_{*,2}^2})\ldots \E(Y_{*,2}^{n_2}|X_{\partial E_{*,n_2}^2})\1_A).
\end{align*}
%
Using $(\ref{mrl:mp1})$ again we get 
\begin{align*}
 &=\E(\E[Y_1'|X_{W_*}]Y_{*,1}^2\E(Y_{*,2}^2|X_{\partial E_{*,2}^2})\ldots \E(Y_{*,2}^{n_2}|X_{\partial E_{*,n_2}^2})\1_A)\\
 &\vdots \\
 &=\E(\E[Y_1'|X_{W_*}]Y_{*,1}^2Y_{*,2}^2\ldots Y_{*,2}^{n_2}\1_A) =\E(\E[Y_1'|X_{W_*}]Y_2'\1_A).
\end{align*}
Finally, conditioning everything under the expectation with respect to $X_{W_*}$ 
we get
$$=\E(\E[Y_1'|X_{W_*}]\E(Y_2'|{X_{W_*}})\1_A),$$
and the claim follows.

Further, by Lemma \ref{lem:sepind}, we have that 
\be X_{\cup_{k=1}^{n_2} E_{*,k}^2} \perp X_{R(\cG)\setminus (\cup_{k=1}^{n_2} R(E_{*,k}^2))}|X_{W_*}.\label{mrl:mp2}\ee

Let $A_H\in \sigma (X_{W_{\cH}})$ and $A_*\in \sigma (X_{W_*})$. Now, since $W_{\cH}$
separates $\cG$ into $\cE_1'\cup\{E_*\}$ and $\cE_2'=\cE_2\cap R(\cH)$ we have :
\begin{align}
 & \E(Y_1Y_2Y_1'Y_2'\1_{A_H}\1_{A_*}) = \E(Y_1Y_2Y_1'\E(Y_2'|X_{W_*})\1_{A_H}\1_{A_*}) \label{mrlalg:1}\\
=& \E(Y_1\E(Y_2|X_{W_{\cH}})Y_1'\E(Y_2'|X_{W_*})\1_{A_H}\1_{A_*})\label{mrlalg:2}\\
=& \E(\E(Y_1Y_1'|X_W)\E(Y_2|X_{W_{\cH}})\E(Y_2'|X_{W_*})\1_{A_H}\1_{A_*})\nonumber \\
=& \E(\E(Y_1Y_1'|X_W)Y_2\E(Y_2'|X_{W_*})\1_{A_H}\1_{A_*}) \nonumber \\
=& \E(\E(Y_1Y_1'|X_W)Y_2Y_2'\1_{A_H}\1_{A_*}) \nonumber \\
=& \E(\E(Y_1Y_1'|X_W)\E(Y_2Y_2'|X_W)\1_{A_H}\1_{A_*}).\nonumber
\end{align}
To get $(\ref{mrlalg:1})$ we use \ref{mrl:mp2}. 
In $(\ref{mrlalg:2})$ we use the fact that $W_{\cH}$ separates $\cE_2'$ 
from the rest of $\cG$, and then the property proven in the previous $\bullet$ case.

Now, by Monotone Class Theorem the claim follows.
%
%
\end{itemize}}}

The following corollary, gives us a connection to the Markov random fields\index{Markov random field (MRF)} and classical
graphical models (see Appendix \S \ref{sec:mrf}).

\pos{For a Markov family $\cM$, let $X$ be a natural $\cM$-process on a TLG$^*$ $\cG=(\cV,\cE)$.
Let $W$ be a finite set of points on $\cG$ such that $\{t\in \cV: d(t)\geq 3\}\subset W$, then $(X(t):t\in W)$
is a random Markov field with a global Markov property. Further, $X_W$ is a random Markov field indexed by the graph $G_W=(W,E_{W})$ where 
$E_W$ contains an edge between $w_1$ and $w_2$ if there is a time path $\pi$
in $\cG^\heartsuit$ between $w_1$ and $w_2$ such that $R(\pi)\cap R(W)=\{w_1,w_2\}$.
\begin{figure}[ht]
\begin{center}
\includegraphics[width=6cm]{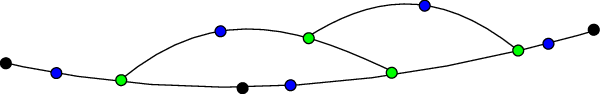}\quad \quad\includegraphics[width=6cm]{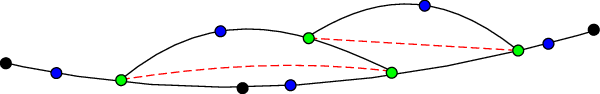}\\
\includegraphics[width=6cm]{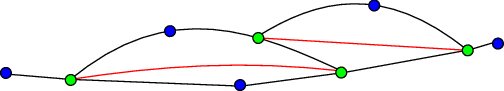}\quad \quad\includegraphics[width=6cm]{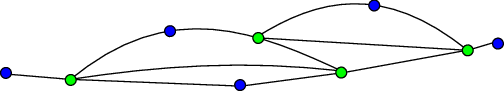}\\
\caption{Getting the MRF: The radnom variables at green and blue points form a Markov random field, where
the underlying graph is given on the last figure. } \label{pic:grph_model}
  \end{center}
\end{figure}

\dok{It is easy to see that $C\subset W$ separates graph $E$ if and only if it separates $\cG^{\heartsuit}$. 
Now it follows that $X_A\perp X_B | X_C$, since $A$ and $B$ are in two different 
components in $\cG^{\heartsuit}$ separated by $C$.}}
\index{Graph-Markovian property!moralized|)}

\subsection*{The constructed process is edge-Markovian}\index{Edge-Markovian property|(}

\defi{\label{def:edg_mrk}We say that the process $X$ on a TLG $\cG=(\cE,\cV)$ is \textbf{edge-Markovian} if
for each $E'=E_{j'k'}\in\cE$ the process $(X(t):t\in E')$ is independent of $(X(t):t\in E,E\in \cE\setminus\{E'\})$
given $X(t_{j'})$ and $X(t_{k'})$.}

\pos{Let $X$ be a natural $\cM$-process on a TLG$^*$ $\cG$. Let $\pi$ be a time-path
between $t_*$
and $t^*$ two points on $\cG$ such that $\pi$ (in the interior) doesn't contain a vertex of degree 3 or more.
Then $X_{\pi}$ and $X_{R(\cG)\setminus R(\pi)}$ are independent given $X(t_*)$ and $X(t^*)$.
\dok{Except the endpoints, the path, can't contain an edge in $\cE^{\heartsuit}\setminus \cE$. Therefore,
endpoints $t_*$ and $t^*$ separate the graph $\cG^{\heartsuit}$ with representations 
of components being $R(\pi)$ and $R(\cG)\setminus R(\pi)$. The calim follows.}}

\teo{The process $X$ defined on $\cG$ defined in \ref{constr} for a Markov family $\cM$ is an edge-Markovian process.}\index{Edge-Markovian property|)}

\subsection{Summary}

Everything we proved so far, can be summarized in the following theorem.

\teo{ \label{cons:sum}
For every TLG$^*$ $\cG$ with finite vertex set $\cV$ and every Markov
family  $\cM$ there exists a hereditary spine-Markovian $\cM$-process $X$ on $\cG$, and the distribution of such
a process is unique. This process also has time-Markovian, cell-Markovian, moralized graph-Markovian and edge-Markovian
properties. Further, if $\cG$ can be constructed from a TLG$^*$ $\cH$, then 
$(X(t):t\in \cH)$ also has these properties.}

\pos{\label{cor:sum_int}Let $X$ be a natural $\cM$-process on a TLG$^*$ $\cG$, where $\cM$ is a Markov family. Then
for $\tau_1\prec \tau_2$
the process $(X(t):t\in\cG[\tau_1,\tau_2])$ has time-Markovian, cell-Markovian, moralized graph-Markovian and edge-Markovian
properties (induced by the structure of $\cG[\tau_1,\tau_2]$).
\dok{We can assume that $\tau_1$ and $\tau_2$ are vertices on $\cG$.
By Theorem \ref{teo:inttgl}  $\cG[\tau_1,\tau_2]$ is a TLG$^*$. Further, 
$$\cM(\cG[\tau_1,\tau_2])=\{\mu_{\sigma}\circ \pi_{[\tau_1,\tau_2]}^{-1}:\sigma \in P_{0\to 1}(\cG), \tau_1,\tau_2\in \sigma\}$$
satisfies (3T) properties.  By Theorem \ref{cor:tlg*cnstint} we can construct first construct a  full time path 
$\sigma$ containing $\tau_1$ and $\tau_2$, and then $\cG[\tau_1,\tau_2]$, and after that the rest of 
$\cG$. It is not hard to see that when we are done constructing $\cG[\tau_1,\tau_2]$ in that TLG$^*$-tower, 
the process restricted to $\cG[\tau_1,\tau_2]$ will be a natural $\cM(\cG[\tau_1,\tau_2])$-process.
}}

\section{Homogeneous Markov family $\cM_{\cP}$}

Let $\cP$ be distribution of a continuous or RCLL Markov process on $[0,1]$. Then we will call
$\cM_{\cP}=\{\mu_{\sigma}=\cP:\sigma\in P_{0\to 1}(\cG)\},$
a \textbf{homogeneous Markov family}. Note that for this family
properties (T2) and (T3) are automatically satisfied.\vspace{0.3cm}

Further, using Theorem \ref{cons:sum}. we have the following fact.

\teo{For every TLG$^*$ $\cG$ with finite vertex set $\cV$ and every  Markov
process $\cP$ there exists a hereditary spine-Markovian $\cM_{\cP}$-process $X$ on $\cG$, 
and the distribution of such
a process is unique. This process also has time-Markovian, cell-Markovian, moralized graph-Markovian, and edge-Markovian
properties.}\vspace{0.2cm}

We will refer to the process $X$ described in the Theorem as the \textbf{natural 
$\cP$-process} on the TLG$^*$ $\cG$.\vspace{0.2cm}


\subsection{The graph-Markovian property doesn't hold}\index{Graph-Markovian property!doesn't hold|(}\label{not_grph}

In paper \cite{tlg1} it was claimed that for the the natural $\cP$-process, 
on what they called NCC graphs, the graph-Markovian property holds.\vspace{0.2cm}

However, the following example shows that this is not true.\vspace{0.2cm}

In our simple model we look at a family of random variables
$\{X_0,X_a,X_b,X_1\}$. Such that 
$(X_0,X_a,X_1)$ is a Markov chain. $(X_0,X_b,X_1)$ is also a Markov
chain independent of the fist one given $(X_0,X_1)$ and has the
same distribution. 

\begin{figure}[ht]
\begin{center}
\psfrag{a}{$\boldsymbol{a}$}
\psfrag{b}{$\boldsymbol{b}$}
\psfrag{1}{$\boldsymbol{1}$}
\psfrag{0}{$\boldsymbol{0}$}

\includegraphics[width=5cm]{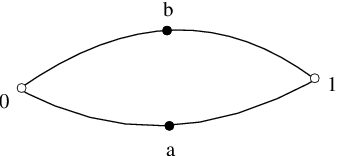}\\
  \end{center}
\end{figure}

We suppose that
the state space $S$ is finite or countable.\vspace{0.2cm}

We set for $x_0,x_a,x_b,x_1\in S$
$$\P(X_0=x_0,X_a=x_a,X_1=x_1)=\P(X_0=x_0,X_b=x_a,X_1=x_1)=p^2_{x_1x_a}p^1_{x_ax_0}p^0_{x_0},$$
with the usual assumptions on initial probabilities $(p^0_s)_{s\in S}$
and transition probabilities $(p^1_{ss'})_{s,s'\in S}$ and $(p^2_{ss'})_{s,s'\in S}$.
Further from the assumption of independence given $(X_0,X_1)$ we have
$$\P(X_a=x_a,X_b=x_b|X_0=x_0,X_1=x_1)=$$
$$\P(X_a=x_a|X_0=x_0,X_1=x_1)\P(X_b=x_b|X_0=x_0,X_1=x_1).$$

Our ultimate goal is to see does 
$$\P(X_1=x_1|X_a=x_a,X_b=x_b,X_0=x_0)\eqno{(*)}$$
depend on $x_0$.
We will first calculate
\begin{align*}
 \P(X_1=x_1,X_0=x_0)&=\sum_{\alpha\in S}\P(X_1=x_1,X_a=\alpha, X_0=x_0)\\
&= \sum_{\alpha\in S}p^2_{x_1\alpha }p^1_{\alpha x_0}p^0_{x_0}.
\end{align*}\vspace{0.2cm}

Next, using the definition of conditional probability and conditional 
independence we calculate 
\begin{align*}
 &\ \P(X_1=x_1,X_a=x_a,X_b=x_b,X_0=x_0)\\
&= \P(X_a=x_a,X_b=x_b|X_1=x_1,X_0=x_0)\P(X_1=x_1,X_0=x_0)\\
&=\P(X_a=x_a|X_1=x_1,X_0=x_0)\P(X_b=x_b|X_1=x_1,X_0=x_0)\P(X_1=x_1,X_0=x_0)\\
&=\frac{\P(X_1=x_1,X_a=x_a,X_0=x_0)\P(X_1=x_1,X_b=x_b,X_0=x_0)}{\P(X_1=x_1,X_0=x_0)}\\
&=\frac{p^2_{x_1x_a }p^1_{x_a x_0}p^0_{x_0}p^2_{x_1x_b }p^1_{x_b x_0}p^0_{x_0}}{\sum_{\alpha\in S}p^2_{x_1\alpha }p^1_{\alpha x_0}p^0_{x_0}}\\
&=p^0_{x_0}\frac{p^2_{x_1x_a }p^1_{x_a x_0}p^2_{x_1x_b }p^1_{x_b x_0}}{\sum_{\alpha\in S}p^2_{x_1\alpha }p^1_{\alpha x_0}}.
\end{align*}


To get $(*)$ we need to calculate
\begin{align*}
  \P(X_a=x_a,X_b=x_b,X_0=x_0)
&=\sum_{\gamma_1\in S}\P(X_1=\gamma_1,X_a=x_a,X_b=x_b,X_0=x_0)\\
&= p^0_{x_0}\sum_{\gamma_1\in S}\frac{p^2_{\gamma_1x_a }p^1_{x_a x_0}p^2_{\gamma_1x_b }p^1_{x_b x_0}}{\sum_{\alpha\in S}p^2_{\gamma_1\alpha }p^1_{\alpha x_0}}\\
\end{align*}
Finally, we have 
\begin{align*}
 &\ \P(X_1=x_1|X_a=x_a,X_b=x_b,X_0=x_0)\\
&=\frac{\P(X_1=x_1,X_a=x_a,X_b=x_b,X_0=x_0)}{\P(X_a=x_a,X_b=x_b,X_0=x_0)}\\
&=p^0_{x_0}\frac{p^2_{x_1x_a }p^1_{x_a x_0}p^2_{x_1x_b }p^1_{x_b x_0}}{\sum_{\alpha\in S}p^2_{x_1\alpha }p^1_{\alpha x_0}}\left(p^0_{x_0}\sum_{\gamma_1\in S}\frac{p^2_{\gamma_1x_a }p^1_{x_a x_0}p^2_{\gamma_1x_b }p^1_{x_b x_0}}{\sum_{\alpha\in S}p^2_{\gamma_1\alpha }p^1_{\alpha x_0}}\right)^{-1}\\
&=\frac{p^2_{x_1x_a }p^2_{x_1x_b }}{\sum_{\alpha\in S}p^2_{x_1\alpha }p^1_{\alpha x_0}}\left(\sum_{\gamma_1\in S}\frac{p^2_{\gamma_1x_a }p^2_{\gamma_1x_b }}{\sum_{\alpha\in S}p^2_{\gamma_1\alpha }p^1_{\alpha x_0}}\right)^{-1}.
\end{align*}

The last shows that $X_1$ given $X_a$, $X_b$, $X_0$ depends on the value
of $X_0$. If the graph-Markovian property holds this should not be so.

Simplifying our model to $S=\{0,1\}$, and setting $p^0_0=p^0_1=1/2$,and $p^1_{10}=p^2_{10}=3/4$,
and $p^1_{11}=p^2_{11}=1/4$, we get that 
 $$\P(X_1=1|X_a=0,X_b=1)=1/2,$$
while
$$\P(X_1=1|X_a=0,X_b=1,X_{0}=0)=3/8.$$
Hence, the graph-Markovian property doesn't hold.\index{Graph-Markovian property!doesn't hold|)}
\subsection{Construction problems on non-TLG$^*$ TLG's}\label{subsec:cnsprb}\index{Process indexed by a TLG!construction problems|(}
Why the construction described in \ref{constr} (on page \pageref{constr}) can't work for all TLG's? 
As an example of Burdzy and Pal 
presented in \cite{tlg1} shows it may not be possible to construct such a
process and have all the properties Markov processes on TLG$^*$'s had.

Let's take a look at the example of a TLG that is not a TLG$^*$  given in Theorem \ref{teo:tlg_exmpl}.(\ref{teo:tlg_exmpl:1}).

$\cG=(\cV,\cE)$, where $\cV=\{t_j=j/5:j=0,1,\ldots,5\}$ and
$$\cE=\{E_{01},E_{02}, E_{14}, E_{13},  E_{23},E_{24},E_{45},E_{35}\}.$$

\begin{figure}[ht]
\begin{center}
\psfrag{0}{$\boldsymbol{t_0}$} \psfrag{1}{$\boldsymbol{t_1}$}
\psfrag{2}{$\boldsymbol{t_2}$}\psfrag{3}{$\boldsymbol{t_4}$}
\psfrag{4}{$\boldsymbol{t_3}$} \psfrag{5}{$\boldsymbol{t_5}$}
\includegraphics[width=7cm]{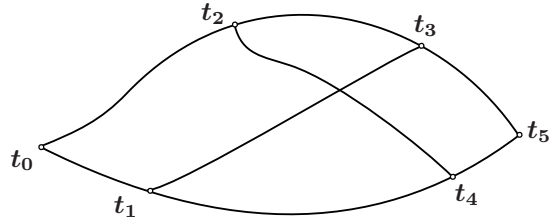}\\
\caption{ Example from Theorem \ref{teo:tlg_exmpl}.(\ref{teo:tlg_exmpl:1}).} \label{pic1_1}
  \end{center}
\end{figure}

Let's take $\cP$ to be Markov process on $[0,1]$.

We will try to construct a process on $\cG$, with a similar approach
as in the construction of Markov processes on TLG$^*$ (see \ref{constr}.)\vspace{0.3cm}

\subsection*{Construction attempt} We first define the process on $\sigma(0,2,4,5)$
 with distribution $\cP$, we construct a $\cP$-Markov bridge
on $\sigma(2,3,5)$ between $(t_2,X(t_2))$ and $(t_5,X(t_5))$ which is independent of the rest of the process 
already defined given $X(t_2)$ and $X(t_5)$.
Further, we construct a $\cP$-Markov bridge  between $(t_0,X(t_0))$ and $(t_5,X(t_5))$  on $\sigma(0,1,4)$ independent of the rest given $X(t_0)$ and $X(t_4)$.
Finally, we construct a $\cP$-Markov bridge  between $(t_1,X(t_1))$ and $(t_3,X(t_3))$  on $\sigma(1,3)$ that is independent of everything
already defined given $X(t_1)$ and $X(t_3)$.\vspace{0.2cm}

The problem in this construction is in the last step. Since, at that time 
$t_1$ and $t_3$ are not connected, the process on 
the full time-path $\sigma(0,1,3,5)$ doesn't have to be $\cP$-distributed.\vspace{0.2cm}

We will prove this when $\cP$ is Brownian motion on $[0,1]$. Then $X_{\sigma(2,3,5)}$
and $X_{\sigma(0,1,4)}$ are Brownian bridges. Using Theorem \ref{bbridg}. we can have
\be X(t_1)=\frac{t_1}{t_4}(X(t_4)-W(t_4))+W(t_1), \label{pokcon1}\ee
\be X(t_3)=\frac{t_5-t_3}{t_5-t_2}(X(t_2)-B(t_2))+B(t_3)+\frac{t_3-t_2}{t_5-t_2}(X(t_5)-B(t_5)),\label{pokcon2}\ee
where $W$, $B$, $X_{\sigma(0,2,4,5)}$ are independent Brownian motions on $[0,1]$.
If $X_{\sigma(0,1,3,5)}$ is Brownian motion on $[0,1]$ then $\E(X(t_1)X(t_3))=t_1=1/5$,
but in our case we get from $(\ref{pokcon1})$ and $(\ref{pokcon2})$:
$$\E(X(t_1)X(t_3))=\frac{1}{3}.$$

\subsection*{Problems with cell-Markovian property} 
The other problem, that might occur, is that the processes we defined 
so far on TLG$^*$'s have the cell-Markovian property (recall Definition \ref{def:cllmk}.), while
on this TLG there might not exist such a process.

Will show this, again, on the example when $\cP$ 
is the distribution of Brownian motion on $[0,1]$.

\prop{\label{prop:inbb}If $Z$ and $Y$ be distributed as Brownian motion on $[s_1,s_2]$ such that $Z(s_j)=Y(s_j)$ for $j=1,2$
and $Z$ and $Y$ are independent given $Y(s_1)$ and $Y(s_2)$. Then for $\tau_1,\tau_2\in [s_1,s_2]$ 
we have
$$\E(Z(\tau_1)Y(\tau_2))=s_1+\frac{(\tau_1-s_1)(\tau_2-s_1)}{(s_1-s_2)}.$$
\dok{We will use the representation given in Theorem \ref{bbridg}. Let $Y$
be Brownian motion on $[s_1,s_2]$, and $(W(t):t\geq 0)$ Brownian motion independent
of $Y$. Then we can take $Z$ to be
$$Z(t)=\frac{s_2-t}{s_2-s_1}(Y(s_1)-W(s_1))+W(t)+\frac{t-s_1}{s_2-s_1}(Y(s_2)-W(s_2)),$$
for $t\in [s_1,s_2]$. Now, we have
\begin{align*}
 \E[Z(\tau_1)Y(\tau_2)]=&\frac{s_2-\tau_1}{s_2-s_1}\E[Y(\tau_2)(Y(s_1)-W(s_1))]+\E(Y(\tau_2)W(\tau_1))\\
&+\frac{\tau_1-s_1}{s_2-s_1}\E[Y(\tau_2)(Y(s_2)-W(s_2))]=\frac{s_2-\tau_1}{s_2-s_1}s_1+\frac{\tau_1-s_1}{s_2-s_1}\tau_2.
\end{align*}
}}

\teo{\label{teo:exncstr}There doesn't exist a process $X$ on $\cG$ such that:
\begin{itemize}
 \item $X$ is cell-Markovian.\index{Cell-Markovian property}
 \item For each full-time $\sigma$ the process $X_{\sigma}$ is distributed as Brownian
motion on $[0,1]$.
\end{itemize}}
\dok{Assume otherwise. Note that cells $(\sigma(2,3,5),\sigma(2,4,5))$ and $(\sigma(1,4,5),\sigma(1,3,5))$
are simple.

Then $X_{\sigma(2,3,5)}$ and $X_{\sigma(2,4,5)}$ are distributed as Brownian motions
on $[t_1,t_5]$, so using the cell-Markovian property of $X$, i.e. the fact that 
$X_{\sigma(2,3,5)}$ and $X_{\sigma(2,4,5)}$ are independent given $X(t_2)$ and $X(t_5)$
from Proposition \ref{prop:inbb}. we have:
$$\E(X(t_3)X(t_4))=\E(X_{\sigma(2,3,5)}(t_3)X_{\sigma(2,4,5)}(t_4))=t_2+\frac{(t_3-t_2)(t_4-t_2)}{(t_5-t_2)}=\frac{8}{15}.$$

For $X_{\sigma(1,4,5)}$ and $X_{\sigma(1,3,5)}$ in a similar way we get:
$$\E(X(t_3)X(t_4))=\E(X_{\sigma(1,3,5)}(t_3)X_{\sigma(1,4,5)}(t_4))=t_1+\frac{(t_3-t_1)(t_4-t_1)}{(t_5-t_1)}=\frac{13}{10}.$$

This shows the claim.}

\index{Process indexed by a TLG!construction problems|)}

\section{Three simple examples}

\begin{figure}[ht]
\begin{center}
\psfrag{a}{$\boldsymbol{E_{01}^1}$}
\psfrag{b}{$\boldsymbol{E_{01}^2}$}
\psfrag{0}{$\boldsymbol{0}$}
\psfrag{1}{$\boldsymbol{1}$}

\includegraphics[width=7cm]{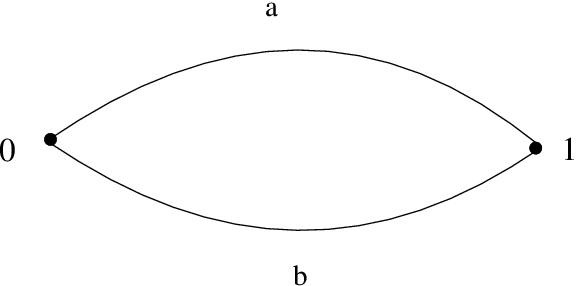}\\
\caption{Graph $\cG$} \label{pic:sl23}
  \end{center}
\end{figure}

Through this section $\cG$ is a graph consisting of one cell (see Figure \ref{pic:sl23}):
$$\cG=(\{0,1\},\{E_{01}^1,E_{01}^2\}).$$
We define three functions $f_1,f_2,f_3:[0,1]\to [0,1]$:
$$f_1(t)=t,\quad f_2(t)=t^2$$
$$f_3(t)=\begin{cases}
          2x, & 0\leq x\leq 1/3;\\
          1-x, & 1/3\leq x\leq 2/3;\\
	  2x-1,& 2/3\leq x\leq 1.
         \end{cases}
$$

\begin{figure}[ht]
\begin{center}
\psfrag{a}{$\boldsymbol{y=f_1(x)}$}
\psfrag{b}{$\boldsymbol{y=f_2(x)}$}
\psfrag{c}{$\boldsymbol{y=f_3(x)}$}

\includegraphics[width=13cm]{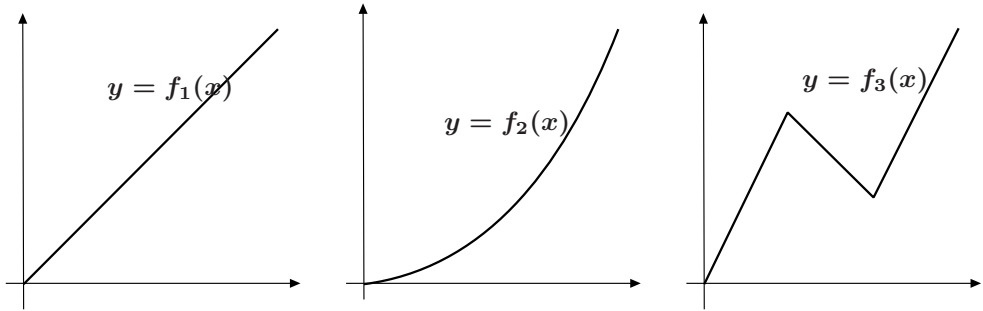}\\
\caption{Graphs of $f_1$, $f_2$ and $f_3$} \label{pic:sl24}
  \end{center}
\end{figure}

Further, let $(B_{t}:t\in [0,1])$ be Brownian motion on $[0,1]$. 
For $k=1,2,3$ we set $\mu_k$ to be the law of $(B_{f_k(t)}:t\in [0,1])$.
Notice that $\mu_k$ for $k=1,2,3$ are all laws of continuous processes.
Also $\mu_1$ and $\mu_2$ are laws of Markov processes, while 
$\mu_3$ is not a law of a Markov process.
Now we set
$$\cM_1:=\{\mu_{E_{01}^1}=\mu_1,\mu_{E_{01}^2}=\mu_1\},$$
$$\cM_2:=\{\mu_{E_{01}^1}=\mu_1,\mu_{E_{01}^2}=\mu_2\},$$
$$\cM_3:=\{\mu_{E_{01}^1}=\mu_1,\mu_{E_{01}^2}=\mu_3\}.$$
Since $f_1(0)=f_2(0)=f_3(0)=0$ and 
$f_1(1)=f_2(1)=f_3(1)=1$,
$\cM_1$, $\cM_2$ and $\cM_3$ satisfy (3T) properties.
Therefore we can construct a natural $\cM_k$-process on $\cG$
for each $k\in\{1,2,3\}$.

Now, note the following:
\begin{itemize}
 \item $\cM_1$ is a homogeneous Markov family.
 \item $\cM_2$ is a Markov family that is not 
homogeneous.
\item $\cM_3$ is not a Markov family.
\end{itemize}


\chapter{Filtrations, martingales and stopping times}\label{sec:3b}

Let's look at a simple example of process on a time-like graph.

$Y$ a value two persons (\textcolor{blue}{1}\&\textcolor{red}{2}) are trying to estimate
based on the information they are getting over time.

\begin{itemize}
 \item The information they collect will be modeled as a filtration 
$$\{\textcolor{blue}{\F^1_t}:t\in [0,1]\}\quad \textrm{and} \quad \{\textcolor{red}{\F^2_t}:t\in [0,1]\}.$$
\item At $t=0$ they start with the same information $\F^1_0=\F^2_0$.
 \item At time $t=1$ everything is known: $\F^1_1=\F^2_1=\F\supset\sigma(Y)$.
 
\end{itemize}
 Set $\textcolor{blue}{X^1_t}=\E(Y|\textcolor{blue}{\F^1_t})$ and $\textcolor{red}{X^2_t}=\E(Y|\textcolor{red}{\F^2_t})$. 

For a TLG $\cG=(\{0,1\},\{E_{01}^1,E_{01}^2\})$, we can define $X=(X(t):t\in \cG)$
to be given by $X_{E_{01}^1}=X^1$ and $X_{E_{01}^2}=X^2$. In this way the process is 
well defined.

\begin{figure}[ht]
\begin{center}

\psfrag{0}{\small $\boldsymbol{0}$} 
\psfrag{1}{\small $\boldsymbol{1}$}
\psfrag{t}{\tiny $\boldsymbol{t}$}
\includegraphics[width=4cm]{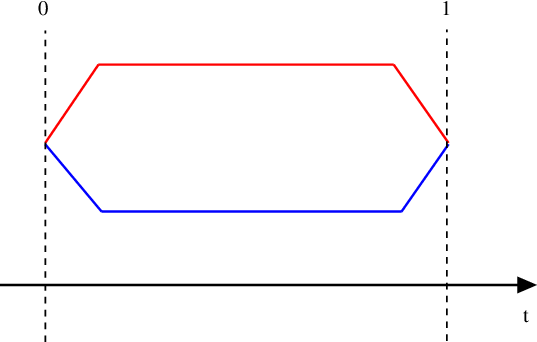}
  \end{center}
\end{figure}
$X$ will be a martingale indexed by a TLG\index{Martingales!indexed by a TLG}, and we will talk more about martingales in this chapter.
We will also show some results on the right-continuity of filtrations, define stopping times
and prove the Optional Sampling Theorem for this class of processes.

\section{Expanding the filtrations}
The following will state some equivalent forms of the time-Markovian 
property.

\pos{ 
Let $X$ be a process on a TLG$^*$ $\cG$. The following are equivalent:
we have:
\begin{enumerate}[(a)]
\item $X$ is time-Markovian on $\cG$
 \item For any point 
$$\F_t=\sigma(X(u): u\preceq t),\ \textrm{and}\ \cG_t=\sigma(X(t):u\succeq t) $$
 are conditionally independent given $X(t)$. ('$\preceq$' is the order induced by $\cG$.)

\item If $Y\in b\cG_t$, then we have 
\be\E(Y|\F_t)=\E(Y|X(t)).\label{t:mrk:eqv}\ee
\end{enumerate}
}

The main result in this section will be to show under which conditions we can expand the $\sigma$-algebra $\F_t$
 so that the relation  (\ref{t:mrk:eqv}) still holds. The main idea is to choose the filtration 
that is right continuous\index{Filtration}.\vspace{0.3cm}

%

\defi{Let $X$ be a process on a TLG $\cG$
\begin{itemize}
\item Set $\F_t^{0}=\sigma(X(u):u \preceq t)$ and $\G_t^{0}=\sigma(X(u):u\succeq t)$. If not specified otherwise 
$\F_t=\F^0_t$ and  $\G_t=\G^0_t$.
 \item For each $\pi\in P_{0\to 1}(\cG)$ and $t\in [0,1]$ we define 
\be \F_{t+}^\pi:=\bigcap_{t\prec s,s\in \pi} \F_s.\label{filt:ft+}\ee
\end{itemize}
}

\defi{For the probability space $(\Omega,\F,\P)$ and $\G$ a sub-$\sigma$-algebra of $\F$
we will denote 
$$\cN^{\P}=\{A\subset \Omega: (\exists B\in \cG) (A\subset B)(\P(B)=0) \}.$$
\be\G^{\P}=\sigma(\G\cup \cN^{\P}). \label{aug_sig}\ee
}
\lem{For $\G^{\P}$ defined by $(\ref{aug_sig})$ the following holds
$$\G^{\P}=\{A\subset \Omega : (\exists B\in \G)(A\triangle B\in \cN^{\P})\}.$$
}
\lem{\label{aug_lem2}For the probability space $(\Omega,\F,\P)$ and $\G_1$ and $\cG_2$ sub-$\sigma$-algebras 
of $\F$ the following are equivalent:
\begin{enumerate}[(i)]
 \item $\cG_1^{\P}=\cG_2^{\P}$;
 \item For each $Y\in L^1(\Omega,\F, \P)$
$$\E(Y|\cG_1)=\E(Y|\cG_2)\ a.s.$$
 \item For each $A\in \F$ 
$$\P(A|\cG_1)=\P(A|\cG_2)\ a.s.$$
\end{enumerate}
\dok{(i)$\Rightarrow$(ii): Let $A\in \G_1$. Since $\G_1\subset \G_1^\P=\G_2^\P$,
there exists $B\in \G_2$ and $N_1$, $N_2$ $\P$-null sets such that $A\cup N_1=B\cup N_2$.
Now, for $Y\in L^1(\Omega,\F,\P)$ we have
\begin{align*}
 \int_A\E(Y|\G_2)\, d\P&=\int_{A\cup N_1}\E(Y|\G_2)\, d\P= \int_{B\cup N_2}\E(Y|\G_2)\, d\P\\
&=\int_{B}\E(Y|\G_2)\, d\P=\int_{B}Y\, d\P=\int_{B\cup N_2}Y\, d\P\\
&=\int_{A\cup N_1}Y\, d\P=\int_{A}Y\, d\P=\int_{A}\E(Y|\G_1)\, d\P
\end{align*}
%
Since this holds for all $A\in \G_1$ the claim follows.

(ii)$\Rightarrow$(iii): This is clear.

(iii)$\Rightarrow$(i): Let $A\in \G_1$, then 
$$\1_A=\P(A|\G_1)=\P(A|\G_2)\quad a.s.$$
Since, $\P(A|\G_2)$ is $\G_2$-measurable, hence
$\1_A$ is $\G_2^{\P}$ measurable. Therefore, $\G_1\subset \G_2^{\P}$, and
we have $\G_1^{\P}\subset \G_2^{\P}$. By symmetry $\G_2^{\P}\subset \G_1^{\P}$, and the claim follows. }}

\teo{\label{teo:rcfiltr}Let $\cM$ be the (3T)-family, and $X$ a natural $\cM$ process on 
a TLG$^*$ $\cG$ such that for each $\pi\in P_{0\to 1}(\cG)$ the process $X_{\pi}$ 
is Markov with respect to the $(\F_{t+}^{\pi}:t\in [0,1])$ (recall $(\ref{filt:ft+})$). 
Then
$$\{\F^{\P}_t:t\in \G\}$$
is a right-continuous filtration\index{Filtration!right-continuous}, that is
$$\F^{\P}_t=\bigcap_{t\prec s}\F^{\P}_s.$$

\dok{Let $(\Omega,\F,\P)$ be the probability space on which $X$ is defined. We pick 
$t\in \cG$. Now, we pick a path $\pi$ that contains $t$, and let $E_{k_1k_2}$ be the
edge that is contained in $\pi$ such that $t_{k_1}\leq t<t_{k_2}$.
To prove that that at $t$ the filtration is right continuous we will restrict
our probability space to $(\Omega,\F',\P'=\P|_{\F'})$ where
$$\F'=\sigma(\F_{t_{k_2}}\cup \G_{t_{k_1}} \cup \cN^\P ).$$.

With $F\in \F_t$ and $G\in \cG_t$ we have 

$$\P(F\cap G|\F^{\pi}_{t+})=\1_{F}\P(G|\F^{\pi}_{t+})=\1_{F}\P(G|X(t))= \1_{F}\P(G|\F^{\pi}_t)=\P(F\cap G|\F^{\pi}_t).$$
Using the monotone class theorem we have that for all $A\in \F'$
\be\P(A|\F^{\pi}_t)=\P(A|\F^{\pi}_{t+}).\label{p|t+}\ee
Since $\cN^{\P'}=\cN^{\P}$, we have by 
Lemma \ref{aug_lem2} (iii) that

$$\F^{\pi,\P}_t=\F^{\pi,\P}_{t+}.$$
Further, note that $\F^{\pi}_{t+}\subset \F^{\pi,\P}_t$.

Now, let $$A\in \bigcap_{t<s}\F^{\pi,\P}_s=\bigcap_{n=1}^{\infty}\F^{\pi,\P}_{t+1/n}.$$
Hence, we have $A\in \F^{\pi,\P}_{t+1/n}$, then there exists $B_n\in \F^{\pi}_{t+1/n}$
such that $A\triangle B_n \in \cN^{\P}$. Set
$$B:=\bigcap_{n=1}^{\infty} \bigcup_{m=n}^{\infty} B_m=\bigcap_{n=M}^{\infty} \bigcup_{m=n}^{\infty} B_m\in \F^{\pi}_{t+1/M},$$
hence $B\in \F^{\pi}_{t+}$, hence $B\in \F^{\pi,\P}_{t}$.  Now, we can show that
$$ B\setminus A\subset\left(\bigcup_{n=1}^{\infty}B_n\right)\setminus A= \bigcup_{n=1}^{\infty}(B_n\setminus A)\in \cN^{\P}.$$
$$A\setminus B=A\cap B^c=A\cap \left(\bigcap_{n=1}^{\infty} \bigcup_{m=n}^{\infty} B_m\right)^c= \bigcup_{n=1}^{\infty} A\cap \left(\bigcap_{m=n}^{\infty} B_m^c\right)\subset$$
$$\subset \bigcup_{n=1}^{\infty} A\cap B_m^c= \bigcup_{n=1}^{\infty} (A\setminus B_m)\in \cN^{\P}.$$
This implies that $A\in \F^{\pi,\P}_t$, and the filtration 
$(\F^{\pi,\P}_s:s\in [0,1])$ is right-continuous at $t$, 
and to prove the claim we should note that $\F^{\pi}_t=\F_t$, 
hence this implies $(\F_s^{\P}:s\in \pi)$ is right-continuous at 
$t$, but since $\pi$ is an arbitrary path that contains $t$ the claim follows, since there is only
finitely many such paths. Hence
$$\F_t^\P=\bigcap_{\pi:\pi\ni t }\F_t^{\pi, \P}=\bigcap_{\pi:\pi\ni t }\bigcap_{t<s}\F_s^{\pi, \P}=\bigcap_{t\prec s}\F_s^{ \P}.$$}}

It turns out that the condition from the previous theorem is
satisfied by the natural Brownian motion. Before we prove that 
we need the following lemma.

\lem{\label{mp:bmtlg*}Let $\G$ be a TLG$^*$ and $X$ the natural Brownian motion on $\cG$. For $t\in \G$ we have that the processes 
$$(X(s):s \preceq t)\quad and \quad (X(s)-X(t):s\succeq t)$$
are independent. 
\dok{Pick a full time-path $\pi\in P_{0\to 1}(\cG)$ such that $t\in \pi$. We pick a TLG$^*$ tower 
$(\cG_j)_{j=0}^n$ where $\cG_0$ has the same representation as $\pi$ and $\cG_n=\cG$.\vspace{0.2cm}

Let $X^j:=(X(s):s\in \cG_j)$ be the natural Brownian motion indexed by $\cG_j$. By induction we will show that
\be (X^j(s):s\in P^j(t))\quad and \quad (X^j(s)-X^j(t):s\in F^j(t))\label{bmtlgmp}\ee

It is clear that $(\ref{bmtlgmp})$ holds for $j=0$. Let's assume that it holds for 
$j=h\geq 0$. Let's show the claim for $j=h+1$. 

If a new edge not in $P^{h+1}(t)$ and not in $F^{h+1}(t)$
has been added to $\cG_h$ to construct $\cG_{h+1}$, then the processes in $(\ref{bmtlgmp})$ are the same for $j=h$ and $j=h+1$,
and the claim follows.

If a new edge $E_{k_1k_2}$ in $F^{h+1}(t)$ has been added to $\cG_{h}$ to obtain 
$\cG_{h+1}$. Then since for $s\in E_{k_1k_2}$ we have 
$$X^{h+1}(s)-X(t)=\frac{t_{k_2}-s}{t_{k_2}-t_{k_1}}(X^h(t_{k_2})-X(t))+\frac{s-t_{k_1}}{t_{k_2}-t_{k_1}}(X^h(t_{k_1})-X(t))+B^{br}_{k_1k_2}(s),$$ 
where $B^{br}_{k_1k_2}$ is a Brownian bridge independent of $X^h$. Hence, both
$(X^h(s)-X(t):s\in F^h(t))$ and $(X_{E_{k_1k_2}}(s)-X(t):s\in E_{k_1k_2})$ are independent pointwise of 
$(X^h(s):s\in P^j(t))$, and $(\ref{bmtlgmp})$ follows for $j=h+1$.}}

\teo{\label{teo:bmfilt}Let $\G$ be a TLG$^*$ and $X$ the natural Brownian motion on $\cG$. For $t\in \G$ and $\pi \in P_{0\to 1}(\cG)$ such that 
$t\in \pi$ we have that 
$$\F^{\pi}_{t+}\quad and \quad (X_{\pi}(s)-X_{\pi}(t):s\geq t))$$
are independent. (See Figure \ref{pic:sl57}.)
\begin{figure}[ht]
\begin{center}
\psfrag{p}{\textcolor{red}{$\boldsymbol{\pi}$}} 
\psfrag{t}{$\boldsymbol{t}$}
\includegraphics[width=10cm]{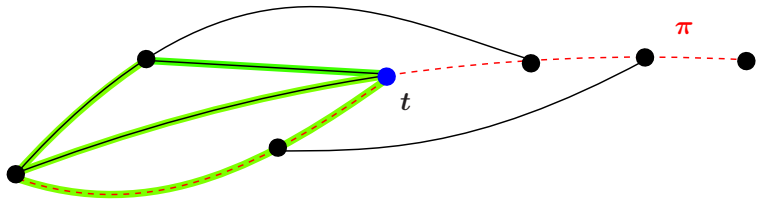}\\
\caption{ Illustration of Theorem \ref{teo:bmfilt}.} \label{pic:sl57}
\end{center}
\end{figure}
\dok{ Let $n\in \N$, and $A\in \F^{\pi}_{t+}$ and $t\prec s_1,s_2,\ldots, s_{n}\in \pi$. For small 
$\varepsilon>0$ we know that $Y:=\1_A\in b\F_{t+\varepsilon/2 }$ and $\Delta^{\varepsilon}:=(X_{\pi}(s_1)-X_{\pi}(t+\varepsilon),\ldots X_{\pi}(s_n)-X_{\pi}(t+\varepsilon))$
are independent. Now using the characteristic functions 
$\varphi_{Y}(t)=\E(\exp(itY))$ and $\varphi_{\Delta^{\varepsilon}}(\mathbf{t})=\E(\exp(i\mathbf{t}\cdot \Delta^{\varepsilon}))$ 
we have 
\be\varphi_{Y,\Delta^{\varepsilon}}(t,\mathbf{t})=\varphi_{Y}(t)\varphi_{\Delta^{\varepsilon}}(\mathbf{t})\label{phi:ind}\ee
Continuity of $X$ gives us $\lim_{\varepsilon \downarrow 0}(Y,\Delta^{\varepsilon})=(Y,\Delta^0)$ a.s. 
Hence, from $(\ref{phi:ind})$ we have
$$\varphi_{Y,\Delta^{0}}(t,\mathbf{t})=\varphi_{Y}(t)\varphi_{\Delta^{0}}(\mathbf{t}).$$
Therefore, $\1_A$ and $(X_{\pi}(s_1)-X_{\pi}(t),\ldots X_{\pi}(s_n)-X_{\pi}(t))$ are independent.}}

\pos{For the natural Brownian motion $X$ on the TLG$^*$ $\G$ the following claims hold:
\begin{enumerate}[(a)]
 \item The filtration $(\F_t^{\P}: t\in \G)$ is right continuous.
 \item $\F_t^{\P}$ and $(X(s)-X(t):s\in F(t))$ are independent.
 \item For $t\prec s$ we have $\E(X(s)|\F_t^{\P})=X(t)$.
 \item For $t\prec s$ and $Y\in b\G_s$ we have $$\E(Y|\F_t^{\P})=\E(Y|X(t)).$$
\end{enumerate}
}

%

%

\section{Markov martingales}

Here we will show that under some conditions we can get a martingale property
for the process defined on a TLG$^*$.

%

\defi{The Markov family of measures 
$$\cM=\{\mu_{\sigma}:\sigma\in P_{0\to 1}(\cG)\}$$
will be called a \textbf{Markov martingale family} if for each $\mu_{\sigma}$-distributed process $(X_{\sigma}(t)\in [0,1])$, we have 
\begin{itemize}
 \item $\E|X_{\sigma}(t)|<\infty$;
 \item $\E(X_{\sigma}(t)| (X_{\sigma}(u):u\in [0,s]))=X_{\sigma}(s)$.
\end{itemize}
}


\teo{\label{teo:mrt}Let $\cM$ be a Markov martingale family, and $X$ an $\cM$-process on a TLG$^*$
$\cG$. Then we have 
\be\E(X(t)|(X(u): u\preceq s))=X(s),\label{rel:mrt}\ee
for all points $s\preceq  t$ in $\cG$.
\dok{First from the time-Markovian property we have that
$$\E(\varphi_M(X(t))|(X(u): u\preceq s))=\E(\varphi_M(X_t)|X_s),$$
where $\varphi_M(x)=\left\{\begin{array}{cl}
                                                x,& |x|<M,\\
                                                M,& |x|\geq M.\\
                                               \end{array}
\right.$. Using the dominated convergence theorem when $M\to \infty $ we have
$$\E(X(t)|(X(u): u\preceq s))=\E(X(t)|X(s)).$$
Now, we pick a full time-path $\sigma$ such that $t$ and $s$ are on it, and we
get 
$$\E(X(t)|(X(u): u\preceq s))=\E(X_{\sigma}(t)|X_{\sigma}(s))=X_{\sigma}(s)=X(s).$$}}

The following is a consequence of Lemma \ref{aug_lem2} (ii).
\pos{Let $\cM$ be a Markov martingale family, and $X$ an $\cM$-process on a TLG$^*$
$\cG$. Then we have 
\be\E(X(t)|\F_s^\P)=X(s),\label{rel:mrt2}\ee
for all points $s\preceq t$ in $\cG$.}

The equality $(\ref{rel:mrt})$ says that $X$ defined in Theorem \ref{teo:mrt}.
is an example of a {\bf martingale indexed by  directed set} $\cG$.
These types of martingales have been investigated and there are a lot of
results including the  optional sampling theorem. We will talk more about this in section \ref{stptm}.

\subsection{Example of glued diffusions}\label{gldiff}

In this subsection we give an example of a general non-homogeneous 
Markov martingale family $\cM$.

\defi{The family of functions
$$f_{\cG}=\{f_{\sigma}:[0,1]\to \R:\sigma\in P_{0\to 1}(\cG)\}$$
is called consistent on the TLG $\cG$ if for
$\sigma_1,\sigma_2\in P_{0\to 1}(\cG)$ 
$$f_{\sigma_1}|_T=f_{\sigma_2}|_T$$
where $T=\{t\ :\ t\in E, E\in \sigma_1 \, \&\, E\in \sigma_2\}$.
}
 \teo{\label{t:cnsdens}Let $F_{\cG}=\{F_{\sigma}:[0,1]\to \R:\sigma\in P_{0\to 1}(\cG)\}$
be a consistent family of absolutely continuous functions. Then there exists
a consistent family $f_{\cG}=\{f_{\sigma}:[0,1]\to \R:\sigma\in P_{0\to 1}(\cG)\}$
of densities of $F_{\cG}$, that is for all $\sigma\in P_{0\to 1}(\cG)$ and
all $t\in [0,1]$
$$F_{\sigma}(t)-F_{\sigma}(0)=\int_{0}^tf_{\sigma}(s)\, ds.$$ 
\dok{Let $\sigma_1$ and $\sigma_2$ be full time-paths. Then
$$T_{12}=\{t:t\in E, E\in \sigma_1\& E\in \sigma_2\}$$
is a finite union of closed segments. For each $a<b$
such that $(a,b)\subset T_{12}$ we have 
$$\int_{a}^{b}f_{\sigma_1}(s)\, ds=\int_{a}^{b}f_{\sigma_2}(s)\, ds,$$
so therefore $f_{\sigma_1}=f_{\sigma_2}$ $\lambda$-almost everywhere 
on $T_{12}$.\vspace{0.3cm}

Assume $\cG=(\cV,\cE)$, for each edge $E_{kj}\in \cE$ 
choose some fixed full time-path $\sigma^*$ containing $E_{jk}$.
For each full time-path $\sigma$ containing that edge we can fix $f_{\sigma}$
on $(t_j,t_k)$, to be some density of the function $t\mapsto F_{\sigma^*}(t)-F_{\sigma^*}(t_j)$
defined on $(t_j,t_k)$.\vspace{0.3cm}

Since there are only at most countably many vertices (in this case
finitely many) the values at the vertices won't influence the values of the 
integrals, hence we can set the values at vertices to be any real numbers. 
Now, we have constructed a consistent family.}}

Let $\cG$ be a TLG$^*$ and $V:R(\cG)\to \R^+$ be a positive function, such that for each
full time-path $\sigma$ the restriction of $V$ along $R(\sigma)$ 
$V_{\sigma}:[0,1]\to \R^+$ is an increasing continuous function.\vspace{0.2cm}

From the theory of functions of bounded variation, we know that
there exists a positive function $f_\sigma$ in $L^1[0,1]$, such that
$$V_{\sigma}(t)=\int_0^tf_{\sigma}(s)\, ds,$$ 
for all $t\in [0,1]$. By Theorem \ref{t:cnsdens}. we 
can assume that
$$\{f_{\sigma}:[0,1]\to \R:\sigma\in P_{0\to 1}(\cG)\} $$
is a consistent family of densities.\vspace{0.3cm}

For $\sigma\in P_{0\to 1}(\cG)$ let $\mu_{\sigma}$ be the distribution of
the process $(N(t):t\in [0,1])$, given by the formula
$$N_{\sigma}(t):=\int_0^t \sqrt{f_{\sigma}(s)}\, dB_s,$$
for $t\in [0,1]$. (This is an Ito integral with respect to the Brownian 
motion $(B_t)$.) This is well defined since $\sqrt{f_{\sigma}}\in L^2[0,1]$.\vspace{0.3cm}

Clearly, $N_\sigma$ is a Markov process with
zero expectation on $[0,1]$. The variance is
$$\E(N_{\sigma}^2(t))=\E\left(\int_0^t \sqrt{f_{\sigma}(s)}\, dB_s\right)^2=\int_0^t f_{\sigma}(s)\, ds=V_{\sigma}(t).$$
We will show that $\{\mu_{\sigma} : \sigma \in P_{0\to 1}(\cG)\}$ is a consistent
family. Again, let $\sigma_1$ and $\sigma_2$ be two full time-paths, and 
$T_{12}$ as before. Let $\tau_1\leq \tau_2$ be from $T_{12}$.
 We have
\begin{align*}
 \E(N_{\sigma_1}(\tau_1)N_{\sigma_1}(\tau_2))&=\E(N_{\sigma_1}(\tau_1)[(N_{\sigma_1}(\tau_2)-N_{\sigma_1}(\tau_1))+N_{\sigma_1}(\tau_1)])\\
&=V_{\sigma_1}(\tau_1)=V_{\sigma_2}(\tau_1)\\
&=\E(N_{\sigma_2}(\tau_1)N_{\sigma_2}(\tau_2)).
\end{align*}
Since, the covariance structure of the Gaussian processes $N_{\sigma_1}$
and $N_{\sigma_2}$ on $T_{12}$ is the same, we have
that the finite dimensional distributions on $T_{12}$ are the same.
Hence, by Kolmogorov's Existence Theorem we have that their
distributions on $T_{12}$ are the same. Therefore 
$\{\mu_{\sigma}:\sigma \in P_{0\to 1}(\cG)\}$ is a
consistent Markov martingale family.\vspace{0.2cm}

\begin{figure}[ht]
\begin{center}
\psfrag{a}{$\boldsymbol{t}$}
\psfrag{b}{$\boldsymbol{t^2}$}
\psfrag{0}{$\boldsymbol{0}$}
\psfrag{1}{$\boldsymbol{1}$}

\includegraphics[width=7cm]{sl23.eps}\\
\caption{Graph $\cG$} \label{pic:sl23a}
  \end{center}
\end{figure}
Let $$\cG=(\{0,1\},\{E^1_{01},E^2_{01}\}).$$
If we define $V$ as 
$$V(t)=\begin{cases}
     t& \textrm{for}\ t\in E^1_{01},\\
     t^2& \textrm{for}\ t\in E^2_{01};\\
    \end{cases}
$$

Brownian motion runs along $E^1_{01}$, while 
$N(t)=\int_0^ts\, dB_s$ runs along $E^2_{01}$. (See Figure \ref{pic:sl23a}.)\vspace{0.2cm}

Glued diffusions have several nice properties. Since along each path the distribution 
is inducing a martingale and a Markov process the whole process is a martingale 
and a Markov process indexed by the underlying TLG$^*$.

Further, we have the following property which is a generalization of the 
Lemma \ref{mp:bmtlg*}. 

\lem{Let $X$ be a natural glued diffusion on a TLG$^*$ $\G$. Then for each 
$t\in \G$ 
$$(X(s) : s\preceq t) \quad \textrm{and} \quad (X(s)-X(t) : t\preceq s)$$
are independent.\dok{Pick a full time-path $\pi\in P_{0\to 1}(\cG)$ such that $t\in \pi$. We pick a TLG$^*$ tower 
$(\cG_j)_{j=0}^n$ where $\cG_0$ has the same representation as $\pi$ and $\cG_n=\cG$.\vspace{0.2cm}

Let $X^j:=(X(s):s\in \cG_j)$ be the natural natural glued diffusion indexed by $\cG_j$. By induction we will show that
\be (X^j(s):s\in P^j(t))\quad and \quad (X^j(s)-X^j(t):s\in F^j(t))\label{bmtlgmp:A}\ee

It is clear that $(\ref{bmtlgmp:A})$ holds for $j=0$. Let's assume that it holds for 
$j=h\geq 0$. Let's show the claim for $j=h+1$. 

If a new edge not in $P^{h+1}(t)$ and not in $F^{h+1}(t)$
has been added to $\cG_h$ to construct $\cG_{h+1}$, then the processes in $(\ref{bmtlgmp:A})$ are the same for $j=h$ and $j=h+1$,
and the claim follows.

If a new edge $E_{k_1k_2}$ in $F^{h+1}(t)$ has been added to $\cG_{h}$ to obtain 
$\cG_{h+1}$. Then since for $s\in E_{k_1k_2}$ we have 
$$X^{h+1}(s)-X(t)=\frac{V(t_{k_2})-V(s)}{V(t_{k_2})-V(t_{k_1})}(X^h(t_{k_2})-X(t))+\frac{V(s)-V(t_{k_1})}{V(t_{k_2})-V(t_{k_1})}(X^h(t_{k_1})-X(t))+N^{t_{k_1}t_{k_2}}_{0,0}(s),$$ 
where $N^{t_{k_1}t_{k_2}}_{0,0}(s)$ (see Corollary \ref{difbridg}.) is a diffusion bridge independent of $X^h$. Hence, both
$(X^h(s)-X(t):s\in F^h(t))$ and $(X_{E_{k_1k_2}}(s)-X(t):s\in E_{k_1k_2})$ are independent pointwise of 
$(X^h(s):s\in P^j(t))$, and $(\ref{bmtlgmp:A})$ follows for $j=h+1$.}}

\teo{Let $\G$ be a TLG$^*$ and $X$ the natural glued diffusion on $\cG$. For $t\in \G$ and $\pi \in P_{0\to 1}(\cG)$ such that 
$t\in \pi$ we have that 
$$\F^{\pi}_{t+}\quad and \quad (X_{\pi}(s)-X_{\pi}(t):s\geq t)$$
are independent.
\dok{Let $n\in \N$, and $A\in \F_{\pi,t+}$ and $t\prec s_1,s_2,\ldots, s_{n}\in \pi$. For small 
$\varepsilon>0$ we know that $Y:=\1_A\in b\F_{t+\varepsilon/2 }$ and $\Delta^{\varepsilon}:=(X_{\pi}(s_1)-X_{\pi}(t+\varepsilon),\ldots, X_{\pi}(s_n)-X_{\pi}(t+\varepsilon))$
are independent. Now using the characteristic functions $\varphi_{Y}(t)=\E(\exp(itY))$ and $\varphi_{\Delta^{\varepsilon}}(\mathbf{t})=\E(\exp(i\mathbf{t}\cdot \Delta^{\varepsilon}))$ 
we have 
\be\varphi_{Y,\Delta^{\varepsilon}}(t,\mathbf{t})=\varphi_{Y}(t)\varphi_{\Delta^{\varepsilon}}(\mathbf{t})\label{phi:ind:A}\ee
Continuity of $X$ gives us $\lim_{\varepsilon \downarrow 0}(Y,\Delta^{\varepsilon})=(Y,\Delta^0)$ a.s. Hence, from $(\ref{phi:ind:A})$ we have
$$\varphi_{Y,\Delta^{0}}(t,\mathbf{t})=\varphi_{Y}(t)\varphi_{\Delta^{0}}(\mathbf{t}).$$
Therefore, $\1_A$ and $(X_{\pi}(s_1)-X_{\pi}(t),\ldots, X_{\pi}(s_n)-X_{\pi}(t))$ are independent.}}

\pos{For the natural glued diffusion $X$ on the TLG$^*$ $\G$ the following claims hold:
\begin{enumerate}[(a)]
 \item The filtration $(\F_t^{\P}: t\in \G)$ is right continuous.
 \item $\F_t^{\P}$ and $(X(s)-X(t):s\in F(t))$ are independent.
 \item For $t\prec s$ we have $\E(X(s)|\F_t^{\P})=X(t)$.
 \item For $t\prec s$ and $Y\in b\G_s$ we have $$\E(Y|\F_t^{\P})=\E(Y|X(t)).$$
\end{enumerate}
}

\section{Optional sampling theorem for martingales indexed by directed sets}\label{stptm}
\index{Martingale indexed by directed set|(}
In his paper \cite{kurtz} Kurtz defined stopping times for martingales on 
{\it directed sets}. The way they are defined, TLG's are directed sets.
We will state some of the results obtained by Kurtz and apply them to the 
processes on TLG's.\vspace{0.3cm}

Let $\mathcal{ S}$ be a {\bf directed set}\index{Directed set} with partial ordering denoted by $t \prec s$. That is, 
$\mathcal{S}$ is partially ordered and for $t_1,t_2\in \mathcal{ S}$ there exists $t_3\in {\mathcal S}$
such that $t_1\prec t_3$ and $t_2\prec t_3$.\par 
\noindent\emph{Remark.} Note that TLG's satisfy this definition.\vspace{0.2cm}

Let $(\Omega,\F,\P)$ be a probability space and let $(\F_t)_{t\in {\cal S}}$ be a filtration 
indexed by ${\cal S}$, that is 
\begin{itemize}
 \item $(\F_t)_{t\in {\cal S}}$ is a family of sub-$\sigma$-algebras of $\F$;
 \item $t\prec s$ implies $\F_t \subset \F_s$.
\end{itemize}

A stochastic process $X$ indexed by ${\cal S}$ is a \textbf{martingale}\index{Martingale indexed by directed set} with respect to
$(\F_t)_{t\in {\cal S}}$ if 
$$\E(X(t)|\F_s)=X(s),$$
for all $s\preceq t$.\vspace{0.2cm}

A ${\cal S}$-valued random variable $T$ is a {\bf stopping time}\index{Stopping times} if $(T\preceq t)\in \F_t$
for all $t\in {\cal S}$.\vspace{0.2cm}

As usual we define 
$$\F_T=\{A\in \F: A\cap (T\preceq t) \in \F_t, \forall t\in {\cal S}\}.$$

The following is the first form of the optional stopping theorem.
\lem{\label{lem:opt0}Let $X(t)$ be martingale and let $T_1\preceq T_2$ be stopping times assuming
countably many values. If there exists a sequence $(t_m)$ in ${\cal S}$ such that 
\be\lim_{m\to \infty}\P(T_2\preceq t_m)=1, \label{con:mc1}\ee
and  
\be\lim_{m\to\infty} \E(|X(t_m)|\1_{(T_2\preceq t_m)^c})=0,\label{con:mc2}\ee
and $\E(|X(T_2)|)<\infty$, then 
$$\E(X(T_2)|\F_{T_1})=X(T_1).$$
}

\noindent\emph{Remark.} In a TLG with a finite number of vertices, we could pick the sequence
$t_m=1$. In that case conditions $(\ref{con:mc1})$ and $(\ref{con:mc2})$ would be automatically 
satisfied.\vspace{0.3cm}

In order to extend the result of Lemma \ref{lem:opt0} to general stopping times 
we need to make some assumptions about the index set ${\cal S}$ and the process $X$.
The assumption we make on ${\cal S}$ is that it is a topological lattice. 

Recall the Definition \ref{def:tplt}. of a topological lattice from Section \ref{tlg*tltt}:\vspace{0.2cm}

\noindent {A Hausdorff space $X$ with some order '$\leq$' is called a \textbf{topological
lattice} if for $x_1,x_2\in X$:
\begin{itemize}
 \item  there exists a unique element $x_1\wedge x_2$ such that
$$\{x\in X:x\leq x_1 \}\cap \{x\in X:x\leq x_2\} = \{x\in X:x\leq x_1\wedge x_2\};$$
\item  there exists a unique element $x_1\vee x_2$ such that
$$\{x\in X:x\geq x_1 \}\cap \{x\in X:x\geq x_2\} = \{x\in X:x\geq x_1\vee x_2\}.$$
\end{itemize}
and $x_1\wedge x_2$ and $x_1\vee x_2$ are continuous mappings of $X\times X$ (with product topology) onto $X$.}\vspace{0.2cm}

If ${\cal S}$ is a topological lattice, note that this implies that the sets 
of the form $[t_1,t_2]=\{t: t_1\preceq t\preceq t_2\}$ (intervals) are closed, and hence Borel measurable.\vspace{0.2cm}

\defi{We will say that a \textbf{topological lattice ${\cal S}$ is separable from above}\index{Topological lattice!separable from above} if there
exists a \textbf{separating sequence} $\{t_k\}\subset {\cal S}$, such that all $t\in {\cal S}$ we have
$$t=\lim_{n\to \infty }t^{(n)}$$
where 
\be t^{(n)}:=\min\{t_k:k\leq n, t_k\succeq t\}.\label{deft(n)}\ee }

In Section \ref{tlg*tltt} (see Theorem \ref{teo:tpltt}) we have shown that TLG$^*$ $\cG$ is
a topological lattice, and clearly we can set $\{t_k\}$ to be the set of points with rational times.


The following is the main result for the martingales on directed sets.

\teo{Let ${\cal S}$ be separable from above with separating set $\{t_k\} $, 
$\F_t=\bigcap_{n=1}^{\infty}\F_{t^{(n)}}$ for all $t$, and let $X(t)$ be a martingale 
satisfying $$\lim_{n\to\infty }X(t^{(n)},\omega)=X(t,\omega),$$
for all $(t,\omega)$ for which the limit exists. Let $T_1\preceq T_2$ be ${\cal S}$-valued 
stopping times. Suppose there exists a sequence $(s_m)$ in $\{t_k\}$ such that
$$\lim_{m\to\infty}\P(T_2\preceq s_m)=1,$$
and 
$$\lim_{m\to \infty }\E(|X(s_m)|\1_{(T_2\preceq s_m)^c})=0,$$
and that $\E(|X(T_1)|)<\infty$. Then 
$$\E(X(T_2)|\F_{T_1})=X(T_1).$$
}

The following theorem will translate the results we have into 
the ones of the process indexed by time-like graphs. 

\teo{\label{teo:OST}Let $\cG$ be a TLG$^*$.\begin{enumerate}[(a)]
      \item Let $X(t)$ be a martingale with respect to the filtration $(\F_t)_{t\in {\cG}}$ and let $T_1\preceq T_2$ be stopping times assuming
countably many values. If $\E(|X(T_2)|)<\infty$ then 
$$\E(X(T_2)|\F_{T_1})=X(T_1).$$
      \item Let $X(t)$ be a RCLL martingale with respect to the filtration $(\F_t)_{t\in {\cG}}$ such that 
\be\F_t=\bigcap_{t\prec s}\F_s.\label{filtcn}\ee

For stopping times $T_1\preceq T_2$, if $\E(|X(T_2)|)<\infty$ then 
$$\E(X(T_2)|\F_{T_1})=X(T_1).$$
     \end{enumerate}
}

The key problem will be choosing a good filtration $(\F_t)_{t\in \cG}$ such that 
the $(\ref{filtcn})$ is satisfied.\index{Martingale indexed by directed set|)}

\section{TLG - valued stopping times}\index{Stopping times|(}

Let's assume that $(\F_t:t\in \G)$ is a right-continuous filtration and $X$ is
an RCLL process adapted to this filtration.\vspace{0.2cm}

First, let's define two random times\index{Stopping times!TLG - valued} that we want to make stopping times.

If $\sigma$ is a path in $\cG$, then clearly
$$H_U^{\sigma}:=\inf\{t\in \sigma : X_{\sigma}(t)\in U\},$$
where $U$ is an opened set. This is a standard 
one-dimensional stopping time. A more interesting example is 
$$T_U^{\sigma}:=\inf\{t\in \sigma : (\exists \tau\preceq t)(  X(\tau) \in U)\}.$$
It is not hard to see that $H_U^{\sigma}\preceq T_U^{\sigma}$.

\lem{\label{lem:TU}$T_U^{\sigma}$ is an $(\F_t)$ stopping time.
\dok{Let $t\in \sigma$, then by right continuity we have
$$(T_U^{\sigma}<t)=\bigcup_{n=1}^{\infty}\bigcup_{s\prec t-\frac{1}{n},s\in \mathbb{Q}}(X(s) \in U) \in \F_t.$$
Where $t-\frac{1}{n}$ is the point on $\sigma$ with that time, and $s\prec t-1/n$, $s\in \mathbb{Q}$ means the
point on TLG $\G$ that is before $t-1/n$ and has rational time. 
If $t\notin \sigma$ then there exists 
$$t^{\sigma}=\max\{s\in \sigma \ : \ s\prec t  \}.$$
Now, from the continuity of the filtration we have 
$$(T_U^{\sigma}\prec t)=(T_U^{\sigma}\preceq t^{\sigma})=\bigcap_{n=k}^{\infty} (T_U^{\sigma}\prec t^{\sigma}+1/n) \in \F_{t^{\sigma}+1/k}, $$
for all $k\in \N$. Therefore, the right-continuity of the filtration implies 
$$(T_U^{\sigma}\prec t)\in \F_{t^{\sigma}}\subset \F_t.$$ 
}}

Let $K$ be a compact set. We define $T^{\sigma}_K$ and $H^{\sigma}_K$ in the same 
way as we did $T^{\sigma}_U$ and $H^{\sigma}_U$. We know from classical Markov processes 
that $H^{\sigma}_K$ is a stopping time. 

\lem{$T_K^{\sigma}$ is an $(\F_t)$ stopping time. 
\dok{Define $U_n=\{x:d(x,K)<1/n\}$. It is clear  that $K=\bigcap_{n=1}^{\infty}\overline{U_n}=\bigcap_{n=1}^{\infty}U_n$, and also it
is clear that $T_{U_n}^{\sigma}\leq T_{U_{n+1}}^{\sigma}\leq T_{K}^{\sigma}$.
Set $T:=\sup_{n}T_{U_n}^{\sigma}=\lim_{n\to\infty} T_{U_n}^{\sigma}$. If $T\geq 1$ then 
clearly $T_K^{\sigma}=T$, on the event $T<1$ we have 
$$\lim_{n\to \infty}X(T_{U_n}^{\sigma})=X(T),$$
but then $X(T) \in \overline{U_n}$, and hence 
$$X(T)\in K.$$
Therefore, $T_{K}^{\sigma}\leq T$, and this implies $T_{K}^{\sigma}=T$. But we know that 
$T$ is an $(\F_t)$ stopping time, and hence so is $T_K^{\sigma}$.}}

Here are is a general result about stopping times.

\prop{\begin{enumerate}[(a)]
       \item If $S$ and $T$ are TLG$^*$ valued stopping times, so is $S\vee T$.
       \item If $(T_n)$ is a sequence of stopping times then $\vee_{n=1}^{\infty} T_n$ is also a stopping time. 
      \end{enumerate}
\dok{We have 
$$(\vee_{n=1}^{\infty} T_n \preceq t)=\bigcap_{n=1}^{\infty}(T_n\preceq t),$$
and the claim follows. The case (a) is proved similarly.}}

\begin{figure}[ht]
\begin{center}
\psfrag{a}{$\boldsymbol{t_0}$} 
\psfrag{b}{$\boldsymbol{t_1}$}
\psfrag{c}{$\boldsymbol{t_2}$}
\psfrag{d}{$\boldsymbol{t_3}$}
\includegraphics[width=10cm]{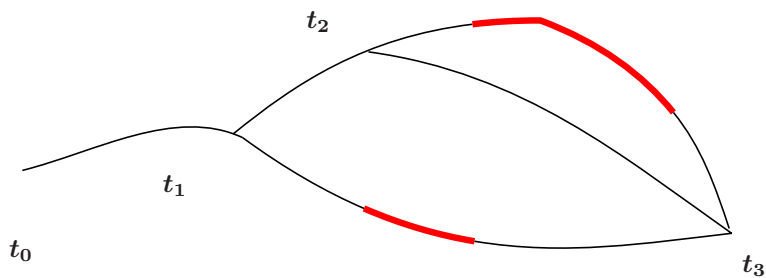}\\
\caption{ The red part of the graph is the part where Brownian bridges are defined,
everywhere else we set the process to be 0.} \label{pic:sl39}
\end{center}
\end{figure}

On the other hand, unlike in the classical case, the \textbf{minimum of two stopping times is not 
a stopping time}. The following example will illustrate that. Let $\cG$ be a TLG$^*$ 
like in Figure \ref{pic:sl39}. where $t_0=0$, $t_1=1/3$, $t_2=1/2$, $t_3=1$ and $\sigma_1$
is the bottom time-path, $\sigma_2$ the middle time-path, and  $\sigma_3$ the upper time-path. 
Let $B_{a,b}^{br}$ represent the Brownian bridge starting at $a$ and ending at $b$, 
and set 
$$X_{\sigma_1}(t)\stackrel{d}{=}\left\{\begin{array}{cc}
                          B^{br}_{1/2,2/3}(t) & t\in [1/2,2/3]\\
			  0	& {\rm otherwise}
                         \end{array}
 \right. \quad \textrm{and}\quad X_{\sigma_3}(t)\stackrel{d}{=}\left\{\begin{array}{cc}
                          B^{br}_{2/3,3/4}(t) & t\in [2/3,3/4]\\
			  0	& {\rm otherwise}
                         \end{array}
 \right. ,$$
and let $X_{\sigma_2}\stackrel{d}{=}0$ (it can  be any other Markov process consistent with the distributions
of $X_{\sigma_1}$ and $X_{\sigma_3}$). 

Let $U=(1,\infty)$, and $T_1=H_U^{\sigma_1}$
and $T_3=H_U^{\sigma_3}$. It is clear that $T_3\wedge T_1$ equals $t_1$ with probability greater than $0$. But the event
$$(T_3\wedge T_1=t_1)$$
depends on events that happen after time $1/3$, and it will not be contained in 
$\F_{t_1}$.{Stopping times|)}

\section{A simple coupling and branching process}
In this section we will describe a simple coupling and branching process.

We are reconstructing the movement of two persons/objects, and we have the following information
\begin{itemize}
 \item 2 persons moving around;
 \item (time $t_0=0$) started at the \textbf{same time from point} $A\in \R^2$;
 \item (time $t_3=1$) stooped  at the \textbf{same time in point} $B\in \R^2$;
 \item we have an {\it additional information} that \textbf{from time
$t_1=1/3$ to time $t_2=2/3$} they were moving \textbf{together}.
\end{itemize}
Note, that we only know that the two persons were together in \textbf{time} interval $[1/3,2/3]$,
but we don't know anything about the \textbf{locations} they visited together!

We will model this as a process on a TLG. Let, $\cG=(\cV,\cE)$ be given by
$$\cV=\{t_0,t_1,t_2,t_3\}, \quad \cE=\{E^1_{01},E^2_{01},E_{12}, E^1_{23},E^2_{23}\}.$$

\begin{figure}[ht]
\begin{center}
\psfrag{0}{$\boldsymbol{t_0}$} 
\psfrag{1}{$\boldsymbol{t_1}$}
\psfrag{2}{$\boldsymbol{t_2}$}
\psfrag{3}{$\boldsymbol{t_3}$}
\includegraphics[width=12cm]{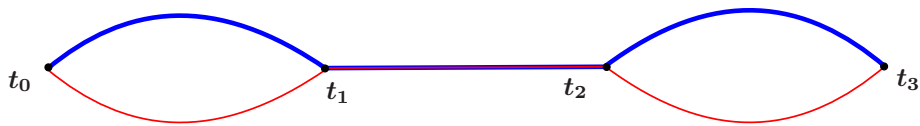}\\
\caption{ The blue path is indexing the movement of the first person, and
red path the movement of the second person.} \label{pic:pr_6}
\end{center}
\end{figure}

Let $\sigma_j=(E_{01}^j,E_{12},E_{23}^j)$ for $j=1,2$. ($\sigma_1$ is the blue path,
and $\sigma_2$ is the red path on Figure \ref{pic:pr_6}.)

Now we set $\cP$ to be the the distribution of the (two-dimensional) Brownian bridge from
$A$ to $B$ on $[0,1]$ with variance $\sigma^2$. We define $X$ on $\cG$ to be a natural $\cM_{\cP}$-process.

\begin{figure}[ht]
\begin{center}
\includegraphics[width=7cm]{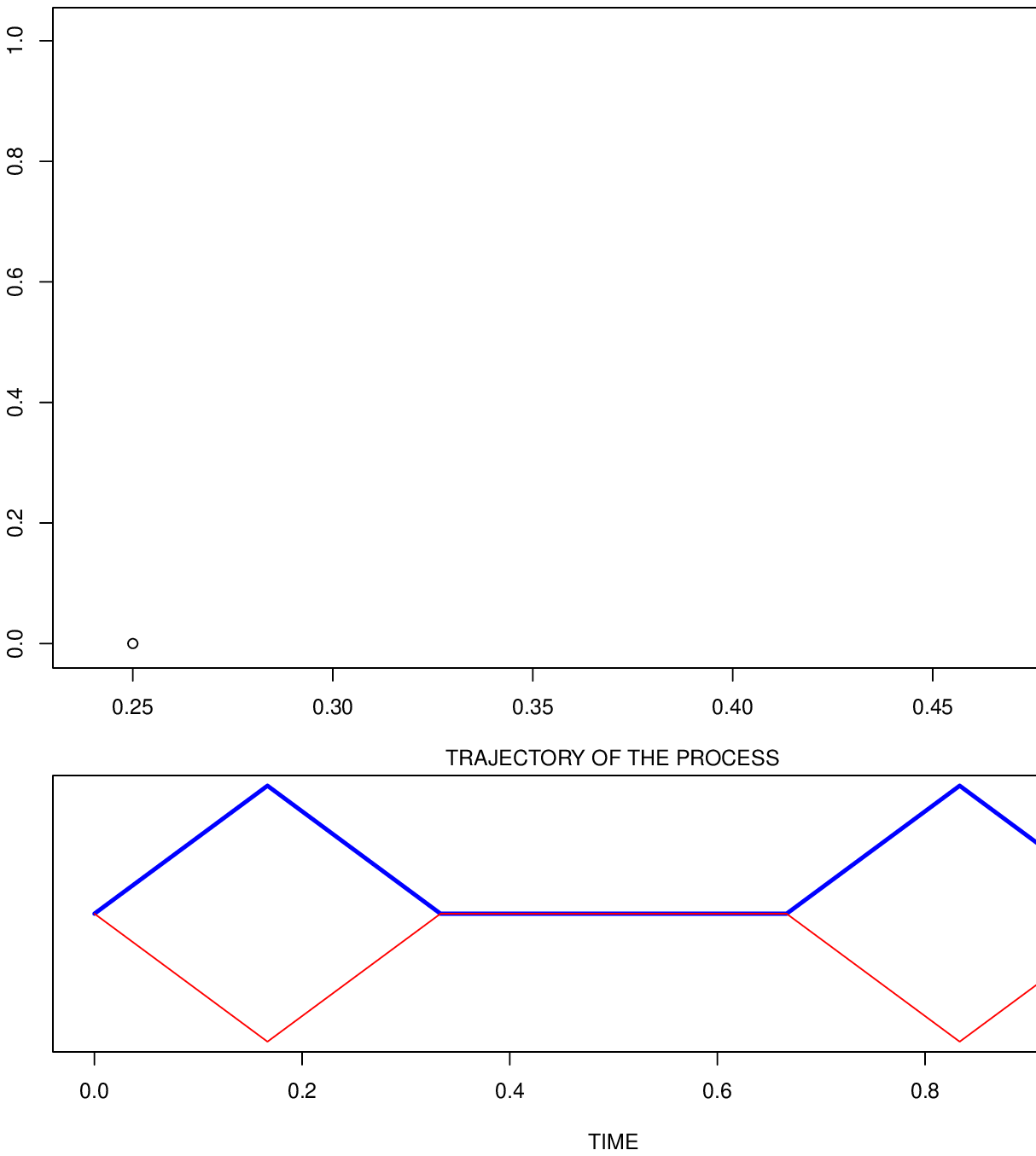}\quad \includegraphics[width=7cm]{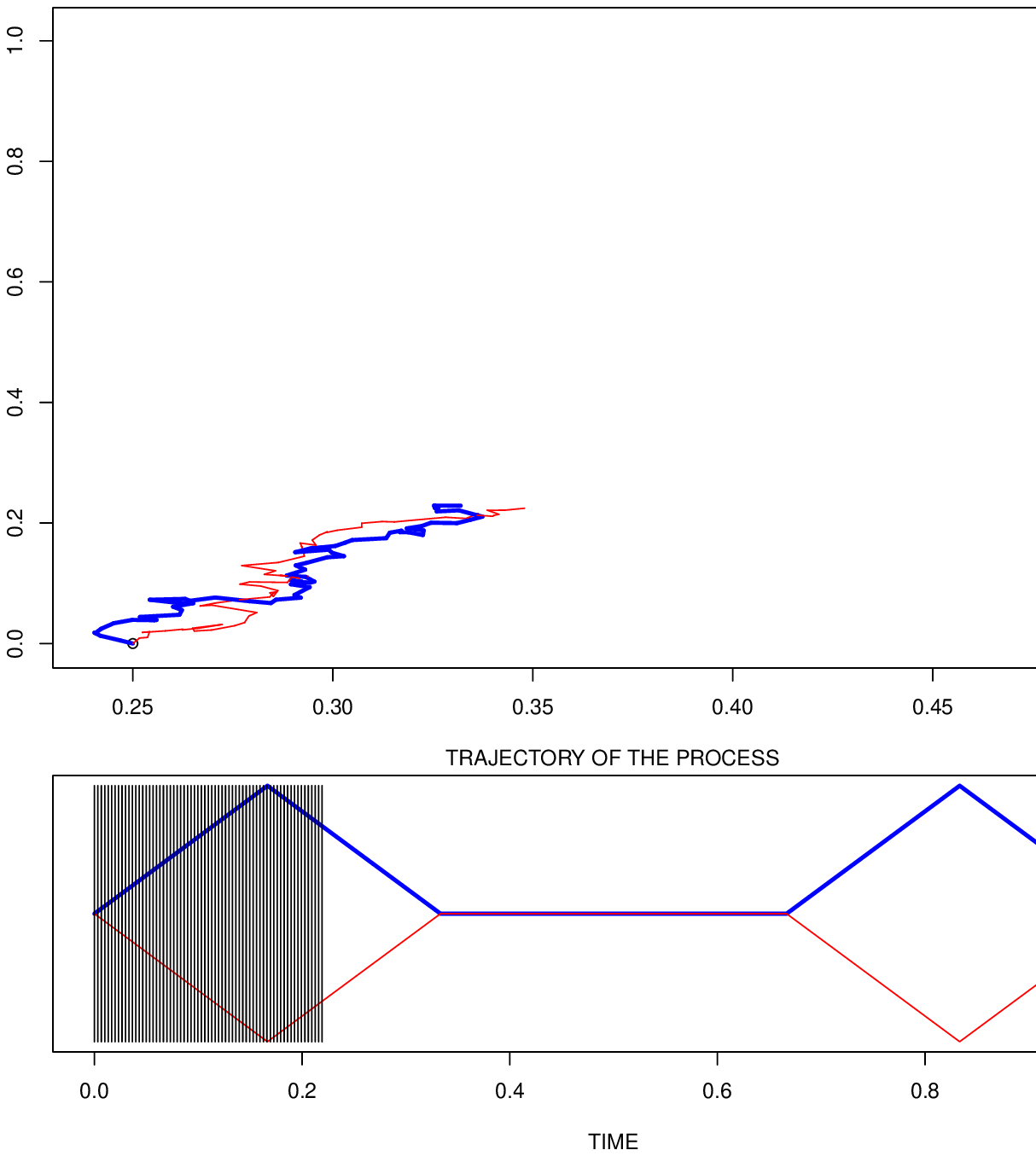}\\
\includegraphics[width=7cm]{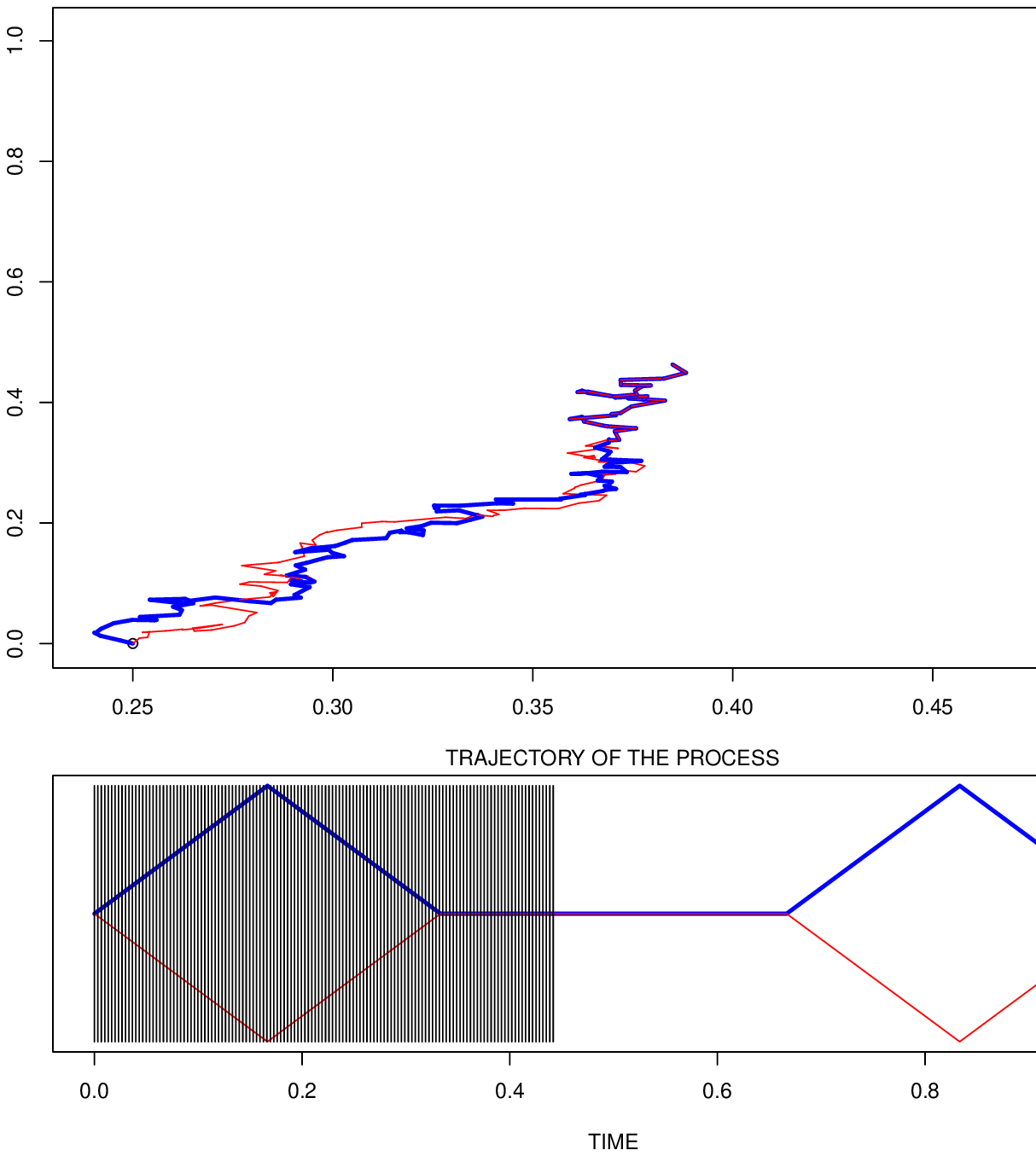}\quad \includegraphics[width=7cm]{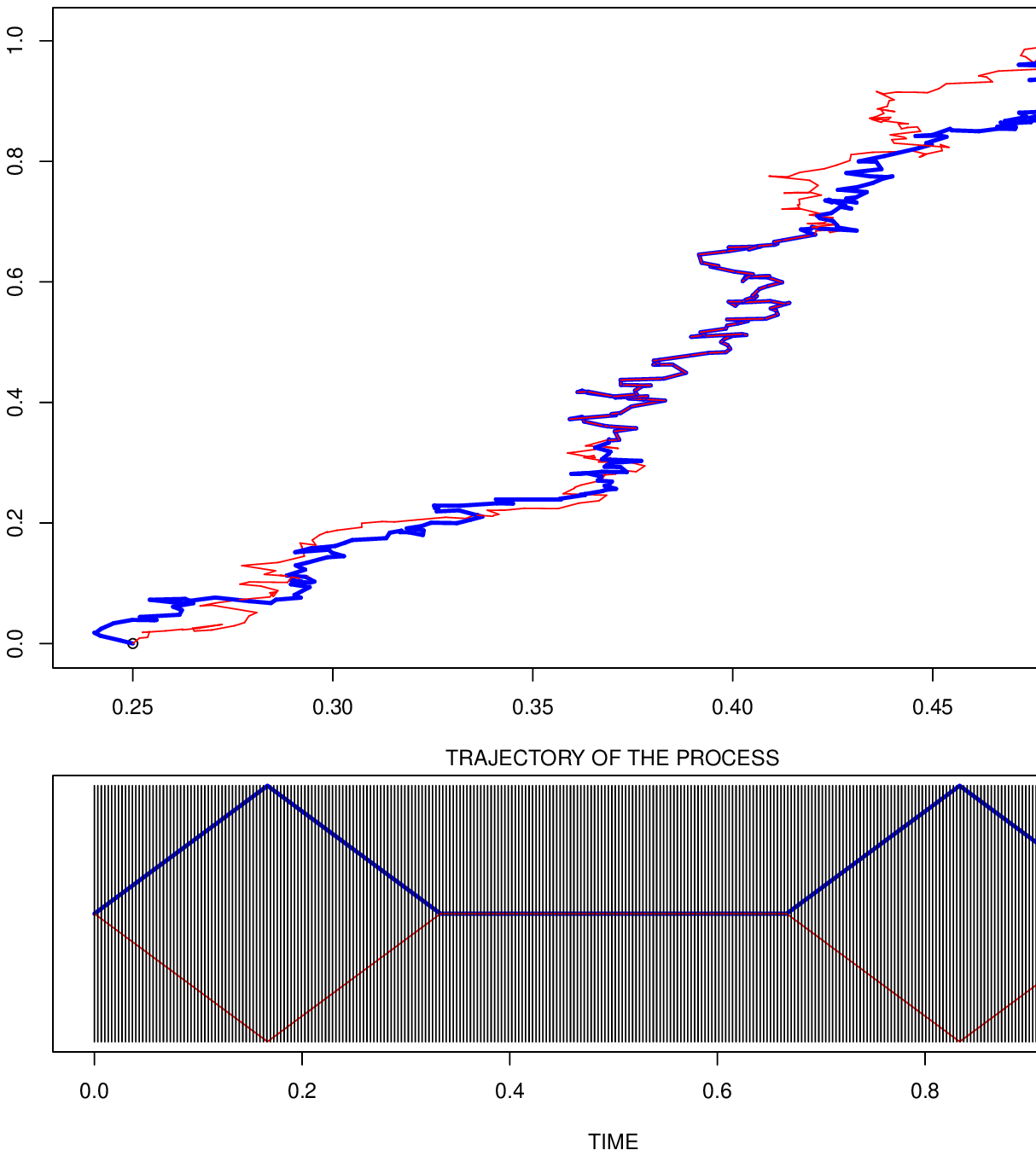}\\
\caption{Simulation of the simple coupling and branching process with $\sigma^2=0.005$.}\label{fig:2prmv}
\end{center}
\end{figure}

Note, that $X_{\sigma_1}$ and $X_{\sigma_2}$ are Brownian bridges from $A$ to $B$ with variance 
$\sigma$, and $X_{\sigma_1}|_{[t_1,t_2]}=X_{\sigma_2}|_{[t_1,t_2]}$. Figure \ref{fig:2prmv} shows a simulation of 
such a process.

Further, for this model we can calculate the expectations
$$\E(X_{\sigma_j}(t))=\left(tx_a+(1-t)x_b,ty_a+(1-t)y_b\right),\quad j=1,2.$$

Also, it is not hard to calculate the covariance structure. The two processes have a known covariance structure 
$${\bf Cov} (X^l_{\sigma_j}(\tau_1),X^l_{\sigma_j}(\tau_2))=\sigma^2\tau_1(1-\tau_2), \quad j,l=1,2,\ \tau_1\leq \tau_2.$$
Since the all full-time paths have the same distribution we have that
for $\tau_1\leq t_2$, and $t_1\leq \tau_2$
$${\bf Cov} (X^l_{\sigma_j}(\tau_1),X^l_{\sigma_i}(\tau_2))={\bf Cov} (X^l_{\sigma_j}(\tau_1),X^l_{\sigma_j}(\tau_2)),\quad i\neq j.$$
The last case is when $\tau_1\leq \tau_2$ are on different sides of a cell:
$${\bf Cov} (X^l_{\sigma_j}(\tau_1),X^l_{\sigma_i}(\tau_2))=\left\{
\begin{array}{cl}
\sigma^2\frac{\tau_1\tau_2(1-t_1)}{t_1} & \tau_1,\tau_2\in [t_0,t_1],\\
\sigma^2\frac{(1-\tau_1)(1-\tau_2)t_2}{1-t_2} & \tau_1,\tau_2\in [t_2,t_3].\\
\end{array}
\right.$$

\part{Natural Brownian motion and the stochastic heat equation}

In this part we study what happens (in some special cases) when the 
the process is indexed by a time-like graph whose representation is dense in 
(a subset of) the $t$-$x$ plane. 

\begin{figure}[ht]
\begin{center}
\psfrag{t}{\small $\boldsymbol{t}$ }
\includegraphics[width=5cm]{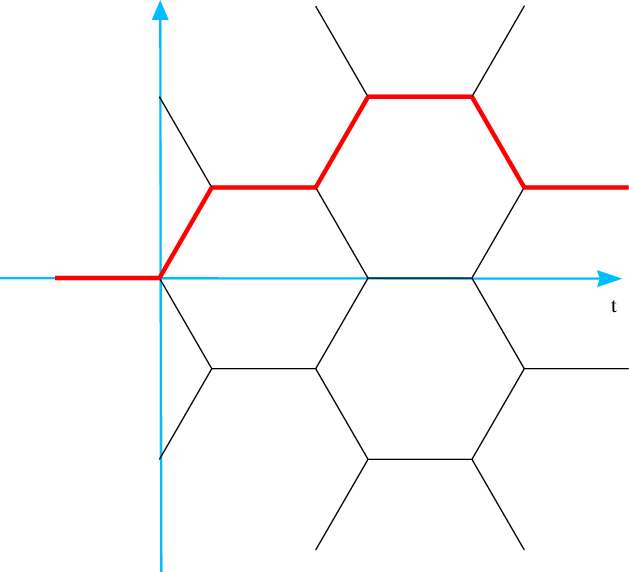}\quad
\psfrag{x}{\small $\boldsymbol{x}$ } 
\psfrag{a}{\tiny \textcolor{red}{$\boldsymbol{n^{-1}}$}} 
\psfrag{b}{\tiny \textcolor{blue}{$\boldsymbol{n^{-\frac{1}{2}-\alpha}}$}} 
\includegraphics[width=10cm]{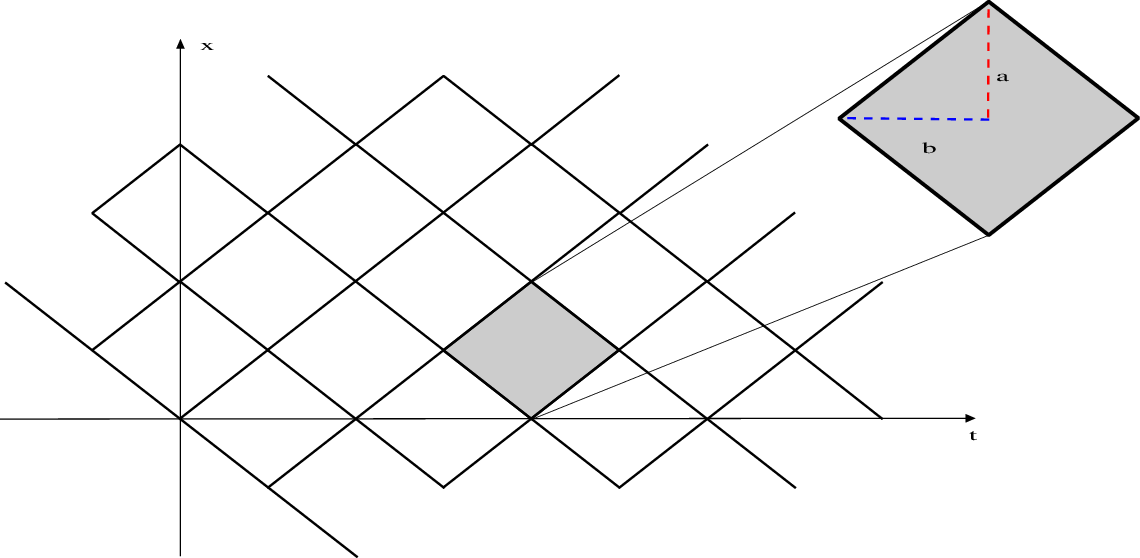}\\
\caption{Index set: Honeycomb graph and $\alpha$-rhombus grid} \label{sl22:1}
  \end{center}
\end{figure}

We will restrict our process to be a natural two-sided 
Brownian motion indexed by the graph whose representation is a rhombus grid. Burdzy and Pal studied 
the same process indexed by a honeycomb graph and found that (under certain scaling) when the 
mesh size goes to zero, the covariance structure is non-trivial (see Theorem 6.1. in \cite{tlg1}). (See Figure \ref{sl22:1}.)

 \begin{figure}[ht]
\begin{center}
\includegraphics[width=7cm]{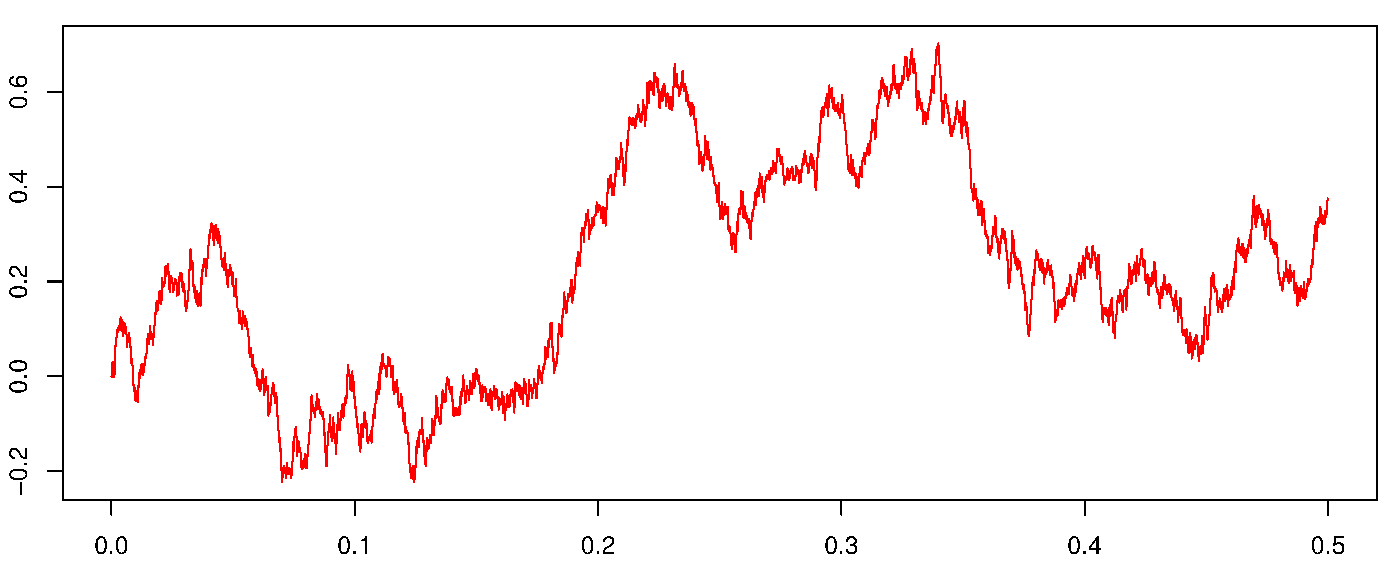}\quad 
\includegraphics[width=7cm]{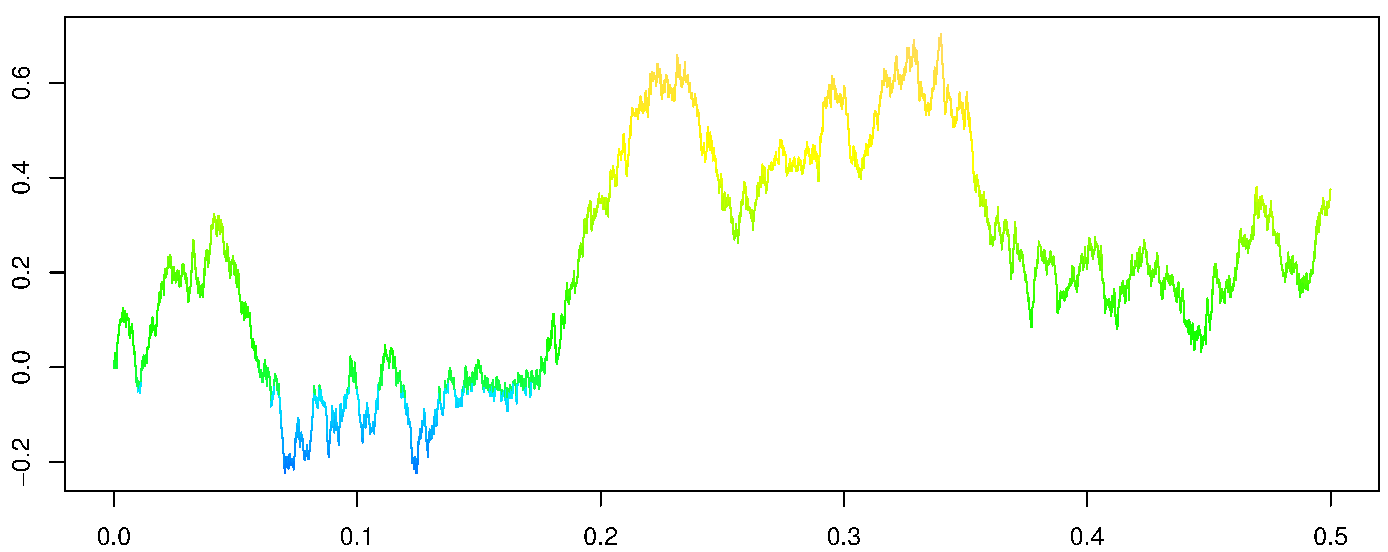}
\caption{We use topographical colors to represent values of the Brownian motion}
\psfrag{x}{\small $\boldsymbol{x}$ } 
\psfrag{a}{$\boldsymbol{\delta}$} 
\psfrag{b}{$\boldsymbol{f(\delta)}$ } 
\psfrag{t}{\small $\boldsymbol{t}$ }
\includegraphics[width=7.5cm]{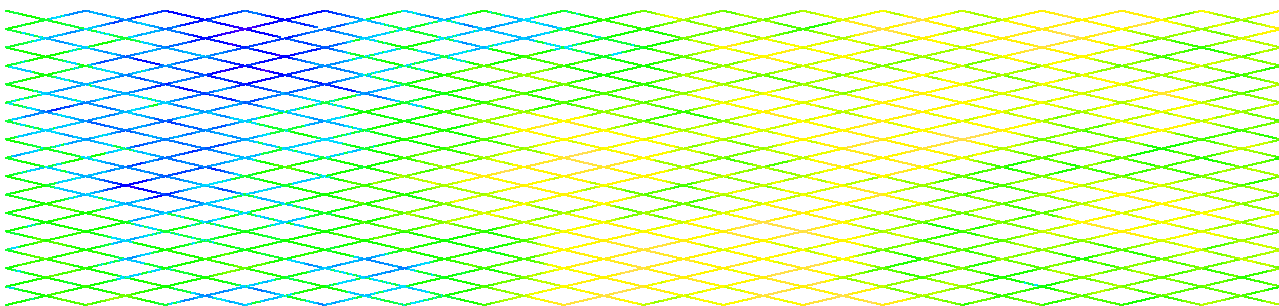}\quad \includegraphics[width=7.5cm]{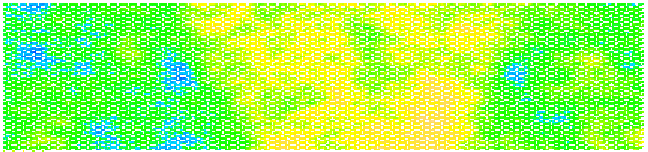}\\[0.2cm]
\includegraphics[width=7.5cm]{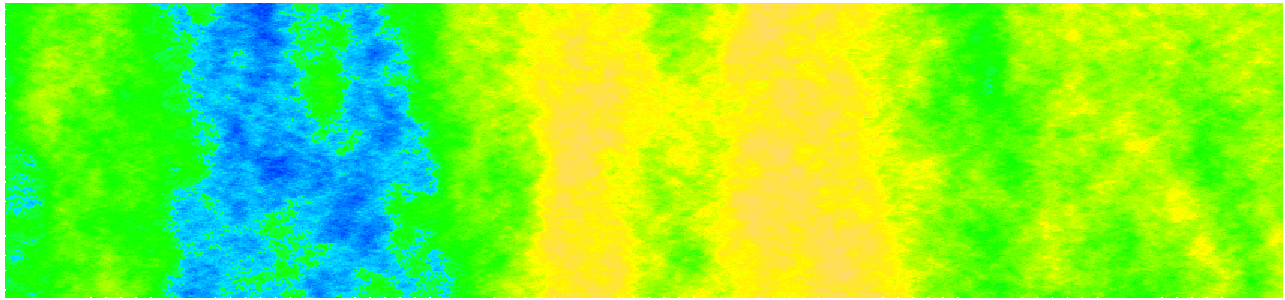}\quad \includegraphics[width=7.5cm]{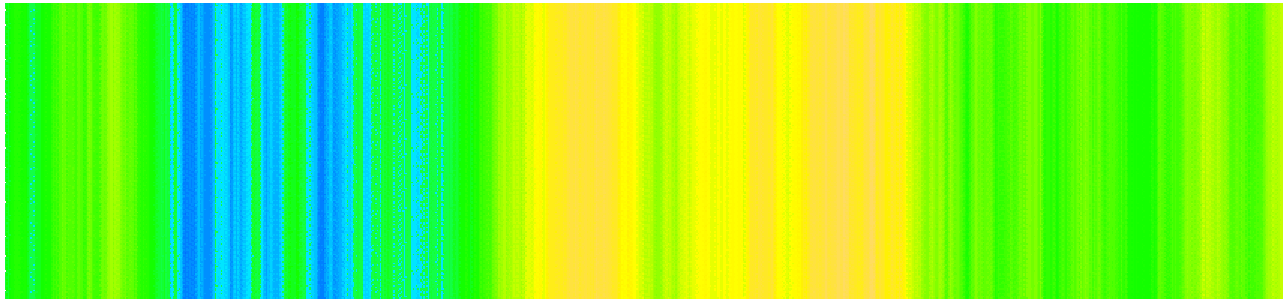}
\caption{Scaled simulation of the process when the ratio is $1/n:1/n$, for $n=32,128,512,\infty$.} \label{sim_net}
  \end{center}
\end{figure}

The images in Figure \ref{sim_net} show what happens when ratio of the half-diagonals is $n^{-1/2-\alpha}:n^{-1}$ for $\alpha>0$.
It turns out, in this case, the process in the limit only depends on the time coordinate ($t$)
and not on the space coordinate ($x$).

\begin{figure}[ht]
\begin{center}
\includegraphics[width=10cm]{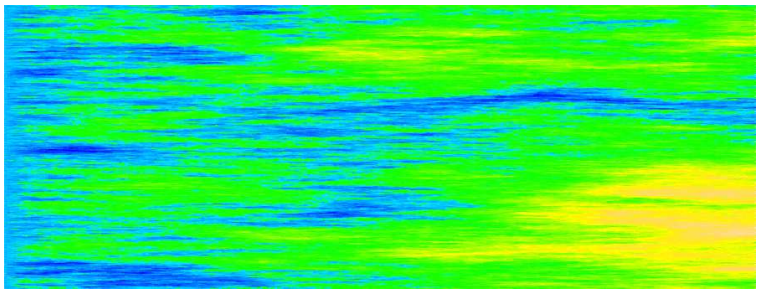}
\caption{Simulation of a natural Brownian motion indexed by a rhombus grid where the ratio
is $n^{-1/2}:n^{-1}$ and $n=1024$.}
\label{heat_eq:im}
\end{center}
\end{figure}

For the limit case $\alpha=0$, however, the simulation (see Figure \ref{heat_eq:im}) indicates that 
the structure of the process in the plane is more complex. It turns out that the process in the limit
is the stochastic heat equation.\vspace{0.2cm}

In this part we prove that these are the limits. 

First, we introduce some results about maximums of Gaussian processes in Chapter \ref{sec:03}.
Then in Chapter \ref{chp:heatrw}, we prove some general results about the (stochastic) heat equation, its 
approximation by Euler's method and the connections to the random walk.\vspace{0.1cm}

After developing those tools, in Chapter \ref{sec:04}. we prove the claims stated in this introduction
in Theorem \ref{teo:eul_grd}.

\chapter{Maximums of Gaussian processes }\label{sec:03}

In this section we will review the
\begin{itemize}
 \item bounds for the second moment of the maximum of a finite sequence of independent Brownian bridges
 \item bounds for the second moment of the maximum of a finite sequence of (not necessarily independent) 
normal random variables;
 \item concentration of the maximum of Gaussian random element in $C(K)$ for some 
compact set $K$. 
\end{itemize}


\section{Sequence of Brownian bridges}\index{Brownian bridge!maximum|(}

For $k =1,2,\ldots$ we will denote $(B^{br}_k(t):t\in[0,1])$ a Brownian bridge starting and 
ending at $0$. (See Definition \ref{bbridg}.)  We are interested in getting some estimation on  moments of $$M_n:=\sup\{|B^{br}_k(t)|:t\in [0,1], k=1,2,\ldots, n\}.$$ 

In order to do this, we will estimate the moments of 
$$M^+_n:=\sup\left\{B^{br}_k(t):t\in [0,1], k=1,2,\ldots, n\right\}$$
$$M^-_n:=\inf\left\{B^{br}_k(t):t\in [0,1], k=1,2,\ldots, n\right\}=-\sup\left\{-B^{br}_k(t):t\in [0,1], k=1,2,\ldots, n\right\}.$$

Since, $-B^{br}_k$ has the same distribution as $B^{br}_k$, it follows that $M^+_n$ has the same
distribution as $-M_n^-$.

Now, since $B^{br}_k(0)=B^{br}_k(1)=0$, $M_n^+>0$ and $M_n^-<0$. Further, it is clear that 
$$M_n=\max\{M_n^+,-M_n^-\}.$$ 

So if we find, a bound on moments of $M_n^+$ we will be able to find a bound
on the moments of $M_n$.

From classical results on boundary crossing probabilities for Brownian motion
(see \cite[Karatzas, Shreve]{karatzas}, page 262-265), we have the following Lemma:

\lem{If $(W_t:t\geq 0)$ is a Brownian motion starting at $0$, then 
\be\P\left(\max_{0\leq t\leq T} W_t\geq \beta|W_T=a \right)=e^{-2\beta(\beta -a)/T} \label{bb:max:1}\ee
for $T>0$ and $\beta >\max\{0,a\}$.}

From the last Lemma we get what we need to calculate $\E(M_n^{+2})$.

\prop{\begin{enumerate}[(a)]
      \item For the Brownian bridge $B^{br}_k$ we have
$$\P\left(\max_{t\in [0,1]} B^{br}_k(t)>\beta \right)=e^{-2\beta^2}.$$
      \item If $(B^{br}_k)$ are independent Brownian bridges, the following equality holds:
\be 4\E(M_n^{+2})=\frac{1}{1}+\frac{1}{2}+\ldots +\frac{1}{n}.\label{bb:max:3}\ee
     \end{enumerate} 
\dok{(a) This follows from (\ref{bb:max:1}) when we set $T=1$, and $a=0$. 
(b) For this we first note that, the independence of the sequence $(B^{br}_k)_{k=1}^n$ 
implies
\begin{align*}
 \P(M_n^+>\beta)& = 1 - P\left(M_n^+\leq \beta\right) = 1- \P\left(\bigcap_{k=1}^n \left\{\max_{t\in [0,1]} B^{br}_k(t) \leq\beta\right\}\right) \\
& = 1- \prod_{k=1}^n\P\left( \max_{t\in [0,1]} B^{br}_k(t) \leq\beta\right)= 1- \prod_{k=1}^n\left[1-\P\left( \max_{t\in [0,1]} B^{br}_k(t) >\beta\right)\right]\\
&\stackrel{(\ref{bb:max:1})}{=} 1-\left(1-e^{-2\beta^2}\right)^n.
\end{align*}
Now, we get 

\be \E(M_n^{+2})= \int_0^\infty\beta \P\left(M_n^+>\beta\right)\, d\beta =\int_0^\infty\beta \left[1-\left(1-e^{-2\beta^2}\right)^n\right]\, d\beta \label{bb:max:2}\ee
Now, we will use some simple algebra,
\begin{align*}
\beta (1-(1-e^{-2\beta^2})^n) &=\beta e^{-2\beta^2 }\frac{1-(1-e^{-2\beta^2})^n}{1-(1-e^{-2\beta^2 })}\\
& = \beta e^{-2\beta^2}\left[1+(1-e^{-2\beta^2})+\ldots +(1-e^{-2\beta^2})^{n-1}\right]. 
\end{align*}
Using Tonelli's Theorem the equality (\ref{bb:max:2}) becomes
$$
 \E(M_n^{+2})=\int_0^\infty \sum_{k=0}^{n-1} \beta e^{-2\beta^2}(1-e^{-2\beta^2})^k=\sum_{k=0}^{n-1} \int_0^\infty \beta e^{-2\beta^2}(1-e^{-2\beta^2})^k= \frac{1}{4}\sum_{k=0}^{n-1}\frac{1}{k+1},
$$
since the expression under the last integral is the derivative of $$-\frac{(1-e^{-2\beta^2})^{k+1}}{4(k+1)}.$$ }}

Recall, that for the harmonic sequence\index{Harmonic sequence} $$H_n=\sum_{k=1}^n\frac{1}{k}$$ we have,
\be\ln n\leq H_n\leq \ln (n+1).\label{bb:max:4}\ee

\pos{\label{cor:bbmax}The following inequalities hold for $M_n$ the maximum of $n$ independent Brownian bridges 
\be\E(M_n)\leq \sqrt{\ln (n+1)}. \label{bb:max:5}\ee
\be \E(M_n^2)\leq \frac{1}{2}\ln (n+1) \label{bb:max:5a}\ee
\dok{ From $(\ref{bb:max:3})$, using $(\ref{bb:max:4})$ we get 
$$\E(M_n^+)\leq\sqrt{E(M_n^{+2})}\leq \frac{1}{2}\sqrt{\ln (n+1)}.$$
Now, using the fact that $\max\{a,b\}\leq a+b$ for $a,b\geq 0$ we have, since
$M_n^+$ and $-M_n^-$ have the same distribution we get
$\E(M_n)\leq \E(M_n^+)+\E(-M_n^-)=2\E(M_n^+)$. The inequality $(\ref{bb:max:5})$ follows. The inequality $(\ref{bb:max:5a})$ follows
in the similar way.}}

\index{Brownian bridge!maximum|)}
\section{Sequence of normal variables}\index{Maximal inequality for normal random variables}

When we have several normal random variables, what can we say
about the expectation of the maximum of these random variables?\vspace{0.2cm}

Let $X_k\sim N(0,\sigma_k^2)$ for $k=1,\ldots,n $, we are interested in
the upper bounds for
$$\E(\max\{|X_1|,\ldots, |X_n|\}),$$
and 
$$\E(\max\{|X_1|^2,\ldots, |X_n|^2\}).$$

Using the result we got for the Brownian bridge we can get the upper bound.

\prop{\label{prop:max_id}For $(X_k)_{k=1}^n$ independent we have   
\be\E(\max\{|X_1|,\ldots, |X_n|\})\leq 2\max_{1\leq k\leq n} \sigma_k \sqrt{\ln(n+1)}\label{bb:max:6}\ee
\be\E(\max\{|X_1|^2,\ldots, |X_n|^2\})\leq 2\left(\max_{1\leq k\leq n} \sigma_k^2\right) \ln(n+1)\label{bb:max:6a}\ee
\dok{Let $\sigma_{max}=\max_{1\leq k\leq n} \sigma_j>0$ and $X_k':=\frac{1}{2\sigma_{max}}X_k$.
Now $X_k'$ are normal random variables with expectation $0$ and variance at most
$1/4$. Hence, for each $k=1,\ldots, n$ there exists a time $t_k$ such that the variance
of $B^{br}_k(t_k)$ is the same as that of $X_k'$, where $(B^{br}_k)_{k=1}^n$ is a sequence
of independent Brownian bridges starting and ending at $0$. Hence, $X_k'$ has the same distribution 
as $B^{br}_k(t_k)$, so the distribution of $\max\{|X_1'|,\ldots, |X_n'|\}$ is the same
as of $\max\{|B^{br}_1(t_1)|,\ldots, |B^{br}_n(t_n)|\}$ and this is less than 
$M_n=\sup\{|B^{br}_k(t)|:t\in [0,1], k=1,\ldots, n\}$. So by $(\ref{bb:max:5})$, we 
have 
$$\E(\max\{|X_1'|,\ldots, |X_n'|\})\leq \E(M_n)\leq \sqrt{\ln (n+1)}.$$ 
Multiplying this with $2\sigma_{max}$ we get (\ref{bb:max:6}). Similar argument 
using inequality (\ref{bb:max:5a}) will give (\ref{bb:max:6a})}}

Now, we will deal with the case when $(X_k)_{k=1}^n$ are not necessarily
independent.
We will do this with the help of a lemma that is due to \v{S}id\'ak (see \cite{sidak}).

\lem{\label{lem:sidak}\index{Sidak@\v{S}id\'ak's inequality}\begin{enumerate}[(a)]
      \item (\v{S}id\'ak 1967.) For positive numbers $c_1,c_2,\ldots, c_n$ $$\P(|X_1|\leq c_1,|X_2|\leq c_2,\ldots |X_n|\leq c_n)\geq \P(|X_1|\leq c_1)\P(|X_2|\leq c_2)\ldots \P(|X_n|\leq c_n).$$
      \item Let $Y_1,\ldots, Y_n$ be independent random variables, such that for each $k=1,2,\ldots,n$
      $Y_k$ and $X_k$ have the same distribution, then
      $$\E(\max\{|X_1|,\ldots, |X_n|\}^p)\leq \E(\max\{|Y_1|,\ldots, |Y_n|\}^p)$$
      for all $p\geq 1$.
     \end{enumerate}
}

The proof of part (a) of this Lemma can be found in \cite{sidak} or \cite{gausproc}.
Part (b) is a direct consequence of part (a).

As a consequence of Lemma \ref{lem:sidak} (part (b) for $p=1$)  and Proposition \ref{prop:max_id}  we get the following theorem.

\teo{\label{teo:max_nnc}For $(X_k\sim N(0,\sigma_k))_{k=1}^n$ (possibly correlated) we have   
\be\E(\max\{|X_1|,\ldots, |X_n|\})\leq 2\max_{1\leq k\leq n} \sigma_k \sqrt{\ln(n+1)}\label{bb:max:7}\ee
\be\E(\max\{|X_1|^2,\ldots, |X_n|^2\})\leq 2\left(\max_{1\leq k\leq n} \sigma_k^2\right) \ln(n+1)\label{bb:max:7a}\ee
}

A much general result (of the same order) can be found in the paper by Chatterjee in \cite{chatterjee}.

\section{Some concentration and convergence results}
We will shortly state some concentration results taken from 
Chapter 3.1. of Talagrand's book \cite{talagrand1}, and apply it
to the convergence of Gaussian processes.\vspace{0.3cm}

Let $B$ be a Banach space\index{Banach space}, and $D$ some countable subset of the unit ball
of the dual space $B'$\index{Banach space!dual space} such that
$$\|x\|=\sup_{f\in D}|f(x)|,$$
for all $x\in B$.

We say that $X$ is a {\bf Gaussian random variable}\index{Gaussian random variable (in Banach space)} in $B$ if 
$f(X)$ is measurable for every $f\in D$ and if every finite linear combination 
$$\sum_i\alpha_i f_i(X),$$
where $\alpha_i\in \R$ and $f_i\in D$, is Gaussian.

Let $X$ be a Gaussian, $M=M(X)$ be the {\bf median} of $\|X\|$, that is $M$ has the property
that $$\P(\|X\|\geq M)\geq 1/2\quad \textrm{and}\quad \P(\|X\|\leq M)\geq 1/2.$$
Further, set the {\bf supremum of weak deviations} to be
$$\sigma=\sigma(X)=\sup_{f\in D}\E[f(X)^2]^{1/2}.$$
The following result is a Lemma 3.1. from \cite{talagrand1}.\index{Gaussian random variable (in Banach space)!concentration inequality}

\lem{\label{lm31:conct}Let $X$ be a Gaussian with median $M=M(X)$ and the supremum of weak deviations $\sigma=\sigma(X)$, then 
$$\P(|\|X\|-M|>t)\leq \exp(-t^2/2\sigma^2).$$ }

\pos{Let $X$ be a Gaussian with $\E[\|X\|^2]<\infty$, then 
\be\P(\|X\|>t)\leq 4\exp\left(\frac{-t^2}{2\E[\|X\|^2]}\right).\label{conc:ineqgs} \ee
\dok{Follows from Lemma \ref{lm31:conct} and the fact that $\sigma^2\leq \E[\|X\|^2]$ 
and $M^2\leq \E[\|X\|^2]$. 
}
}

Let $Y=(Y_1,Y_2,\ldots, Y_d)$ be a Gaussian vector with expectation vector 0, then since $B=\R^d$ is the Banach
space with usual norm, and the set of projections $D=\{\pi_k:k\in \{1,2,\ldots, d\}\}$
is a subset of $B'$ we have that 
\be\P(\max_{1\leq k\leq n}|Y_k|>t)\leq 4\exp\left(\frac{-t^2}{2\E[\max_{1\leq k\leq n}|Y_k|^2]}\right).\label{eq:cncmg}\ee

It will not always be easy to get an estimate for $\E[\max_{1\leq k\leq n}|Y_k|^2]$, but
when we do the inequality $(\ref{eq:cncmg})$  will tell us a lot.  

\teo{\label{teo:cncrt}Let $Y^k$ be a sequence of Gaussian vectors (not necessarily of the same 
dimension) with expectation vector 0 on the same probability space, such that
\be\E[\|Y^k\|_{\infty}^2]\leq \frac{C}{k^\alpha},\label{eq:cnvrt}\ee
for some $C>0$ and $\alpha>0$. Then 
$$\|Y^k\|_{\infty}\to 0\quad a.s.$$
\dok{It is easy to show that for sufficiently large $k$ we have
$$\exp\left(-\frac{t^2k^\alpha}{2C}\right)\leq \frac{1}{k^2}.$$
Now, using $(\ref{eq:cncmg})$ and $(\ref{eq:cncmg})$ the previous inequality 
for sufficiently large $k$ gives 
$$\P(\|Y^k\|_{\infty}>t)\leq \frac{1}{k^2},$$
and hence 
$$\sum_{k=1}^\infty \P(\|Y^k\|_{\infty}>t)<\infty.$$
Since this holds for all $t>0$, $\|Y^k\|_{\infty}\to 0$ a.s.} }

\pos{\label{cor:cncrt}Let $Y^k$ be a sequence of Gaussian vectors (not necessarily of the same 
dimension) with expectation vector 0 on the same probability space, such that
\be\E[\|Y^k\|_{\infty}^2]\leq \frac{C}{k^\alpha},\label{eq:cnvrt2}\ee
for some $C>0$ and $\alpha>0$. Then for $0<\beta<\alpha/2 $
$$k^\beta\|Y^k\|_{\infty}\to 0\quad a.s.$$
\dok{We see that for $Z^k:=k^\beta Y^k$ we have 
$$\E[\|Z^k\|_{\infty}^2]\leq \frac{C}{k^{\alpha-2\beta}}.$$
Hence, since $\alpha-2\beta>0$ by Theorem \ref{teo:cncrt} we have $\|Z^k\|_{\infty}\to 0$ a.s.}}

We can get similar results for continuous Gaussian fields. 
\teo{\label{teo:cnv0}Let $(\Omega,\F,\P)$ be a probability space, $K\subset \R^n$ be a compact
set, and $X_n:K\times \Omega\to \R$ have the following properties:
\begin{enumerate}[(1)]
 \item For each $x\in K$ $X_n(x)$ is a Gaussian random variable.
 \item For each $\omega \in \Omega$ $x\mapsto X_n(x,\omega)$ is a continuous function.
\end{enumerate}

 Then if  
$$\E[\|X_n\|_{\infty}^2]\leq \frac{C}{n^\alpha},$$
we have \be \|X_n\|_{\infty}\to 0\quad a.s.\ee
Further, for $0<\beta<\alpha/2$ we have 
\be n^{\beta}\|X_n\|_{\infty}\to 0\quad a.s.\ee
\dok{Let $B=C(K)$ with the usual $\|\cdot\|_{\infty}$ supremum norm, and set
$D=\{\pi_q\ :\ q\in \mathbb{Q}^n\cap K\}$. Now it is clear that $X_n$ is a Gaussian 
random variable in $B$, and since $X_n$ is continuous we have $\|X_n\|_{\infty}=\sup_{q\in \mathbb{Q}^n\cap K}|\pi_q (X_n)|$, we
have from $(\ref{conc:ineqgs})$ that 
$$\P(\|X_n\|_{\infty}>t)\leq \exp\left(-\frac{t^2k^\alpha}{2C}\right).$$
Using the same technique as in proof of Theorem \ref{teo:cncrt} we have the desired results.  }}

\chapter{Random walk and stochastic heat equation reviewed}\label{chp:heatrw}

\section{Modification of the Local Limit Theorem}\index{Local Limit Theorem|(}
In the rest of this chapter $(S_n)$ will denote the \textbf{simple random walk}\index{Simple random walk|textbf}\index{Simple random walk|seealso{Local Limit Theorem}}, where
$S_n=X_1+X_2+\ldots + X_n$, $(X_k)$ are i.i.d. and $\P(X_1=\pm 1)=1/2$.

First we introduce some notation. For a simple random walk $(S_n)$ we set
$$p_n^k(x)=\P\left(\frac{S_k}{\sqrt{n}}=x\right),\quad \textrm{for}\quad  x\in {\cal L}_n^k:=\{(k+2z)/\sqrt{n}:z\in \mathbb{Z}\},$$
and 
$$\rho_n^k(x)=\frac{1}{\sigma^k_n\sqrt{2\pi }}\exp\left(-\frac{x^2}{2(\sigma^k_n)^2}\right),$$
where $(\sigma_n^k)^2=\frac{n}{k}$. The main result of this section is Theorem \ref{thm:llt} which gives  
the bound  on the difference of $p_n^k$ and $\rho_n^k$.

We will need the following two lemmas. The first lemma is a consequence of the inversion formula for characteristic functions. (See \cite{durrett}.)

\lem{\label{lem:prchrf}If $Y$ is a random variable with $\P(Y\in a+\theta \mathbb{Z})=1$, and $\psi(t)=\E(e^{itY})$
is its characteristic function, then
$$\P(Y=x)=\frac{1}{2\pi/\theta}\int_{-\pi/\theta}^{\pi/\theta}e^{-itx}\psi(t)\, dt.$$
}

The second lemma is a consequence of the Stirling formula\index{Stirling formula}.
\lem{\label{lem:strlng}For $k\in \N$ set
$$I_k:=\int_0^{\pi/2} \cos^k(x)\, dx,$$
 there exists a $C>0$ such that
\be \left|\sqrt{k}I_k-\sqrt{\frac{\pi}{2}}\right|\leq \frac{C}{k}.\label{eq:lmllt1}\ee
\dok{From integration by parts we have 
$$I_k=\frac{k-1}{k}I_{k-2},$$
and further we can calculate $I_1=1$ and $I_2=\frac{\pi}{4}$.
Now, this recursion gives us
$$I_{2k+1}=\frac{2k}{2k+1}\cdot\frac{2k-2}{2k-1}\cdots \frac{2}{3}\cdot I_1=\frac{2^{2k}(k!)^2}{(2k+1)!},$$
$$I_{2k}=\frac{2k-1}{2k}\cdot \frac{2k-3}{2k-2}\cdots \frac{3}{4}\cdot I_2=\frac{(2k)!}{2^{2k}(k!)}\cdot \frac{\pi}{2}.$$
Using Stirling's Formula (see for example Gamelin \cite{gamelin} page 368), i.e. the fact that 
$$n!=\left(\frac{n}{e}\right)^n\sqrt{2n\pi} \exp\left(\frac{1}{12n}+O\left(\frac{1}{n^3}\right)\right),$$
we have
$$2^k\sqrt{2k+1}(k!)^2=\frac{(2k)^{2k+1}}{e^{2k}} \pi\sqrt{2k+1}\exp\left(\frac{1}{6k}+O\left(\frac{1}{8k^3}\right)\right),$$
$$(2k+1)!=\left(\frac{2k+1}{e}\right)^{2k+1}\sqrt{2\pi}\sqrt{2k+1}\exp\left(\frac{1}{12(2k+1)}+O\left(\frac{1}{8k^3}\right)\right),$$
and therefore
$$\sqrt{2k+1}I_{2k+1}=\left(1-\frac{1}{2k+1}\right)^{2k+1}e\sqrt{\frac{\pi}{2}}\exp(\frac{1}{6k}-\frac{1}{12(2k+1)}+O\left(\frac{1}{k^3}\right)).$$
Now, 
\begin{align}
\nonumber & (2k+1)|\sqrt{2k+1}I_{2k+1}-\sqrt{\pi/2}|\\
\nonumber=&(2k+1)\sqrt{\frac{\pi}{2}}\left|\left(1-\frac{1}{2k+1}\right)^{2k+1}e(1+\frac{1}{6k}-\frac{1}{12(2k+1)}+O\left(\frac{1}{k^2}\right))-1\right|\\
\leq &(2k+1)\sqrt{\frac{\pi}{2}}e\left|\left(1-\frac{1}{2k+1}\right)^{2k+1}-e^{-1}\right|+\sqrt{\frac{\pi}{2}}e\left|\frac{2k+1}{6k}-\frac{1}{12}+O\left(\frac{1}{k}\right)\right| \label{eq:instr1}
\end{align}
The second absolute value is clearly bounded. For the first absolute value we use the well-known fact that 
if $|u|,|z|\leq 1$ then for $m\in \N$ we have $|u^m-z^m|\leq m|u-z|$. So, by setting $m=2k+1$,
$u=1-(2k+1)^{-1}$ and $z=e^{-(2k+1)^{-1}}$ we have
$$\left|\left(1-\frac{1}{2k+1}\right)^{2k+1}-e^{-1}\right|\leq (2k+1)\left|1-\frac{1}{2k+1}-e^{-1/(2k+1)}\right|$$
$$\leq (2k+1)\frac{1}{2(2k+1)^2}=\frac{1}{2(2k+1)},$$
where the last inequality follows from the Taylor's Theorem. Hence, the first absolute
value in $(\ref{eq:instr1})$ is also bounded.\vspace{0.3cm}

Using the same methods we get the same result for $(I_{2k})$. }}
\teo{\label{thm:llt}There exists a $C>0$ such that for any $\beta (n)$ we have
$$\sup_{\beta(n)\leq k } \sup_{x\in {\cal L}_n^k}\left|\frac{n^{1/2}}{2}p_n^k(x)-\rho_n^k(x)\right|\leq \frac{C}{\pi}\sqrt{\frac{n}{\beta(n)^3}}, $$
for all $n$.
\dok{Using Lemma \ref{lem:prchrf}. for $\theta =2/\sqrt{n}$ and function
$$\psi_k(t)=\E\left[\exp(\frac{itS_k}{\sqrt{n}})\right]=\varphi^k\left(\frac{t}{\sqrt{n}}\right),$$
we have that 
\be\frac{n^{1/2}}{2}p_n^k(x)=\frac{1}{2\pi}\int_{-\pi\sqrt{n}/2}^{\pi\sqrt{n}/2}e^{-itx}\varphi^k(t/\sqrt{n})\, dt.\label{eq:llt1}\ee
The inversion formula gives that
\be\rho_n^k(x)=\frac{1}{2\pi}\int_{\R}e^{-itx}\exp(-(\sigma_{n}^k)^2t^2/2)\, dt.\label{eq:llt2}\ee 
From $(\ref{eq:llt1})$ and $(\ref{eq:llt1})$ we have
\begin{align*}
 \left|\frac{n^{1/2}}{2}p_n^k(x)-\rho_n^k(x) \right| & \leq \frac{1}{2\pi}\int_{-\pi\sqrt{n}/2}^{\pi\sqrt{n}/2}|\varphi^k(t/\sqrt{n})-\exp(-(\sigma_{n}^k)^2t^2/2)|\, dt\\
&+\frac{1}{\pi}\int_{\pi\sqrt{n}/2}^{\infty}\exp(-(\sigma_{n}^k)^2t^2/2)\, dt.
\end{align*}\index{Simple random walk|)}
%
First, note that the right side the inequality doesn't depend on $x$.  Now by substituting
$u=t\sqrt{k/n}$ in both integrals, we get that the right side of the inequality is 
\be\frac{\sqrt{n/k}}{\pi}\left[\frac{1}{2}\int_{-\pi\sqrt{k}/2}^{\pi\sqrt{k}/2}|\varphi^k(u/\sqrt{k})-\exp(-u^2/2)|\, du+\int_{\pi\sqrt{k}/2}^{\infty}\exp(-u^2/2)\, du\right].\label{eq:llt3}\ee
For the first integral in $(\ref{eq:llt3})$, first note that $\varphi(t)=\E(e^{itX_1})=\cos t$, so 
since the function under the integral is even we have 
$$\frac{1}{2}\int_{-\pi\sqrt{k}/2}^{\pi\sqrt{k}/2}|\varphi^k(u/\sqrt{k})-\exp(-u^2/2)|\, du=\int_{0}^{\pi\sqrt{k}/2}|\varphi^k(u/\sqrt{k})-\exp(-u^2/2)|\, du.$$
Further, it is not hard to show that $e^{-\frac{x^2}{2}}\geq \cos x $ for $x\in [-\pi/2,\pi/2]$. So,
$e^{-\frac{x^2}{2k}}\geq \cos(x/\sqrt{k})$ for $x$ in the bounds of the integral, and therefore
$$\int_{0}^{\pi\sqrt{k}/2}|\varphi^k(u/\sqrt{k})-\exp(-u^2/2)|\, du=
\int_{0}^{\pi\sqrt{k}/2}\exp(-u^2/2)-\varphi^k(u/\sqrt{k})\, du,$$
and now right-side of $(\ref{eq:llt3})$ becomes 
\begin{align*}
 &\frac{\sqrt{n/k}}{\pi}\left[\int_0^{\infty} \exp(-u^2/2)\, du -\int_{0}^{\pi\sqrt{k}/2} \varphi^k(u/\sqrt{k})\, du\right]\\
 =&\frac{\sqrt{n/k}}{\pi}\left[\sqrt{\frac{\pi}{2}}-\int_{0}^{\pi\sqrt{k}/2} \cos^k(u/\sqrt{k})\, du\right]\\
=&\frac{\sqrt{n/k}}{\pi}\underbrace{\left[\sqrt{\frac{\pi}{2}}-\sqrt{k}\int_{0}^{\pi/2} \cos^k(u)\, du\right]}_{\stackrel{(\ref{eq:lmllt1})}{\leq} \frac{C}{k}}.
\end{align*}
From Lemma \ref{lem:strlng} we have that $(\ref{eq:llt3})$ is less than 
$$\frac{C}{\pi}\sqrt{\frac{n}{k^3}}$$
 

}}
\pos{If $\lim_{n\to \infty}\frac{n}{\beta(n)^3}=0$, then 

$$\lim_{n\to \infty}\sup_{\beta(n)\leq k } \sup_{x\in {\cal L}_n^k}\left|\frac{n^{1/2}}{2}p_n^k(x)-\rho_n^k(x)\right|=0.$$
Specially, in the case when $\beta(n)=n$, we have 
$$\sup_{x\in {\cal L}_n^n}\left|\frac{n^{1/2}}{2}p_n^n(x)-\rho_n^n(x)\right|\leq\frac{C}{\pi n}\to 0,$$
as $n\to \infty$.}
\index{Local Limit Theorem|)}

\section{Approximations of the classical heat equation solution}\index{Heat equation|(}

In this section we will review the one-dimensional heat equation
(mostly classical results that can be found in books
that deal with connections to PDEs like Karatzas and Shreve \cite{karatzas},
and some books on classical PDEs like Folland \cite{inpde})
and develop more general results that will later help us.\vspace{0.2cm}

In this section we use the usual space-time ($x$-$t$) coordinate system. We 
are considering the classical initial value problem
\be \left\{\begin{array}{cl} \partial_t w=\frac{1}{2}\partial_{xx} w&\quad {\rm on}\quad  \R\times (0,\infty),\\
    w(0,x)=f(x) & \quad {\rm for}\quad  x\in\R.
    \end{array}\right.\label{pde1}
\ee

If we assume that $f:\R\to \R$ is a Borel measurable function satisfying
\be\int_{-\infty}^{\infty}e^{-ax^2}|f(x)|\, dx<\infty \label{cond1}\ee
for some $a>0$. Then the solution exists.

\teo{\label{pde:t1}If the condition $(\ref{cond1})$ is satisfied, then 
\be w(t,x):=\E(f(x+W_t))=\int_{-\infty}^{\infty}\frac{1}{\sqrt{2\pi t}}f(y)\exp\left(\frac{-(y-x)^2}{2t}\right)\, dy,\label{solpde}\ee
for $0<t<\frac{1}{2a}$ and $x\in \R$ is the solution to the initial value $(\ref{pde1})$. This solution\index{Heat equation!solution} has
derivatives of all orders. Furthermore, if $f$ is continuous at $x$, then
\be \lim_{(t,y)\to (0,x)}w(t,y)=f(x). \ee }
\dok{This follows from the fact that the so called  Gaussian kernel 
$$K_t(x):= \frac{1}{\sqrt{2\pi t}}\exp\left(\frac{-x^2}{2t}\right),$$
satisfies the heat equation. (This can be checked by a direct calculation.) 
The rest follows from the dominated convergence theorem.  }

The main question that will be of interest to us is: if $f$ satisfies $(\ref{cond1})$
and it is continuous,
 for a simple random walk $(S_n)$\index{Simple random walk|(} starting from zero is
\be \E\left[f\left(\frac{S_{\mic{nt}}}{\sqrt{n}}+x\right)\right]\to w(t,x), \label{htapx}\ee
where $u$ is given by $(\ref{solpde})$ and how strong is this convergence. \vspace{0.1cm}

It is clear from the definition of convergence in distribution and the Donsker's
theorem\index{Donsker's Theorem} that this convergence holds if $f$ is bounded. We will show that this holds 
for a much wider set of functions.

\lem{\label{lem:hfdmd}\begin{enumerate}[(a)] 
      \item {\sc (Hoeffding's Inequality)}\index{Hoeffding's Inequality} For $y\geq 0$ we have 
$$\P\left(\left|\frac{S_n}{\sqrt{n}}\right|\geq y\right)\leq 2e^{-y^2/2},$$
where $(S_n)$ is a simple random walk.
      \item  If $\tau>0$ then  for
all $ t\leq \tau $ 
$$\P\left(\left|\frac{S_{\mic{nt}}}{\sqrt{n}}\right|\geq y\right)\leq 2e^{-y^2/(2\tau)}$$
for all $y\geq 0$.
     \end{enumerate}
\dok{(a) 
This is a well known inequality. For the proof see, for example, \cite{hoeffding} or \cite{cnctrinq}.
(b) For $\mic{nt}=0$ the claim is clear. Otherwise, we have
$$\P\left(\left|\frac{S_{\mic{nt}}}{\sqrt{n}}\right|\geq y\right)=\P\left(\left|\frac{S_{\mic{nt}}}{\sqrt{\mic{nt}}}\right|\geq \frac{y\sqrt{n}}{\sqrt{\mic{nt}}}\right).$$
Now, since 
$$\frac{y\sqrt{n}}{\sqrt{\mic{nt}}}\geq \frac{y\sqrt{n}}{\sqrt{nt}}=\frac{y}{\sqrt{t}}\geq \frac{y}{\sqrt{\tau}},$$
we have 
$$\P\left(\left|\frac{S_{\mic{nt}}}{\sqrt{\mic{nt}}}\right|\geq \frac{y\sqrt{n}}{\sqrt{\mic{nt}}}\right)\leq \P\left(\left|\frac{S_{\mic{nt}}}{\sqrt{\mic{nt}}}\right|\geq \frac{y}{\sqrt{\tau}} \right),$$
and the claim follows from part (a).
}}

Define $B^n$ to be the linear interpolation of $t\mapsto\dfrac{S_{\mic{nt}}}{\sqrt{\mic{nt}}}$,
that is 
$$B^n(t):=\frac{S_{\mic{nt}}}{\sqrt{n}}+(nt-\mic{nt})\left(\frac{S_{\mic{nt}+1}}{\sqrt{n}}-\frac{S_{\mic{nt}}}{\sqrt{n}}\right).$$

\lem{\label{cnt:lem}For any $f$ continuous, $a<b$ real numbers, and $\varepsilon>0$ we have
$$\E\left[f\left(\frac{S_{\mic{nt}}}{\sqrt{n}}+x\right)g^{\varepsilon}_{a,b}\left(\frac{S_{\mic{nt}}}{\sqrt{n}}\right)\right]\to \E\left[f(W_t+x)g^{\varepsilon}_{a,b}(W_t)\right]$$
uniformly on compact sets in $(t,x)$, where  
$$g^{\varepsilon}_{a,b}(x)=\left\{\begin{array}{cc}
                         1 & x\in [a,b],\\
			\frac{x-a+\varepsilon}{\varepsilon}& x\in [a-\varepsilon, a],\\
			\frac{x-b-\varepsilon}{-\varepsilon}& x\in [b,b+\varepsilon],\\
			0 & x\notin [a-\varepsilon,b+\varepsilon].
                        \end{array}\right.
$$
\dok{Let $K\subset \R^+\times\R$ be a compact set and define $K_t:=\pi_t(K)$ and $K_x:=\pi_x(K)$.
They are also compact. Hence, the function $h:\R\times K_x\to \R$
given by $h(u,x):=f(u+x)g_{a,b}^{\varepsilon}(u)$ is a continuous function supported on a compact set
(which is a subset of $[a-\varepsilon,b+\varepsilon]\times K_x$).
Now, since $K_t$ is compact, there exists $T>0$ such that $K_t\subset [0,T]$. By Donsker's Theorem
we know that $B^n\stackrel{d}{\to}W$ on $[0,T]$, hence by Skorohod's Representation Theorem\index{Skorohod's Representation Theorem} there exists
a probability space $(\tilde{\Omega},\tilde{\F},\tilde{\P})$ with random elements $\tilde{B}^n\stackrel{d}{=}B^n$ and $\tilde{W}\stackrel{d}{=}W$ 
such that $$\|\tilde{B}^n(\omega)-\tilde{W}(\omega)\|=\sup_{t\in [0,T]}|\tilde{B}^n(t)(\omega)-\tilde{W}(t)(\omega)|\to 0,$$
for all $\omega\in \tilde{\Omega}$. Note that if we define 
$$\tilde{S}^n_t:=\tilde{B}^n_{\mic{nt}/n},$$
$\tilde{S}^n$ has the same distribution as $S_{\mic{nt}}/\sqrt{n}$. Further,
it is clear that 
$$\|\tilde{S}^n-B^n\|\leq \frac{1}{\sqrt{n}}.$$ 
Therefore $\|\tilde{S}^n-\tilde{W}^n\|\to 0$. Now since $h\in C_c(\R^2)$ it is uniformly continuous function,
and therefore 
$$\sup_{(t,x)\in K}|h(\tilde{S}^n_t,x)-h(\tilde{W}(t),x)|\to 0,$$
Now
\begin{align*}
 & \left|\E\left[f\left(\frac{S_{\mic{nt}}}{\sqrt{n}}+x\right)g^{\varepsilon}_{a,b}\left(\frac{S_{\mic{nt}}}{\sqrt{n}}\right)\right]- \E\left[f(W_t+x)g^{\varepsilon}_{a,b}(W_t)\right]\right|\\
&=\left|\tilde{\E}[f(\tilde{S}^n_t+x)g^{\varepsilon}_{a,b}(\tilde{S}^n_t))-f(\tilde{W}_t+x)g^{\varepsilon}_{a,b}(\tilde{W}_t))]\right|\\
&\leq \tilde{\E}\left[\sup_{(t,x)\in K}\left|h(\tilde{S}^n_t,x)-h(\tilde{W}(t),x)\right|\right].
\end{align*}
The convergence follows from the dominated convergence theorem.
}}

In order to get $(\ref{htapx})$ we have to make some \textit{mild} assumptions on $f$. 

\lem{\label{lem:htco}Let $f:\R\to\R$ be a continuous function such that there exist $C>0$ and a locally integrable $g:\R^+\to \R^+$ with the property
\be |f(x)|^2\leq C+ \int_0^{|x|}g(y)\, dy,\label{cond1:ht}\ee
for all $x\in \R$, and
\be\int_0^{\infty}g(y)e^{-y^2/(2\tau)}\, dy<\infty,\label{cond2:ht}\ee
for some $\tau>0$. Then there exists $M>0$ (that depends on $C$, $g$ and $\tau$) such that for all $t\leq \tau$
and all $n$
\be\E[|f(S_{\mic{nt}}/\sqrt{n})|^2]<M.\label{2mbd:ht} \ee
Further, 
\be\int_{-\infty }^{\infty}|f(x)|e^{-x^2/(2\tau )}\, dx<\infty .\label{1mbd:ht}\ee
\dok{We set $G(x):=\int_0^{|x|}g(y)\, dy$. Now, 
\begin{align*}
\E[|f(S_{\mic{nt}}/\sqrt{n})|^2]&\leq C+\E(G(|S_{\mic{nt}}/\sqrt{n}|))\\
&\leq C+ \E( \int_0^\infty g(y)\1_{(y\leq |S_{\mic{nt}}/\sqrt{n}|) }\, dy )\\
&=C+  \int_0^\infty g(y)\E(\1_{(y\leq |S_{\mic{nt}}/\sqrt{n}|) }\, dy\\
&=C+  \int_0^\infty g(y)\P( |S_{\mic{nt}}/\sqrt{n}|\geq y) \, dy\\
&=C+\int_0^{\infty}g(y)\P( |S_{\mic{nt}}/\sqrt{n}|\geq y) \, dy\\
&=C+2\int_0^{\infty}g(y)e^{-y^2/(2\tau)}\, dy=:M.
\end{align*}
For $(\ref{1mbd:ht})$ we first show a similar results using the same arguments. Let $X\sim N(0,\tau)$.
Then $\P(|X|>x)\leq 2e^{-x^2/(2\tau)}$. Now,
\begin{align*}
 \int_{-\infty }^{\infty}|f(x)|^2e^{-x^2/(2\tau )}\, dx &= \sqrt{2\pi\cdot \tau }\E(|f(X)|^2)\\
&\leq \sqrt{2\pi\tau }(C+\E(G(X)))\\
&=C\sqrt{2\pi\tau  }+\sqrt{2\pi\tau }\E(G(X))\\
&=C\sqrt{2\pi\tau  }+\sqrt{2\pi\tau  }\int_0^{\infty}g(y)\P( |X|\geq y) \, dy\\
&=C\sqrt{2\pi\tau  }+2\sqrt{2\pi\tau  }\int_0^{\infty}g(y)e^{-y^2/(2\tau)}\, dy\\
&=M \sqrt{2\pi\tau  } .
\end{align*}
Now, it is clear from Cauchy-Schwarz inequality we have that
\begin{align*}
\int_{-\infty }^{\infty}|f(x)|e^{-x^2/(2\tau )}\, dx &\leq \left(\int_{-\infty }^{\infty}e^{-x^2/(2\tau )}\, dx\right)^{1/2} \left(\int_{-\infty }^{\infty}|f(x)|^2e^{-x^2/(2\tau )}\, dx\right)^{1/2}\\
&= \sqrt{2\pi\tau M}.
\end{align*}

}}

\noindent \emph{Remark.} The conditions $(\ref{cond1:ht})$ and $(\ref{cond2:ht})$ given by the previous lemma are satisfied by a wide 
family of functions. For instance, if for $\alpha\geq 1$ we have
$$\limsup_{|y|\to\infty}\frac{|f(y)|}{|y|^\alpha}=:L<\infty.$$
Then there exists a $C>0$ such that
$$|f(y)|\leq C+L|y|^\alpha ,$$
for all $y\in \R$. Now, the function $g(y):=\alpha y^{\alpha-1}$ satisfies 
$(\ref{cond1:ht})$. Further, since the normal distribution has all
the $\alpha$-moments for $\alpha\geq 1$, $g$ satisfies $(\ref{cond2:ht})$.

\teo{\label{teo:aprxclhe}Let $f:\R\to\R$ be a continuous function and  $a<b$ finite real numbers, such that there exist $C>0$ and a locally integrable $g:\R^+\to \R^+$ with the property
$$\sup_{r\in [a,b]}(f(x+r))^2\leq C+ \int_0^{|x|}g(y)\, dy,$$ 
for all $x\in \R$, and
$$\int_0^{\infty}g(y)e^{-y^2/(2\tau)}\, dy<\infty,$$
for some $\tau>0$. Then for all $r\in [a,b]$
\be\E\left[f\left(\frac{S_{\mic{nt}}}{\sqrt{n}}+r\right)\right]\to w(t,r),\label{conv:uni}\ee
as $n\to \infty$ where $t<\tau$ and $w$ is the solution to the initial value problem $(\ref{pde1})$ given
by 
$$w(t,x)=\E(f(x+W_t)).$$
Further, the convergence in $(\ref{conv:uni})$ is uniform on $[0,\tau)\times [a,b]$
\dok{From Lemma \ref{lem:htco}. and Theorem \ref{pde:t1}. we know that $w(\cdot,\cdot)$
is the solution to $(\ref{pde1})$. From Lemma \ref{lem:htco}. applied on the function 
$f(\cdot+r)$
we know that for each $K>0$ 
\begin{align*}
 \left|\E\left[f\left(\frac{S_{\mic{nt}}}{\sqrt{n}}+r\right)\left(1-g^\varepsilon_{-K,K}(\frac{S_{\mic{nt}}}{\sqrt{n}})\right)\right]\right|&\leq \sqrt{\E\left[f\left(\frac{S_{\mic{nt}}}{\sqrt{n}}+r\right)^2\right]\E\left[\left(1-g^\varepsilon_{-K,K}\left(\frac{S_{\mic{nt}}}{\sqrt{n}}\right)\right)^2\right]}\\
&\leq  \sqrt{\E\left[f\left(\frac{S_{\mic{nt}}}{\sqrt{n}}+r\right)^2\right]\E\left[\1_{(|S_{\mic{nt}}/\sqrt{n}|>K)}\right]}\\
&\leq \sqrt{ M \P\left(|S_{\mic{nt}}/\sqrt{n}|>K\right)}\\
&\leq \sqrt{2M\exp\left(-\frac{K^2}{2\tau }\right)}.
\end{align*}
Pick $\varepsilon>0$, then there exists $K>0$ such that 
$$\left|\E\left[f\left(\frac{S_{\mic{nt}}}{\sqrt{n}}+r\right)\left(1-g^\varepsilon_{-K,K}\left(\frac{S_{\mic{nt}}}{\sqrt{n}}\right)\right)\right]\right|<\varepsilon/3$$
for all $n$ and in the same way
$$|\E[f(x+W_t)(1-g^\varepsilon_{-K,K}(W_t)]|<\varepsilon/3$$. Finally, using Lemma \ref{cnt:lem}, $a=-K$, $b=-K$ we
have that for sufficiently large  $n$
$$\left|\E\left[f\left(\frac{S_{\mic{nt}}}{\sqrt{n}}+r\right)\right]- w(t,r)\right|<\varepsilon.$$}}

\subsection{The case when $\alpha>0$}
In this subsection we will show that for $\alpha>0$
\be \E\left[f\left(\frac{S_{\mic{nt}}}{n^{1/2+\alpha}}+x\right)\right]\to f(x), \label{htapx:1}\ee
and uniformly for $(t,x)$ over a compact set.

\lem{If $\tau>0$ then  for
all $ t\leq \tau $ 
$$\P\left(\left|\frac{S_{\mic{nt}}}{n^{1/2+\alpha}}\right|\geq y\right)\leq 2e^{-y^2/(2\tau)}$$
for all $y\geq 0$.
\dok{It is not hard to show that
$$\left\{\left|\frac{S_{\mic{nt}}}{n^{1/2+\alpha}}\right|\geq y\right\}\subset \left\{\left|\frac{S_{\mic{nt}}}{n^{1/2}}\right|\geq y\right\},$$
and the claim now follows from the result of Lemma \ref{lem:hfdmd}. part (b).}}

Using exactly the same argumentation we get a version of Lemma \ref{lem:htco}:
\lem{\label{lem:htco:2}Let $f:\R\to\R$ be a continuous function such that there exist $C>0$ and a locally integrable $g:\R^+\to \R^+$ with the property
\be |f(x)|^2\leq C+ \int_0^{|x|}g(y)\, dy,\label{cond1:ht:2}\ee
for all $x\in \R$, and
\be\int_0^{\infty}g(y)e^{-y^2/(2\tau)}\, dy<\infty,\label{cond2:ht:2}\ee
for some $\tau>0$. Then there exists $M>0$ (that depends on $C$, $g$ and $\tau$) such that for all $t\leq \tau$
and all $n$
\be\E[|f(S_{\mic{nt}}/n^{1/2+\alpha})|^2]<M.\label{2mbd:ht:2} \ee}

Now, under similar conditions as in Theorem \ref{teo:aprxclhe}, we have:
\teo{\label{teo:aprxclhe:2}Let $f:\R\to\R$ be a continuous function and  $a<b$ finite real numbers, such that there exist $C>0$ and a locally integrable $g:\R^+\to \R^+$ with the property
$$\sup_{r\in [a,b]}(f(x+r))^2\leq C+ \int_0^{|x|}g(y)\, dy,$$ 
for all $x\in \R$, and
$$\int_0^{\infty}g(y)e^{-y^2/(2\tau)}\, dy<\infty,$$
for some $\tau>0$. Then for all $r\in [a,b]$
\be\E\left[f\left(\frac{S_{\mic{nt}}}{n^{1/2+\alpha}}+r\right)\right]\to f(r),\label{conv:uni:2}\ee
as $n\to \infty$ where $t<\tau$ 
Further, the convergence in $(\ref{conv:uni:2})$ is uniform on $[0,\tau)\times [a,b]$.
\dok{Let $\varepsilon>0$. The function $f$ on $[a-1,b+1]$ is uniformly continuous, and hence there 
exists $\delta\in (0,1)$ such that for all $y,y'\in [a-1,b+1]$ if $|y-y'|<\delta$ then $|f(y)-f(y')|<\varepsilon$.
Now,
$$\left|\E\left[f\left(\frac{S_{\mic{nt}}}{n^{1/2+\alpha}}+r\right)\right]- f(r)\right|\leq \E\left|f\left(\frac{S_{\mic{nt}}}{n^{1/2+\alpha}}+r\right)- f(r)\right|$$
$$\leq \E\left[\underbrace{\left|f\left(\frac{S_{\mic{nt}}}{n^{1/2+\alpha}}+r\right)- f(r)\right|}_{\leq \varepsilon}\1(|S_{\mic{nt}}/n^{1/2+\alpha}|<\delta)\right]$$ $$+\E\left[\left|f\left(\frac{S_{\mic{nt}}}{n^{1/2+\alpha}}+r\right)- f(r)\right|\1(|S_{\mic{nt}}/n^{1/2+\alpha}|\geq \delta)\right], $$
by uniform continuity of $f$, and triangle inequality we get
$$\leq \varepsilon +\E\left[\left|f\left(\frac{S_{\mic{nt}}}{n^{1/2+\alpha}}+r\right)\right|\1\left(\left|\frac{S_{\mic{nt}}}{n^{1/2+\alpha}}\right|\geq \delta\right)\right]+|f(r)|\P\left(\left|\frac{S_{\mic{nt}}}{n^{1/2+\alpha}}\right|\geq \delta\right).$$
By Cauchy-Schwarz we get 
$$\leq \varepsilon +\E\left[\left|f\left(\frac{S_{\mic{nt}}}{n^{1/2+\alpha}}+r\right)\right|^2\right]^{1/2}\P\left(\left|\frac{S_{\mic{nt}}}{n^{1/2}}\right|\geq n^{\alpha}\delta\right)^{1/2}+\left(\max_{y\in [a,b]}|f(y)|\right)\P\left(\left|\frac{S_{\mic{nt}}}{n^{1/2}}\right|\geq n^{\alpha}\delta\right).$$
Using Lemma \ref{lem:htco:2} for the function $f(\cdot+r)$ on the expectation, and Lemma \ref{lem:hfdmd}. part (b) on the probabilities, we get 
$$\leq \varepsilon +M^{1/2}\sqrt{2}\exp\left(\frac{-n^{2\alpha}\delta^2}{4}\right)+\left(\max_{y\in [a,b]}|f(y)|\right)\cdot 2\exp\left(\frac{-n^{2\alpha}\delta^2}{2}\right).$$
Note that the bound doesn't depend on $t$ or $r$, and we have 
$$\limsup_{n\to\infty}\sup_{t,r}\left|\E\left[f\left(\frac{S_{\mic{nt}}}{n^{1/2+\alpha}}+r\right)\right]- f(r)\right|\leq \varepsilon.$$
Since, $\varepsilon>0$ is arbitrary the claim follows.}}

\subsection{Summary}
\teo{\label{teo:aprxclhe:3}Let $f:\R\to\R$ be a continuous function and  $a<b$ finite real numbers, such that there exist $C>0$ and a locally integrable $g:\R^+\to \R^+$ with the property
$$\sup_{r\in [a,b]}(f(x+r))^2\leq C+ \int_0^{|x|}g(y)\, dy,$$ 
for all $x\in \R$, and
$$\int_0^{\infty}g(y)e^{-y^2/(2\tau)}\, dy<\infty,$$
for some $\tau>0$. Then for all $r\in [a,b]$
\be\E\left[f\left(\frac{S_{\mic{nt}}}{n^{1/2+\alpha}}+r\right)\right]\to w_{\alpha}(t,r),\label{conv:uni:3}\ee
as $n\to \infty$ where $t<\tau$ and $w_{\alpha}$ is the solution to the initial value problem given
by 
$$\left\{\begin{array}{l}
\partial_t w_{\alpha}=\left\{\begin{array}{cl}
                                     0, & \alpha>0\\
				  \frac{1}{2}\partial_{xx}w_{\alpha} & \alpha =0
                                    \end{array}\right.\\
w_{\alpha}(0,x)=f(x)
   
  \end{array}\right. .$$
Further, the convergence in $(\ref{conv:uni:3})$ is uniform on $[0,\tau)\times [a,b]$.}

\index{Simple random walk|)}
\section{Euler method for the stochastic heat equation}\label{sec:eumth}

Let $u$ be the solution to the heat equation 
\be \partial_tu=\beta\partial_{xx}u+f.\label{es:j1}\ee
Now, we discretize this equation at the point $(t,x)$

$$u_t(t,x)\approx \frac{u(t+\Delta t,x)-u(t,x)}{\Delta t},$$
$$u_{xx}(t,x)\approx \frac{u(t,x+\Delta x)-2u(t,x)+u(t,x-\Delta x)}{(\Delta x)^2},$$
where $\Delta t$ and $\Delta x$ are small and positive. So equation $(\ref{es:j1})$ becomes
$$u(t+\Delta t,x)\approx\frac{\beta\Delta t}{(\Delta x)^2}(u(t,x+\Delta x)+u(t,x-\Delta x))+\left(1-2\frac{\beta\Delta t}{(\Delta x)^2}\right)u(t,x)+\Delta t f(t,x).$$

\begin{figure}[ht]
\begin{center}
\psfrag{a}{$\boldsymbol{t_0}$}
\psfrag{b}{$\boldsymbol{t_1}$}
\psfrag{c}{$\boldsymbol{t_2}$}
\psfrag{d}{$\boldsymbol{t_3}$}
\psfrag{e}{$\boldsymbol{t_4}$}
\psfrag{v}{$\boldsymbol{\vdots}$}
\psfrag{l}{$\boldsymbol{\ldots}$}

\psfrag{0}{$\boldsymbol{x_0}$}
\psfrag{1}{$\boldsymbol{x_1}$}
\psfrag{2}{$\boldsymbol{x_2}$}
\psfrag{3}{$\boldsymbol{x_3}$}
\psfrag{4}{$\boldsymbol{x_4}$}

\includegraphics[width=7cm]{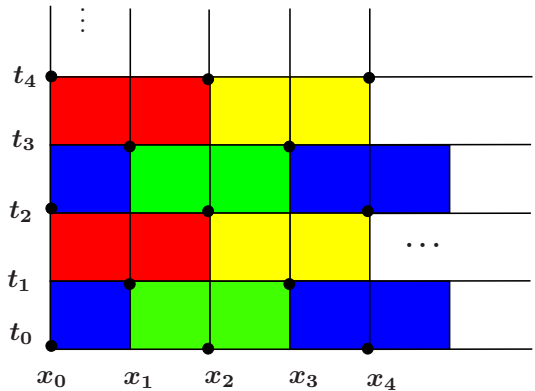}\\
\caption{Euler method} \label{pic:sl29}
  \end{center}
\end{figure}

Now, if we set $t_k=k\Delta t$, $x_k=k\Delta x$, for $k=0, 1, 2, \ldots$, and if 
we replace $u(t_j,x_k)$ by $U^j_k$ and $f(t_j,x_k)$ by $f^j_k$, we get
$$U^{j+1}_{k}=s(U^j_{k+1}+U^j_{k-1})+(1-2s)U^j_{k}+(\Delta t) f^j_k,$$
where $s=\frac{\beta \Delta t}{(\Delta x)^2}$.\vspace{0.3cm}

Further, if the equation $(\ref{es:j1})$ has 
\begin{enumerate}[(a)]
 \item initial condition $u(x,0)=g(x)$;
 \item 0-boundary condition $u(0,t)=0$;
\end{enumerate}
then we set
\begin{enumerate}[(a)]
 \item initial condition $U^0_k=g(x_k)$;
 \item 0-boundary condition $U^j_0=0$.
\end{enumerate}

The given scheme is called the {\bf explicit Euler method}\index{Heat equation!Euler method|(} for the one dimensional
heat equation. It is stable if $s\leq 1/2$. (See \cite{lattice}.)\vspace{0.2cm}

We are interested for the Euler method in the case of the stochastic heat equation
$$\partial_t v=\frac{1}{2}\partial_{xx}v+\W, $$\index{Heat equation!stochastic|(}
with initial and boundary value conditions $v(0,x)=0$ and $v(t,0)=0$.
We will look at the method when $\Delta x=n^{-1/2}$ and $\Delta t=1/n$. In this case
$s=\frac{1}{2}$
and Euler method looks like this
\be V^{j+1}_{k}=\frac{1}{2}(V^j_{k+1}+V^j_{k-1})+\frac{\sqrt{n}}{2} \W(R_{jk}),\quad k\geq 1, j\geq 0\label{esh:eq1}\ee
where $R_{jk}=[x_{k-1},x_{k+1}]\times [t_{j},t_{j+1}]$. We will work with the case
when $V^j_0=0$ and $V^0_k=0$. (Initial and boundary value conditions are $0$.)

\begin{figure}[ht]
\begin{center}
\psfrag{a}{$\boldsymbol{x_{k-1}}$}
\psfrag{b}{$\boldsymbol{x_k}$}
\psfrag{c}{$\boldsymbol{x_{k+1}}$}
\psfrag{d}{$\boldsymbol{t_{j+1}}$}
\psfrag{e}{$\boldsymbol{t_{j}}$}

\includegraphics[width=4cm]{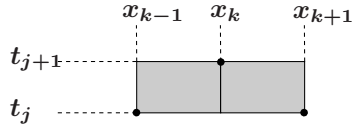}\\
\caption{Rectangle $R_{jk}$} \label{pic:sl30}
  \end{center}
\end{figure}

It is not hard to get the following result.\index{Simple random walk|(}
\lem{The solution to the difference equation $(\ref{esh:eq1})$ with initial and boundary condition $0$
is given by
\begin{align}
& V^j_{k}\nonumber\\
=&\frac{\sqrt{n}}{2}\sum_{j'=1}^j\sum_{k'=1}^{k+j} \left(\P\left(\frac{1}{\sqrt{n}}S_{n(t_j-t_{j'})}+x_{k}=x_{k'}\right)-\P\left(\frac{1}{\sqrt{n}}S_{n(t_j-t_{j'})}+x_{k}=-x_{k'}\right)\right)\W(R_{j'-1,k'}) \label{rel:rek1}\\
=&\frac{\sqrt{n}}{2}\sum_{j'=1}^{j}\sum_{k'=1}^{k+j} \left(\P(S_{j-j'}=k'-k)-\P(S_{j-j'}=-k-k')\right)\W(R_{j'-1,k'}). \label{rel:rek2}
\end{align}
%
}

\noindent\emph{Remark.} Note that $j-j'$ and $k-k'$ have to be of both either even or odd. Otherwise,
$\P(S_{j-j'}=k'-k)=\P(S_{j-j'}=-k-k')=0$. So in the upper formula the only use the rectangles 
$R_{j'-1,k'}$ where $j-j' \equiv k-k'$ (mod 2).
\vspace{0.3cm}

Having in mind this remark we will focus on lattice points 
that are in the same class as $(0,0)$:
$${\cal L}^{n}_0=\{(t_j,x_k): k,j\geq 0,\ j \equiv k \ (\textrm{mod} \ 2)\}$$
(Points in ${\cal L}^{n}_0$ are marked with $\bullet$ on Figure \ref{pic:sl29}.) Note that if $(t_j,x_k)\in {\cal L}_0^{n}$ then if $(t_{j'},x_{k'})\notin {\cal L}_0^{n}$ then
$$\P\left(\frac{1}{\sqrt{n}}S_{n(t_j-t_{j'})}+x_{k}=x_{k'}\right)-\P\left(\frac{1}{\sqrt{n}}S_{n(t_j-t_{j'})}+x_{k}=-x_{k'}\right)=0.$$

Also, note that if $|k'-k|>j-j'$ then the previous equality also holds. \vspace{0.2cm}

Our aim is to show that for a compact set $K$ when $n\to \infty$
to discover the rate of convergence to 0 of
$$\sup\{|V^j_k-v(t_j,x_k)|:(t_j,x_k)\in K\cap {\cal L}_0^{n}\},$$ 
where $v$ is the solution to the stochastic heat equation 
with $0$-boundary and $0$-initial condition. Recall (see \S\ref{sec:shee}), that the solution to the homogeneous 
stochastic heat equation with $0$ initial and boundary conditions is
$$v(t,x)=\int_0^t\frac{1}{\sqrt{2\pi (t-s)}}\int_{\R^+}\left(e^{-\frac{|x-y|^2}{2(t-s)}}-e^{-\frac{|x+y|^2}{2(t-s)}}\right)\W(ds,dy).$$

Note that for $(t_j,x_k)\in {\cal L}^n_0$, $j,k>0$ we have 
$$V_k^j=\sum_{j'=1}^j\sum_{(t_{j'},x_{k'})\in {\cal L}^{n}_0}\frac{\sqrt{n}}{2} \left(\P\left(\frac{1}{\sqrt{n}}S_{n(t_j-t_{j'})}+x_{k}=x_{k'}\right)-\P\left(\frac{1}{\sqrt{n}}S_{n(t_j-t_{j'})}+x_{k}=-x_{k'}\right)\right)\W(R_{j'-1,k'})$$
where $R_{j,0}=[t_j,t_{j+1}]\times [x_0,x_1]$. (Notice, that in the case when $k'=0$, the expression in the sum
is equal to 0.)

Note that (see Figure \ref{pic:sl29})
$$\bigcup_{j'=1}^j\bigcup_{(t_{j'},x_{k'})\in {\cal L}_0^n} R_{j'-1,k'}=[0,t_j]\times [0,\infty],$$
therefore 
\begin{align*}
 & V_j^k-v(t_j,x_k)\\
=&\sum_{j'=1}^j\sum_{(t_{j'},x_{k'})\in {\cal L}^{n}_0} \\
& \frac{\sqrt{n}}{2} \left(\P(\frac{1}{\sqrt{n}}S_{n(t_j-t_{j'})}+x_{k}=x_{k'})-\P(\frac{1}{\sqrt{n}}S_{n(t_j-t_{j'})}+x_{k}=-x_{k'})\right)\W(R_{j'-1,k'})-\\
 &\int_{R_{j'-1,k'}}\frac{1}{\sqrt{2\pi (t_j-s)}}\left(e^{-\frac{|x_k-y|^2}{2(t_j-s)}}-e^{-\frac{|x_k+y|^2}{2(t_j-s)}}\right)\W(ds,dy).
\end{align*}

Since $\E(V_j^k-v(t_j,x_k))=0$, for the variance we have
\begin{align}
 & \E[(V_j^k-v(t_j,x_k))^2] \nonumber\\
=&\sum_{j'=1}^j\sum_{(t_{j'},x_{k'})\in {\cal L}^{n}_0}\nonumber\\
& \int_{R_{j'-1,k'}}\left[\frac{\P(\frac{1}{\sqrt{n}}S_{n(t_j-t_{j'})}+x_{k}=x_{k'})-\P(\frac{1}{\sqrt{n}}S_{n(t_j-t_{j'})}+x_{k}=-x_{k'})}{2/\sqrt{n}}\right.\nonumber\\
&\left. -\frac{1}{\sqrt{2\pi (t_j-s)}}\left(e^{-\frac{|x_k-y|^2}{2(t_j-s)}}-e^{-\frac{|x_k+y|^2}{2(t_j-s)}}\right)\right]^2dy\, ds\label{eq:esh3}
\end{align} 
\index{Heat equation!Euler method|)}

\subsection{Convergence of the Euler Method}\index{Heat equation!Euler method|(}

The main result of this subsection is the following,
\teo{\label{thm:cnvemt}There exists $\gamma>0$ and $\Gamma>0$ such that
$$\E[(v(t_j,x_k)-V_k^j)^2]\leq \frac{\Gamma}{n^{\gamma}}.$$
}
This rate will help us to prove the convergence of the Euler method.\vspace{0.2cm}

Recall from $(\ref{eq:esh3})$ that
\begin{align*}
 & \E[(V_j^k-v(t_j,x_k))^2] \\
=&\sum_{j'=1}^j\sum_{(t_{j'},x_{k'})\in {\cal L}^{n}_0}\\
& \int_{R_{j'-1,k'}}\left[\frac{\sqrt{n}}{2}\left(\P\left(\frac{S_{j-j'}}{\sqrt{n}}=x_{k'}-x_{k}\right)-\P\left(\frac{S_{j-j'}}{\sqrt{n}}=-x_{k'}-x_{k}\right)\right)\right.\\
&\left. -\frac{1}{\sqrt{2\pi (t_j-s)}}\left(e^{-\frac{|x_k-y|^2}{2(t_j-s)}}-e^{-\frac{|x_k+y|^2}{2(t_j-s)}}\right)\right]^2dy\, ds
\end{align*} 
We will decompose our sum over $j$'s into two sums:
\be\sum_{j'=1}^j\ \ =\sum_{j-j'<n^{1/3+\alpha}}+\sum_{j-j'\geq n^{1/3+\alpha}}
\label{sumdec}
\ee
where $\alpha>0$ is some small positive number to be determined later.

\begin{align}
&\sum_{j-j'<n^{1/3+\alpha}}\sum_{(t_{j'},x_{k'})\in {\cal L}^{n}_0} \nonumber\\
& \int_{R_{j'-1,k'}}\left[\frac{\sqrt{n}}{2}\left(\P\left(\frac{S_{j-j'}}{\sqrt{n}}=x_{k'}-x_{k}\right)-\P\left(\frac{S_{j-j'}}{\sqrt{n}}=-x_{k'}-x_{k}\right)\right)\right.\nonumber\\
&\left. -\frac{1}{\sqrt{2\pi (t_j-s)}}\left(e^{-\frac{|x_k-y|^2}{2(t_j-s)}}-e^{-\frac{|x_k+y|^2}{2(t_j-s)}}\right)\right]^2dy\, ds\label{esh:sum1}.
\end{align} 
From the inequality $(a+b+c+d)^2\leq 4a^2+4b^2+4c^2+4d^2$, so we have that the sum $(\ref{esh:sum1})$ is less
than: 
\begin{align}
&\sum_{j-j'<n^{1/3+\alpha}}\sum_{(t_{j'},x_{k'})\in {\cal L}^{n}_0} \nonumber\\
& \left[\frac{2}{\sqrt{n}}\left(\P\left(\frac{S_{j-j'}}{\sqrt{n}}=x_{k'}-x_{k}\right)^2+\P\left(\frac{S_{j-j'}}{\sqrt{n}}=-x_{k'}-x_{k}\right)^2\right)\right.\nonumber\\
&\left. +\int_{R_{j'-1,k'}}\frac{2}{\pi (t_j-s)}\left(e^{-\frac{|x_k-y|^2}{(t_j-s)}}+e^{-\frac{|x_k+y|^2}{(t_j-s)}}\right)\right]dy\, ds\label{esh:sum2}.
\end{align}
It is not hard to see that using sub-additivity we have
$$\sum_{(t_{j'},x_{k'})\in {\cal L}^{n}_0}\P\left(\frac{S_{j-j'}}{\sqrt{n}}=\pm x_{k'}- x_{k}\right)^2\leq \sum_{(t_{j'},x_{k'})\in {\cal L}^{n}_0}\P\left(\frac{S_{j-j'}}{\sqrt{n}}=\pm x_{k'}- x_{k}\right)\leq 1.$$
Therefore 

\begin{align}
&\sum_{0\leq j-j'<n^{1/3+\alpha}}\sum_{(t_{j'},x_{k'})\in {\cal L}^{n}_0} \left[\frac{2}{\sqrt{n}}\left(\P\left(\frac{S_{j-j'}}{\sqrt{n}}=x_{k'}-x_{k}\right)^2+\P\left(\frac{S_{j-j'}}{\sqrt{n}}=-x_{k'}-x_{k}\right)^2\right)\right]\nonumber\\
\leq & \frac{4}{\sqrt{n}}n^{1/3+\alpha}=4n^{\alpha-1/6}.\label{esh:sum2a}
\end{align}

Further 
\begin{align}
 &\sum_{j-j'<n^{1/3+\alpha}}\sum_{(t_{j'},x_{k'})\in {\cal L}^{n}_0} \int_{R_{j'-1,k'}}\frac{2}{\pi (t_j-s)}e^{-\frac{|x_k\pm y|^2}{(t_j-s)}}dy\, ds\nonumber\\
 =& \int_{t_{j*}}^{t_j} \int_0^\infty \frac{2}{\pi (t_j-s)}e^{-\frac{|x_k\pm y|^2}{(t_j-s)}}dy\, ds
\leq   \int_{t_{j*}}^{t_j} \int_{-\infty}^\infty \frac{2}{\pi (t_j-s)}e^{-\frac{|x_k\pm y|^2}{(t_j-s)}}dy\, ds\nonumber\\
=& \int_{t_{j*}}^{t_j} \frac{2}{\sqrt{\pi (t_j-s)}}\, ds=\frac{4}{\sqrt{\pi}}\sqrt{t_j-t_{j*}}=\frac{8}{\sqrt{\pi}}\sqrt{\frac{j-j^*}{n}},\label{esh:sum3}
\end{align}
where $j^*+1=\min\{j':j'\geq 1,\ 0\leq j-j'\leq n^{1/3+\alpha} \}$, so the sum 
$(\ref{esh:sum3})$ is less than 
\be \frac{8}{\sqrt{\pi}}\sqrt{\frac{n^{1/3+\alpha}+1}{n}}=\frac{8}{\sqrt{\pi}}\sqrt{n^{\alpha - 2/3}+n^{-1}}.\label{esh:bnd1}\ee
Hence, from $(\ref{esh:sum2a})$  and $(\ref{esh:bnd1})$ the sum $(\ref{esh:sum1})$ is bounded by
\be 4n^{\alpha-1/6}+\frac{8}{\sqrt{\pi}}\sqrt{n^{\alpha - 2/3}+n^{-1}}.\label{esh:bnd1a}\ee

In order to estimate $\sum_{j-j'\geq n^{\alpha+1/3}}$ - part of the sum $(\ref{sumdec})$ we first need to do some
estimates on the gradient of the function 
$$F(t,x)=\frac{1}{\sqrt{2\pi (t_j-t)}}\exp\left(-\frac{(x_j-x)^2}{2(t_j-t)}\right)$$
for $(t,x)\in [0,t_j)\times \R$.

\lem{\begin{enumerate}[(a)]
      \item For fixed $t<t_j$ we have
\be\sup_{x\in \R}\|\nabla_{t,x}F(t,x)\|^2\leq \max\left\{\frac{1}{8\pi(t_j-t)^3},\frac{1}{2e(t_j-t)^2},\frac{(1-(t_j-t))e^{-3+4(t_j-t)}}{2\pi(t_j-t)^3}\right\}.\label{eq:grdgs}\ee
\item For $A>0$ there exists a constant $C_A$ (depending only on $A$) such that for
$0\leq t< t_j\leq A$ we have 
$$\sup_{x\in \R}\|\nabla_{t,x}F(t,x)\| \leq \frac{C_A}{(t_j-t)^{3/2}}.$$
     \end{enumerate}
\dok{(a) By doing taking derivatives we have 
$$D^2(t,x)=\|\nabla_{t,x}F(t,x)\|^2 = (\partial_t F(t,x))^2+(\partial_x F(t,x))^2=$$
$$=\exp\left(-\frac{(x-x_k)^2}{t_j-t}\right)\left[\frac{1}{8\pi}\left(\frac{1}{(t_j-t)^{3/2}}-\frac{(x-x_k)^2}{(t_j-t)^{5/2}}\right)^2+\frac{1}{2\pi}\cdot\frac{(x-x_k)^2}{(t_j-t)^3} \right]$$
It is clear that $\lim_{x\to \pm \infty}D^2(t,x)=0$, so there exists a maximum,
and it is obtained at the zeros of $\partial_x (D^2(t,x))=$
$$=-\frac{(x-x_k)((x-x_k)^2-(t_j-t))((x-x_k)^2-(3(t_j-t)-4(t_j-t)^2)}{4\pi(t_j-t)^6}\exp\left(-\frac{(x-x_k)^2}{t_j-t}\right).$$
If we set $x-x_k=0$ we get  $D^2(t,x)=(8\pi(t_j-t)^3)^{-1}$; for $(x-x_k)^2=(t_j-t)$
we have $D^2(t,x)=(2e(t_j-t)^2)^{-1}$; for $(x-x_k)^2=3(t_j-t)-4(t_j-t)^2$ (note that this
may not be solvable) we have $D^2(t,x)=\frac{(1-(t_j-t))e^{-3+4(t_j-t)}}{2\pi(t_j-t)^3}$. If we can solve the equation in the last case then
we have an equality in $(\ref{eq:grdgs})$, otherwise we have an inequality.\vspace{0.2cm}

(b) Since $(t,t_j)\mapsto (1-(t_j-t))e^{-3+4(t_j-t)}$ obtains a maximum $M_A$
on the compact set $[0,A]^2$, we have
$$(t_j-t)^3\sup_{x\in \R}\|\nabla_{t,x}F(t,x)\|^2\leq \max\{\frac{1}{8\pi}, \frac{\overbrace{t_j-t}^{\leq 2A}}{2e}, M_A\}.$$
}}

\pos{For $0\leq t< t_j\leq A$ and $0<L\leq U$ we have
$$\sup_{L\leq t_j-t\leq U}\ \sup_{x\in \R}\|\nabla_{t,x}F(t,x)\|\leq \frac{C_A}{L^{3/2}}.$$}

We now have everything we need to estimate $\sum_{j-j'\geq n^{1/3+\alpha}}$ - part
of the sum:
\begin{align}
&\sum_{j-j'\geq n^{1/3+\alpha}}\sum_{(t_{j'},x_{k'})\in {\cal L}^{n}_0} \nonumber\\
& \int_{R_{j'-1,k'}}\left[\frac{\sqrt{n}}{2}\left(\P\left(\frac{S_{j-j'}}{\sqrt{n}}=x_{k'}-x_{k}\right)-\P\left(\frac{S_{j-j'}}{\sqrt{n}}=-x_{k'}-x_{k}\right)\right)\right.\nonumber\\
&\left. -\frac{1}{\sqrt{2\pi (t_j-s)}}\left(e^{-\frac{|x_k-y|^2}{2(t_j-s)}}-e^{-\frac{|x_k+y|^2}{2(t_j-s)}}\right)\right]^2dy\, ds\label{esh:sum4}
\end{align}   
We first give an upper bound for
\begin{align}
 &\left|\frac{\sqrt{n}}{2}\left(\P\left(\frac{S_{j-j'}}{\sqrt{n}}=x_{k'}-x_{k}\right)-\P\left(\frac{S_{j-j'}}{\sqrt{n}}=-x_{k'}-x_{k}\right)\right)\right.\nonumber \\
&\left. -\frac{1}{\sqrt{2\pi (t_j-s)}}\left(e^{-\frac{|x_k-y|^2}{2(t_j-s)}}-e^{-\frac{|x_k+y|^2}{2(t_j-s)}}\right)\right| \label{esh:ex1}
\end{align}
where $(s,y)\in R_{j'-1,k'}$. By triangle inequality, expression $(\ref{esh:ex1})$ is less or equal to
\begin{align}
 &\left|\frac{\sqrt{n}}{2}\P\left(\frac{S_{j-j'}}{\sqrt{n}}=x_{k'}-x_{k}\right)-\frac{1}{\sqrt{2\pi (t_j-t_{j'})}}e^{-\frac{|x_k-x_{k'}|^2}{2 (t_j-t_{j'})}}\right|\nonumber \\
+&\left|-\frac{\sqrt{n}}{2}\P\left(\frac{S_{j-j'}}{\sqrt{n}}=-x_{k'}-x_{k}\right)+\frac{1}{\sqrt{2\pi (t_j-t_{j'})}}e^{-\frac{|x_k+x_{k'}|^2}{2 (t_j-t_{j'})}}\right|\nonumber \\
+&\left|\frac{1}{\sqrt{2\pi (t_j-t_{j'})}}e^{-\frac{|x_k-x_{k'}|^2}{2 (t_j-t_{j'})}} -\frac{1}{\sqrt{2\pi (t_j-s)}}e^{-\frac{|x_k-y|^2}{2(t_j-s)}}\right|\nonumber\\
+&\left|-\frac{1}{\sqrt{2\pi (t_j-t_{j'})}}e^{-\frac{|x_k+x_{k'}|^2}{2 (t_j-t_{j'})}}+\frac{1}{\sqrt{2\pi (t_j-s)}}e^{-\frac{|x_k+y|^2}{2(t_j-s)}}\right|. \label{esh:ex2}
\end{align}
The first two terms in $(\ref{esh:ex2})$, by Theorem \ref{thm:llt}, are less than $\frac{C}{n^{3\alpha}}$. By mean-value theorem, the last two terms in $(\ref{esh:ex2})$ are
less than 
$$\sup_{(s,y)\in R_{j'-1,k'}}\|\nabla_{t,x}F(s,y)\| \sqrt{(t_j-s)^2+(x_k-y)^2},$$
and by the definition of $R_{j'-1,k'}=[t_{j'-1},t_{j'}]\times [x_{(k'-1)\wedge 0},x_{k'+1}]$, this is less than
$$\left(\sup_{t_j-{t_{j'}}\leq t_j-s\leq t_j-t_{j'-1}}\|\nabla_{t,x}F(s,y)\|\right) \sqrt{(t_j-s)^2+(x_k-y)^2}\leq \frac{C_A}{(t_j-t_{j'})^{3/2}}\sqrt{\frac{1}{n^2}+\frac{1}{n}}$$
$$=\frac{C_A}{(j-j')^{3/2}}\sqrt{n+n^2}.$$
Now, for $j-j'\geq n^{1/3+\alpha}$ we have 
\begin{align}
&\sum_{(t_{j'},x_{k'})\in {\cal L}^{n}_0}  \int_{R_{j'-1,k'}}\left[\frac{\sqrt{n}}{2}\left(\P\left(\frac{S_{j-j'}}{\sqrt{n}}=x_{k'}-x_{k}\right)-\P\left(\frac{S_{j-j'}}{\sqrt{n}}=-x_{k'}-x_{k}\right)\right)\right.\nonumber\\
&\left. -\frac{1}{\sqrt{2\pi (t_j-s)}}\left(e^{-\frac{|x_k-y|^2}{2(t_j-s)}}-e^{-\frac{|x_k+y|^2}{2(t_j-s)}}\right)\right]^2dy\, ds \nonumber\\
\leq &\left(\frac{C_A}{(j-j')^{3/2}}\sqrt{n+n^2}+\frac{C}{n^{3\alpha}}\right) \sum_{(t_{j'},x_{k'})\in {\cal L}^{n}_0}\nonumber\\
&\int_{R_{j'-1,k'}}\left|\frac{\sqrt{n}}{2}\left(\P\left(\frac{S_{j-j'}}{\sqrt{n}}=x_{k'}-x_{k}\right)-\P\left(\frac{S_{j-j'}}{\sqrt{n}}=-x_{k'}-x_{k}\right)\right)\right.\nonumber\\
&\left. -\frac{1}{\sqrt{2\pi (t_j-s)}}\left(e^{-\frac{|x_k-y|^2}{2(t_j-s)}}-e^{-\frac{|x_k+y|^2}{2(t_j-s)}}\right)\right|dy\, ds \label{esh:je2}
\end{align}
First note, that we have 
\begin{align}
&\sum_{(t_{j'},x_{k'})\in {\cal L}^{n}_0}\int_{R_{j'-1,k'}}\left|\frac{\sqrt{n}}{2}\P\left(\frac{S_{j-j'}}{\sqrt{n}}=\pm x_{k'}-x_{k}\right)\right|\nonumber \\
\leq & \frac{1}{n}\sum_{(t_{j'},x_{k'})\in {\cal L}^{n}_0}\P\left(\frac{S_{j-j'}}{\sqrt{n}}=\pm x_{k'}-x_{k}\right)\leq  \frac{1}{n}.
\end{align}
Further, 
\begin{align}
& \sum_{(t_{j'},x_{k'})\in {\cal L}^{n}_0}\int_{R_{j'-1,k'}}\left|\frac{1}{\sqrt{2\pi (t_j-s)}}\left(e^{-\frac{|x_k-y|^2}{2(t_j-s)}}-e^{-\frac{|x_k+y|^2}{2(t_j-s)}}\right)\right|dy\, ds\nonumber\\
\leq & \int_{t_{j'-1}}^{t_{j'}}\int_0^{\infty}\left|\frac{1}{\sqrt{2\pi (t_j-s)}}\left(e^{-\frac{|x_k-y|^2}{2(t_j-s)}}-e^{-\frac{|x_k+y|^2}{2(t_j-s)}}\right)\right|dy\, ds\nonumber\\
\leq &\int_{t_{j'-1}}^{t_{j'}}\int_{-\infty}^{\infty}\frac{1}{\sqrt{2\pi (t_j-s)}}e^{-\frac{|x_k-y|^2}{2(t_j-s)}} \, dy\, ds\nonumber\\
= & t_{j'}-t_{j'-1}=\frac{1}{n}.
\end{align}

Using triangle inequality the expression $(\ref{esh:je2})$ is bounded by
$$\left(\frac{C_A}{(j-j')^{3/2}}\sqrt{n+n^2}+\frac{C}{n^{3\alpha}}\right)\frac{3}{n}.$$
Finally, the $\sum_{j-j'\geq n^{\alpha +1/3}}$-part of the sum is less than
$$\sum_{j-j'\geq n^{\alpha +1/3}}\left(\frac{3C_A}{(j-j')^{3/2}}\sqrt{\frac{1}{n}+1}+\frac{3C}{n^{3\alpha+1}}\right).$$
Since the sum goes over $j'$ with the property $n^{\alpha +1/3}\leq j-j'\leq \mic{nA}\leq nA$ we have 
$$\leq \int_{n^{\alpha +1/3}-1}^{\infty}\frac{1}{h^{3/2}}\, dh+ \frac{3A}{n^{3\alpha}}=$$
\be (n^{\alpha +1/3}-1)^{-1/2}+\frac{3A}{n^{3\alpha}}. \ee

Now, from $(\ref{esh:bnd1a})$ and the last bound we have that 
$$\E[(V_j^k-v(t_j,x_k))^2]\leq 4n^{\alpha-1/6}+\frac{8}{\sqrt{\pi}}\sqrt{n^{\alpha - 2/3}+n^{-1}}+(n^{\alpha +1/3}-1)^{-1/2}+\frac{3}{n^{3\alpha}}$$
Set $\gamma:=\min\{1/6-\alpha, 3\alpha\}$, where $\alpha >0$ such that $\gamma>0$. Then there exists
$\Gamma>0$ (that depends on $\gamma$) such that 
$$\E[(V_j^k-v(t_j,x_k))^2]\leq \frac{\Gamma}{n^\gamma}.$$
This discussion proves Theorem \ref{thm:cnvemt}.\vspace{0.2cm}

Now from Theorem \ref{teo:max_nnc} (inequality $(\ref{bb:max:7a})$) we know that 
\be\E\left[\sup_{(t_j,x_k)\in {\cal L}_0^n\cap K}|V_j^k-v(t_j,x_k)|^2\right]\leq 2\frac{\Gamma}{n^{\gamma}}\ln (n^{3/2}AB+1).\label{esh:in10}\ee
since $|{\cal L}_0^n\cap K|\leq \mic{nA}\mic{\sqrt{n}B}\leq n^{3/2}AB$.\vspace{0.2cm}

The following shows th convergence of the Euler method.
\pos{Using the same notation as before we have,
 \be\sup_{(t_j,x_k)\in {\cal L}_0^n\cap K}|V_j^k-v(t_j,x_k)|\to 0\quad \textrm{as}\ \ n\to \infty\quad a.s. \label{esh:cnv1}\ee
Further for $\beta <\gamma/2$, 
$$n^{\beta}\sup_{(t_j,x_k)\in {\cal L}_0^n\cap K}|V_j^k-v(t_j,x_k)|\to 0\quad \textrm{as}\ \ n\to \infty\quad a.s.$$
\dok{Since $\{V_j^k-v(t_j,x_k):(t_j,x_k)\in {\cal L}_0^n\cap K\}$ is a family of Gaussian 
random variables, by Theorem \ref{teo:cncrt} the inequality $(\ref{esh:in10})$ implies (\ref{esh:cnv1}). 

The second inequality follows from the fact that  for every $\varepsilon\in (0,\gamma)$ there exists $C>0$ such that
$$\frac{\Gamma}{n^{\gamma}}\ln (n^{3/2}AB)\leq \frac{C}{n^{\gamma-\varepsilon}}.$$
By using Corollary \ref{cor:cncrt}. we get the desired result.}}

\index{Euler method|see{Heat equation}}

\section{Convergence of interpolation of the Euler method}\label{cvnesh}

We know the values $V_j^k$ at $(t_j,x_k)\in {\cal L}_0^n$ and we
want to approximate the solutions to the heat equation on the rest of the plane.

We are doing the interpolation\index{Interpolation} in the following way:

\begin{itemize}
 \item We do a linear interpolation between points 
$(t_j,x_k)$ and $(t_{j+1},x_{k+1})$ for all $(t_j,x_k)\in {\cal L}_0^n$.
\item  We do a linear interpolation between points 
$(t_j,x_k)$ and $(t_{j-1},x_{k+1})$ whenever $(t_j,x_k)$, $(t_{j-1},x_{k+1})\in {\cal L}_0^n$.
\item We set all values on $x$ and $y$ axis to be $0$.

\item Finally, each point $(t,x)$ is linearly approximated by the values 
$(t,x_-)$ and $(t,x_+)$ the closest points previously defined with respect to the $x$-coordinate.
\end{itemize}

\begin{figure}[ht]
\begin{center}
\psfrag{0}{\small $\boldsymbol{0}$\ }
\psfrag{P}{\small $\boldsymbol{(t,x)}$\ } \psfrag{Q}{\small $\boldsymbol{(t,x_-)}$}
\psfrag{R}{\small $\boldsymbol{(t,x_+)}$}
 \includegraphics[width=6cm]{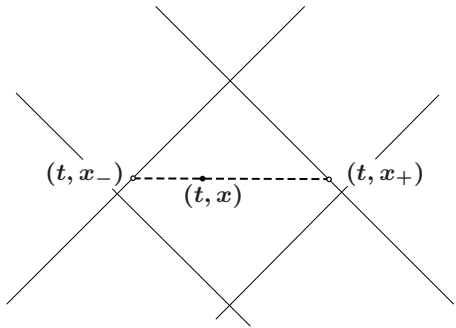} \\
\caption{Interpolation}\label{sl17:2}
  \end{center}
\end{figure}

In this way we obtain the approximation $V_n(t,x)$ of the stochastic heat equation 
on $\R_+^2$, and we want to show convergence to $u$ on compact sets, where
\be  v(t,x)=\int_0^t\frac{1}{\sqrt{2\pi (t-s)}}\int_{\R^+}\left(e^{-\frac{|x-y|^2}{2(t-s)}}-e^{-\frac{|x+y|^2}{2(t-s)}}\right)\W(ds,dy).\ee

We will show that for a compact set $K\subset \R_+^2$ we have
$$\sup_{(t,x)\in K}|V_n(t,x)-v(t,x)|\to 0.$$

\teo{\label{teo:intBM}For a compact set $K\subset \R_+^2$ we have 
$$\lim_{n\to \infty}\sup_{(t,x)\in K}|V_n(t,x)-v(t,x)|= 0\quad a.s. $$
\dok{Pick $K$, and then pick $K_{ab}=[0,a]\times [0,b]$, such that 
$$\sup\{ x: (t,x)\in K\}<a,$$
and 
$$\sup\{ t: (t,x)\in K\}<b.$$ 
For large $n$, the points ${\cal L}_0^n\cap K_{ab}$ will be enough
to calculate the value of $V_n$ for all points on $K$.

$V_n$ at point $(t,x)$ can be written as a convex combination of the values 
of the four points in ${\cal L}_0^n$ that make the rhombus in which the point is.
Therefore 
$$V_n(t,x)=\alpha_1 V_j^k+\alpha_2V_j^{k+2}+\alpha_3V_{j+1}^{k+1}+\alpha_4V_{j-1}^{k+1}, $$
where $\alpha_1+\alpha_2+\alpha_3 +\alpha_4=1$ (one or more of $\alpha_i$'s will be $0$).

Now,
\begin{align*}
 V_n(t,x)-v(t,x)&=\alpha_1 (V_j^k-v(t_j,x_k))+\alpha_2(V_j^{k+2}-v(t_j,x_{k+2}))+\alpha_3(V_{j+1}^{k+1}-v(t_j,x_{k+2}))\\
&+\alpha_4(V_{j-1}^{k+1}-v(t_{j-1},x_{k+1}))\\
&+\alpha_1 (v(t_j,x_k)-v(t,x))+\alpha_2(v(t_j,x_{k+2})-v(t,x))\\
&+\alpha_3(v(t_j,x_{k+2})-v(t,x))+\alpha_4(v(t_{j-1},x_{k+1})-v(t,x)).
\end{align*}
%
Not that $v$ is a continuous function, and $K_{ab}$ a compact set. Therefore, $u$
is uniformly continuous on $K_{ab}$. The distance between $(t,x)$ and points $(t_j,x_k)$,
$(t_{j-1},x_{k+1})$, $(t_{j+1},x_{k+1})$ and $(t_j,x_{k+2})$ goes to $0$ uniformly. So, by uniform
continuity for any $\varepsilon>0$ we have $|v(t_*,x_{*})-v(t,x)|<\varepsilon$,
when the distance between $(t_*,x_*)$ and $(t,x)$ is less than some $\delta$.

Hence, for a large $n$

$$\sup_{(t,x)\in K}|V_n(t,x)-v(t,x)|\leq \sup_{{\cal L}_0^n\cap K_{ab}}|V_j^k-v(t_j,x_k)|+\varepsilon.$$
When $n\to \infty$ we have, by $(\ref{esh:cnv1})$
$$\limsup_{n\to \infty }\sup_{(t,x)\in K}|V_n(t,x)-v(t,x)|\leq \varepsilon.$$
Finally, since $\varepsilon$ is an arbitrary positive number, the claim follows.}}

We have shown that the Euler method converges uniformly on compact subsets $\R_+^2$ to the stochastic heat equation.

\prop{\label{prop:eumtsum}For each $n$ let $(E_{jk}^n)$ be an i.i.d. sequence 
of $N(0,\frac{1}{\sqrt{2n}})$, and 
$$Y_{k}^{j+1}=\frac{1}{2}(Y_{k+1}^j+Y_{k-1}^j)+E_{jk}, \quad j\geq 0,k\geq 1,$$
with $Y^0_{k}=0$ and $Y^j_0=0$. Then the described interpolation $\widetilde{Y}_n(t,x)$
converges in distribution to the solution of the stochastic heat equation 
$$\left\{\begin{array}{l}
          v_t=\frac{1}{2}v_{xx}+\W \\
	  v(0,\cdot)=0, v(\cdot,0)=0
         \end{array} \right.$$ 
}\index{Heat equation!stochastic|)}

\subsection{Euler method with weaker noise}

We finish the study of Euler method by looking at the case when the noise is weak, so that
in the limit it has no effect. What happens 
if $E_{jk}^n$ would be distributed as $N(0,\frac{1}{\sqrt{2}n^{1/2+\alpha}})$
in Proposition \ref{prop:eumtsum}? If we have noise with slightly lower variance, 
would we still have convergence. It turns out we would and that convergence would be to $0$.

\lem{\label{lem:eushto0}Let $\alpha>0$. For each $n$ let $(E_{jk}^n: k\equiv j+1\ (\textrm{mod}\ 2))$ be an i.i.d. sequence 
of $N(0,\frac{1}{\sqrt{2}n^{1/2+\alpha}})$, and 
\be Y_{k}^{j+1}=\frac{1}{2}(Y_{k+1}^j+Y_{k-1}^j)+E_{jk}, \quad j\geq 0,k\geq 1,\label{eq:heeulal}\ee
with $Y^0_{k}=0$ and $Y^j_0=0$. Then for all $0<a<1+2\alpha$,$b>0$ and $A,B>0$ there exist $\Gamma>0$ and $\gamma>0$

such that 
\be\E\left[\sup_{j\leq An^a,k\leq Bn^b}|Y^j_{k}|^2\right] \leq \frac{\Gamma}{n^{\gamma}}.\label{esh:ineq0}\ee
\dok{It can be shown that 
$$Y^j_k=\sum_{j'=1}^{j}\sum_{k'=1}^{k+j} \left(\P(S_{j-j'}=k'-k)-\P(S_{j-j'}=-k-k')\right)E_{j'-1,k'}^n$$
is a solution to $(\ref{eq:heeulal})$, where $(S_j)$ is a simple random walk. 
From the fact that $E_{jk}^n$ are i.i.d. we have 
$$\E(Y^j_k)^2=\sum_{j'=1}^{j}\sum_{k'=1}^{k+j} \left(\P(S_{j-j'}=k'-k)-\P(S_{j-j'}=-k-k')\right)^2\frac{1}{\sqrt{2}n^{1/2+\alpha}}.$$
Since $|k'-k|<|k'+k|$,
we have $0\leq \P(S_{j-j'}=k'-k)-\P(S_{j-j'}=-k-k')\leq \P(S_{j-j'}=k'-k)$, and so 
$$\E(Y^j_k)^2\leq \sum_{j'=1}^{j}\sum_{k'=1}^{k+j} \P(S_{j-j'}=k'-k)^2\frac{1}{\sqrt{2}n^{1/2+\alpha}}.$$
Now, it follows from the properties of the random walk that 
$$\P(S_{j-j'}=k'-k)^2=\P(S_{2(j-j')}-S_{(j-j')}=-(k'-k))\P(S_{j-j'}=k'-k)=$$
$$=\P(S_{2(j-j')}-S_{(j-j')}=-(k'-k),S_{j-j'}=k'-k)=\P(S_{2(j-j')}=0,S_{j-j'}=k'-k).$$
Furthermore, 
$$\sum_{k'=1}^{k+j} \P(S_{j-j'}=k'-k)^2 = \sum_{k'=1}^{k+j} \P(S_{2(j-j')}=0,S_{j-j'}=k'-k)$$
$$ \leq \sum_{k'}\P(S_{2(j-j')}=0,S_{j-j'}=k'-k) =\P(S_{2(j-j')}=0).$$ 
If $j\leq An^a$ by Stirling's formula we have
$$\E(Y^j_k)^2\leq \frac{1}{\sqrt{2}n^{1/2+\alpha}}\sum_{j'=1}^{An^a}\P(S_{2(j-j')}=0)\sim \frac{1}{\sqrt{2}n^{1/2+\alpha}}\sum_{j'=1}^{An^a}\frac{1}{\sqrt{\pi j}}.$$
The last sum can be bounded by $1+\int_1^{An^a}\frac{1}{\sqrt{t}}\, dt=2\sqrt{An^a}$. Hence, there exists $C>0$ 
such that
$$\E(Y^j_k)^2\leq C\frac{2\sqrt{An^a}}{\sqrt{2}n^{1/2+\alpha}} =\frac{C\sqrt{2A}}{n^{1/2+\alpha-a/2}}.$$
Now, by $(\ref{bb:max:7a})$ we have 
$$\E\left[\sup_{j\leq An^a,k\leq Bn^b}|Y_{jk}|^2\right] \leq \frac{C\sqrt{2A}}{n^{1/2+\alpha-a/2}} \ln(ABn^{ab}+1).$$
Now, for any $\gamma \in (0,1/2+\alpha-a/2)$ there exists $\Gamma$ such that $(\ref{esh:ineq0})$ holds.}}

\prop{\label{prop:eushto0}Let $Y^j_k$ be as in the previous Lemma. $t_j=jn^{-1}$ and $x_k=\frac{k}{n^{1/2+\alpha}}$, and set ${\cal L}^n_0=\{(t_j,x_k): k\equiv j \ (\textrm{mod}\ 2)\}$, if we define $V_n(t_j,x_k)$ for $(t_j,x_k)\in {\cal L}^n_0$ to have a value 
$Y^j_k$, and do the interpolation described in \S\ref{cvnesh}, $V_n$ converges in distribution 
to $0$. \dok{Let $K$ be a compact set, there exists $A>0,B>0$ such that $K\subset [0,A]\times [0,B]$,  the value $\max_{(t,x)\in K}V_n(t,x)$ 
is obtained at some point ${\cal L}^n_0\cap [0,A]\times [0,B]$. Now, from Lemma \ref{lem:eushto0}. and Theorem \ref{teo:cnv0}. the claim follows.}}
\section{Euler method with initial value condition and no external noise }
We saw what happens with Euler scheme when the external source
is replaced by noise.

In this section we shortly comment what happens with the Euler scheme
for the heat equation with no external force and a non-zero initial value
function:
\be W^{j+1}_{k}=\frac{1}{2}(W^j_{k+1}+W^j_{k-1}),\quad k\geq 1, j\geq 0\label{esh:eq2}\ee 
We will work with the case
when $W^j_0=0$ and $W^0_k=g(x_k)$. 

In order to simplify our work we will set  $W^0_{-k}:=-W^0_k$. In this way we don't have 
to think about the boundary condition $W^j_{0}=0$, because it is easy to show that
the scheme

\be\begin{cases} 
W^{j+1}_{k}=\frac{1}{2}(W^j_{k+1}+W^j_{k-1}), & k\in \mathbb{Z}, j\geq 0 \\  
W^0_{k}=\tilde{g}(x_k) & k\in \mathbb{Z}.  
   \end{cases}
\label{esh0n}\ee

Where $\tilde{g}(x)=g(x)$ for $x\geq 0$ and $\tilde{g}(x)=-g(-x)$ for $x\leq 0$.

With this setup the following lemma holds.

\lem{For the scheme given by (\ref{esh0n}) we have 
$$W^j_k=\E\left[\tilde{g}\left(\frac{S_j}{\sqrt{n}}+\frac{k}{\sqrt{n}}\right)\right].$$}

Under some mild properties on $\tilde{g}$ we have that the Euler method converges to
the heat equation

\be\begin{cases} 
\partial_t w=\partial_{xx}w & \textrm{on}\ \R\times (0,\infty),  \\  
w(x,0)=\tilde{g}(x) & x\in \R.  
   \end{cases}
\label{heg}\ee
 
It is well-known that $u$ restricted to $\R^+\times (0,\infty)$ is the solution to
\be\begin{cases} 
\partial_t w=\partial_{xx}w & \textrm{on}\ \R\times (0,\infty),  \\  
w(x,0)=g(x) & x\in \R,  \\
w(0,t)=0 & t\geq 0.
   \end{cases}
\label{heg0}\ee

\lem{\label{lem:cnin}If $g:\R^+\to \R$ is continuous and we have 
\be\limsup_{x\to\infty}\frac{|g(x)|}{|x|}<\infty,\label{he:con1}\ee
then for a compact set $K\subset \R^+\times [0,\infty)$ we have 
$$\lim_{n\to \infty }\sup_{(x_k,t_j)\in {\cal L}^n_0\cap K}|w(x_k,t_j)-W^j_k|=0,$$ 
where $W^j_k$ is the solution to $(\ref{esh:eq2})$ and $w$ to $(\ref{heg0})$.
\dok{First, note that there exists $L,C>0$ such that 
$$|g(x)|^2\leq L|x|^2+C.$$
Further we have $|\tilde{g}(x)|^2\leq L|x|^2+C$, and now for any $a<b$ we have
$$\sup_{r\in [a,b]}|\tilde{g}(x+r)|^2\leq 2L\max\{|a|^2, |b|^2\}+2L|x|^2+C=2L\max\{|a|^2, |b|^2\}+C+\int_0^{|x|}4Ly\, dy,$$
for all $x$. Also, we have 
$$\int_0^\infty Lye^{-\frac{y^2}{2\tau}}\, dy<\infty$$
for all $\tau>0$. For a compact set $K\subset \R \times [0,\infty) $,there exists 
$\tau>0$ and $a<b$ such that $K\subset [a,b]\times [0,\tau)$, hence Theorem \ref{teo:aprxclhe:3}
implies that
$$\lim_{n\to\infty}\sup_{(x_k,t_j)\in {\cal L}^n_0\cap K}|W^j_k-w(x_j,t_j)|= 0.$$
Where $W^j_{k}$ is the solution to $(\ref{esh0n})$ and $w$ to $(\ref{heg0})$.\vspace{0.2cm}

It is now clear that the result follows.}}

We can interpolate $\{W^j_k\ :\ (x_k,t_j)\in {\cal L}^n_0\}$ as in sections 
\ref{cvnesh} and \ref{ntwbb} to obtain a function $W_n(x,t)$. The following can be shown using
the same techniques.
\teo{Let $W_n$ be the interpolation  described in \ref{cvnesh} of the Euler method, and $w$ the solution to the equation 
$(\ref{heg0})$ where $g$ satisfies $(\ref{he:con1})$. Then we have 
$$\lim_{n\to\infty}W_n=w,$$
uniformly on compact sets.} 

\pos{\label{pos:bmheq}If we set $g(x)=B(x)$ where $(B(x):x\geq 0)$ is Brownian motion, then for $W_n$ be the described interpolation of the Euler method, and $w$ the solution to the equation 
$(\ref{heg0})$ we get  
$$\lim_{n\to\infty}W_n=w,$$
uniformly on compact sets almost surely.
\dok{Follows from the strong law of large numbers for the Brownian motion since
$$\lim_{x\to \infty}\frac{|B(x)|}{|x|}=0\ a.s.$$
Now we can apply Lemma \ref{lem:cnin}, and all the results after.}}

\index{Heat equation!Euler method|)}

\chapter{Limit of the natural Brownian motion on a rhombus grid}\label{sec:04}
\section{Natural Brownian motion on a rhombus grid}
We will investigate the process on a rhombus grid where
the ratio of diagonals depends on $n$. A natural question is what happens
when $n\to \infty$.

\begin{figure}[ht]
\begin{center}
\psfrag{x}{\small $\boldsymbol{x}$ } 
\psfrag{a}{\textcolor{red}{$\boldsymbol{n^{-1}}$}} 
\psfrag{b}{\textcolor{blue}{$\boldsymbol{n^{-1/2-\alpha}}$}} 
\psfrag{t}{\small $\boldsymbol{t}$ }
\includegraphics[width=15cm]{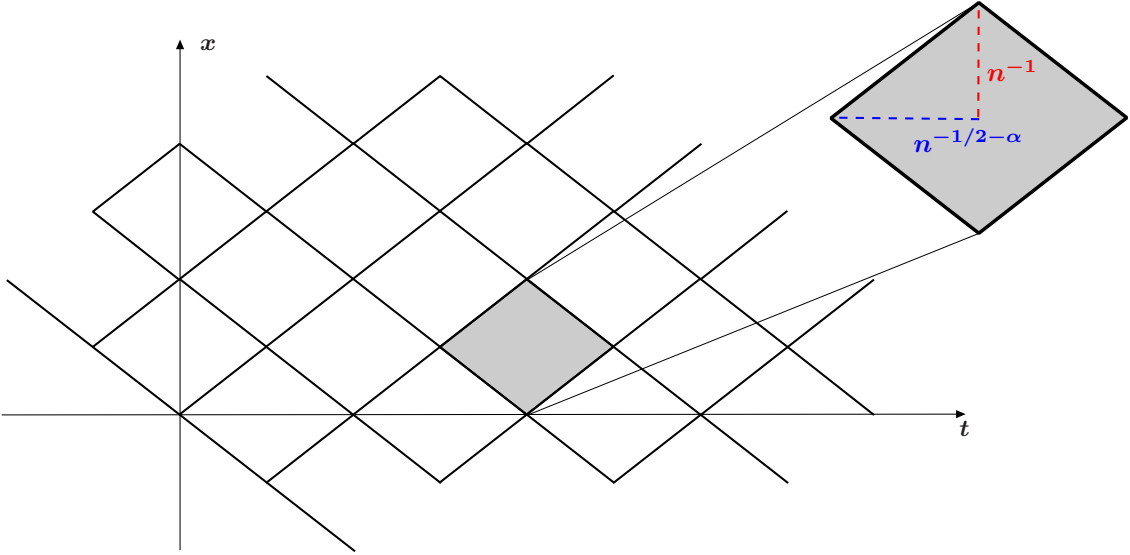}\\
\caption{$\alpha$-rhombus grid} \label{sl22}
  \end{center}
\end{figure}

\defi{We will call the TLG$^*$ $\cG$  the $(\alpha,n)$-rhombus grid if the
plane is divided into congruent rhombuses, diagonals of which are parallel to the
$x$ and $t$ axis, the  length of the  half-diagonal parallel to $x$ is $\frac{1}{n}$ and the length of the half-diagonal 
parallel to $t$
is $\frac{1}{n^{1/2+\alpha}}$, and there is a rhombus that has vertex $(0,0)$. (See Figure \ref{sl22}.)}\vspace{0.2cm}

\noindent\emph{Remark.} In our representation  the vertices are represented by the set, where $t_j=jn^{-1/2-\alpha}$ and 
$x_k=k/n$:
$${\cal L}^{\alpha,n}_0=\{(t_j,x_k): k,j\geq 0,\ j \equiv k \ (\textrm{mod} \ 2)\}$$

Using Theorem \ref{teo:burdzy-pal}. we construct a natural two-sided Brownian motion $X_{\alpha,n}$
on this grid. Recall, that if $(W_1(t):t\geq 0)$ and $(W_2(t):t\geq 0)$ are two independent
Brownian motions, then
$$B(t):=\begin{cases}
         W_1(t), & t\geq 0\\
         W_2(-t), & t< 0\\
        \end{cases}
$$
is a two-sided Brownian motion. It is not hard to check that this is 
a Markov process on $T=\R$. Further, covariance of this process 
is 
$$C_B(t,s)=\frac{1}{2}(|t|+|s|-|t-s|).$$

The following result will be useful.

\lem{\label{ind:posneg}The processes $(X_{\alpha,n}(\bar{t}):t\geq 0)$ and  $(X_{\alpha,n}(\bar{t}):t\leq 0)$ are independent.
\dok{Let $Y^+_{\alpha,n}$ be a natural $\cP_+$-Markovian process on a $(\alpha,n)$-grid,
where $\cP_+$ is the distribution of
$$B^0_+(t):=\begin{cases}
         W_1(t), & t\geq 0\\
         0, & t< 0\\
        \end{cases}.$$ 
In the same way we can construct $Y^-_{\alpha,n}$ as a natural $\cP_-$-Markovian process on a $(\alpha,n)$-grid,
where $\cP_-$ is the distribution of
$$B^0_-(t):=\begin{cases}
         0, & t\geq 0\\
         W_2(-t), & t< 0\\
        \end{cases}.$$
We can construct $Y^+_{\alpha,n}$ and $Y^-_{\alpha,n}$ such that they are independent
and on the same space and using the same TLG$^*$-towers.
Then, it is not hard to see that $Y^+_{\alpha,n}+Y^-_{\alpha,n}$ in
each member of a TLG$^*$-tower has the same distribution as $X_{\alpha,n}$
on this TLG$^*$. Therefore, the distribution of $X_{\alpha,n}$ and  
$Y^+_{\alpha,n}+Y^-_{\alpha,n}$ are the same by Theorem \ref{teo:burdzy-pal}. Furthermore, 
$((X_{\alpha,n}(\bar{t}):t\geq 0),(X_{\alpha,n}(\bar{t}):t\leq 0))$
are distributed as $(Y^+_{\alpha,n},Y^-_{\alpha,n})$.
}}

Due to the last lemma, we can focus on what happens with the process $X_{\alpha,n}(x,t)$ for $t\geq 0$.

The final distribution of the process, by Theorem \ref{teo:burdzy-pal}. doesn't depend 
on the way we construct the process. We fix a construction that we will refer to.\vspace{0.2cm}

For our construction we need:

\begin{itemize}\index{Brownian bridge|(}
 \item two-sided Brownian motion $(B(t):t\in \R)$;
 \item for  $j\neq -1 $, $k\in \Z$: $(B^{br}_{jk}(t):t\in [t_j,t_{j+2}])$ be a collection
of Brownian bridges ($n\in \N$);
 \item for $j=-1 $, $k\in \Z$: $(B^{br}_{jk-}(t):t\in [t_j,t_{j+1}])$, $(B^{br}_{jk+}(t):t\in [t_{j+1},t_{j+2}])$ be a collection
of Brownian bridges ($n\in \N$);
\end{itemize}
all of these things are independent.

\begin{description}
 \item[Step 0] We run the two-sided Brownian motion on the time-path $\sigma$ that is going through 
$(-\infty,\infty)\times [0,\frac{1}{n}]$ (this will be our spine), that is we define $X_{\sigma}(t)=B(t)$.
(See Figure \ref{sl21ff}.)

\begin{figure}[ht]
\begin{center}
\psfrag{x}{\small $\boldsymbol{x}$ } 
\psfrag{t}{\small $\boldsymbol{t}$ } 
\psfrag{0}{\small $\boldsymbol{0}$ } 
\psfrag{1}{\small $\boldsymbol{1/n}$ }    
\psfrag{s}{\small $\textcolor{red}{\boldsymbol{\sigma}}$ }    
\includegraphics[width=12cm]{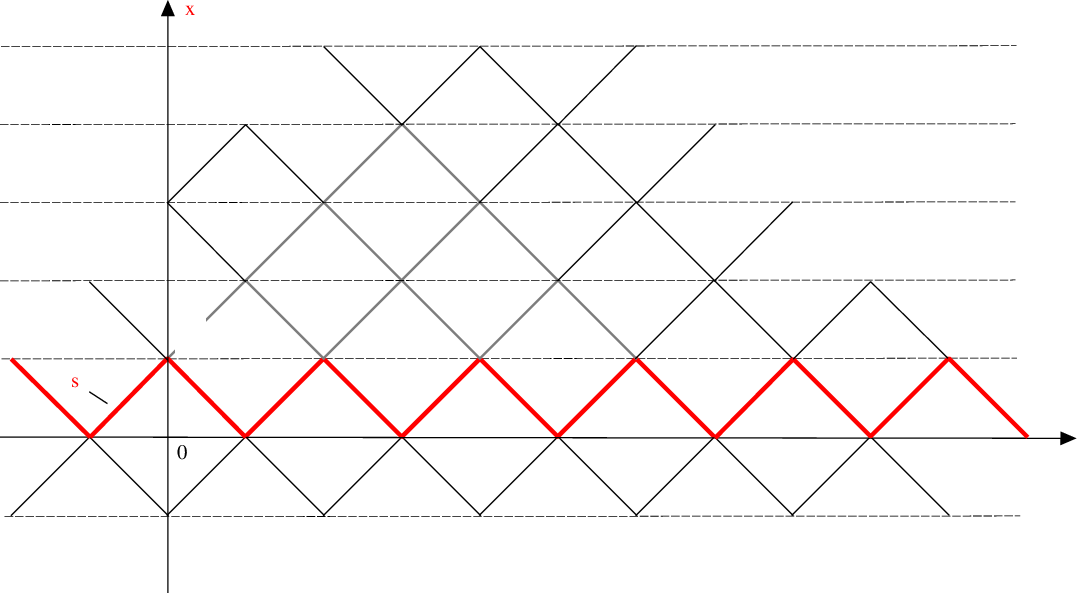}\\
\caption{} \label{sl21ff}
  \end{center}
\end{figure}

\item[Step 1]  Now if we defined the process at points $(t_j,x_k)$ and $(t_{j+2},x_{k})$
then we define the process on the time-path $\pi_+$ (if $x_k>0 $) 
$(t_j,x_k)-(t_{j+1},x_{k+1})- (t_{j+2},x_{k})$ or time-path $\pi_-$ $(t_j,x_k)-(t_{j-1},x_{k-1})- (t_{j+2},x_{k})$
(if $x_k\leq 0$)
by setting $X_{\pi_{\pm}}$ to be:
\begin{itemize}
 \item if $j=-1$ two Brownian bridges (on the intervals $[t_j,t_{j+1}]$ and $[t_{j+1},t_{j+2}]$, where
the value of the process at $t_j$ is $X_{\alpha,n}(t_j,x_k)$, $t_{j+1}$ is $0$ and $t_{j+2}$ is $X_{\alpha,n}(t_{j+2},x_k)$;
 \item a Brownian bridge at times $t_j$ and $t_{j+2}$ between 
values $X_{\alpha,n}(t_j,x_k)$ and $X_{\alpha,n}(t_{j+2},x_k)$.
\end{itemize}
\begin{figure}[ht]
\begin{center}
\psfrag{A}{\small $\boldsymbol{(t_j,x_k)}$\ }
\psfrag{B}{\small $\boldsymbol{(t_{j+2},x_k)}$\ } 
\psfrag{C}{\small $\boldsymbol{(t_{j+1},x_{k+1})}$}
\includegraphics[width=5cm]{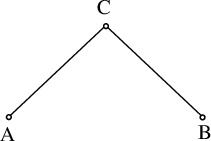}
  \end{center}
\end{figure} 
  Specially, if the path if $t_j\geq 0$
then $X_{\pi_{\pm}}$ will be of the form
\be X_{\pi_{\pm}}(t):= \frac{t_{j+2}-t}{t_{j+2}-t_j}X_{\alpha,n}(t_j,x_k)+\frac{t-t_{j}}{t_{j+2}-t_{j}}X_{\alpha,n}(t_{j+2},x_k)+B^{br}_{jk}(t)
\label{eq:cons_grid}
\ee
where $B^{br}_{jk}$ is a Brownian bridge that has value $0$ at times $t_j$ and $t_{j+2}$,
and independent of the other Brownian bridges.

If $j=-1$ ($t_{j+1}=0$), then $X_{\pi_{\pm}}$ on $[t_{j+1},t_{j+2}]$ is equal 
\be X_{\pi_{\pm}}(t):= \frac{t}{t_{j+2}}X_{\alpha,n}(t_{j+2},x_k)+B^{br}_{jk+}(t)
\label{eq:cons_grid2}\ee
where $B^{br\pm}_{jk}(t)$ is a Brownian bridge with value $0$ at times $t_{j+1}$ and $t_{j+2}$.
\item[Step 2] We repeat {\bf Step 1} in such a way that at every point in ${\cal L}^{\alpha,n}_0$
the process will be eventually defined.
\end{description}

\noindent\emph{Remark.}
Along every path from $-\infty$ to $+\infty$ we have a two-sided Brownian motion.\vspace{0.2cm}

Having in mind Lemma \ref{ind:posneg}, we will focus our attention to the process $X_{\alpha,n}$
defined in the first quadrant. The convergence of the process in other quadrants  can be shown  in a similar way.\vspace{0.2cm}

The most important thing to note from the construction of the process, that
if from the equation $(\ref{eq:cons_grid})$ is that when we set $t=t_{j+1}$ we get:
\be X_{\alpha,n}(t_{j+1},x_{k+1})=\frac{1}{2}X_{\alpha,n}(t_j,x_k)+\frac{1}{2}X_{\alpha,n}(t_{j+2},x_k)+E_{j+1,k+ 1}, \label{eq:eush}\ee
where 
$$E_{j+1,k+1}=B^{br}_{jk}(t_{j+1})\stackrel{d}{=}N(0,2^{-1/2}n^{-1/2-\alpha}),$$
 for $j\geq 0, k\geq 1$ such that $(t_j,x_k)\in {\cal L}^{\alpha,n}$. This is a form of 
the discrete stochastic heat equation\index{Heat equation!discrete}\index{Heat equation!stochastic} (see \cite{rw_heat}) with random external source.

We discussed the convergence of these equations in Chapter \ref{chp:heatrw} (see \S\ref{sec:eumth}), that is what happens when $n \to \infty$.

\subsection{Interpolation} \label{subsec:intbnet}

\begin{figure}[ht]
\begin{center}
\psfrag{0}{\small $\boldsymbol{0}$\ }
\psfrag{P}{\small $\boldsymbol{(t,x)}$\ } \psfrag{Q}{\small $\boldsymbol{(t_-,x)}$}
\psfrag{R}{\small $\boldsymbol{(t_+,x)}$}
\includegraphics[width=5cm]{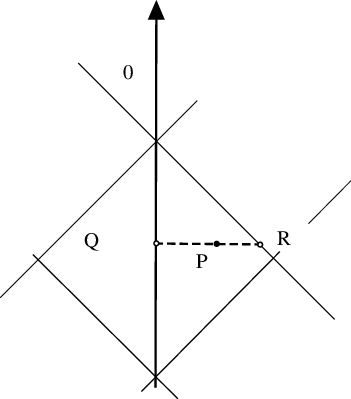}\ \ \includegraphics[width=6cm]{sl17} \\
\caption{Interpolation of the process}\label{sl17}
  \end{center}
\end{figure}

Now our process is defined on the representation of the whole $(\alpha,n)$-rhombus grid, and we will extend
the definition of the process on the whole plane (see Figure \ref{sl17}):
\begin{itemize}
 \item $X_{\alpha,n}(0,x)=0$ (the process on the $x$-axis is $0$);
 \item by interpolation we will extend the definition of our process on the whole plane:
$$X_{\alpha,n}(t,x):=\frac{t_+-t}{t_+-t_-}X_{\alpha,n}(t_-,x)+\frac{t-t_-}{t_+-t_-}X_{\alpha,n}(t_+,x),$$
where $(t_+,x)$ and $(t_-,x)$ are points on the representation of the graph or on the $x$-axis that
are the closest to $(t,x)$. 

\end{itemize}

In further text we will denote the interpolated process as $(Y_{\alpha,n}(t,x) :t\in \R,x\in \R)$.
Note that this is a continuous Gaussian process on $\R^2$. 

We will study two cases $\alpha=0$ and $\alpha>0$, and how does $Y_{\alpha,n}$ 
behaves as $n\to\infty$.

\section{Network of Brownian bridges}\label{ntwbb}
The final result that we need to show that the Brownian motion on a rhombus grid
converges is the fact that a network of Brownian bridges will converge to $0$ on compact sets.\vspace{0.2cm}


%

\teo{\label{teo:nbrb}Let $K$ be a compact subset of $\R_+^2$, then 
$$Z_n(K)=\sup\{\max|B_{jk(+)}^{br}|: (t_j,x_k)\in K\}\stackrel{L^2}{\to} 0.$$
\index{Brownian bridge!maximum}
\dok{Pick $K$, and then pick $K_{ab}=[0,a]\times [0,b]$, such that 
$$\sup\{ x: (t,x)\in K\}<a,$$
and 
$$\sup\{ t: (t,x)\in K\}<b.$$ 
For each point in ${\cal L}_0^n$ there are at most 2 bridges going out of this point
(in the direction of time). Hence on $K_{ab}$ there are at most $an^{1/2+\alpha} \cdot  bn=abn^{3/2+\alpha}$ such bridges.
The bridges on $K_{ab}$ define the $Z_n$ on $K$. We have at most $abn^{3/2+\alpha}$ on intervals of length
$[0,n^{-\frac{1}{2}-\alpha}]$. Hence, by Corollary \ref{cor:bbmax} (inequality $(\ref{bb:max:5a})$) and the fact that for $(B^{br}(t):t\in [0,1])$ Brownian bridge on $[0,1]$
$(n^{-\frac{1}{4}-\alpha/2}B^{br}(n^{1/2+\alpha}t)\ :\ t\in [0,n^{-1/2-\alpha}])$ is the Brownian bridge on $[0,n^{-1/2-\alpha}]$, the second moment of their maximum is bounded by
\be\frac{1}{2n^{\frac{1}{2}+\alpha}}\ln (abn^{3/2+\alpha}+1).\label{ntbbb:e1}\ee
Since the maximum is obtained in the points where the Brownian bridges have been 
defined, the claim follows. }}

From the rate of convergence in $(\ref{ntbbb:e1})$, and Theorem \ref{teo:cnv0} we get the following
result.

\pos{\label{pos:nbb}Let $K$ be a compact subset of $\R_+^2$, then 
\be \lim_{n\to \infty}Z_n(K)= 0 \quad a.s.\label{bbnet:cnv}\ee
Further, for $\beta <\frac{1}{4}+\alpha/2$ we have
$$\lim_{n\to \infty}n^{\beta}Z_n(K)= 0$$. }

\section{The main result}

The process that we will be more interested is the interpolation (similar to the one described in 
\S\ref{subsec:intbnet}) between values of $X_{\alpha,n}$ at points in ${\cal L}^{\alpha,n}\cap \R_+^2$:

\begin{itemize}
 \item $X_{\alpha,n}(0,x)=0$ (the process on the $x$-axis is $0$);
 \item $(t_j,x_k)$ and $(t_{j+1},x_{k\pm 1})$ we interpolate between the values at these points;
 \item by interpolation we will extend the definition of our process on the whole plane:
\be X_{\alpha,n}(t,x):=\frac{t_+-t}{t_+-t_-}X_{\alpha,n}(t_-,x)+\frac{t-t_-}{t_+-t_-}X_{\alpha,n}(t_+,x),\label{eq:bnnet}\ee
where $(t_+,x)$ and $(t_-,x)$ are points on the representation of the graph or on the $x$-axis that
are the closest to $(t,x)$. 
\end{itemize}

We will call this process $\widetilde{Y}_{\alpha,n}$.

From the construction of $Y_{\alpha,n}$ 
and $\widetilde{Y}_{\alpha,n}$ it is not hard to see that 
for each rhombus $\Diamond$
the value 
$$\max_{(t,x)\in \Diamond}|\widetilde{Y}_{\alpha,n}(t,x)-Y_{\alpha,n}(t,x)|,$$
due to linear interpolation, is obtained on $\partial \Diamond$. That means we can
focus on the process $|\widetilde{Y}_{\alpha,n}(t,x)-Y_{\alpha,n}(t,x)|$ on the representation 
of the $(\alpha,n)$-rhombus grid. 

\lem{The process on the path $(t_j,x_k)-(t_{j+1},x_{k\pm1})-(t_{j+1},x_{k+1})$
is bounded by
$$|\widetilde{Y}_{\alpha,n}(t,x)-Y_{\alpha,n}(t,x)|\leq 2\max_{t\in [t_j,t_{j+2}]}|B_{jk}^{br}(t)|$$
\dok{From $(\ref{eq:bnnet})$ and $(\ref{eq:eush})$ we have:
$$Y_{\alpha,n}(t,x)=\frac{t_{j+2}-t}{t_{j+2}-t_j}X_{\alpha,n}(t_j,x_k)+\frac{t-t_j}{t_{j+2}-t_j}X_{\alpha,n}(t_{j+2},x_k)+B^{br}_{jk}(t)$$
$$\widetilde{Y}_{\alpha,n}(t,x)=\frac{t_{j+2}-t}{t_{j+2}-t_j}X_{\alpha,n}(t_j,x_k)+\frac{t-t_j}{t_{j+2}-t_j}X_{\alpha,n}(t_{j+2},x_k)+\alpha(t)E_{j+1,k+1},$$
where $|\alpha(t)|<1$ obtained by interpolation between values at points $(t_j,x_k)$ and $(t_{j+1},x_{k+1})$ (if $t\in [t_j,t_{j+1}]$)
or $(t_{j+1},x_{k+1})$ or $(t_{j+2},x_k)$ (for $t\in [t_{j+1},t_{j+2}]$). 
Since $E_{j+1,k+1}=B_{jk}(t_{j+1})$, the claim follows.}}

\lem{\label{lem:shed0}For a compact set $K\subset \R_+^2$ we have 
$$\sup_{(t,x)\in K}|\widetilde{Y}_{\alpha,n}(t,x)-Y_{\alpha,n}(t,x)|\to 0, \ a.s. $$
\dok{There exists a compact set $\hat{K}$ such that all the rhombi whose interior 
intersects $K$, are contained in $\hat{K}$. Now, 
$$\sup_{(t,x)\in K}|\widetilde{Y}_{\alpha,n}(t,x)-Y_{\alpha,n}(t,x)|\leq Z_n(\hat{K}),$$
and by Corollary \ref{pos:nbb}. the claim follows.}}
\index{Brownian bridge|)}


%



\prop{\label{prop:shelim}The process $(\widetilde{Y}_{\alpha,n}(t,x):(t,x)\in \R_+^2)$ converges to 
$u$, where $u$ is the solution to the stochastic heat equation

$$ \begin{array}{cl} \partial_x u=\left\{\begin{array}{cl}
				      \frac{1}{2}\partial_{tt} u+\W & \alpha=0,\\
				      0 & \alpha>0,\\
                                   \end{array}\right. 
&\quad {\rm on}\quad  \R_+^2,\\
    u(0,t)=B(t) & \quad {\rm for}\quad  t\in\R.\\
    u(x,0)=0 & \quad {\rm for}\quad  x\in \R.
    \end{array}
$$

\dok{We will write the process $\widetilde{Y}_{\alpha,n}=\widetilde{Y}_{\alpha,n}^1+\widetilde{Y}_{\alpha,n}^2$, where 
$$\begin{array}{c}
   \widetilde{Y}^{1}_{\alpha,n}(t_{j+1},x_{k+1})=\\ \frac{1}{2}\widetilde{Y}^1_{\alpha,n}(t_j,x_k)+\frac{1}{2}\widetilde{Y}^1_{\alpha,n}(t_{j+2},x_k)+E_{j+1,k+ 1}\\
     \widetilde{Y}^{1}_{\alpha,n}(0,x_k)= 0, \widetilde{Y}^{1}_{\alpha,n}(t_j,0)= 0
  \end{array}
\quad
\begin{array}{c}
   \widetilde{Y}^{2}_{\alpha,n}(t_{j+1},x_{k+1})=\\ \frac{1}{2}\widetilde{Y}^2_{\alpha,n}(t_j,x_k)+\frac{1}{2}\widetilde{Y}^2_{\alpha,n}(t_{j+2},x_k)\\
     \widetilde{Y}^{2}_{\alpha,n}(0,x_k)= 0, \widetilde{Y}^{2}_{\alpha,n}(t_j,0)= B(t_j)
  \end{array}
$$
Now, by Propositions \ref{prop:eumtsum}. and \ref{prop:eushto0}. $\widetilde{Y}^1\stackrel{d}{\to}u^1$ where 
$$\left\{\begin{array}{l}
          u^1_x=\left\{\begin{array}{cl}
               \frac{1}{2}u^1_{tt}+\W & \alpha=0\\
		0 & \alpha >0\\
              \end{array}\right. \\
	  u^1(0,\cdot)=0, u^1(\cdot,0)=0
         \end{array} \right. .$$
By Corollary \ref{pos:bmheq}. $\widetilde{Y}^2\stackrel{d}{\to}u^2$ where 
$$\left\{\begin{array}{l}
          u^2_x=\left\{\begin{array}{cl}
               \frac{1}{2}u^2_{tt} & \alpha=0\\
		0 & \alpha >0\\
              \end{array}\right. \\
	  u^2(0,\cdot)=0, u^2(\cdot,0)=B(\cdot)
         \end{array} \right. .$$

Since $u=u^1+u^2$, the claim follows.}}

Therefore by previous results we have the following theorem.

\teo{\label{teo:eul_grd}$Y_{\alpha,n}$ the interpolated natural two-sided Brownian motion on the $(\alpha,n)$-rhombus lattice
converges in distribution to $u$ as $n\to\infty$, where $u$ is the solution 
to following stochastic heat equation 
\be \begin{array}{cl} \partial_x u=\left\{\begin{array}{cl}
				      \frac{1}{2}\partial_{tt} u+\W & \alpha=0,\\
				      0 & \alpha>0,\\
                                   \end{array}\right. 
&\quad {\rm on}\quad  (\R\setminus\{0\})^2,\\
    u(0,t)=B(t) & \quad {\rm for}\quad  t\in\R.\\
    u(x,0)=0 & \quad {\rm for}\quad  x\in \R.
    \end{array}\label{shefull}
\ee
and $t\mapsto B(t)$ is a two-sided Brownian motion independent of $(\W(A):A\in \cB(\R^2))$.
\dok{We will show the claim on $\R_+^2$, the other quadrants are shown in the same way.
From $Y_{\alpha,n}=\widetilde{Y}_{\alpha,n}+(Y_{\alpha,n}-\widetilde{Y}_{\alpha,n})$, Proposition \ref{prop:shelim}. and
Lemma \ref{lem:shed0}. we have $Y_{\alpha,n}\stackrel{d}{\to} u+0=u$. }
}

\index{Heat equation|)}

\part{Processes on general and random time-like graphs}The TLG's defined so far (see Chapter \ref{sec:01}) have only one beginning and one end (usually denoted by $0$ and $1$).

In applications and theory of classical graphical models an important role belongs
to processes indexed by trees. This includes one of the most widely used models - hidden Markov model.

\begin{figure}[ht]
\begin{center}
 
\includegraphics[width=10cm]{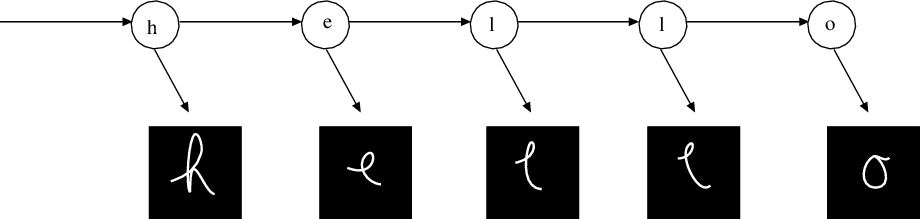}\\
Hidden Markov model in image (letter) analysis.

\end{center}
\end{figure} 

Hidden Markov model\index{Hidden Markov model} is also a collection of Markov processes combined together. 
(See for example \S6.2.3.1 in \cite{prgphmdl}.)\vspace{0.3cm}

Another model that is of wide interest is the branching Brownian motion. 
This is a similar model, but underlying graph is a random tree.\vspace{0.3cm}

Having this in mind, it is natural to ask can we have more than one beginning and more than 
one end. Could we define a process indexed by a TLG with a structure of a tree? 
We will show that this can be so in Chapter 9, and that there is a natural 
embedding into the existing family of TLG's. This embedding will help us 
define processes on a generalized family of TLG$^*$'s in Chapter \ref{chp:prnstlg}.

Later, in Chapter \ref{chp:gw}. we will be able to randomize the underlying graph, and see
how it is connected to the branching Markov processes.

\chapter{Non-simple TLG's}\label{sec:7}

\index{Time-like graph (TLG)|(}
\section{New definitions}

The TLG's defined in Chapter \ref{sec:01}, from now on, we will call {\bf simple TLG's}\index{Time-like graph (TLG)!simple|textbf}.

\defi{A graph $\cG=(\cV,\cE)$ will be called a  \textbf{time-like graph (TLG)}\index{Time-like graph (TLG)|textbf} if its sets 
of vertices $\cV$ and edges $\cE$ satisfy the following properties. 
\begin{enumerate}[(i)]
 \item Let $A,B>0$. The set $\cV$ contains at least two elements, $\cV=\{t_0,t_1,\ldots,t_N\}$,
where for $k=1,2,\ldots,N-1$, $$A\leq t_k\leq t_{k+1}\leq B.$$
 \item An edge between $t_j$ and $t_k$ will be denoted $E_{jk}$. We assume that there 
is no edge between $t_j$ and $t_k$ if $t_j=t_k$. $E_{jk}$ indicates that $t_j<t_k$. (We use 
$E_{jk}^1$, $E_{jk}^2$,\ldots if there is more than one edge connecting $t_k$ and $t_j$.)
\item We assume that  all vertices
have a finite non-zero degree. 
\end{enumerate}
We will call $\cG$ the \textbf{unit} TLG if $A=0$ and $B=1$.}\vspace{0.2cm}

\begin{figure}[ht]
\begin{center}
\psfrag{a}{$\boldsymbol{t_k}$}
\psfrag{b}{$\boldsymbol{t_j}$}

\includegraphics[width=7cm]{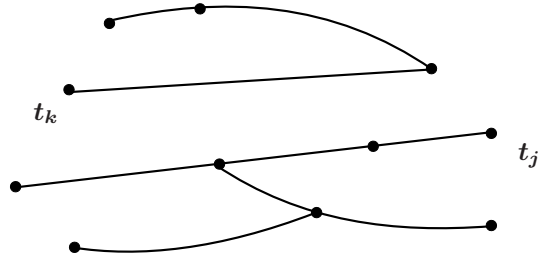}\\
\caption{TLG $\cG$ with entrance vertex $t_k$ and exit vertex $t_j$} \label{pic:sl25a}
  \end{center}
\end{figure}

\noindent \emph{Remarks.} \begin{enumerate}[(a)]
 \item Notice that in the new definition there are no longer unique vertices with times $A$ and $B$.
 \item We dropped part $(iv)$ of the original definition and added an assumption
in $(iii)$ that all vertices are of non-zero degree.
\item  Notice, that this definition no longer
guaranties that the graph is connected. (See Figure \ref{pic:sl25a}.)
\end{enumerate}

Again, as in Chapter \ref{sec:01}, we will restrict our attention to unit TLG's and 
prove all the claims for them.

\defi{\begin{enumerate}[(a)]
       \item A vertex $t_k$ that is not connected to any other vertex with time less than $t_k$
will be called an \textbf{entrance (vertex)}\index{Time-like graph (TLG)!entrance (vertex)}. We will denote the set of entrance vertices by $En(\cG)$.
      \item A vertex $t_j$ that is not connected to any other vertex with time greater than $t_k$
will be called an \textbf{exit (vertex)}\index{Time-like graph (TLG)!exit (vertex)}. We will denote the set of exit vertices by $Ex(\cG)$.
      \end{enumerate}
}\vspace{0.2cm}

The definitions of the \textbf{representation}, {\bf time-paths}, {\bf (simple) cells} remain the same as before (see 
Definitions \ref{def:rptlg}, \ref{def:tmpth}, \ref{defct}).\vspace{0.2cm}

We no longer have the full-time path as in the case of simple TLG's but instead we
define the full path\index{Time-like graph (TLG)!full path}.

\defi{\label{def:full-path}A time-path $\sigma$ is called a \textbf{full path} if it starts with an entrance vertex
and ends with an exit vertex. We denote the set of full paths by $P(\cG)$, while the full paths
starting at $t_k\in En(\cG)$ and ending at $t_j\in Ex(\cG)$ we will denote by
$P_{t_k\to t_j}(\cG)$.}\par
\noindent\emph{Remark.} Note that it can be $P_{t_k\to t_j}(\cG)=\emptyset$ (see Figure \ref{pic:sl25a}.) and further
$$P(\cG)=\bigcup_{t_k\in En(\cG)}\bigcup_{t_j\in Ex(\cG)}P_{t_k\to t_j}(\cG).$$

\section{Embedding TLG's into simple TLG's}\index{Time-like graph (TLG)!embedding|(} 
Although it seems that TLG's are much more general objects than
simple TLG's, there is a natural embedding that will enable us to 
use most of the results that we had for simple TLG's. As a result we will be 
able to construct processes under similar conditions as we did on
simple TLG's.

\subsection*{Minimal embedding}\index{Time-like graph (TLG)!embedding!minimal}
The first embedding  will use the minimal number of edges to 
embed the (unit) TLG into a simple TLG.\vspace{0.1cm}

Procedure is the following:\vspace{0.2cm}

Let $\cG=(\cV,\cE)$ be a TLG. 
\begin{itemize}
 \item Set $t_{-\infty}=-1$, $t_{\infty}=2$.
 \item For all $t_k\in En(\cG)$ we denote $E_{-\infty k}$ and edge between $t_{-\infty}$ and $t_k$, and
for all $t_j\in Ex(\cG)$ we denote $E_{j\infty }$ and edge between $t_j$ and $t_{\infty}$.
\item Set 
$$\cV^{\#}=\cV\cup \{t_{-\infty},t_{\infty}\},$$
and 
$$\cE^{\#}=\cE \cup\{E_{-\infty k}:t_k\in En(\cG) \} \cup \{E_{j\infty } :t_j\in Ex(\cG)\}.$$
\end{itemize}
\begin{figure}[ht]
\begin{center}
\psfrag{a}{$\boldsymbol{-1}$}
\psfrag{b}{$\boldsymbol{2}$}
\psfrag{0}{$\boldsymbol{0}$}
\psfrag{1}{$\boldsymbol{1}$}

\includegraphics[width=14cm]{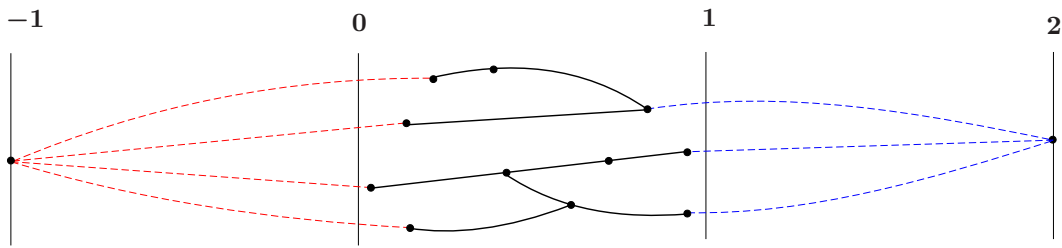}\\
\caption{Minimal embedding of the TLG $\cG$ (from Figure \ref{pic:sl25a}) into a simple TLG.} \label{pic:sl25b}
  \end{center}
\end{figure}

The transformation that defines the embedding has some nice properties.

\prop{\label{prop:tlg**map}The following claims hold:
\begin{enumerate}[(i)]
 \item $\cG'$ is a simple TLG.
\item $\cG\mapsto \cG'$ is an injective map. 
\item $\cG\mapsto \cG'$ preserves the connectedness by time-paths.
\item $\cG\mapsto \cG'$ preserves the order induced by $\G$ and $\G'$, i.e. 
$$t \stackrel{\G}{\prec} s \quad \Leftrightarrow \quad t'\stackrel{\G'}{\prec} s'$$. 
\end{enumerate}
\dok{(i) Follows form Definition \ref{def:tlg}. (ii) If we have $\cG'$, we can delete the edges connected to $t_{-\infty}$ and $t_{+\infty}$
and get $\cG$. (iii) From the definition of the mapping it is clear 
that if $t$ and $s$ are connected by time-path then $t'$ and $s'$ will also be connected. 
If $t$ and $s$ are not connected by a time-path
neither will $t'$ and $s'$ be connected by a time path, since all the new edges added 
include vertices $0$ and $1$ in $\G'$. (iv) This follows from (iii) and the 
the fact that time remains the same. }}

\subsection*{Maximal embedding} \index{Time-like graph (TLG)!embedding!maximal}
The embedding will add a number of edges to 
embed the TLG into a simple TLG.\vspace{0.1cm}

Procedure is the following:\vspace{0.2cm}

Let $\cG=(\cV,\cE)$ be a TLG. 
\begin{itemize}
 \item Set $t_{-\infty}=-1$, $t_{\infty}=2$.
 \item For all $t_k\in \cV$ we denote $E_{-\infty k}$ and edge between $t_{-\infty}$ and $t_k$, and
 $E_{k\infty }$ an edge between $t_k$ and $t_{\infty}$.
\item Set 
$$\cV^{\#}=\cV\cup \{t_{-\infty},t_{\infty}\},$$
and 
$$\cE^{\#}=\cE \cup\{E_{-\infty k}, E_{k\infty }:t_k\in \cV \} 
$$
\end{itemize}

\begin{figure}[ht]
\begin{center}
\psfrag{a}{$\boldsymbol{-1}$}
\psfrag{b}{$\boldsymbol{2}$}
\psfrag{0}{$\boldsymbol{0}$}
\psfrag{1}{$\boldsymbol{1}$}

\includegraphics[width=14cm]{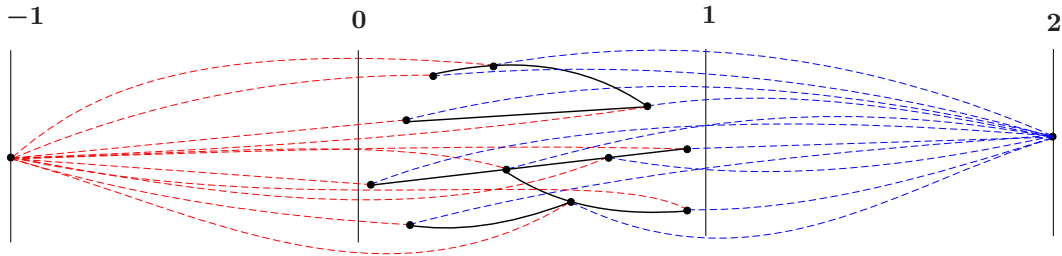}\\
\caption{Maximal embedding of the TLG $\cG$ (from Figure \ref{pic:sl25a}) into a simple TLG.
} \label{pic:sl25c}
  \end{center}
\end{figure}

The transformation that defines the embedding has some nice properties.

\prop{\label{prop:tlg***map}The following claims hold:
\begin{enumerate}[(i)]
 \item $\cG''$ is a simple TLG.
\item $\cG\mapsto \cG''$ is an injective map. 
\item $\cG\mapsto \cG''$ preserves the connectedness by time-paths.
\item $\cG\mapsto \cG''$ preserves the order induced by $\G$ and $\G'$, i.e. 
$$t \stackrel{\G}{\prec} s \quad \Leftrightarrow \quad t''\stackrel{\G''}{\prec} s''$$. 
\end{enumerate}

\subsection*{Remark on the embeddings} We will use both embeddings of a TLG $\cG$ for several reasons. 
It is easier to draw and see properties of $\cG'$ than of $\cG''$. On the other hand, for the many of the 
proofs that we have to do $\cG''$ will be much better to use.
\dok{The proof is similar to the proof of Proposition \ref{prop:tlg**map}.}}
\index{Time-like graph (TLG)!embedding|)}

\index{TLG$^{**}$ family|(}
\section{TLG$^{**}$ family}
As we have already seen in \S\ref{subsec:cnsprb} we might have problems to define 
a process with natural properties on some TLG's. In this section we introduce the family TLG$^{**}$,
similar to the family TLG$^{*}$ that we had defined for simple TLG's.\vspace{0.3cm}

We will describe the family of TLG graphs that is generated from a minimal graph 
by adding vertices, adding edges between vertices connected by a time-path
and adding edges between a new vertex and a vertex already on the graph.

\defi{\label{def:tlg**}The \textbf{TLG$^{**}$-family}\index{TLG$^{**}$ family|textbf} is given in the following inductive way.
\begin{enumerate}[(i)]
 \item The minimal graph $\cG=(\cV,\cE)$, with $\cV=\{t_0,t_N\}$ ($t_0<t_N$) and $\cE=\{E_{0N}\}$ is a 
TLG$^{**}$.

\begin{figure}[ht]
\begin{center}
\psfrag{0}{$\boldsymbol{t_0}$}
\psfrag{1}{$\boldsymbol{t_1}$}
\psfrag{E}{$\boldsymbol{E_{01}}$}
\includegraphics[height=0.5cm]{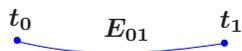}\\
\caption{A minimal graph}
\end{center}
\end{figure}

\item \label{def:tlg**:ii} Let $\cG_1=(\cV_1,\cE_1)$ be a TLG$^*$, where $\cV_1=\{t_0,t_2,\ldots,t_N\}$.
\begin{enumerate}[(1)]
 \item \label{def:tlg**:1} If $\tau_k\in [0,1]$ (not a vertex), and for some $E_{k_1k_2}\in \cE$ and
$t_{k_1}< \tau_k< t_{k_2}$ then set
$$\cV_2:=\cV_1\cup\{\tau_k\}\quad \textrm{and}\quad \cE_2:=\cE_1\cup \{E_{k_1k},E_{kk_2} \}\setminus \{E_{k_1k_2}\}.$$
$\cG_2:=(\cV_2,\cE_2)$ is also a TLG$^{**}$.

\begin{figure}[ht]
\begin{center}
\psfrag{a}{$\boldsymbol{t_{k_1}}$}
\psfrag{b}{$\boldsymbol{t_{k_2}}$}
\psfrag{c}{$\boldsymbol{\tau_{k}}$}
\psfrag{E}{$\boldsymbol{E_{k_1k_2}}$}
\psfrag{1}{$\boldsymbol{E_{k_1k}}$}
\psfrag{2}{$\boldsymbol{E_{kk_2}}$}
\includegraphics[height=2cm]{pr_13a.eps}\\
\caption{Adding a vertex}
\end{center}
\end{figure}
\item \label{def:tlg**:2} If $\tau_k\in [0,1]$ (not a vertex), and for some $\tau_k< t_{k_2}$ then set
$$\cV_2:=\cV_1\cup\{\tau_k\}\quad \textrm{and}\quad \cE_2:=\cE_1\cup \{E_{kk_2} \}.$$
$\cG_2:=(\cV_2,\cE_2)$ is also a TLG$^{**}$.

\begin{figure}[ht]
\begin{center}

\psfrag{a}{$\boldsymbol{t_{k_2}}$}
\psfrag{b}{$\boldsymbol{\tau_{k}}$}
\psfrag{E}{$\boldsymbol{E_{kk_2}}$}
\includegraphics[height=2cm]{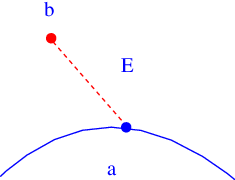}
\psfrag{a}{$\boldsymbol{t_{k_1}}$}
\psfrag{b}{$\boldsymbol{\tau_{k}}$}
\psfrag{E}{$\boldsymbol{E_{k_1k}}$}
\quad\quad\quad\includegraphics[height=2cm]{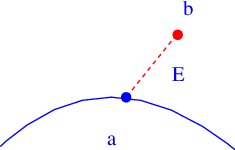}\\
\caption{Adding the edge and the vertex}
\end{center}
\end{figure}

\item \label{def:tlg**:3} If $\tau_k\in [0,1]$ (not a vertex), and for some $t_{k_1}< \tau_k$ then set
$$\cV_2:=\cV_1\cup\{\tau_k\}\quad \textrm{and}\quad \cE_2:=\cE_1\cup \{E_{k_1k} \}.$$
$\cG_2:=(\cV_2,\cE_2)$ is also a TLG$^{**}$.

\item Let $t_j,t_k\in \cV_1$ such that $t_j<t_k$, and assume that there exists a time-path $\sigma(j,\ldots,k)$ 
between these vertices. Then set
$$\cV_2:=\cV_1\quad \textrm{and}\quad \cE_2:=\cE_1\cup \{E_{jk}^* \}.$$
$\cG_2:=(\cV_2,\cE_2)$ is also a TLG$^{**}$. ($E^{*}_{jk}$ is an new edge (not in $\cE_1$).)
\begin{figure}[ht]
\begin{center}
\includegraphics[height=1.3cm]{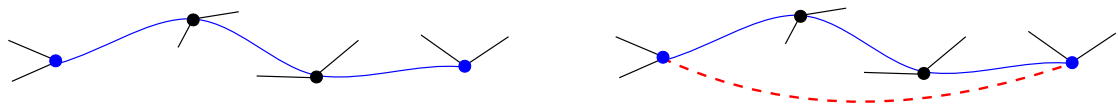}\\
\caption{Adding the edge \textcolor{red}{$E_{jk}^*$}}
\end{center}
\end{figure}
\end{enumerate}
\item If $\cG_1=(\cV_1,\cE_1)$ and $\cG_2=(\cV_2,\cE_2)$ are two disjoint TLG$^{**}$ their union is also a
a TLG$^{**}$.
\item \label{def:tlg**:iii} We will say that $(\cG_j)_{1\leq j\leq k}$ is a \textbf{tower of TLG$^{**}$'s} if 
for $j>1$, $\cG_j$ is constructed from $\cG_{j-1}$ as in (ii). 

\end{enumerate}

}
%
%
\noindent \emph{Remarks.} (1) Applying only the procedure in (ii)  will clearly give us a 
a connected TLG$^{**}$. (2) Any connected component of a TLG$^{**}$ can be obtained
only by using step (ii). (3) It can be easily seen that a TLG$^{**}$ is a TLG. (4) It is
clear that a TLG$^*$ is a TLG$^{**}$.\vspace{0.3cm}

\lem{\label{lem:**|*}Let $\cG$ be a TLG. If $\cG''$ is a TLG$^*$ then $\cG$
is a TLG$^{**}$.
\dok{If $\cG''$ is a TLG$^*$, then there exists a TLG$^*$-tower $(\cH_j)_{j=0}^n$
such that $\cH_n=\cG''$ and $\cH_0$ contains an edge in $\cG$. Now, we construct a 
TLG$^{**}$-tower $(\cG_j)_{j=0}^{m_n}$ from the tower $(\cH)_{j=0}^n$. 
Let $\cG_0$ be the minimal graph that is contained in $R(\cH_0)$. 
If $\cH_{j}$ is obtained by $\cH_{j-1}$ by
\begin{itemize}
 \item adding a new vertex, then let $\cG_{j}$ be obtained from 
$\cG_{j-1}$ by adding a new vertex (step (1));
 \item adding a new edge contained in $R(\cG)$, then let $\cG_{j}$ be obtained from 
$\cG_{j-1}$ by adding a new edge (step (4));
 \item adding a new edge partially contained in $R(\cG)$, then let $\cG_{j}$ be obtained from 
$\cG_{j-1}$ by adding a new edge with a new vertex (steps (2) or (3));
 \item adding a new edge not contained in $R(\cG)$, then let $\cG_{j}=\cG_{j-1}$.
\end{itemize}
Let's assume that $(\cG_j)$ doesn't have repeating TLG's. In order to show that it is 
a TLG$^{**}$-tower, we need to check that each time we add an edge (step (4)) the two endpoints are connected.
This is clear from the fact that one endpoint of all other edges in $\cH_j$ not in $\cG_j$
is in the set $\{t_{-\infty},t_{+\infty}\}$. So the two points on $\cG_j$ are 
connected by a time path in $\cH_j$ only if they are connected by a time path in 
$\cG_j$.}}

\teo{\label{thm:tlg**|*}Let $\cG$ be a TLG, then $\cG$ is a TLG$^{**}$ if and only if its embedding $\cG''$ is a TLG$^*$.
\dok{Let $n$ be the number of vertices and edges of $\cG$. For $n=3$ we have a minimal graph
and the claim is clear. Let's assume that the claim holds for $n\geq 3$, and show that the claim is true for $n+1$.\vspace{0.2cm}

($\Rightarrow$): Let $\cG^\#$ be a TLG$^{**}$ such that we can construct 
$\cG$ using steps $(1)-(4)$ from Definition \ref{def:tlg**}. Then $\cG''$ 
can be constructed from $(\cG^\#)''$ in several steps from Definition \ref{def:tlg*}. 

($\Leftarrow$): See Lemma \ref{lem:**|*}.}}

From the previous proof we get the following fact.
\pos{\label{pos:tlg*t**}If $(\cH_j)$ is a TLG$^{**}$-tower, then $(\cH_j'')$ is a subsequence 
of a TLG$^*$-tower.}

\teo{\label{teo:T*/**}The following statements are equivalent:
\begin{enumerate}[(a)]
 \item $\cG$ is a TLG$^{**}$.
 \item $\cG'$ is a TLG$^{*}$.
 \item $\cG''$ is a TLG$^{*}$.
\end{enumerate}
\dok{(a) $\Leftrightarrow$ (c): Follows from Theorem \ref{thm:tlg**|*}. 

(b) $\Rightarrow$ (c): Every $t_k'$ in $\cG'$ is connected to $t_{-\infty}$ and $t_{\infty}$ by a time-path. 
Therefore, 
we can add an edge to $\cG'$ between $t_{-\infty}$ and $t_k'$, and an edge between $t_k'$ and $t_{\infty}$. 
Hence, we can construct $\cG''$ from $\cG'$ by adding edges. Hence, $\cG''$ 
is a TLG$^*$.

(c) $\Rightarrow$ (b): Let $t_k$ be a vertex that is not an entrance, then the 
if we remove the edge $E_{-\infty,k}$ from $\cG''$ we get a TLG$^*$ by Corollary \ref{pos:sped}.
The same holds if $t_k$ is not an exit for the edge $E_{k,\infty}$. Doing this until all such edges 
are removed gives us $\cG'$, that will, by repeated use of Corollary \ref{pos:sped}, be a TLG$^*$.  
}}
\pos{\label{pos:g'g''}If $\cG''$ is a TLG$^*$, then there is a TLG$^*$-tower $(\cH_j)_{j=1}^n$ such that 
$\cH_1=\cG'$ and $\cH_n=\cG''$.}

The order '$\preceq$' between the points is defined in the same way as in Chapter \ref{sec:01}. See
Definition \ref{def:ptorder}.\index{Time-like graph (TLG)!order}\index{Time-like graph (TLG)!point}

\lem{\label{lem:tlg**min}For points $t_1$ and $t_2$ on a TLG$^{**}$ $\G$  
\begin{itemize}
 \item  there exists  a point $t_1\wedge t_2$ on $\G\cup \{-1\}$ such that
$$\{t\in \G:t\preceq t_1 \}\cap \{t\in \G:t\preceq t_2\} = \{t\in \G:t\preceq t_1\wedge t_2\};$$
\item there exists a point $t_1\vee t_2$ on $\G\cup \{2\}$ such that
$$\{t\in \G:t\succeq t_1 \}\cap \{t\in \G:t\succeq t_2\} = \{t\in \G:t\succeq t_1\vee t_2\};$$
\end{itemize}
in the sense that if we have an empty set on one side we define $t_1\wedge t_2=-1$ in the first case, and $t_1\vee t_2=2$
in the second case.
\dok{By Proposition \ref{prop:tlg**map}, 
we will have $t_1'\wedge t_2'=(t_1\wedge t_2)'$ and $(t_1\vee t_2)'=t_1'\vee t_2'$. Since $\G$ is a TLG$^*$,
$t_1'\wedge t_2'$ and $t_1'\vee t_2'$ exists and can obtain one of the values in $\cV'\cup\{0,1\}$. 
Since, the transformation is injective so are $t_1\wedge t_2$ and $t_1\vee t_2$.}}

We know from Theorem \ref{thm:tlg*1}, that all planar simple TLG's are TLG$^*$. Unfortunately,
the same is not true for TLG$^{**}$'s.

\prop{The following statements hold:
\begin{enumerate}[(a)]
 \item If $\cG$ is a planar TLG\index{Time-like graph (TLG)!planar} its embedding $\cG'$ doesn't have to be a planar TLG.
\item If $\cG$ is a planar TLG$^{**}$ its embedding $\cG'$ doesn't have to be a planar TLG$^*$.
\item There exists a planar TLG that is not a TLG$^{**}$.
\end{enumerate}
\dok{ (a) See Figure \ref{pic:sl26}.
\begin{figure}[ht]
\begin{center}
\includegraphics[width=15cm]{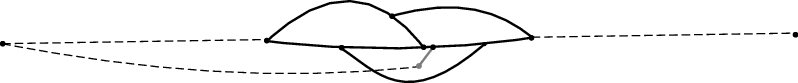}\\
\caption{TLG $\cG$ is planar (full lines), but its embedding $\cG'$ is not.} \label{pic:sl26}
  \end{center}
\end{figure}

(b) The black part of the graph $\cG$ (in Figure \ref{pic:sl26}) is a (planar) simple TLG, so it
is a TLG$^*$. Hence, we can first construct the black part, and then add 
the gray vertex and the gray edge connecting it to the rest of the graph. So, $\cG$
is a TLG$^*$.

(c) See the graph in Figure \ref{pic:sl27}. This is not a TLG$^{**}$, because 
$t_3\wedge t_4$ is not defined, and by 
Lemma \ref{lem:tlg**min} this should be defined in the case of a TLG$^{**}$.

\begin{figure}[ht]
\begin{center}
\psfrag{a}{$\boldsymbol{t_1}$}
\psfrag{e}{$\boldsymbol{t_2}$}
\psfrag{b}{$\boldsymbol{t_3}$}
\psfrag{c}{$\boldsymbol{t_4}$}
\psfrag{d}{$\boldsymbol{t_5}$}
\includegraphics[width=8cm]{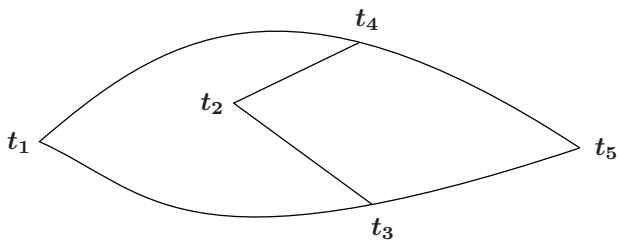}\\
\caption{A planar TLG that is not a TLG$^{**}$.} \label{pic:sl27}
  \end{center}
\end{figure}

}}

There are two important cases of planar TLG's that are planar TLG$^{**}$'s,
and we will encounter with them in the future.

\prop{\label{tlg**_ex}\begin{enumerate}[(a)]
       \item A planar TLG $\cG$ such that all vertices in $En(\cG)$ and all vertices in $Ex(\cG)$
have the same time component is a TLG$^{**}$.
	\item TLG that has the structure of a tree\index{Time-like graph (TLG)!tree} is also a TLG$^{**}$
      \end{enumerate}
\dok{(a) The proof follows from the fact that $\cG'$ the embedding of $\cG$
is a planar TLG, therefore a TLG$^*$, and by Theorem \ref{thm:tlg**|*} $\cG$ is a TLG$^{**}$.
(b) Follows by induction on the number of edges.}}
\index{TLG$^{**}$ family|)}
\index{Time-like graph (TLG)|)}

\chapter{Processes on non-simple TLG's}\label{chp:prnstlg}

\section{Processes on TLG$^{**}$}

Idea of the construction is the similar to the one that we had in the case
of simple TLG's (as described in Section \ref{cdst_pth}):
\begin{itemize}
 \item We take a family $\cM$ of measures $\mu_{\sigma}$ on full paths $P(\cG)$
with certain properties.
 \item Using these properties we create a (3T)-family $\cM'$ of measures 
$\mu_{\sigma'}$ on full-time paths of the embedding $\cG'$.
\item We create a natural $\cM'$-process on $\cG'$, and from that process
we create the process on $X$ on $\cG$.
\end{itemize}
We could do the same approach for $\cG''$, and we will briefly discuss it.

We need the version of the consistent family of measures along full paths.

\defi{Let $\cG$ be a TLG, for a family of distributions
$$\cM=\{\mu_{\sigma}:\sigma\in H\subset P(\cG)\}$$
where if $\sigma$ is a full path from $t_k$ to $t_j$
then $\mu$ is a distribution of a stochastic process on $[t_k,t_j]$,  we say
that it is \textbf{consistent} if for $\sigma_1,\sigma_2 \in H$
$$\mu_{\sigma_1}\circ\pi_T^{-1}=\mu_{\sigma_2}\circ\pi_T^{-1},$$
where $T=\{t:t\in E,E\in\sigma_1\ \&\ E\in\sigma_2\}$. 
}\index{Process indexed by a TLG!consistent family of distributions along time-paths}\vspace{0.2cm}

We also need a notion of the half-cell\index{Time-like graph (TLG)!half-cell} that didn't exist for simple TLG's.

\defi{Let $\cG=(\cE, \cV)$ be a TLG.

\begin{enumerate}[(a)]
       \item We say that time paths $\sigma_1$ and $\sigma_2$ in $\cG$ 
 starting at $t_{k_1}\in En(\cG)$ and respectively at $t_{k_2}\in En(\cG)$ and both ending at $t_m$
which is their only common vertex,
form a \textbf{right half-cell} $(\sigma_1,\sigma_2)$. 

\begin{figure}[ht]
\begin{center}
\psfrag{A}{$\boldsymbol{t_{k_1}}$}
\psfrag{B}{$\boldsymbol{t_{k_2}}$}
\psfrag{C}{$\boldsymbol{t_m}$}
\psfrag{D}{$\boldsymbol{t_{k}}$}
\psfrag{F}{$\boldsymbol{t_{m_1}}$}
\psfrag{G}{$\boldsymbol{t_{m_2}}$}
\includegraphics[width=10cm]{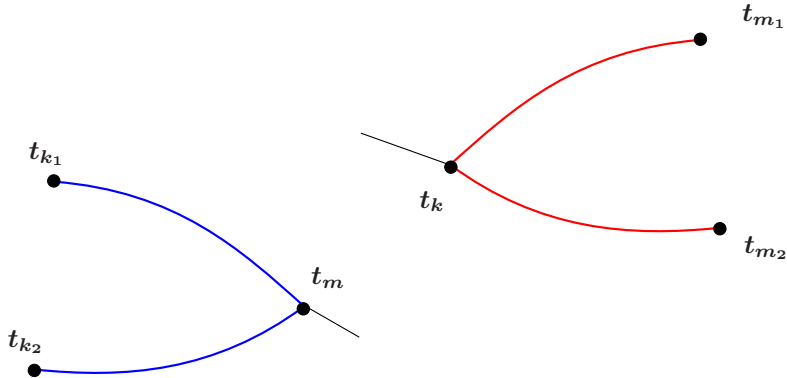}\\
\caption{Right and left half-cells.} \label{pic:sl36}
  \end{center}
\end{figure}

      \item We say that time paths $\sigma_1$ and $\sigma_2$ in $\cG$ 
both starting at $t_{k}$ which is their only common vertex, and ending at $t_{m_1}\in Ex(\cG)$ and respectively $t_{m_2}\in Ex(\cG)$,
form a \textbf{left half-cell} $(\sigma_1,\sigma_2)$. 
\item A half-cell $(\sigma_1,\sigma_2)$ is called \textbf{simple} if there is no time-path
connecting vertex on $\sigma_1$ and a vertex on $\sigma_2$ (both must be different from the connecting vertex $t_m$). 
      \end{enumerate}}
\noindent\emph{Remarks.} Note that a half-cell in $\cG$ will be embedded into a cell in
the embedding $\cG'$.\vspace{0.3cm}
\subsection{Conditions}\label{subsec:cond}\index{Process indexed by a TLG!construction|(}
We will show that an $\cM$-process\index{Process indexed by a TLG!$\cM$-process} exists if the following conditions are satisfied:
\begin{description}
  \item[T1'] $\cG=(\cV,\cE)$ is a TLG$^{**}$.
 \item[T2'] $\cM$ is a consistent family of measures that induce continuous or RCLL processes.
 \item[T3'] Let $\pi$ be a full-path in $P(\cG)$ and $t_k\in \cV$ a vertex on that path. Then 
$(X_{\pi}(t):t\leq t_k)$ and $(X_{\pi}(t):t\geq t_k)$ are independent given $X(t_k)$.
%
\end{description}

\defi{The family $\cM=\cM(\cG)=\{\mu_{\sigma}:\sigma\in P(\cG)\}$ satisfying properties
(T1'), (T2') and (T3') is called the \textbf{(3T') family}.\index{Process indexed by a TLG!construction!(3T') conditions}}

\prop{\label{prop:3t'3t}If $\cG$ is a TLG$^*$ and $\cM$ a (3T') family, then $\cM$ is a (3T) 
family on $\cG$.
\dok{In this case we only need to check the (T3) property of $\cM$. Let $\pi$
be a path that contains $t_{*}$ and $t^*$ endpoints of a simple cell. Let $A_{*}\in \sigma(X_{\pi}(t):t\leq t_*)$
$A_*^*\in \sigma(X_{\pi}(t):t_*\leq t\leq t^*)$ and $A^*\in \sigma(X_{\pi}(t):t^*\leq t)$, while $B_*\in \sigma(X(t_*))$
and $B^*\in\sigma(X(t^*))$. Now we have 
\begin{align*}
 &\E(\P(A_*\cap A_*^*\cap A^*|X(t_1),X(t_2))\1_{B_*}\1_{B^*})\\
&=\E(\E(\1_{A_*}\1_{ A_*^*}\1_{ A^*}|X(t_1),X(t_2))\1_{B_*}\1_{B^*})=\E(\1_{A_*}\1_{ A_*^*}\1_{ A^*}\1_{B_*}\1_{B^*})\\
&=\E(\E(\1_{A_*}|X(t_*))\1_{ A_*^*}\1_{ A^*}\1_{B_*}\1_{B^*})=\E(\E(\1_{A_*}|X(t_*))\1_{ A_*^*}\E(\1_{ A^*}|X(t^*))\1_{B_*}\1_{B^*})\\
&=\E(\E(\1_{A_*}|X(t_*))\E(\1_{ A_*^*}|X(t_*),X(t^*))\E(\1_{ A^*}|X(t^*))\1_{B_*}\1_{B^*})\\
&=\E(\1_{A_*}\P( A_*^*|X(t_*),X(t^*))\1_{ A^*}\1_{B_*}\1_{B^*})\\
&=\E(\P( A_*^*|X(t_*),X(t^*))\E(\1_{A_*}\1_{ A^*}|X(t_*),X(t^*))\1_{B_*}\1_{B^*})\\
&=\E(\P( A_*^*|X(t_*),X(t^*))\P(A_*\cap A^*|X(t_*),X(t^*))\1_{B_*}\1_{B^*}).
\end{align*}
The claim now follows from the Monotone Class Theorem.
}}
\noindent\emph{Remark.} The converse of of the statement of the previous proposition is not true. Take for example a 
non-Markovian process on  the graph $\cG=(\{t_0=1,t_1=1/2,t_2=1\}, \{E_{01},E_{12}\})$,
such that $X(0)$ and $X(1)$ are not independent given $X(1/2)$.
\subsection{Construction}


Let $\cM$ be a (3T') family on a TLG$^{**}$ $\cG$.\vspace{0.2cm}

Let $\cG''$ be the embedding of $\cG$ into simple TLG's. Now for each
time-path $\sigma$ in $\cG$  there exists a full-time path $\sigma'$
in $\cG''$ that corresponds to $\sigma$.\vspace{0.2cm}

\begin{figure}[ht]
\begin{center}
\psfrag{0}{$\boldsymbol{-1}$}
\psfrag{1}{$\boldsymbol{2}$}
\psfrag{k}{$\boldsymbol{t_k}$}
\psfrag{j}{$\boldsymbol{t_j}$}
\psfrag{Y}{\textcolor{blue}{$\boldsymbol{Y_{\sigma}(t)}$}}
\includegraphics[width=10cm]{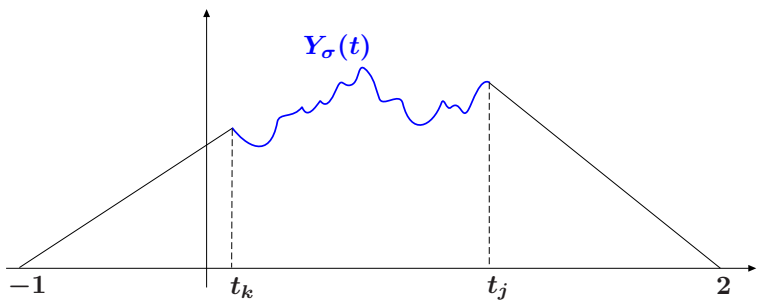}\\
\caption{Construction of $Y_{\sigma'}$.} \label{pic:sl28}
  \end{center}
\end{figure}

If $\sigma$ starts at $t_k$ and ends at $t_j$, then we can define 
a process $(Y_{\sigma}(t):t\in [t_k,t_j])$ whose distribution is $\mu_{\sigma}$. We will define 
$Y_{\sigma'}$ by interpolating $Y_{\sigma}$ on the whole interval $[0,1]$ (see Figure \ref{pic:sl28}
for illustration):
\be Y_{\sigma'}(t)=\left\{\begin{array}{cl}
               \dfrac{1+t}{1+t_k}Y_{\sigma}(t_k) & \textrm{if}\ t\leq t_k\\
		Y_{\sigma}(t) &  \textrm{if}\ t\in [t_k,t_j]\\
               \dfrac{2-t}{2-t_j}Y_{\sigma}(t_j) & \textrm{if}\ t\geq t_j\\
              \end{array}
\right. \ee
Note that if $Y_{\sigma}$ is continuous or RCLL so is $Y_{\sigma'}$.\vspace{0.2cm}

Now, we define $\mu_{\sigma'}$ to be the distribution of $(Y_{\sigma'}(t):t\in [0,1])$, and set
$$\cM'=\cM'(\cG')=\{\mu_{\sigma'}:\sigma'\in P_{0\to 1}(\cG')\},$$
or 
$$\cM''=\cM'(\cG'')=\{\mu_{\sigma'}:\sigma'\in P_{0\to 1}(\cG')\},$$
\teo{$\cM'$ is a (3T')-family on $\cG'$ and $\cM''$ is a (3T')-family on $\cG''$.
\dok{The proof is similar $\cM'$ and $\cM''$. Properties (T1') and (T2') are 
clearly satisfied. While the property (T3') follows from the construction and 
the (T3') property of $\cM$.
%
%
}}
\pos{$\cM'$ is a (3T) family on $\cG'$ and $\cM''$ is a (3T) family on $\cG''$.
\dok{By Theorem \ref{teo:T*/**} both $\cG'$ and $\cG''$ are TLG$^*$'s. The claim follows from Proposition \ref{prop:3t'3t}. }}

We can construct the process on a TLG$^{**}$ similar to the construction of processes 
on TLG$^*$'s (see Definition \ref{def:constr}).

\defi{\label{defi:cnstr_ptg}Let $(\cG_l)_{0\leq l\leq n}$ be a TLG$^{**}$-tower where $\cG_0$
is a minimal graph $\cV_0=\{t_0,t_N\}$, $\cE_0=\{E_{0N}\}$ and $\cG_n=\cG$. Further let 
$\cM$ be a family of distributions satisfying (3T') conditions.
\begin{itemize}
 \item On $\cG_0$ we define a process $X^0$ with $\mu_{E_{0n}}$ distribution.
 \item If we have already defined $X^l$ on $\cG_l$ (for some $l<n$), then
we define $X^{l+1}$ on $\cG_{l+1}$ in the following way depending how we
constructed $\cG_{l+1}$ from $\cG_l$ (recall part (\ref{def:tlg**:ii}) of Definition \ref{def:tlg**}.).
\begin{enumerate}[(1)]
 \item In the construction a new vertex $\tau_l\in [0,1]\setminus \cV_{l}$ was added to graph $\cG_l$, by subdividing
some $E_{jk}$ such that $t_{j}<\tau_l<t_{k}$, into $E_{jl}$ and $E_{lk}$ to get $G_{l+1}$.  In this case, the two graphs $\cG_l$ and $\cG_{l+1}$ have a common representation,
$R(\cG_l)=R(\cG_{l+1})$, and we can define $X^{l+1}$ on $\cG_{l+1}$ to have the same values on this
representation as $X^{l}$.

\item In the construction a new vertex $\tau_k$ and a new edge between the vertex $t_j<\tau_k$ in $\cV_{l}$ and $\tau_k$, was added to get $\cG_{l+1}$. 
So, $\cG_{l+1}$ has a new edge $E^*=E^*_{jk}$ and a new vertex $\tau_k$. 
Let $Z_j=X^l(t_{j})$. 

Now we pick a full-time path $\sigma$ that contains $E^*$. Now we define 
$\mu_*(\cdot |x)$ to be the conditional probability of the process with
the distribution $\mu_{\sigma}\circ \pi_{[t_j,\tau_k]}^{-1}$ conditioned to
have the value $x$ at $t_j$. So we construct the process
$X^{l+1}$ in such a way that $X^{l+1}$ on $R(\cG_l)$ is equal to $X^l$
and $X^{l+1}_{E^*}$ is the process with distribution $\mu_*(\cdot |Z_j)$
and otherwise independent of $X^l$ given $Z_j$.

\item In the construction a new vertex $\tau_k$ and a new edge between the vertex $t_m>\tau_k$ in $\cV_{l}$ and $\tau_k$, was added to get $\cG_{l+1}$. 
So, $\cG_{l+1}$ has a new edge $E^*=E^*_{jk}$ and a new vertex $\tau_m$. 
Let $Z_m=X^l(t_{m})$. 

Now we pick a full-time path $\sigma$ that contains $E^*$. Now we define 
$\mu_*(\cdot |y)$ to be the conditional probability of the process with
the distribution $\mu_{\sigma}\circ \pi_{[\tau_k,t_m]}^{-1}$ conditioned to
have the value $y$ at $t_m$. So we construct the process
$X^{l+1}$ in such a way that $X^{l+1}$ on $R(\cG_l)$ is equal to $X^l$
and $X^{l+1}_{E^*}$ is the process with distribution $\mu_*(\cdot |Z_m)$
and otherwise independent of $X^l$ given $Z_m$.

\item In the construction a new edge between two vertices $t_j<t_k$ in $\cV_{l}$ that are connected by a time path
in $\cG_l$ , was added to get $\cG_{l+1}$. So, $\cG_{l+1}$ has a new edge $E^*=E^*_{jk}$. 
Let $Z_j=X^l(t_{j})$ and $Z_k=X^l(t_{k})$. 

Now we pick a full-time path $\sigma$ that contains $E^*$. Now we define 
$\mu_*(\cdot |x,y)$ to be the conditional probability of the process with
the distribution $\mu_{\sigma}\circ \pi_{[t_j,t_k]}^{-1}$ conditioned to
have values $x$ at $t_j$ and $y$ at $t_k$. So we construct the process
$X^{l+1}$ in such a way that $X^{l+1}$ on $R(\cG_l)$ is equal to $X^l$
and $X^{l+1}_{E^*}$ is the process with distribution $\mu_*(\cdot |Z_j,Z_k)$
and otherwise independent of $X^l$ given $Z_j$ and $Z_k$.

\end{enumerate}

\end{itemize}
Since $n$ is finite this procedure will end and we will have a process $X=X^n$ defined on $\cG$.}

We define:
\begin{itemize}
\item $X$ on $\cG$ to be defined as in Definition \ref{defi:cnstr_ptg};
 \item $X'$ to be a natural $\cM'$-process on $\cG'$ (in the sense of Definition \ref{def:nmprc});
 \item $X''$ to be a natural $\cM''$-process on $\cG''$ (in the sense of Definition \ref{def:nmprc}).
\end{itemize}

\teo{\label{teo:3ver:cstr}The following processes have the same distribution on $\cG$:
\begin{enumerate}[(a)]
 \item $(X(t):t\in \cG)$;
 \item $(X'(t):t\in \cG)$;
 \item $(X''(t):t\in \cG)$.
\end{enumerate}

\dok{To show that $X''$ and $X'$ have the same distribution on $\cG$ we will show that
they have the same distribution on $\cG'$.It is known by Corollary \ref{pos:g'g''}. that there exists 
a TLG$^*$-tower that starts with $\cG'$ and ends with $\cG''$. Now, by Lemma \ref{lem:inr:A13}, $\cM''(\cG')$ is a (3T)-family, and by definition
it coincides with $\cM'$. Since, the construction of a process on TLG$^*$ doesn't depend on the order (Theorem \ref{teo:uniq}), we can first
construct $X''$ on $\cG'$ and later on the rest of $\cG''$. Hence, $X''$ 
on $\cG'$ will be a natural $\cM'$-process, so $(X'(t):t\in \cG')$ and $(X''(t):t\in \cG')$ 
have the same distribution. Therefore, the restriction of these two processes to $\cG$
is also the same.\vspace{0.2cm} 

It remains to show that $X$ and $X_{\cG}''$ have the same distribution.
Let $(\cH)_{j=0}^n$ be a TLG$^{**}$-tower, we will show that the construction 
of $X$ on $\cG$ can embedded int the construction of $X''$ on $\cG''$. For $j=0$ 
it is clear that $X_{\cH_0}$ and $X_{\cH_0''}''$ have the same distribution on $\cH_0$. 
Let's assume for $j\geq 0$ $X_{\cH_j}$ and $X_{\cH_j''}''$ have the same 
distribution on $\cH_j$, and prove it for $j+1$. We have the following cases to consider:
\begin{itemize}
 \item A new vertex has been added to $\cH_j$ to obtain $\cH_{j+1}$. In this case the claim follows clearly.
 \item A new vertex $\tau$ and an edge $E_*$ between that vertex and existing one has been added. In this case the claim follows
from the fact that in $\cH_j''$ we are adding:
\begin{itemize}
 \item a new edge $E_*''$ between $-1$ or $2$ and an vertex $t_l$ on $\cH_j$;
 \item a vertex $\tau$ on that edge;
 \item we are adding an edge between $\tau$ and between the vertex $-1$ or $2$ to which it is not connected.
\end{itemize}
 Since, the  $X''(-1)=X''(2)=0$ the distribution of the process on $E_*''$ whose 
representation is in $R(\cH_{j+1})$ is given and depends only on the value $X(t_l)$. Hence, the claim follows.  
 \item A new edge has been added to $\cH_j$ to obtain $\cH_{j+1}$. In this case the claim follows from
the fact that the distribution of the process on the new edge is given 
and depends only on the values of the process on $\cH_j$ (for both $X$ and $X''$).
\end{itemize}

}}
 
\teo{For a (3T')-family $\cM$ the constructed process $X$ on a TLG$^{**}$ $\cG$
will always have the same distribution.
\dok{By Theorem \ref{teo:3ver:cstr}. we can embed the constructed process into a natural
$\cM'$ process on $\cG'$ and this process has a unique distribution.} }
\index{Process indexed by a TLG!construction|)}
%

\section{Properties of constructed processes}
We know, from Chapters \ref{sec:02} and \ref{sec:3a}, that for the process $X'$ many interesting properties
hold. Many of these properties have their natural analogous for the process $X$.

We will show that for $X$ the following properties hold:
\begin{itemize}
 \item $X$ is an $\cM$-process;
 \item $X$ is a spine-Markovian process;
 \item $X$ is a hereditary spine-Markovian process;
 \item $X$ is a cell-Markovian process.
\end{itemize}
Additionally if $\cM$ is a Markov family of measures we have
\begin{itemize}
 \item $X$ is moralized graph-Markovian;
 \item $X$ is time-Markovian;
 \item $X$ is edge-Markovian.
\end{itemize}
All these properties are (slightly generalized) versions of the
properties we had defined for simple TLG's.

\subsection{$X$ is an $\cM$-process}\index{Process indexed by a TLG!$\cM$-process}
It is
easy to see $X_{\sigma}\sim \mu_{\sigma}$, that is, $X$ is a $\cM$-process.
(This is the same as defined in Subsection \ref{condit}, on page \pageref{def:Mprc}.)
\subsection{$X$ is a spine-Markovian process}\index{Process indexed by a TLG!spine-Markovian property}
We will first define the spine-Markovian property.
\defi{Let $\sigma$ be any
full path  in the TLG $\cG=(\cV,\cE)$.
Let $\cG_-$ be a subgraph (not necessarily a TLG) of $\cG$ whose representation is a connected 
component of $R(\cG)\setminus R(\sigma)$. Let $W$ be the set
of vertices - \textbf{roots} connecting $\cG_-$ to $\sigma$ and let
$\cG_+$ denote the graph represented by $R(\cG)\setminus R(\cG_-)$.

We say that the process $X$ on a TLG $\cG$ is  \textbf{spine-Markovian} if for each such $\sigma$ and $\cG_-$ the processes $(X(t):t\in \cG_-)$ and $(X(t):t\in \cG_+)$ given $(X(t):t\in W)$ are independent.}

\prop{The constructed process $X$ is a \textbf{spine-Markovian} process on $\cG$.
\dok{Let $\sigma$ be the full path, and $\sigma'$ the corresponding full-time path in the embedding
$\cG'$. If $\cG_{-}$ is as in the definition, this is a connected graph and is a connected
component of $R(\cG')\setminus R(\sigma')$. We set $\cG_{+}$ and $\cG_{+}'$ to be graphs that have the 
representation, respectively $R(\cG)\setminus R(\cG_-)$ and $R(\cG')\setminus R(\cG_-)$.
The roots $W'$ of $\cG'$ include all the roots $W$ of $\cG$ and maybe $-1$ and $2$.
Since, $X'(-1)=X'(2)=0$, we have $\sigma(X_{W})=\sigma(X'_{W'})=\sigma(X'_{W'\setminus\{-1,2\}})$. Therefore,
since $X'$ is spine Markovian, and $\sigma(X(t):t\in \cG_{-})\subset \sigma(X'(t):t\in \cG_{-})$,
$\sigma(X(t):t\in \cG_{+})\subset \sigma(X'(t):t\in \cG'_{+})$
the spine-Markovian property for $X$ follows.  }}

\subsection{Hereditary spine-Markovian property}\index{Process indexed by a TLG!spine-Markovian property!hereditary}
Recall, Definition \ref{def:S*}. of $S^*(\cG)$. 
\defi{\label{def:S**}For a TLG$^{**}$ $\cG$ we define $S^{**}(\cG)$ to be the set of all TLG$^{**}$'s 
$\cH$ such that there exists a TLG$^{**}$-tower $(\cK_k)_{k=0}^n$ that starts with 
$\cK_0=\cH$ and ends with $\cK_n=\cG$.}

\defi{\label{def:her_mp**}The process $(X(t):t\in \cG)$ has a \textbf{hereditary spine-Markovian property} 
if $(X(t):t\in \cH)$ is a spine-Markovian process
for each $\cH\in S^{**}(\cG)$. }

\pos{Let $\cG$ be a TLG$^{**}$ and $\cG''$ its embedding, if TLG$^{**}$
$\cH$ is in $S^{**}(\cG)$ then $\cH''$ is in $S^{*}(\cG'')$.
\dok{Let $(\cK_l)_{l=1}^m$ TLG$^{**}$-tower such that $\cK_1=\cH$, and $\cK_m=\cG$. 
By Corollary \ref{pos:tlg*t**}. $\cK_1''$, $\cK_2''$, \ldots $\cK_m''$ are one after 
another in a TLG$^*$-tower, and the claim follows. }}

\prop{The constructed process $X$ is hereditary spine-Markovian.
\dok{Let $\cH\in S^{**}(\cG)$. Then $\cH''$ is in $S^{*}(\cG'')$.

By Theorem \ref{prop:hrdspm}., $X'$ on $\cG''$ is hereditary spine-Markovian, $X'$ is spine Markovian on $\cH''$,
and therefore, $X$ is spine Markovian on $\cH$. }}

\subsection{Cell-Markovian property}\index{Process indexed by a TLG!cell-Markovian property}

A cell will remain truly simple, as in Definition \ref{defi:msct}. We need to extend our definition
to half-cells. 

\defi{\begin{enumerate}[(a)]
       \item A right half-cell $(\sigma_1,\sigma_2)$ ending at $t_m$ is called \textbf{truly simple}\index{Time-like graph (TLG)!half-cell!truly simple} if there is no path 
	$\{t\in \cG: t\prec t_m \}$ that starts on on one side of the cell and ends on the other.
	       \item A left half-cell $(\sigma_1,\sigma_2)$ starting at $t_k$ is called \textbf{truly simple} if there is no path 
	$\{t\in \cG: t_k\prec t \}$ that starts on on one side of the cell and ends on the other.

      \end{enumerate}
 }

\lem{\label{lem:tshcl}A truly simple half-cell in $\cG$ is a part of a truly simple cell in $\cG'$.
\dok{We will prove the claim for the right half-cell, the proof for the left half-cell 
is similar. Let $\sigma_j'$ be the path consisting including $t_{-\infty}$ 
and $\sigma_j$, for $j=1,2$. Now, sigma $(\sigma_1,\sigma_2)$ is a cell. If there exists 
a path in $\cG[t_{-\infty},t_m]\setminus \{t_{-\infty},t_m\}$ connecting vertices on $\sigma_1$
and $\sigma_2$, then these vertices are in $\cG$. Further, since the path can't go
through $t_{-\infty}$, the path it self is in $\cG$. Hence, $(\sigma_1,\sigma_2)$
is not a truly simple half-cell. }}

\defi{\label{def:cllmk:2}We will say that a process $X$ on a TLG $\cG$ is \textbf{cell-Markovian} if for 
\begin{enumerate}[(a)]
 
 \item any truly simple cell $(\sigma_1,\sigma_2)$
starting at $t_*$ and ending at $t^*$ the processes $X_{\sigma_1}$ and $X_{\sigma_2}$
are conditionally independent, given the values $X(t_*)$ and $X(t^*)$;
 \item any truly simple right half-cell $(\sigma_1,\sigma_2)$
ending at $t^*$ the processes $X_{\sigma_1}$ and $X_{\sigma_2}$
are conditionally independent, given the value of $X(t^*)$;
 \item any truly simple left half-cell $(\sigma_1,\sigma_2)$
starting at $t_*$ the processes $X_{\sigma_1}$ and $X_{\sigma_2}$
are conditionally independent, given the value $X(t_*)$.
\end{enumerate}

}

\defi{We will say that a process $X$ on a TLG $\cG$ is \textbf{strong cell-Markovian} if 
it is cell-Markovian and for 
\begin{enumerate}[(a)]
 \item any truly simple cell $(\sigma_1,\sigma_2)$
starting at $t_*$ and ending at $t^*$ the processes 
$(X(t):t_*\preceq t\preceq t^*)$ and $(X(t): t^*\preceq t\ \textrm{or}\ t\preceq t^*)$ are independent, given the values $X(t_*)$ and $X(t^*)$;
 \item any truly simple right half-cell $(\sigma_1,\sigma_2)$
 ending at $t^*$ the processes  
$(X(t):t\prec t^*)$ and $(X(t): t^*\preceq t)$ are independent, given the value $X(t^*)$;
 \item any truly simple left-cell $(\sigma_1,\sigma_2)$
starting at $t_*$ the processes 
$(X(t):t^*\prec t )$ and $(X(t):t\preceq t^* )$ are independent, given the value $X(t_*)$.
\end{enumerate}

}

\prop{The constructed process $X$ on $\cG$ is strong cell-Markovian.
\dok{A simple cell in $\cG$ is clearly a simple cell 
in  $\cG'$, and by Lemma \ref{lem:tshcl}. a truly simple half-cell 
is a part of a truly simple cell in $\cG'$. By Theorem \ref{thm:strcell}. $X'$ (on $\cG'$) is strong cell-Markovian
(in the sense of the Definition \ref{def:scllmk}.), and
all the claims now follow.\vspace{0.2cm}


}}

\subsection{Distribution uniqueness}\index{Process indexed by a TLG!uniqueness of distribution}
\prop{A hereditary spine-Markovian $\cM$-process (satisfying (3T') properties) on a TLG$^{**}$ $\cG$ has a unique 
distribution.
\dok{This is a consequence of the unique distribution of $\cM'$-process on a TLG$^*$ $\cG'$.
(See Theorem \ref{teo:uniq}.) }
}

\section{Properties for Markov family $\cM$}

Again, $\cM$ is called a Markov family, if all the measures in $\cM$
are distributions of Markov processes.

\lem{If $\cM$ is a Markov family, so is $\cM'$ and $\cM''$.
\dok{For $\sigma\in P(\cG)$ if $\mu_{\sigma}$ is the distribution of a Markov
process $Y_{\sigma}$, then the process $Y_{\sigma'}$ is also a Markov process,
and hence $\mu_{\sigma'}$ is a distribution of a Markov process.}}

\subsection{Moralized graph-Markovian property}\index{Graph-Markovian property!moralized}
The definition of moralized graph-Markovian property is the same
as in Definition \ref{def:mrlgrpMrk}.

\lem{\label{lem:mrlgmp}The constructed process $X$  on $\cG$  for a Markov family $\cM$ is a moralized graph-Markovian process.
\dok{Let $\cE_1$ and $\cE_2$ be two components of $\cG$ connected
through points $W$, and let $W$ separate $\cE_1$ and $\cE_2$ 
in $(\cG)^{\heartsuit}$. $\cG'$ we will get new edges connecting $t_{-\infty}$
and $t_{+\infty}$, so all the new cells (that are not in $\cG$) will have 
one endpoint in $\{t_{-\infty},t_{+\infty}\}$. If $\cE_1$ and $\cE_2$ were separated 
by $W$ in 
$\cG^{\heartsuit}$, they will be separated in $(\cG')^{\heartsuit}$ by 
$W\cup \{t_{-\infty},t_{+\infty}\}$. Since $X(t_{-\infty})=X(t_{+\infty})=0$,
$\sigma(X_{W\cup \{t_{-\infty},t_{+\infty}\}})=\sigma(X_{W})$. Now, by Theorem \ref{thm:mrl_gp_mk},
$X'_{\cE_1}=X_{\cE_1}$ and $X'_{\cE_2}=X_{\cE_2}$ are independent given $\sigma(X_{W\cup \{t_{-\infty},t_{+\infty}\}})$.}}

\subsection{Time-Markovian property}\index{Time-Markovian property}
The definition of time-Markovian property is the same
as in Section \ref{grph_tm_mrk}. (see Definition \ref{def:t-mrkov}.).

\lem{The constructed process $X$  on $\cG$  for a Markov family $\cM$ is a time-Markovian process.
\dok{Let $t$ be a point in $\cG$. By construction of $X$
we have that
$$\F_t=\sigma\{X(u):u\in \cG, u\preceq t\}\subset \F'_t=\sigma\{X'(u):u\in \cG',u\preceq t\}, $$
$$\cH_t=\sigma\{X(u):u\in \cG, u\succeq  t\}\subset \cH_t'=\sigma\{X'(u):u\in \cG', u\succeq  t\}. $$
(Actually equalities hold in both expressions.) Since $\cM'$ is a Markov family, $X'$ is 
a time-Markovian process. Therefore, $\F_t'$ and $\cH_t'$ are independent given $X'(t)=X(t)$,
but then also $\F_t$ and $\cH_t$ are independent given $X(t)$.
}}

\subsection{Edge-Markovian property}\index{Edge-Markovian property}
The definition of edge-Markovian processes remains the same 
(see Definition \ref{def:edg_mrk}.).

\prop{The constructed process $X$ on $\cG$ for a Markov family $\cM$ is an edge-Markovian process.
\dok{$E$ be an arbitrary edge in $\cG$. Since $\cM'$ is a Markov family, 
$X'$ is edge Markovian, so since $\sigma(X'_{E})=\sigma(X_{E})$ and $\sigma(X(t):t\in \cG,t\notin E)\subset \sigma(X'(t):t\in \cG',t\notin E')$
are independent given the values at the endpoints of $E$,
$X_E$ is independent of $(X(t):t\in \cG,t\notin E)$ given the values at the endpoints of $E$.}}

\section{Processes on time-like trees}\index{Time-like tree (TLT)|(}\label{sec:ptlt}

Among all graphs trees have a special place. Processes on trees have been widely studied 
and used. For examples see Markov chains indexed trees\index{Markov chains indexed by trees} (\cite[Benjamini, Peres]{mchidxtr}), 
branching Markov processes\index{Branching Markov processes} (where the underlying tree is random), hidden Markov models\index{Hidden Markov models}, \ldots

In this section we will look at the properties processes on trees have. We start by defining time-like trees.

\begin{figure}[ht]
\begin{center}
\psfrag{A}{$\boldsymbol{t_{*}}$}
\psfrag{B}{$\boldsymbol{t^{*}}$}
\includegraphics[width=10cm]{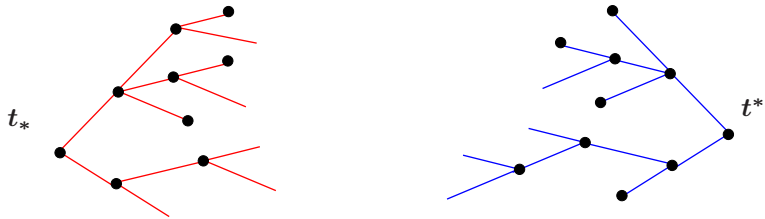}\\
\caption{Forward and backward trees.} \label{pic:sl37}
  \end{center}
\end{figure}

\defi{\label{def:tlt}\begin{enumerate}[(a)]
 \item A \textbf{time-like tree} (TLT)\index{Time-like tree (TLT)|textbf} is a TLG with no cells.
       \item A \textbf{forward time-like tree}\index{Time-like tree (TLT)!forward} $\cT$ is a TLT with exactly one entrance. The entrance vertex we will call the \textbf{root} of the forward time-like tree.
	\item A \textbf{backward time-like tree}\index{Time-like tree (TLT)!backward} $\cT$ is a TLT with exactly one exit
. The exit vertex we will call the \textbf{root} of the backward time-like tree.
      \end{enumerate}}

\noindent\emph{Remark.} Since all vertices are connected to the root, it is clear that $\cT$
is connected as a graph.\vspace{0.3cm}

We know from Theorem \ref{tlg**_ex} (b) that a time-like tree is a TLG$^{**}$, and 
further the following holds:
\lem{\label{lem:tltsub}Let $\cT$ be a TLT, and $\cT'$ be a TLG that is a connected sub-graph of $\cT$.
Then $\cT'$ is a TLT, and $\cT'$ and $\cT$ are elements of the TLG$^{**}$-tower.
\dok{Let $n$ be the difference between the number of edges $\cT$ and $\cT'$ have. For $n=0$ the claim 
is clear. Let's assume the claim holds for $n\geq 0$ and prove it for $n+1$.
Pick a leaf $t_m$ on $\cT$ not in $\cT'$, and an edge $E$ that that is connected to it. 
Now, let $\cT''$ be $\cT$ without $t_m$ and $E$. $\cT''$ is a TLT, and further we can construct 
since the difference between the edges of $\cT''$ and $\cT'$ is $n$, we can construct $\cT''$
from $\cT'$. Hence, they are in some TLG$^{**}$-tower. It is clear that $\cT''$ and $\cT$
are in some TLG$^{**}$-tower. The claim now follows.   }}

For a (3T') family $\cM$ on $\cT$ we can construct a natural $\cM$-process on 
$\cT$. By changing time to each vertex from $t_k$ into $\tilde{t}=1-t$ 
we can transform a backward graph into a forward graph, and in the same way 
transform the process on a backward time-like tree into a process on a forward
time-like tree. Everything we prove for processes on forward TLT's will in a
similar way hold for backward TLT's.

\teo{\label{tree:markov} If $\cM$ is a (3T') family on a TLT $\cT$ and $t_k\in \cV$ is a vertex and $X$ a  natural $\cM$-process on $\cT$ then
\begin{enumerate}[(a)]
 \item the closures of connected components of $R(\cT)\setminus R(t_k)$ are representations of several 
time-like trees $\cT_1$, $\cT_2$,\ldots , $\cT_m$;
 \item the processes $X_{\cT_1}$, \ldots, $X_{\cT_m}$ are independent given the value of $X(t_k)$ . 
\end{enumerate}


\dok{(a) Each of the components is a TLG without any cells. Hence, every component is
a TLT.

\begin{figure}[ht]
\begin{center}
\definecolor{orange}{rgb}{1,0.5,0}
\psfrag{1}{$\boldsymbol{\cT_1}$}
\psfrag{2}{\textcolor{orange}{$\boldsymbol{\cT_2}$}}
\psfrag{3}{\textcolor{blue}{$\boldsymbol{\cT_3}$}}
\psfrag{4}{\textcolor{green}{$\boldsymbol{\cT_4}$}}
\psfrag{k}{\textcolor{red}{$\boldsymbol{t_k}$}}
\includegraphics[width=6cm]{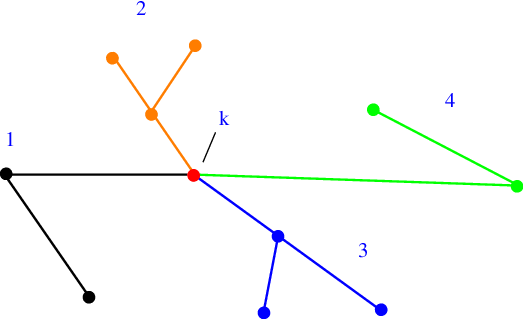}\\
\caption{} \label{pic:sl54}
  \end{center}
\end{figure}

(b) For each $\cT_l$ there is a full path such with no edges in $\cT_l$. Using the spine-Markovian property
we know that $X_{\cT_l}$ is independent of the process on the rest of the graph 
given the value of $X(t_k)$ (since $t_k$ is the only root). Now, applying this fact
several times for $A_l\in \sigma(X_{\cT_l})$ for $l=1,\ldots, m$, we have 
\begin{align*}
 & \E(\P(A_1\cap\ldots \cap A_m|X(t_k))\1_B) = \E(\1_{A_1}\ldots \1_{ A_m}\1_B)\\
 & =\E(\E(\1_{A_1}|X(t_k))\ldots \1_{ A_m}\1_B)=\ldots = \E(\E(\1_{A_1}|X(t_k))\ldots \E(\1_{ A_m}|X(t_k))\1_B)\\
 & = \E(\P(A_1|X(t_k))\ldots \P(A_m|X(t_k))\1_B).
\end{align*}
%
for arbitrary $B\in \sigma(X(t_k))$. Hence the claim follows. }}

\pos{If $\cM$ is a (3T') family on a forward TLT $\cT$ and $t_k$ is a vertex of
degree  at least 3, then the natural $\cM$-process $X$ on $\cT$ will have the property
that given process $X$ on $\cT_{t_k}^+=\{s\in \cT: t_k\preceq s\}$ is independent of the process on the rest of
$\cT$ given $X(t_k)$.}\vspace{0.2cm}

The graph-Markovian property\index{Graph-Markovian property!for time-like trees} was introduced in Definition \ref{def:grpMrk}, and it was shown in
Subsection \ref{not_grph}. that this property doesn't have to hold on TLG$^*$'s. This 
property was replaced by the moralized graph-Markovian property on TLG$^*$'s (see Definition \ref{def:mrlgrpMrk}),
and in Theorem \ref{lem:mrlgmp} it was shown to also hold for natural $\cM$-processes on TLG$^{**}$'s when $\cM$ 
is a Markov family. 

\teo{If $\cM$ is a (3T') Markov family on a TLT $\cT$ the process will have the 
graph-Markovian property.
\dok{By Theorem \ref{lem:mrlgmp}, we know that every natural $\cM$ process on $\cT$ is a moralized graph-Markovian process.
Since $\cT$ has no cells, the claim follows.}}

\pos{If $\cM$ is a (3T') Markov family on a TLT $\cT$ and $\tau_1^t,\ldots , \tau_n^t$ are all the points on $\cT$
with time $t$, 
then the natural $\cM$-process $X$ on $\cT$ will have the property
that 
$$\F^t_{\leftarrow}=\sigma(X(s): s\leq t)\quad and \quad \F^t_{\rightarrow}=\sigma(X(s): s\geq t)$$
are independent given $X(\tau_1^t)$, \ldots, $X(\tau_n^t)$.
\begin{figure}[ht]
\begin{center}
\psfrag{1}{$\boldsymbol{\tau_1^t}$}
\psfrag{2}{$\boldsymbol{\tau_2^t}$}
\psfrag{3}{$\boldsymbol{\tau_3^t}$}
\psfrag{4}{$\boldsymbol{\tau_4^t}$}
\psfrag{5}{$\boldsymbol{\tau_5^t}$}
\psfrag{6}{$\boldsymbol{\tau_6^t}$}
\psfrag{t}{$\boldsymbol{t}$}
\includegraphics[width=12cm]{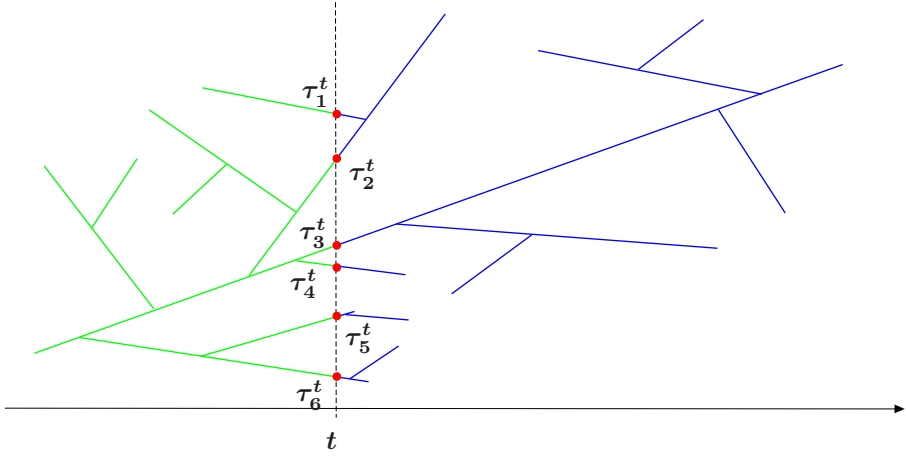}\\
\caption{The process \textcolor{green}{before} and \textcolor{blue}{after} time $t$.} \label{pic:sl55}
  \end{center}
\end{figure}
\dok{

The points $\tau_1^t$, \ldots, $\tau_m^t$ separate the graph into two parts
$\{s\in \cG: s\leq t\}$ and $\{s\in \cG: s\geq t\}$, and the claim follows by
graph-Markovian property. }}\noindent

\noindent\emph{Remark.} The previous corollary states that the process 
$(\tilde{X}(t)=(X(s):s\in R(\cG)\cap (\{t\}\times \R^2)):t\geq 0)$ is a Markov process.\vspace{0.2cm}

The following lemma states that the spine-Markovian property\index{Process indexed by a TLG!spine-Markovian property!for time-like trees} and hereditary spine-Markovian
properties are equivalent on time-like trees. (Note that we didn't have this result 
for TLG$^*$'s.)

\lem{\label{lm:tree_hspn}If $\cT$ is a time-like tree, and $X$ a process indexed by $\cT$ then the following claims 
are equivalent:
\begin{enumerate}[(a)]
 \item $X$ is a spine-Markovian process;
 \item $X$ is a hereditary spine-Markovian process. 
\end{enumerate}
\dok{Clearly (b) implies (a). Now, let's prove that (a) implies (b). Let $(\cG_k)_{k=1}^n$
be a TLG$^{**}$ tower leading towards the construction of $\cT$. Note that since
each $\cG_k$ is a connected subgraph of $\cT$, it is also a tree.

\begin{figure}[ht]
\begin{center}
\psfrag{p}{$\textcolor{blue}{\boldsymbol{\pi'}}$}
\psfrag{q}{$\boldsymbol{\pi}$}
\includegraphics[width=6cm]{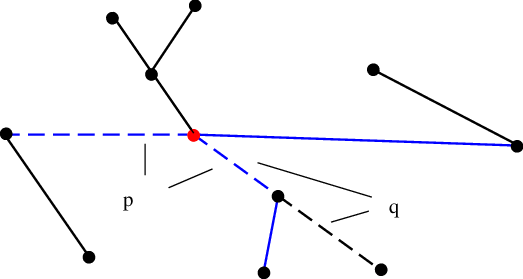}\\
\caption{\textcolor{blue}{$\cG_k$}, the spine \textcolor{blue}{$\boldsymbol{\pi'}$} and the root \textcolor{red}{$\bullet$}.} \label{pic:sl53}
  \end{center}
\end{figure}
  
If we pick a $\pi'$ full-path in $\cG_k$, then there is a full-path $\pi$ in $\cT$ such that
$R(\pi')\subset R(\pi)$. (See Figure \ref{pic:sl53}.) But the representation of roots of $\cT$ with respect to
$\pi$ will contain the representations of roots of $\cG_k$ with respect to 
$\pi'$. Since the roots decompose the graph into disjoint components the claim follows.
}}
\index{Time-like tree (TLT)|)}

\chapter{Galton-Watson time-like trees and the Branching Markov processes}\label{chp:gw}
\section{TLG's with an infinite number of vertices}

We will allow $t_0$ and $t_N$ to take values in $\R$.\index{Time-like graph (TLG)!with infinitely many vertices|(}

\defi{\label{def:intlg**}\begin{enumerate}[(i)]
       \item Suppose that the vertex set of a graph $\cG=(\cV,\cE)$ is infinite. We will call
$\cG$ a time-like graph (TLG) if it satisfies the following conditions.
\begin{enumerate}[(a)]
 \item There is a sequence of TLG's $\cG_n=(\cV_n,\cE_n)$ with finite vertex
set $\cV_n$, $n\geq 1$, and for some representations of $\cG_n$'s and $G$ we have
$$\bigcup_{n=1}^{\infty}R(\cG_n)=R(\cG).$$ 
\item \label{def:intlg**:b} The graph $\cG$ is locally finite, i.e. it has a representation $R(\cG)$
such that for any compact $K\subset \R^3$ a finite number of edges intersects $K$.
\end{enumerate}
\item A TLG $\cG$ with infinite vertex set will be called an TLG$^{**}$ if it satisfies
the following conditions.
\begin{enumerate}[(a)]
  \item We can choose a sequence of TLG$^{**}$'s $\cG_n$ in (i).
(In the sense of the Definition \ref{def:tlg**}.(\ref{def:tlg**:iii}), i.e. $(\cG_j)_{1\leq j\leq n}$
is a tower of TLG$^{**}$'s for all $n$.)
  \item Let $\cV_n=\{t_{0,n},t_{1,n},\ldots,t_{N_n,n}\}$. The initial 
vertices $t_{0,n}\in \cV_n$ and $t_{N_n,n}\in \cV_n$ are the same 
for all $\cG_n$, i.e. for all $1\leq m\leq n$
$$t_{0,n}\leq t_{0,m}\quad \textrm{and}\quad t_{N_n,n}\geq t_{N_m,m}.$$

\end{enumerate}
      \end{enumerate}
 }

%

The following lemma will be useful for the construction of processes. (It is a version of the Lemma \ref{lem:subgphfcv}.
for TLG$^{**}$'s.)

\lem{\label{lem:subgphfcv2}Let $(\cG_n)$ and $(\cG_n')$ be two TLG$^{**}$-towers that lead to the construction 
of $\cG$. Let $\cH$ be a sub-graph (not necessarily a TLG$^*$) of some $\cG_{n_0}$. Then there exists
$\cG_{n_1}'$ such that $R(\cH)\subset R(\cG_{n_1}')$ and all the vertices of $\cH$ 
are contained in $\cG_{n_1}'$.
\dok{Since $\cG$ is locally finite, there are finitely many vertices with representation 
on $R(\cH)$, also these vertices are of finite degree. For each such vertex $v$, by same argument,
there has to be $\cG_{n_v}'$ such $v$ in $\cG_{n_v}'$ has that degree. Now if $n_1$
is the maximum of $n_v$ over each such vertex $v$ the claim follows. }}

The definition of (forward/backward) time-like trees is the same as in Definition \ref{def:tlt}.

\prop{Time-like tree $\cT$ with infinite number of vertices is a TLG$^{**}$.
\dok{
%
%
Pick a vertex $t_k$, and let $K_n$ be a set of compact sets such that 
$$\bigcup_{n=1}^{\infty}K_n= \R^3.$$
It is clear that the connected component of $R(\cT)\cap K_n$ that contains $t_k$ 
is a tree, and we set $\cT_n$ to be the time-like tree such that $R(\cT_{n-1})\subset R(\cT_n)\subset R(\cT)\cap K_n$ 
and the number of $\cT_n$ is as large as possible. By Lemma \ref{lem:tltsub}, $\cT_n$ can be constructed from 
$\cT_{n-1}$. So $(\cT_n)$ is a subsequence of some TLG$^{**}$-tower $(\cH_n)$. }}


\index{Time-like graph (TLG)!with infinitely many vertices|)}

\section{Galton -- Watson time-like tree}\label{GWTLT}\index{Time-like tree (TLT)!Galton-Watson|see{Galton-Watson time-like tree}}\index{Crump - Mode - Jagers trees!Galton-Watson|see{Galton-Watson time-like tree}}\index{Galton-Watson time-like tree|(}
We will encode a continuous version of Galton-Watson process into a (forward) time-like tree.
The idea is to use the setup in the Crump -  Mode - Jagers model\index{Crump - Mode - Jagers trees} (see Section \ref{CMJtree}.).\vspace{0.3cm}

Let $I=\{\emptyset\}\cup\bigcup_{n=1}^{\infty} \N^n$, and we interpret that 
$(x,j)\in I$, $j\in \N$ is a child of $x\in I$.\vspace{0.3cm}

First, lets make some assumptions:
\begin{itemize}
\item Let $(\lambda_x : x\in I)$ be a collection of exponential random variables
with parameter $V$. (Lifetime of an individual.)
 \item Let $(R_x : x\in I)$ be a collection of random variables with distribution given by the
generating function 
$$\Phi (s)=\sum_{k=0}^\infty a_ks^k, \quad \Phi(1)=1.$$
\item $(\lambda_x , R_x)_{x\in I}$ is an i.i.d. sequence.

\end{itemize}
In our model at the end of its lifetime, the individual gets divided 
into nonnegative number of new individuals (0, 1, 2,\ \ldots ), so we define the reproduction 
function to be 
$$\xi_x(t)=R_x\1_{(t\geq \lambda_x)}.$$

Recall, that we defined with $\tau_x$ the birth time of $x$, with
$\tau_{\emptyset} =0$, $\tau_{(x',i)}=\tau_{x'}+\inf\{u: \xi_{x'}(u)\geq i\}$.

We define the Galton - Watson tree in the following way:
\begin{itemize}
 \item Let $\cV_0=\{t_{-1}=0, t_{0}=\lambda_{\emptyset}\}$, and $\cE_1=\{E_{-1, 0}\}$. Set $\cG_0=(\cV_0,\cE_0)$.
 \item Define $l_{Ex}^+(\cG_k)$ all the $x$ labels of exit vertices in $\cG_k$ such that $\xi_x>0$.
 \item Now set 

$$\cV_k:=\cV_{k-1}\cup \bigcup_{x'\in l_{Ex}^+(\cG_{k-1})} \{t_{(x',j)}:=\tau_{(x',j)}: j=1\ldots \xi_{x'}\},$$
$$\cE_k:=\cE_{k-1}\cup \bigcup_{x'\in l_{Ex}^+(\cG_{k-1})} \{E_{x',(x',j)}: j=1\ldots \xi_{x'}\},$$
\be \cG_k :=(\cV_k,\cE_k) \label{df:gk:tlg**}.\ee
\item Now, set 
$$\cV:= \bigcup_{k=1}^{\infty} \cV_k, \quad \cE=\bigcup_{k=1}^{\infty} \cE_k.$$
\end{itemize}
We define $\cG=(\cV,\cE)$ as the {\bf Galton -- Watson time-like tree}.

%

\lem{\label{cor:tree:fin} For all $t\geq 0$ we have $\E(\xi_x(t))\leq \E(R_x)$. 

Therefore, if $\E(R_x)<\infty$, almost surely
for all $t\geq 0$ the number of vertices from $\cV$ with time at most $t$ is finite. 
\dok{Since $\xi_x\leq R_x$ the first claim follows. 
For the second claim,  first note
$\E(\xi(t))<\infty$. Further $(0\geq \lambda_x)=(0=\lambda_x)$, and this is a set of probability $0$, hence
$\E(\xi(0))=0<1$. Therefore, by Theorem \ref{teo:fin_br_pr}. the set of vertices
with time label at most $t$ is finite.  }
}

\teo{If $\E(R_x)<\infty$, the Galton - Watson time-like tree is a TLG$^{**}$. 
Specially, it is a forward time-like tree.
\dok{It is clear that $(\cG_k)$ from $(\ref{df:gk:tlg**})$ is the TLG$^{**}$-tower that leads 
towards the construction of $\cG$. Further, any representation is locally finite, since any compact set 
$K$ will have a finite time component, i.e. it will be contained in $[0,T]\times \R^2$, and by Lemma \ref{cor:tree:fin}. it can contain finitely many points 
finitely many edges. We know by Lemma \ref{cor:tree:fin} that the number of vertices 
whose representation is in $[0,T]\times \R^2$ is finite a.s., and also since $R_x$ is finite a.s. we have 
that number of edges intersecting $K$ is finite. }}\index{Galton-Watson time-like tree|)}

\section{Processes on TLG$^{**}$'s with infinite number of vertices}\label{sec:tlg**incst}

\subsection{Construction}\index{Process indexed by a TLG!with infinite number of vertices|(}
\index{Process indexed by a TLG!with infinite number of vertices!construction}

Let $\cG=(\cV,\cE)$ a TLG$^{**}$ such that $\cV$ is infinite. According to the definition,
there exists a tower of TLG$^{**}$'s $\cG_n=(\cV_n,\cE_n)$, $n\geq 1$, such that 
$\cV_n$ is finite, where $\cV=\bigcup_{n\geq 1}\cV_n$.\vspace{0.3cm}

Let 
\begin{equation}
\cM=\{\mu_{\sigma}: \sigma\in P(\cG)\} \label{kf:2} 
\end{equation}

be a family of distributions of processes along
full-time paths in $\cG$ satisfying conditions (T'1)-(T'3) given in Subsection \ref{subsec:cond}. 
\vspace{0.3cm}

Since $$\cM(\cG_n)=\{\mu_{\sigma}: \sigma\in P(\cG_n)\}$$ 
is well-defined, and we can show similarly as in Lemma \ref{lem:inr:A13} 
that $\cM(\cG_n)$ satisfies (T'1)-(T'3), we can define a hereditary spine-Markovian process $X^n$ on $\cG_n$, such that
for each $\sigma\in P(\cG_n)$ the process $X^n_{\sigma}$ has the distribution $\mu_{\sigma}$. 
Further, 
the restriction of this process to $\cG_k$ ($k\leq n$)  has the same distribution as
the $\cM(\cG_k)$-process $X^k$ defined on $\cG_k$ in the similar manner.\vspace{0.3cm}

Now, Kolomogorov's consistency theorem shows, that there
exists a process $X$ on $\cG$ such that the restriction of $X$ to any $\cG_k$
has same distribution as $X^k$. Note, that since each $\sigma\in P(\cG)$ is in some of the
$\cG_k$'s we have $X_\sigma$ has the distribution $\mu_{\sigma}$.

\subsection{Uniqueness of distribution}\index{Process indexed by a TLG!with infinite number of vertices!uniqueness of distribution}
Using a similar approach as in \S\ref{subsec:undistr} we will get 
that the distribution of the process $X$ doesn't depend on the choice
of the TLG$^{**}$-tower $(\cG_n)$.

\lem{\label{lem:2twtlg**}Let $\cG$ be a TLG$^{**}$ with infinitely many vertices, $(\cG_j^1)$ and $(\cG_j^2)$ two TLG$^{**}$-towers 
that construct $\cG$ and $X^1$ and $X^2$ the natural $\cM$-processes constructed using these two towers.  
The distribution of the processes $X^1$ and $X^2$ restricted on $\cG_{k}^1$ is the same for all $k$.

\dok{We first prove the claim when the vertices of $\cG$ have only real values.
By Lemma \ref{lem:subgphfcv2}. we can choose $k_1$, and $l_1$  in such that 
$$R(\cG_{k}^1)\subset R(\cG_{l_1}^2)\subset R(\cG_{k_1}^1),$$
where $\cV_{\cG_{k}^1}\subset \cV_{\cG_{l_1}^2}\subset \cV_{\cG_{k_1}^1}$.
Now, we look at the embeddings $(\cG_{k}^1)''$, $(\cG_{l_1}^2)''$ and $(\cG_{k_1}^1)''$.
We will have the same relationships, and by Lemma \ref{lem:2seqcrt}, we know that $(\cG_{k}^1)''$ and $(\cG_{l_1}^1)''$ 
are in some TLG$^{*}$-tower.
Now, by  Theorem \ref{teo:uniq}. and Theorem \ref{teo:3ver:cstr}. the result follows.}}

\teo{Let $\cG=(\cV,\cE)$ be a TLG$^{**}$'s with infinitely many vertices in $\cV$,
and let $X^1$ and $X^2$ be two $\M$-processes constructed using the TLG$^{**}$-towers 
$(\cG_n^1)$ and $(\cG_n^2)$, then  $X^1$ and $X^2$ have the same distribution.
\dok{Let $t_1,\ldots, t_m$ be the points on $\cG$ with finite time. Then, by Lemma \ref{lem:subgphfcv2}., there exists
$\cG^1_k$ that contains all of these points. By Lemma \ref{lem:2twtlg**} it follows,
that $X^1$ and $X^2$ have the same distribution on $\cG_k^1$. Specially, 
$(X^1(t_1),\ldots, X^1(t_m))$ and $(X^2(t_1),\ldots, X^2(t_m))$ have the same distribution. Now, by 
Kolomogorov's Consistency Theorem the claim follows. }}

\pos{\label{cor:ind_cnstr}The distribution of the process $X$ on $\cG$ doesn't depend on the choice 
of the TLG$^{**}$-tower $(\cG_j)$ that constructs $\cG$. }

\defi{\label{def:ntr_Minf}We call the constructed process $X$ the natural $\cM$-process on the TLG$^{**}$ $\cG$.}

\index{Process indexed by a TLG!with infinite number of vertices|)}
\section{Natural $\cP$-Markov process}

First, let's define the natural $\cP$-Markov process.

\defi{Let $\cG$ be a TLG$^{**}$ and $\cP$ a distribution of a Markov process on $[0,\infty)$, 
then \textbf{natural $\cP$-Markov process} on $\cG$ 
is a stochastic process $X$ indexed by $\cG$ such that the distribution of $X$ 
along each path $\pi$ from any point $t_j$ to any  
other point $t_k$ is distributed as a $\cP$-Markov process along $[t_j,t_k]$, and satisfies (3T') conditions.
This induces a (3T') family $\cM_{\cP}$, and the natural $\cP$-Markov process on $\cG$
is the natural $\cM_{\cP}$-process on $\cG$ (see Definition \ref{def:ntr_Minf}.)}\vspace{0.2cm}

The following was shown in Section \ref{sec:tlg**incst}.
\teo{For any distribution $\cP$ of a Markov process on $[0,\infty)$ 
and any TLG$^{**}$ $\cG$ whose time components are all greater or equal to $0$, there
exists a natural  $\cP$-Markov process.}

\section{Branching $\cP$-Markov process}\label{sec:mbmpp}\index{Branching Markov process|(}
Idea of this section is to construct a natural $\cP$-Markov process  on a random Galton - Watson tree,
where $\cP$
is a distribution of an RCLL or continuous process.
We will also show its connection to the branching $\cP$-Markov process. Specially, to show that in the case
when $\cP$ is the distribution of the Brownian motion, that we have the branching Brownian motion.\vspace{0.2cm}

Basically, we first construct a Galton -- Watson tree, and then  on that tree we 
construct the $\cP$-Markov process indexed by it.

\begin{itemize}
 \item Based on the construction in Section  \ref{GWTLT} construct a Galton-Watson time-like
tree $\cT$.
 \item Construct a natural $\cP$-Markov process on $\cT$ whose values are independent of $\cT$.
\end{itemize}

Note, that the probability space on which we live can be written as 
$$\left[\prod_{x\in I}(\R\times \N_0, \cB(\R)\times \cP(\N_0))\right]\times \left[\prod_{x\in I}(D[0,\infty), \cB(D[0,\infty)))\right]
$$
This is a product of countably many Borel spaces, and therefore it is a Borel space.
The first part of the product encodes the tree, while the second part is used to construct 
the process on the tree.

\subsection*{Construction of the tree}

As discussed in Section \ref{GWTLT}. the sequence $(\lambda_x, \xi_x)_{x\in I}$ encodes the whole tree, and 
from there we can get the time $\tau_x$ of birth of each individual $x\in I$. 
(Recall, that $\lambda_x$ is the lifetime of $x$ and $\xi_x$ is the number of children.)

If $\tau_x=\infty$ 
then $x$ was never born. Since the sequence was i.i.d. we can construct a 
probability measure on 
$$(\Omega_{\cT},\F_{\cT})=\prod_{x\in I}(\R\times \N_0, \cB(\R)\times \cP(\N)).$$
We know that $\cT$ is a time-like tree a.s.

\subsection*{Construction of the process}

We will construct a probability on the space 
$$(\Omega,\F)= (\Omega_{\cT},\F_{\cT})\times\prod_{x\in I}(D[0,\infty), \cB(D[0,\infty))).$$
For each element $((\lambda_x, \xi_x)_{x\in I},(f_x)_{x\in I})$:

\begin{itemize}
 \item $(\lambda_x, \xi_x)_{x\in I}$ is distributed as Galton-Watson time-like tree
 \item $f_x|_{[\tau_x,\tau_x+\lambda_x)}$ represents the space
position of $x$ during its lifetime
 \item $f_{x}|_{\R\setminus[\tau_x,\tau_x+\lambda_x)}=\Delta$ for all $h\geq 0$ (represents cemetary). 
%
 \item If $\tau_{(x,j)}<\infty$ then $f_{(x,j)}(\tau_{(x,j)})=f_x((\tau_{x}+\lambda_x)^-)$ almost surely for all $x\in I$ and $j\in \N$ (last position of the 
parent, is the first position of the child).
\end{itemize}
Specially, if $\tau_x=\infty$ then 
\begin{itemize}
 \item $f_{x}(h)=\Delta$ for all $h\geq 0$ (never born, remains on cemetary).
\end{itemize}

Let's make some assumptions on the distribution $\cP$ and introduce some notation.
Let $(X(t):t\geq 0)$ be a $\cP$-distributed process:
\begin{itemize}
 \item by $\cP_{\tau}^x$ we are denoting the 
distribution of the process $(X(\tau+t):t\geq 0)$ conditioned on the event $X_{\tau}=x$.
\end{itemize}
 
We will assume the following on $(\cP_{\tau}^x\ :\ \tau\geq 0,x\in \R)$
for all $A\in {\cal B}(D[0,\infty))$ the map 
$$(\tau,x)\mapsto \cP_{\tau}^x(A)$$
is a measurable function. This clearly holds in the case of many time-homogeneuos 
Markov process (e.g. Brownian motion or Levy processes).\vspace{0.2cm}

We do the following construction, based on first child - next sibling\index{Firs child - next sibling} idea from 
computer science.

\begin{center}
\begin{algorithm}[H]
\SetKwIF{If}{ElseIf}{Else}{loop}{}{else if}{else}{endif}
$A_0=\{\emptyset\}$\;
$k=0$\;

\If{}{$k=k+1$\;
\For{$x\in A_{k-1}$}{
add to $A_k$ first child and next sibling of $x$\; }
}
\caption{First child - next sibling search of the plane tree}\label{alg:1cnsb}
\end{algorithm}
\end{center}

We now order the $I$ in a sequence $(x_n)$, such that we first all 
the elements of $A_0$ appear, then of all the elements of $A_1$ appear, 
then of $A_2$ \ldots \vspace{0.3cm}

Now $(\Omega_{\cT},\F_{\cT})$, $(\Omega_{x_1},\F_{x_1})$, $(\Omega_{x_2},\F_{x_2})$, 
\ldots is a sequence of measurable spaces, and we have the following probability measures on them:
\begin{itemize}
 \item On $(\Omega_{\cT},\F_{\cT})$ we define $\P_{\cT}$ as explained in the previous subsection;
 \item On $(\Omega_{x_1},\F_{x_1})$ we define $\P_{x_1}^{\lambda_{x_1}}$ as the distribution of the process $(Y(t):t\geq 0)$
where
$$Y(t)=\left\{\begin{array}{cc}
              X(t), & t< \lambda_{x_1};\\
	      \Delta, &t\geq  \lambda_{x_1};
             \end{array}
\right. $$
where the distribution of $(X(t):t\geq 0)$ is $\cP$. 
 \item On $(\Omega_{x_j},\F_{x_j})$ we define $\P_{x_j}^{\tau_{x_{j}},\lambda_{x_j},f_{x_{j'}}}$
to be the distribution of the process $(Y(t):t\geq 0)$ given by
$$Y(t)=\left\{\begin{array}{cc}
	      \Delta & t< \tau_{x_j}\\
              X(t), & \tau_{x_j} \leq t<\tau_{x_j}+ \lambda_{x_j};\\
	      \Delta, &t\geq  \tau_{x_j}+\lambda_{x_j};
             \end{array}
\right.$$
where $(X(t):t\geq 0)$ is distributed as
$\cP_{\tau_{x_j}}^{f_{x_{j'}}(\tau_{x_{j}})}$ where 
$f_{x_{j'}}\in \Omega_{x_{j'}}$, $x_{j'}$ is the parent of $x_j$, and we can show 
that $\tau_{x_{j}}$ is a measurable function on
$(\Omega_{\cT},\F_{\cT})$. Therefore, since $j'<j$, for $A\in \F_{x_j}$ 
$$\P_{x_j}^{\tau_{x_{j}},\lambda_{x_j},f_{x_{j'}}}(A)$$
is $\prod_{k=1}^{j-1}(\Omega_{x_j},\F_{x_j})$-measurable.
\end{itemize}

Now we can define a product probability on $(\Omega,\F)$  using Theorem  \ref{cond_seq:inf}. 

Now, for $\omega=(\omega_{\cT},(f_{x_j}))\in \Omega$. $\cT(\omega)$ is represented by
$\omega_{\cT}$, and for $E_{x_{j'}x_{j}}$ an edge in $\cT(\omega)$ we define 
$$X_{E_{x_{j'}x_{j}}}(\omega)(t)=f_{x_{j}}(t)$$ for $\tau_{x_{j'}}\leq t< \tau_{x_{j}}$.

\subsubsection*{Properties of the construction}

\teo{The probability measure is well-defined, that is $\P$ doesn't depend 
on the choice of $(x_n)$ as along as: 
\begin{enumerate}[(1)]
 \item $x_0=\emptyset$;
 \item $\{x_n:n\in \N_0\}=\{x:x\in I\}$;
 \item For each $j\geq 1$ there exists $j'<j$ such that $x_{j'}$ is a parent of $x_{j}$.
\end{enumerate}
\dok{Conditioned on $(\cT=T)$ the constructed process can be mapped into a construction of 
a natural $\cP$-process on a TLT $T$. The distribution of the process by Corollary \ref{cor:ind_cnstr} doesn't depend 
on the constrcution, hence the probability measure is well-defined.}}

\pos{The distribution of constructed process conditioned that the underlying tree
$\cT=T$ is a natural $\cP$-process on $T$.}

\teo{If $\cP$ is a distribution of a Markov process, for constructed process $(\cT,X)$
the process 
\be Y(t)=\{X(\tau):\bar{\tau}\in R(\cT)\cap (\{t\}\times \R^2)\}\label{eq:brPrc}\ee
is a Branching $\cP$-Markov process.
\dok{Follows from description stated in \S\ref{sec:mbmpp}. }}

\pos{If $\cP$ is a distribution of Brownian motion, then the process given by $(\ref{eq:brPrc})$ is the branching Brownian motion\index{Branching Brownian motion}.}\index{Branching Markov process|)}

\part*{Open questions and appendix}

\chapter{Open questions}
In this chapter we will state some open problems that could be of interest 
for further research.
\section{Construction of process on all TLG's}
As it was pointed out by Burdzy and Pal in \cite{tlg1} (and in \S\ref{subsec:cnsprb} of this paper), it is not possible 
to construct a natural Markov process on every  TLG.\vspace{0.2cm}

\begin{figure}[ht]
\begin{center}

\includegraphics[width=7cm]{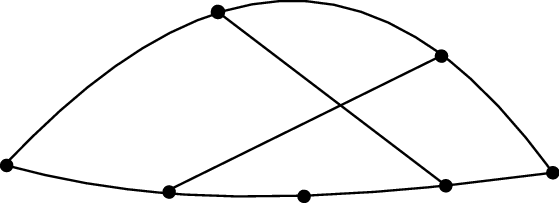}\quad \includegraphics[width=6.8cm]{tlgn_hded.eps}\\[0.1cm]
\includegraphics[width=7cm]{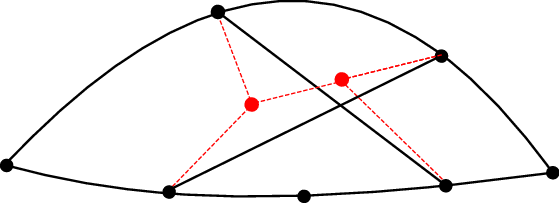}\quad \includegraphics[width=6.8cm]{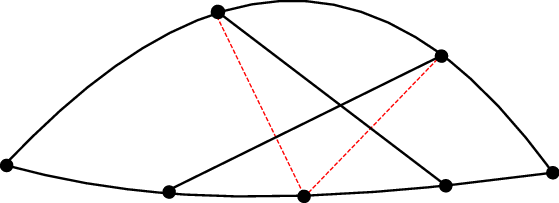}\\
\caption{ Example from Theorem \ref{teo:exncstr}. and different embeddings into a TLG$^*$.} \label{pic1_2}
  \end{center}
\end{figure}

Theorem \ref{teo:exncstr}. shows that a Brownian motion with the cell-Markovian property  indexed by the TLG $\cG$ 
given on the first image of the Figure \ref{pic1_2} does not exist.

We know, by discussion in Section \ref{constr_cnsfm}, that it is possible to construct a Markov process on 
a TLG$^*$. We could try to embed $\cG$ into some TLG$^*$ $\cH$, define a natural Brownian motion $X$ on $\cH$
and then restrict $X$ to $\cG$ (i.e. set $X_{\cG}=(X(t):t\in \cG)$).\vspace{0.2cm}

It is possible to embed any TLG into a TLG$^*$.

\teo{Let $\cG=(\cG,\cV)$ be a (unit) TLG, then there exists a TLG$^*$ $\cH$ that is a sup-graph of $\cG$.
\dok{Let $\tau_1$, \ldots, $\tau_m$ be times of vertices of $\cV$. Now, we construct $\cV_{\cH}$ that 
contains $\cV$ and vertices
$t^*_{1/2}$,$t^*_{3/2}$,\ldots, $t^*_{m+1/2}$ with times $\tau_{1/2}=-1$, $\tau_{3/2}=\frac{\tau_1+\tau_2}{2}$, \ldots, $\tau_{m-1/2}=\frac{\tau_{m-1}+\tau_{m}}{2}$, $\tau_{m+1/2}=2$. 
Now, we set $\cE_{\cH}^0$ is constructed in such a way that $t_k\in \cV$ with time $\tau_j$ the edge
\begin{itemize}
 \item $E_{j-1/2,k}$ between $t^*_{j-1/2}$ and $t_k$ is in $\cE_{\cH}^0$;
 \item $E_{k,j+1/2}$ between $t_k$ and $t^*_{j+1/2}$ is in $\cE_{\cH}^0$.
\end{itemize}
It is not hard to see that $\cH_0=(\cV_{\cH},\cE_{\cH}^0)$ is a planar simple TLG, therefore by Theorem \ref{thm:tlg*1}.
a TLG$^*$. 

\begin{figure}[ht]
\begin{center}
\includegraphics[width=6cm]{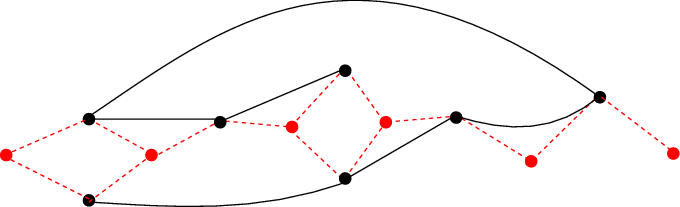}\\
\caption{\textcolor{red}{$\cH_0$} (induced by dashed \textcolor{red}{edges}) is planar. } \label{sl58}
  \end{center}
\end{figure} 

Further, every two vertices $t_j$ and $t_k$ are connected by a time-path in  $\cH_0$. 
Hence, we can add one by one edge from $\cE$ to $\cH_0$, and $\cH=(\cV_{\cH},\cE_{\cH}^0\cup\cE)$ is a TLG$^*$. }}

It is clear that the distribution of a Brownian motion on  $\cG$
will depend on the embedding $\cH$.

\begin{itemize}
 \item For a given (simple) TLG $\cG$, under what conditions on the distributions along time-paths 
can we construct a process on $\cG$? 
 \item Are there examples of distributions of (Markov) processes along time-paths 
for which this is not possible?
 \item Is there a way of getting the uniqueness of distribution of $X$ on $\cG$?
 \item What properties will the constructed process have?
\end{itemize}

\section{Reconstruction of TLG's based on the process}

As we saw in the previous section, the fact that the underlying graph is not a
TLG$^*$ or TLG$^{**}$ does not have to prevent us from defining a process on it.\vspace{0.2cm}

It could be that a part of the graph and a part of the process on that graph is 
hidden from us.\vspace{0.2cm}

Suppose $X$ is a natural $\cM$-process on a TLG$^*$ $\cH$ where $\cM$ is a family 
of distributions of Gaussian Markov processes. Let $\cG$ be a TLG such that  $R(\cG)\subset R(\cH)$.\vspace{0.2cm}

\begin{itemize}
 \item If we know how the graph $\cG$ looks like and we know the distribution of $X_{\cG}=(X(t):t\in \cG)$,
how much can we say about $\cH$?
 \item What if we don't know the distribution of the process $X$ on the whole $\cG$,
but only on the part of it?
 \item Could we use any of this on the branching Markov process (specially on 
branching Brownian motion)?
\end{itemize}

In classical graphical models problems of hidden (latent) variables have been studied
(see Chapter 20. in \cite{prgphmdl} or \S17.4 in \cite{stat_learning}). 
One of the strong 
tools in solving the problems could be the moralized graph-Markovian\index{Graph-Markovian property!moralized} property, which 
enables us to project a process on a TLG into Markov random field (MRF) \index{Markov random field (MRF)}.  We could 
use some of the properties of MRF's to detect hidden parts of the graph.\vspace{0.2cm}

There is an interesting criteria for finding edges in a Gaussian MRF.\index{Markov random field (MRF)!Gaussian} Before we state 
that we will need the following lemma.

\lem{If $X=(X_1,\ldots, X_n)$ is a Gaussian random vector with positive definite covariance matrix $\Sigma$\index{Covariance matrix}, then 
$X_i\perp X_j|(X_k:k\in \{1,2,\ldots , n\}\setminus \{i,j\}) $ if and only if 
$\Sigma_{ij}^{-1}=0$.}

\prop{\label{prop:conedge}Let $G=(V,E)$ be a undirected graph and $X=(X_v:v\in V)$ a Gaussian Markov random field. Let $K$ be 
the positive definite covariance matrix of $X$. If $\{u,v\}\notin E$ then $K^{-1}_{u,v}=0$.}

For proof of these claims see Chapter 7. in \cite{prgphmdl} or Chapter 5. in \cite{graph_models}.
With these results we can show the following.

\prop{Let $X$ be a natural $\cM$-process on a unit TLG$^*$  $\cH$, where $\cM$ is a family of Gaussian Markov 
processes. Assume we know the distribution of  $X_{R(\cG)}$ where  $\cG=(\{0,1\},\{E_{01}^1,E_{01}^2\})$. 
If $R(\cG)$ is a representation of a truly simple cell in $\cH$ then for the covariance matrix $K(t_1,t_2)$ of the vector 
$$X=(X(0),X(t_1),X(t_2),X(1))$$
we have $K(t_1,t_2)^{-1}_{2,3}=0$ for every point $t_1\in E_{01}^1$ and $t_2\in E_{01}^2$ with times 
in the interval $(0,1)$.
\begin{figure}[ht]
\begin{center}
\psfrag{0}{$\boldsymbol{0}$}\psfrag{1}{$\boldsymbol{1}$}
\psfrag{a}{$\boldsymbol{t_1}$}\psfrag{b}{$\boldsymbol{t_2}$}
\includegraphics[width=6cm]{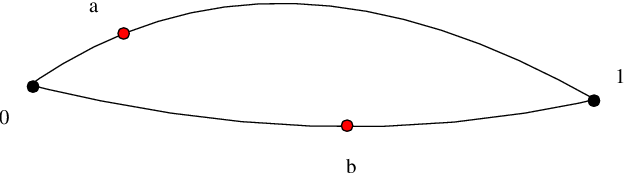}\quad \includegraphics[width=6cm]{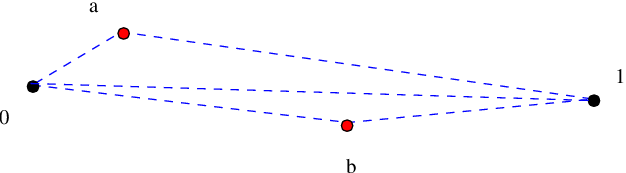}\\
\caption{$\cG$ and the induced graphical model } \label{sl56}
  \end{center}
\end{figure} 
\dok{If $\cG$ is a representation of truly simple cell, then by the Corollary \ref{cor:sum_int}, the strong cell-Markovian property (Theorem \ref{thm:strcell}) and moralized graph-Markovian property
(Theorem \ref{thm:mrl_gp_mk}.)
we know that $X$ can be represented as a graphical model. In this representation there will be no edge 
between $t_1$ and $t_2$, and by Proposition \ref{prop:conedge}. the claim follows.}} 

Making some natural conditions on the distributions on the family $\cM$ and using
the variable elimination algorithm (see Chapter 9. in \cite{prgphmdl}) in  for MRF's we could try to get the converse 
of the statement.

\begin{itemize}
 \item If $\cG$ is not the representation of a truly simple which paths 
can we detect?
\end{itemize}

\section{Strong Markov property, parametrization, evolution over time,\ldots}
In Chapter \ref{sec:3b} we defined stopping times and proved the Optional Sampling Theorem. 
We also proved the time-Markovian property\index{Strong Markov property}, and the the following question naturally follows.
\begin{itemize}
 \item Do we have a version of the strong Markov property for a natural $\cM$-process, where $\cM$ is a Markov family?
\end{itemize}\vspace{0.2cm}

Parametrizng the process in suitable way and calculating probabilities is always
a challenge.
\begin{itemize}
 \item Is there a convenient way to parametrize the family along time-paths of a TLG $\cG$?
 \item Is there a procedure how to calculate finite dimensional distributions of the 
process on the TLG $\cG$?
 \item Is there a procedure how to calculate finite dimensional distributions
conditioned that we know some values of the 
process on the TLG $\cG$?
\end{itemize}\vspace{0.2cm}

We could evolve the process on a graph $\cG$ over time, and maybe even make 
the graph evolve over time.
\begin{itemize}
 \item Could we define a process $(X^\tau:\tau\geq 0)$ such that 
$X^\tau=(X^{\tau}(t):t\in \cG)$ is a process indexed by a TLG $\cG$?
\item Could we define a process $(X^\tau:\tau\geq 0)$ such that 
$X^\tau=(X^{\tau}(t):t\in \cG(\tau))$ is a process indexed by a TLG $\cG(\tau)$?
\end{itemize}

We saw one way to randomize the underlying graph in Chapter \ref{chp:gw}, we could 
try to randomize the underlying graph in a different way.

\begin{itemize}
 \item Let $\cG$ be infinite TLG, suppose we run site or bond percolation\index{Percolation} on $\cG$,
and then on the connected component we define a Markov process. What properties
will the process have?
\end{itemize}

\appendix 
{\small 
\chapter{Independence and processes}

\section{Conditional independence and expectations}\index{Conditional independence|(}
The results in this section are taken from Section 21.5. in \cite{modernapp}. 

We will often use conditional independence, so we need to define it.

\defi{Let $(\Omega,\F,\P)$ be a probability space and $\F_1$, $\F_2$, and 
${\cal G}$ sub-$\sigma$-fields of $\F$. The $\sigma$-fields $\F_1$ and $\F_2$ are \textbf{conditionally 
independent}\index{Conditional independence|textbf} given ${\cal G}$ if 
$$\P(A_1\cap A_2| {\cal G})=\P(A_1| {\cal G}) \P(A_2| {\cal G})\ \ \ a.s.$$
for all $A_1\in \F_1$ and $A_2\in \F_2$.}

\prop{Let $(\Omega,\F,\P)$ be a probability space and $\F_1$, $\F_2$, and 
${\cal G}$ sub-$\sigma$-fields of $\F$, and suppose that $\F_2\subset {\cal G}$. Then 
$\F_1$ and $\F_2$ are conditionally independent.}

\prop{Let $\cG$, $\cH$, and $\cK$ be $\sigma$-fields of events in a probability
space. If $\cG$ and $\cH$ are conditionally independent given $\cK$, then $\cG$
and $\sigma(\cH,\cK)$ are conditionally independent  given $\cK$.}

\prop{Let $\cG$ and $\cH$ be two $\sigma$-fields of events in a probability space,
and let $\cG_1$ and $\cH_1$ be sub-$\sigma$ fields of $\cG$ and $\cH$, receptively.
Suppose that $\cG$ and $\cH$ are independent. Then $\cG$ and $\cH$ are conditionally
independent given $\sigma(\cG_1,\cH_1)$.}
\index{Conditional independence|)}

\subsection*{Conditional expectations}
\index{Conditional expectation|(}
\prop{Let $X$ be $(\Psi,{\cal H})$-valued random variable on a probability
space $(\Omega,\F,\P)$ and suppose that a conditional distribution $Z$ of $X$ given ${\cal G}$
exists where ${\cal G}$ is a sub-$\sigma$-field of ${\cal F}$. Let 
$f$ denote a $\overline{\R}$-valued function on $(\Psi,{\cal H})$. Then
$$\E(f(X)|{\cal G})=\int_{\Psi}f(x)Z(dx) \quad a.s.$$ }

\prop{For $i=1,2$, let $X_i$ be a $(\Psi_i,{\cal H}_i)$-valued random variable on a probability
space $(\Omega,\F,\P)$ and let ${\cal G}$ be a sub-$\sigma$-field of ${\cal F}$, such that
$X_2$ is measurable with respect to ${\cal G}$. Suppose that each $(\Psi_i,{\cal H}_i)$ is a Borel space.
Let $f$ be a measurable $\overline{\R}$-valued function defined on $(\Psi_1 ,{\cal H}_1)\times (\Psi_2 ,{\cal H}_2)$. If $Q_1$ 
is the distribution of $X_1$, then
$$\E(f(X_1,X_2)|{\cal G})(\omega)=\int_{\Psi_1}f(x,X_2(\omega))Q_1(dx|{\cal G})(\omega) \ \ \ a.s.$$
in the sense that the set of $\omega$ such that one side exist but the other does not 
is a null event. }\index{Conditional expectation|)}

\section{Construction of a conditional sequence}

\lem{\label{cond_seq:1}Let $(\Psi_0,\cG_0)$ and $(\Psi_1,\cG_1)$ be two measurable spaces, let $R_0$
denote the probability measure on $(\Psi_0,\cG_0)$, and let $x_0\mapsto R_1(x_0,\cdot)$
be a random distribution on $(\Psi_1,\cG_1)$ whose domain is the probability space
$(\Psi_0,\cG_0,R_0)$. Then there is a unique distribution $Q$ on $(\Psi_0\times\Psi_1 ,\cG_0\times \cG_1)$
such that if $X=(X_0,X_1)$ is any valued $\Psi_0\times\Psi_1$-valued random variable 
having distribution $Q$, then $R_0$ is the distribution of $X_0$ and $R_1$ is a conditional 
distribution of $X_1$ given $\sigma(X_0)$. Moreover $Q$ is given by
$$Q(A)=\int_{\Psi_0}\int_{\Psi_1}\1_A(x_0,x_1)R_1(x_0,\,dx_1)R_0(dx_0)$$
for $A\in \cG_0\times \cG_1$.}

\teo{\label{cond_seq:2}{\sc (Conditional Fubini)}\index{Conditional Fubini's Theorem} Let $(\Psi_0,\cG_0)$ and $(\Psi_1,\cG_1)$ be two measurable spaces and let 
$$(\Omega, \F)=(\Psi_0,\cG_0)\times (\Psi_1,\cG_1).$$ 
Let $R_0$, $R_1$, and $Q$ be as in Lemma \ref{cond_seq:1}. If $f$ is and $\bar{\R}$-valued measurable function 
defined on $(\Omega,\F,Q)$ whose integral with respect to $Q$ exists, then the function 
$$x_0\mapsto \int_{\Psi_1}f(x_0,x_1)R_1(x_0,dx_1)$$
is an $R_0$-almost surely defined $\cG_0$-measurable function, and 
$$\int_{\Omega}f\, dQ=\int_{\Psi_0}\int_{\Psi_1}f(x_0,x_1) R_1(x_0,dx_1)R_0(dx_0).$$}

\teo{\label{cond_seq:inf} Let $(\Psi_n,\cG_n)_{n\geq 0}$ be a sequence of measurable spaces. Let $R_0$
be a probability measure on $\cG_0$, and for each $n\geq 0$, let $R_{n+1}$ be a measurable function from 
$(\Psi_0,\cG_0)\times \ldots \times (\Psi_n,\cG_n)$ to the measurable space of probability measures 
on $(\Psi_{n+1},\cG_{n+1})$. Then there exists a probability space $(\Omega,\F,\P)$
and a random sequence $(X_k:k=0,\ldots )$ defined on the space such that the distribution
of $X_0$ is $R_0$, and for $n\geq 0$, conditional distribution of $X_{n+1}$
given $\sigma(X_0,\ldots, X_n)$ is given by 
$$\omega \mapsto R_{n+1}(X_0(\omega),X_1(\omega), \ldots , X_n(\omega), \cdot).$$
The distribution of $X$ is uniquely determined by the relations 
$$
 \P((X_0,\ldots,X_n)\in A_n)
= \int_{\Psi_0}\ldots \int_{\Psi_n} \1_A(x_0,\ldots, x_n)R_n((x_0,\ldots, x_{n-1}),dx_n)\ldots R_0(dx_0),
$$
$n\in \N_0$ and $A_n\in \cG_0\times \ldots \times \cG_n$.

}
\section{Markov and Brownian bridges}\index{Markov bridge|(}\label{mk_bd}

The best way to describe a Markov bridge $(Y_t)$ is as a Markov process on the time
interval $[s,u]$ conditioned that we know the value of the process at times $s$ and $u$.\vspace{0.3cm}

In oder to construct such a process we need to see what is happening with a Markov process
when we condition it on the outside of that interval. Here we will prove a slight generalization of the result stated in \cite{mark_bridge}.
In this section we are working on a probability space $(\Omega,\F,\P)$
until we extend it later. The proof of the following theorem can be found in 
\cite{mark_bridge}.

\teo{({\sc Two-sided Markovian property})\label{m_most:1}Let $(X_t)_{t\in T}$ be a Markov process with respect to the filtration $(\F_t)_{t\in T}$, and let $\G_t=\sigma\{X_u:u\geq t\}$. For $s< u$ in $T$ and $T'\subset T\cap [s,u]$, if $Y$
is a bounded $\sigma\{X_t:t \in T'\}$-measurable random variable then
$$\E(Y|X_s,X_u)= \E(Y|\F_s\vee \G_u)\quad a.s.$$ 
}
\pos{\label{m_most:2}Let $(X_t)_{t\in T}$ be a Markov process with respect to the filtration $(\F_t)_{t\in T}$, and let $\G_t=\sigma\{X_u:u\geq t\}$. For $s< u$ and $t\in [s,u]$ in $T$, if $f$ is a bounded $\R$-valued measurable function on the state 
space, then
$$\E(f(X_t)|X_s,X_u)= \E(f(X_t)|\F_s\vee \G_u)\quad a.s.$$ }

If $T\subset \R$ be a closed finite interval, and $(X_t)_{t\in T}$ is RCLL (or continuous) process with real values. 
Then $X$ can be viewed as a random map
into a Borel space $(\Sigma, \cS)$ consisting of all $x\in \R^T$, such that $t\mapsto x_t$ is RCLL (or continuous)
with the usual Skorohod (or uniform) topology. (See \cite{billing} for more on this.) 
Under those conditions, since the space of RCLL functions on  a compact set
 is  a Borel space, we can define a conditional
probability $\mu(\omega,\cdot ) $ for $\omega \in \Omega$ such that

\be \mu(\omega,H)=\P(X^{-1}(H)|X_u,X_s)(\omega),\label{mb:jd1}\ee
where $H$ is an element in the $\sigma$-algebra of that Borel space,
for $\P$-almost all $\omega\in \Omega$.

Specially, since coordinate projection $\pi_t:\R^T\to \R$ are measurable, $\mu$ we get the conditional
distribution for each $X_u$:
$$\mu(\omega, \pi_u^{-1}(A))=\P(X_t\in A|X_u,X_s).$$

A property of this random measure.
\prop{\label{mp:promb}For $u\in \{s,t\} $, we have
$$\mu(\cdot , \pi_u^{-1}(A))=\delta_{X_u}(A).$$
\dok{Since $\1_A(X_u)$ is $\F_s\vee \G_t$-measurable, from Corollary \ref{m_most:2} we have 
\begin{multline*}
 \mu(\cdot , \pi_u^{-1}(A))=\P(X_u\in A|X_s,X_t)=\E(\1_A(X_u)|X_s,X_t)=\\ =\E(\1_A(X_u)|\F_s\vee \G_t)=\1_A(X_u)=\delta_{X_u}(A).
\end{multline*}
%
}}

 $\P$-almost all $\omega \in \Omega$ the measure $\mu(\omega,\cdot)$ on $(\Sigma,\cS)$ defines a random map
$Y$ such that $Y_u=X_u(\omega)$ $\mu(\omega,\cdot)$-a.s. for $u\in [0,s]\cup[t,\infty)$.

\subsection*{The construction}

We will focus on RCLL (or continuous) Markov process $(X(t) : t\in [0,1])$
with distribution $D$.\vspace{0.3cm}

\defi{For times $t_1<t_2$ in $[0,1]$ we say that a  process
$(Y(t) : t\in [t_1,t_2])$ is a \textbf{Markov bridge}\index{Markov bridge|textbf} between $(t_1,y_{t_1})$ and $(t_2,y_{t_2})$ on some probability space if :
\begin{itemize}
 \item $Y_{t_1}=y_{t_1}$ and $Y_{t_2}=y_{t_2}$;
 \item The distribution of $(Y_t:t\in [t_1,t_2])$ 
is the same as $(X(t):t\in [t_1,t_2])$ given $(X(t_1)=y_{t_2},X(t_2)=y_{t_2})$.
\end{itemize}

}

\teo{A Markov bridge between $(t_1,X_{t_1})$ and $(t_2,X_{t_2})$ exists, for $D$-almost all values of $(X_{t_1},X_{t_2})$. 
\dok{The process $(X_t:t\in [t_1,t_2])$ is still Markov and RCLL (or continuous). 
Now, from the previous discussion (see (\ref{mb:jd1})) and since the space of RCLL functions on  a compact set
 is  a Borel space, there
exists $R$ such that
$$R(X(t_1),X(t_2))(\cdot )=\P(X\in\cdot|X(t_1),X(t_2)).$$
Now the measure $B\mapsto R(X(t_1),X(t_2))(B)$ defines a  process $Y$ on $D[t_1,t_2]$ (or $C[0,1]$).
From the Proposition \ref{mp:promb}. we get that $Y(t_1)=X(t_1)$ and $Y(t_2)=X(t_2)$ $R(X(t_1),X(t_2))$ - a.s.}} 

Often we will have a probability space a Markov process $X$ and maybe some
other process $Y$ on that space, and we will need to extend that process
to get construct an additional Markov bridge of the process $X$.

\teo{\label{2mrk_pro}Let $(\Omega_0,\F_0,\P_0)$ be a probability space, $(X:t\in [0,1])$ a RCLL (or continuous) Markov process with distribution $D$,
and $Z$ some other random element on that space. Assume $Q$ is the law of the Markov process
on $[t_1,t_2]$, where $D\circ \pi_{t_1,t_2}^{-1}=Q\circ \pi_{t_1,t_2}^{-1}$. 
Then for all $t_1<t_2$ in $[0,1]$ there exist a probability
space $(\Omega,\F,\P)$ with a process $(\hat{X}(t):t\in [0,1])$, random element $\hat{Z}$, 
and a $Q$ - Markov bridge $(Y(t):t\in [t_1,t_2])$ between $(t_1,\hat{X}(t_1))$ and $(t_2,\hat{X}(t_2))$  such that:
\begin{itemize}
 \item The joint distribution of $(X,Z)$ is the same as of $(\hat{X},\hat{Z})$;
\item $(\hat{X},\hat{Z})$ and $Y$ are conditionally independent given $(\hat{X}(t_1),\hat{X}(t_2))$.
\end{itemize}
The process $\tilde{X}$ given by $\tilde{X}(t)=\hat{X}(t)$ for $t\in [0,t_1]\cup [t_2,1]$
and $\tilde{X}(t)=Y(t)$ for $t\in (t_1,t_2)$ is a Markov process.  
Further, if $D$ on $[t_1,t_2]$ is distributed  as $Q$ then $\tilde{X}$ has the
same distribution ($D$) as $X$. }

\dok{We construct a Markov bridge and the space $(\Omega,\F,\P)$ using Lemma \ref{cond_seq:1}.
Let's prove that the process $\tilde{X}$ is Markov. Pick $u\in [0,1]$ and with $A_{lk}$ we denote a set in $\sigma(X_t:t\in [l,k])$. 

If $u\in  (t_1,t_2)$, and let $B_u\in \sigma(\tilde{X}_u)$ 
then when we condition on $\tilde{X}(t_1)$ and $\tilde{X}(t_2)$ from the construction we have
\begin{align*}
&\E(\1_{A_{0t_1}}\1_{A_{t_1u}}\1_{A_{ut_2}}\1_{A_{t_21}}\1_{B_u})\\
 &=\E(\E(\1_{A_{0t_1}}\1_{A_{t_1u}}\1_{A_{ut_2}}\1_{A_{t_21}}\1_{B_u}|\tilde{X}(t_1),\tilde{X}(t_2)))\\
 &=\E(\1_{A_{0t_1}}\1_{A_{t_21}}\E(\1_{A_{t_1u}}\1_{B_u}\1_{A_{ut_2}}|\tilde{X}(t_1),\tilde{X}(t_2)))
\end{align*}
Now using the Markov property of the process $\hat{X}$ we have
\begin{align*}
 &=\E(\E(\1_{A_{0t_1}}\1_{A_{t_21}}\E(\1_{A_{t_1u}}\1_{B_u}\1_{A_{ut_2}}|\tilde{X}(t_1),\tilde{X}(t_2))|X(t_1)))\\
  &=\E(\E(\1_{A_{0t_1}}|X(t_1))\1_{A_{t_21}}\E(\1_{A_{t_1u}}\1_{B_u}\1_{A_{ut_2}}|\tilde{X}(t_1),\tilde{X}(t_2)))\\
&=\E(\E(\E(\1_{A_{0t_1}}|X(t_1))\1_{A_{t_21}}\E(\1_{A_{t_1u}}\1_{B_u}\1_{A_{ut_2}}|\tilde{X}(t_1),\tilde{X}(t_2))|X(t_2)))\\
&=\E(\E(\1_{A_{0t_1}}|X(t_1))\E(\1_{A_{t_21}}|X(t_2))\E(\1_{A_{t_1u}}\1_{B_u}\1_{A_{ut_2}}|\tilde{X}(t_1),\tilde{X}(t_2)))
\end{align*}
Now, again using the properties of the conditional expectation we have 
\begin{align*}
 &=\E(\E(\E(\1_{A_{0t_1}}|X(t_1))\E(\1_{A_{t_21}}|X(t_2))\1_{A_{t_1u}}\1_{B_u}\1_{A_{ut_2}}|\tilde{X}(t_1),\tilde{X}(t_2)))\\
 &=\E(\E(\1_{A_{0t_1}}|X(t_1))\E(\1_{A_{t_21}}|X(t_2))\1_{A_{t_1u}}\1_{B_u}\1_{A_{ut_2}}).
\end{align*}
Since $(\tilde{X}(t):t\in [t_1,t_2])$
is a $Q$-Markov process, conditioning on $\tilde{X}(u)$ we get
\begin{align*}
 &=\E(\E(\E(\1_{A_{0t_1}}|X(t_1))\E(\1_{A_{t_21}}|X(t_2))\1_{A_{t_1u}}\1_{B_u}\1_{A_{ut_2}}|\tilde{X}(u)))\\
&=\E(\E(\E(\1_{A_{0t_1}}|X(t_1))\E(\1_{A_{t_21}}|X(t_2))\1_{A_{t_1u}}\1_{A_{ut_2}}|\tilde{X}(u))\1_{B_u})\\
&=\E(\E(\E(\1_{A_{0t_1}}|X(t_1))\1_{A_{t_1u}}|\tilde{X}(u))\E(\E(\1_{A_{t_21}}|X(t_2))\1_{A_{ut_2}}|\tilde{X}(u))\1_{B_u})\\
&=\E(\E(\1_{A_{0t_1}}|X(t_1))\1_{A_{t_1u}}\E(\E(\1_{A_{t_21}}|X(t_2))\1_{A_{ut_2}}|\tilde{X}(u))\1_{B_u})
\end{align*}
We again condition on $X(t_1)$ and $X(t_2)$ and we get
\begin{align*}
 &=\E(\E(\E(\1_{A_{0t_1}}|X(t_1))\1_{A_{t_1u}}\E(\E(\1_{A_{t_21}}|X(t_2))\1_{A_{ut_2}}|\tilde{X}(u))\1_{B_u}|X(t_1),X(t_2)))\\
 &=\E(\E(\1_{A_{0t_1}}|X(t_1))\E(\1_{A_{t_1u}}\E(\E(\1_{A_{t_21}}|X(t_2))\1_{A_{ut_2}}|\tilde{X}(u))\1_{B_u}|X(t_1),X(t_2)))
\end{align*}
Now, using Markov property of the process $\hat{X}$, and later the construction we get
\begin{align*}
  &=\E(\1_{A_{0t_1}}\E(\1_{A_{t_1u}}\E(\E(\1_{A_{t_21}}|X(t_2))\1_{A_{ut_2}}|\tilde{X}(u))\1_{B_u}|X(t_1),X(t_2)))\\
 &=\E(\1_{A_{0t_1}}\1_{A_{t_1u}}\E(\E(\1_{A_{t_21}}|X(t_2))\1_{A_{ut_2}}|\tilde{X}(u))\1_{B_u})
\end{align*}
Now we again condition everything on $\tilde{X}(u)$ and using properties of the conditional
expectation we get:
\begin{align*}
 &=\E(\E(\1_{A_{0t_1}}\1_{A_{t_1u}}|\tilde{X}(u))\E(\E(\1_{A_{t_21}}|X(t_2))\1_{A_{ut_2}}|\tilde{X}(u))\1_{B_u})\\
&=\E(\E(\E(\1_{A_{0t_1}}\1_{A_{t_1u}}|\tilde{X}(u))\E(\1_{A_{t_21}}|X(t_2))\1_{A_{ut_2}}\1_{B_u}|\tilde{X}(u)))\\
&=\E(\E(\1_{A_{0t_1}}\1_{A_{t_1u}}|\tilde{X}(u))\E(\1_{A_{t_21}}|X(t_2))\1_{A_{ut_2}}\1_{B_u})
\end{align*}
Again conditioning on $X(t_1)$ and $X(t_2)$, and using Markov property of $\hat{X}$, and the construction we get
\begin{align*}
&=\E(\E(\E(\1_{A_{0t_1}}\1_{A_{t_1u}}|\tilde{X}(u))\1_{A_{ut_2}}\1_{B_u}|X(t_1),X(t_2))\E(\1_{A_{t_21}}|X(t_2)))\\
&= \E(\E(\E(\1_{A_{0t_1}}\1_{A_{t_1u}}|\tilde{X}(u))\1_{A_{ut_2}}\1_{B_u}|X(t_1),X(t_2))\1_{A_{t_21}})\\
&= \E(\E(\1_{A_{0t_1}}\1_{A_{t_1u}}|\tilde{X}(u))\1_{A_{ut_2}}\1_{B_u}\1_{A_{t_21}})
\end{align*}
Finally, conditioning on $\tilde{X}(u)$ we get
\begin{align*}
 &= \E(\E[\E(\1_{A_{0t_1}}\1_{A_{t_1u}}|\tilde{X}(u))\1_{A_{ut_2}}\1_{B_u}\1_{A_{t_21}}|\tilde{X}(u)])\\
 &= \E(\E(\1_{A_{0t_1}}\1_{A_{t_1u}}|\tilde{X}(u))\E[\1_{A_{ut_2}}\1_{A_{t_21}}|\tilde{X}(u)]\1_{B_u}).
\end{align*}
This proves, using monotone class theorem  that $(\tilde{X}(t):t\leq u)$ and $(\tilde{X}(t):t\geq u)$ are conditionally independent
given $\tilde{X}(u)$.
 
When $u\in [0,t_1]\cup [t_2,1]$ this can be shown in a similar way.}
\subsection*{Brownian bridge}\index{Brownian bridge|(}

Brownian bridges are Markov bridges when the given Markov process is Brownian
motion.

The following representation holds.
\teo{\label{bbridg}For $0<t_1<t_2$ the process $(B^{br}(t):t\geq 0)$ given by
$$B^{br}(t)=\frac{t_2-t}{t_2-t_1}(x_1-W_{t_1})+W_t+\frac{t-t_1}{t_2-t_1}(x_2-W_{t_2}),$$
where $(W_t:t\geq 0)$ is Brownian motion has the same distribution 
as a Brownian bridge conditioned at times $t_1$ and $t_2$ to have values
$x_1$ and $x_2$. }

\pos{\label{difbridg}Let $(N(t):t\in [0,T])$ be given for each $t$ by the Ito integral 
$$N(t)=\int_0^t f(s)\, dB_s.$$
For $0\leq t_1<t_2\leq T$ the distribution of the process $N$ conditioned 
at times $t_1$ and $t_2$ to have values
$x_1$ and $x_2$ is the same as that of
$$N^{t_1,t_2}_{x_1,x_2}(t)=\frac{V(t_2)-V(t)}{V(t_2)-V(t_1)}(x_1-W_{V(t_1)})+W_{V(t)}+\frac{V(t)-V(t_1)}{V(t_2)-V(t_1)}(x_2-W_{V(t_2)}),$$
where $(W_t:t\geq 0)$ is Brownian motion and $V(t)=\int_0^t (f(s))^2\, ds$.}

\index{Brownian bridge|)}\index{Markov bridge|)}
\section{Markov random fields}\label{sec:mrf}

Let $G=(V,E)$ be a simple undirected graph, where $V$ is a finite set
of vertices and $E$ is a set of edges. We are looking a process $(X_v:v\in V)$.

\defi{\label{def:gr_mpr}The process $(X_v:v\in V)$ has a
\begin{enumerate}[(a)]
       \item \textbf{pairwise Markov property}\index{Pairwise Markov property} if for all $v,u\in V$ such that $\{u,v\}\notin E$
we have 
$$X_v\perp X_u | X_{V\setminus\{v,u\}};$$

\item \textbf{local Markov property}\index{Local Markov property} if for all $v\in V$ 
$$X_v \perp X_{V\setminus \{v\}}|X_{\{u: \{u,v\}\in E\})};$$
\item  \textbf{global Markov property}\index{Global Markov property|textbf} 
if  for every $A$, $B$ and $C$ subsets of $V$ such that $C$ separates 
$A$ and $B$, we have
$$X_A \perp X_B|X_C.$$
      \end{enumerate}
}

\defi{We say that the process $(X_v:v\in V)$ is a \textbf{Markov random field} (MRF) \index{Markov random field (MRF)|textbf} 
if it satisfies one of the three properties (a), (b) or (c) in Definition \ref{def:gr_mpr}.}

\lem{The global Markov property implies local Markov property, and the local Markov property 
implies the pairwise Markov property.}

If the random vector $(X_v : v\in V)$ has a positive density then 
we have several interesting results. (For more details see \cite{prgphmdl}.)

\teo{Let $X=(X_v:v\in V)$ have a positive density function $f$.
Then global, local, and pairwise Markov properties are equivalent.}

The following theorem was proven in an unpublished paper by Hammeresley and Clifford. There have been several
proofs published obtained in different ways, see for example  \cite[Grimmett (1973)]{grimmett} or 
\cite[Clifford (1990)]{clifford}.

\teo{(Hammeresley-Clifford, 1971)\index{Markov random field (MRF)!Hammeresley-Clifford Theorem} Let $X=(X_v:v\in V)$ be a continuous or discrete random vector 
with a positive density function $f$. $X$ is a Markov random field if and only if 
$f$ is of the form 
$$f(x)=\frac{1}{Z}\prod_{C\in {\cal C}(G)}\phi_C(x_C),$$
where ${\cal C}(G)$ is the set of all maximal cliques in $G$.}
\section{White noise}\label{sec:whtnos}\index{White noise|(}
In this section we define the one dimensional white noise on $\R^n$\index{White noise|textbf}.
This is a mean-zero Gaussian process indexed by Borel $\sigma$-algebra on $\R^n$
($\cB(\R^n)$), i.e. 
$$(\mathbb{W}(A)\,:\, A\in \cB(\R^n)),$$
with the covariance function
\be\Sigma (A,B)=\lambda (A\cap B),\label{wn1}\ee
where $\lambda$ is the Lebesgue measure, and $A,B\in \cB(\R^n)$.

\lem{The function $\Sigma :\cB(\R^n)\times \cB(\R^n)\to \R$ given by $(\ref{wn1})$ is symmetric and positive definite.}

By Kolmogorov's Consistency Theorem, the process $\mathbb{W}$ exists, and has the following properties:

\teo{Let $\mathbb{W}=(\mathbb{W}(A)\,:\, A\in \cB(\R^n))$ be the white noise on $\R^n$.
\begin{enumerate}[(a)]
 \item For all disjoint $A,B\in \cB(\R^n)$, $\W(A)$ and $\W(B)$ are independent.
\item For all $A,B\in \cB(\R^n)$, $\W(A\cup B)=\W(A)+\W(B)-\W(A\cap B)$ a.s.
\item If $A_1,A_2,\ldots \in \cB(\R^n)$ are disjoint and $\sum_{i=1}^{\infty}\lambda(A_i)<\infty$, then
a.s. 
$$\W\left(\bigcup_{i=1}^\infty A_i\right)=\sum_{i=1}^{\infty}\W(A_i).$$
\end{enumerate}
}

Although $\W$ is not a measure, it has enough properties (see details in Khoshnevisan)
that for $h\in L^2(\lambda)$ we can define the \textbf{Wiener integral}\index{Wiener integral|textbf}
$$W(h)=\int h(s)\W(ds).$$ 

The stochastic process $(W(h):h\in L^2(\lambda))$ is called the isonormal process.

\teo{The isonormal process $(W(h):h\in L^2(\lambda))$ is a mean zero Gaussian process
indexed by $L^2(\lambda)$ such that for all $h_1,h_2\in L^2(\lambda)$,
$$\E(W(h_1)W(h_2))=\int h_1h_2\, d\lambda.$$
Moreover, for every $\alpha,\beta \in \R$ and $f,g\in L^2(\lambda)$ 
$$W(\alpha f+\beta g)=\alpha W( f)+\beta W(g), \quad a.s. $$}
\index{White noise|)}
\section{The stochastic heat equation}\label{sec:shee}\index{Heat equation!stochastic|(}
\index{Stochastic heat equation|see{Heat equation}}

The usual heat equation is the initial value problem
\be \begin{array}{cl} \partial_t u=c\partial_{xx} u+f&\quad {\rm on}\quad  (0,\infty)\times \R,\\
    u(0,x)=g(x) & \quad {\rm for}\quad  x\in\R.
    \end{array}\label{pde0}
\ee

Under mild assumptions (see \cite[Folland]{inpde}) it is well known that the following is a solution to $(\ref{pde0})$:
\be u(t,x)=\frac{1}{2\sqrt{\pi c t}}\int_{\R}e^{-\frac{|x-y|^2}{4c t}}g(y)\, dy+\int_0^t\frac{1}{2\sqrt{\pi c(t-s)}}\int_{\R}e^{-\frac{|x-y|^2}{4c(t-s)}}f(s,y)\, dy\, ds\ee
 
The idea of the stochastic heat equation is to replace the external force $f$,
with random noise, in our case the white noise $\W$. So the stochastic heat equation will be given by
\be \begin{array}{cl} \partial_t u=c\partial_{xx} u+\sigma\W&\quad {\rm on}\quad  (0,\infty)\times \R,\\
    u(0,x)=g(x) & \quad {\rm for}\quad  x\in\R.
    \end{array}\label{spde0}
\ee
where $\sigma:\R^+\times \R\to \R$ is a {\it nice} function.
The so called {\it mild} solution to $(\ref{spde0})$ is 
\be u(t,x)=\frac{1}{2\sqrt{\pi c t}}\int_{\R}e^{-\frac{|x-y|^2}{4c t}}g(y)\, dy+\int_0^t\frac{1}{2\sqrt{\pi c(t-s)}}\int_{\R}e^{-\frac{|x-y|^2}{4c(t-s)}}\sigma(s,y)\W(ds,dy).\ee
\vspace{0.2cm}

We will state the results adapted from \cite{shiga} (see the appendix of the paper). 
The case that will interest us is the case when we have a boundary condition, and 
$c=\frac{1}{2}$ and $\sigma=1$ are constants: 
\be \begin{array}{cl} \partial_t u=\frac{1}{2}\partial_{xx} u+\W&\quad {\rm on}\quad  (0,\infty)\times \R^+,\\
    u(0,x)=g(x) & \quad {\rm for}\quad  x\in\R^+.\\
    u(t,0)=0 & \quad {\rm for}\quad  t\geq 0.
    \end{array}\label{spde01}
\ee

We need to define precisely what the solution of this equation is, and 
when it is unique (and in what sense). The following definition and results
have been taken from \cite[Section 3 \& 4]{inv_meas}, where more general
result were obtained and by modification of results from \cite{shiga}. \vspace{0.3cm}

First we will define a space of $C_{tem}(\R^+)$, and we will require that
for all $t\geq 0$ the function $u(t,\cdot) \in C_{tem}(\R^+)$.
\defi{We denote by $C_{tem}(\R^+)$ the family of all continuous functions $f:\R^+\to
\R$ satisfying 
$$\|f\|_{(-\lambda)}=\sup_{x\in \R^+}|e^{-\lambda|x|}f(x)|<\infty,$$
for all $\lambda >0$.}

\defi{\begin{enumerate}[(a)]
       \item We call a random function $\{u=u(t,x),t\geq 0, x\in \R^+\}$ a \textbf{weak solution}\index{Heat equation!stochastic!weak solution}
of the SPDE (\ref{spde01}) with an initial value $u_0\in C_{tem}(\R^+)$ if it is $(\F_t)$
adapted and has the following two conditions:
\begin{itemize}
 \item $u\in C([0,\infty), C_{tem}(\R^+))$, a.s.
 \item For every $\varphi\in C_{c}^{\infty}(\R^+)$ such that $\varphi(0)=0$, the following is
satisfied:
\begin{multline*}
 \int_{\R^+} u(t,x)\varphi(x)\, dx=\int_{\R^+}u_0(x)\varphi(x)\, dx+\\ 
+\frac{1}{2}\int_0^t\int_{\R^+}u(s,x)\varphi''(x)\, dx\, ds+\int_0^t\int_{\R^+}\varphi(x) \W(ds,dx)
\end{multline*}

\end{itemize}
\item We call $u$ under the same assumptions a \textbf{mild solution}\index{Heat equation!stochastic!mild solution} if the following holds
\begin{multline}
 u(t,x)=\frac{1}{\sqrt{2\pi  t}}\int_{\R^+}\left(e^{-\frac{|x-y|^2}{2 t}}-e^{-\frac{|x+y|^2}{2 t}}\right)g(y)\, dy+\\ +\int_0^t\frac{1}{\sqrt{2\pi (t-s)}}\int_{\R^+}\left(e^{-\frac{|x-y|^2}{2(t-s)}}-e^{-\frac{|x+y|^2}{2(t-s)}}\right)\W(ds,dy).
\label{wkslhe}
\end{multline}
($u$ is a $C_{tem}$-version of the integral on the right.)
\item We say that the pathwise uniqueness of the weak solution of the SPDE (\ref{spde01}) holds if for 
arbitrary two weak solutions $u^{(1)}$ and $u^{(2)}$ of the SPDE (\ref{spde01}) with the respect to the same 
filtration $(\Omega, \F,(\F_t),\P)$ and the same noise $\W$ we have 
$$\bigcap_{t\geq 0}\{u^{(1)}(t,\cdot)\neq u^{(2)}(t,\cdot)\}\subset N,$$
where $N\in \F$ such that $\P(N)=0$.
      \end{enumerate}
}
In order to show that a $C_{tem}$-version of (\ref{wkslhe}) exists we will need
the following results.

\lem{If $\phi:\R^+\times \R^+\to\R$ is in $L^2$, for each $p>0$ there exists a constant 
$C_p>0$ such that 
\be\E\left[\left(\int_0^t\int_{\R^+}\phi(s,x)\W(ds,dx)\right)^{2p}\right]\leq C_p\left(\int_0^t\int_{\R^+}\phi(s,x)^2ds\, dx\right)^{p}\ee
}
\lem{(i) There exists a constant $C>0$ such that
$$\int^{t\vee t'}_0\int_{\R}(G(t-s,x,y)-G(t'-s,x',y))^2\, ds\, dy\leq C(|t-t'|^{1/2}+|x-x'|) $$
for $t,t'\geq 0$ and $x,x'\in \R$, where $G(t,x,y)=(2\pi t)^{-1/2}\exp(-(x-y)^2/(2t))$ 
for $t> 0$ and
$G(t,x,y)=0$ if $t\leq 0$.

(ii) For every $\lambda \in \R$ and $T>0$ 
$$\sup_{0\leq t\leq T}\sup_{x\in \R}e^{-\lambda |x|}\int_{\R}G(t,x,y)e^{\lambda |y|}dy<\infty.$$}
\teo{\label{spde:sol}If $g\in C_{tem}$ a.s., the following claims are true:
\begin{enumerate}[(a)]
 \item The SPDE (\ref{spde01}) has a at most one pathwise unique weak solution.
 \item If $u$ is a mild solution to the SPDE (\ref{spde01})
then it is also a weak solution.
\end{enumerate}
}

\lem{Brownian motion is in $C_{tem}$ a.s.
\dok{The claim follows from the strong law of large numbers for the Brownian motion, that is if
$\lambda>0$ then 
$$\lim_{x\to\infty}e^{-\lambda x}W_x=\lim_{x\to\infty} (xe^{-\lambda x})\frac{W_x}{x}=0\cdot 0=0.$$}}\vspace{0.3cm}

The following result is Lemma 4.4. from \cite{inv_meas}.
\teo{\label{spde:inv}The Wiener measure is an invariant measure for the SPDE (\ref{spde01}),
i.e. if $g$ is Brownian motion, then for each $x\in \R$ the process $t\mapsto u(x,t)$ is also Brownian motion.}
\index{Heat equation!stochastic|)}

\section{Crump - Mode - Jagers trees}\label{CMJtree}\index{Crump - Mode - Jagers trees|(}
Here we present an introduction to Crump -Mode - Jagers model which we will later mention in 
the context of time-like trees. We will use the notation given by Dawson in \cite{stoch_pop_sys}.%
\vspace{0.3cm}

First some notation. We define $I=\{\emptyset\}\cup\bigcup_{n=1}^{\infty} \N^n$.
Given $u=(u_1,\ldots, u_m), v=(v_1,\ldots, v_n)\in I$ we denote the composition by
$uv:=(u_1,\ldots, u_m,v_1,\ldots, v_n)$.

\defi{A \textbf{plane rooted tree}\index{Plane rooted tree} $T$ with root $\emptyset$ is a subset of $I$ such that:
\begin{enumerate}
 \item $\emptyset \in T$,
 \item If  $v=uw\in T$ for some $u\in I$ and $w\in I$, then $u\in T$.
 \item For every $u\in T$, there exists a number $k_u(T)\geq 0$, such that 
$uj\in T$ if and only if $1\leq j\leq k_{u}(T)$.
\end{enumerate}
Set $\mathbb{T}$ to be the \textbf{set 
of all plane rooted trees}. For $u\in T$ define the \textbf{level of the vertex} to be 
$|u|=|(u_1,\ldots, u_m)|=m$.}\vspace{0.2cm}

A plane tree $T$ can be given a structure of a graph in which $uw\in T$ is   
\textbf{descendant} of $u$. Specially, $(u)(j)\in T$ is the \textbf{child} of $u$. \vspace{0.2cm}

Consider the following process: For each individual $x\in I$
\begin{itemize}
 \item  We denote his \textbf{birth time} $\tau_x$.
 \item \textbf{Lifetime} $\lambda_x$.
\item  Point process $\xi_x$ denoting \textbf{reproduction} function. ($\xi_x(t)$ is the number of offsprings 
produced by individual $x$ born at 0 during $[0,t]$. )
\item Assume that the pairs $(\lambda_x,\xi_x)$ are i.i.d.
\item Assume $\P(\xi_x(\lambda_x,\infty)=0)=1$. (Offsprings can't be produced after 
$x$ is no longer alive.) 
\end{itemize}

\begin{figure}[ht]
\begin{center}
\psfrag{0}{$\tau_{x'}$} \psfrag{1}{$\tau_{(x',1)}$}
\psfrag{2}{$\tau_{(x',2)}$}\psfrag{3}{$\tau_{x'}+\lambda_{x'}$}
\psfrag{4}{$\tau_{(x',1,2)}$} \psfrag{5}{$\tau_{(x',1)}+\lambda_{(x',1)}$}
\includegraphics[width=14cm]{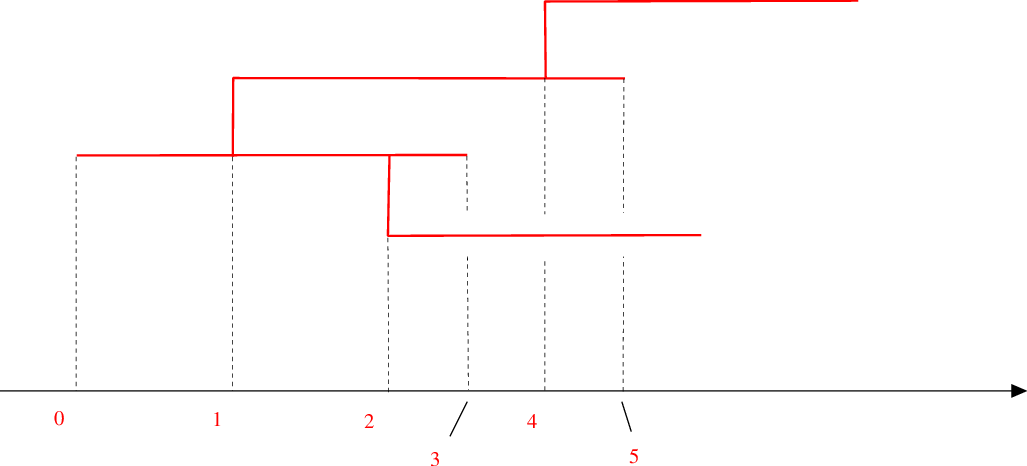}
\caption{Crump-Mode-Jagrers tree}
\end{center}
\end{figure}


The probability space that we are working in is 
$$(\Omega,\F,\P)=\prod_{x\in I}(\Omega_x,\F_x,\P_x),$$
where each $(\Omega_x,\F_x,\P_x)$ supports $(\lambda_x,\xi_x)$.\vspace{0.3cm}

We can determine the birth times $\{\tau_x:x\in I\}$ as follows,
$$\tau_{\emptyset} =0,$$ $$\tau_{(x',i)}=\tau_{x'}+\inf\{u: \xi_{x'}(u)\geq i\}.$$

The natural question that one may many individuals were born in the the time 
period $[0,t]$. Is that number even finite? We will introduce some results 
on this.

Set $\mu(t):=\E(\xi(t))$, and we define
$$T_t=\sum_{x\in I}\1_{(\tau_x\leq t)},$$
to be the number of individuals born up to time $t$. The following two results are form \cite{jagers} (Theorem 6.2.1. and Theorem 6.2.2. pages 126-127).

\teo{If $\mu(0)>1$, then for all $t\geq 0$, $\P(T_t=\infty)>0$.}
\teo{\label{teo:fin_br_pr}If $\mu(0)<1$ and $\mu(t)$ is finite for some $t>0$, then
$$\P(\forall t: T_t<\infty)=1.$$}

\index{Crump - Mode - Jagers trees|)}

\section{Branching Markov processes and branching Brownian motion}

The following is a definition given in \cite{etheridge} of the branching Brownian motion.

\defi{\textbf{Branching Brownian motion}\index{Branching Brownian motion|textbf} has three ingredients:

\begin{itemize}
 \item {\sc The spatial motion:} During its lifetime, each individual in the population moves around in $\R^d$
(independently of all other individuals) according to a Brownian motion.
 \item {\sc The branching rate $V$:} Each individual has an exponentially 
distributed lifetime with parameter $V$.
 \item {\sc The branching mechanism $\Phi$:} When it dies, and individual leaves behind (at the location
where it died) a random number of offsprings with probability generating function $\Phi(s)=\sum_{k=0}^\infty p_ks^k$.
Conditional on their time and place of birth, offsprings evolve independently of each other 
(in the same way as their parent).
\end{itemize}
}

We could have defined any Markov process on any Polish space $E$ to evolve in the same way, and in that 
case this would be the \textbf{branching Markov process}\index{Branching Markov process|textbf}.

For more details on the definition of the branching Markov process see \cite{brmp1} and \cite{brmp2}.
}
\backmatter
\chapter*{Acknowledgments}
The author would like to thank his advisor Krzysztof Burdzy for his advice and discussion during the work on this topic.

The author was partially supported by the:
\begin{itemize}
 \item NSF Grant DMS-1206276;
 \item MZOS grant 037-0372790-2799 of the Republic of Croatia;
 \item Croatian Science Foundation grant 3526.
\end{itemize}

%
%
%
%
%


\printindex
\end{document}